\pdfoutput=1
\documentclass[11pt,nodots,reqno,eqnsec]{style/ocelot}

\usepackage{bigints}
\usepackage{bookmark}
\usepackage{caption}
\usepackage{graphicx}
\usepackage{pdfpages}
\usepackage{subcaption}
\usepackage{ytableau}
\usetikzlibrary{arrows.meta, decorations.markings, calc, positioning, decorations.pathmorphing, decorations.shapes, shapes.geometric, calligraphy, patterns}
\tikzset{->-/.style={decoration={
  markings,
  mark=at position .5 with {\arrow{>}}},postaction={decorate}}}
\tikzset{-<-/.style={decoration={
  markings,
  mark=at position .5 with {\arrow{<}}},postaction={decorate}}}
\tikzset{-->-/.style={decoration={
  markings,
  mark=at position .75 with {\arrow[scale=0.75]{>}}},postaction={decorate}}}
\tikzset{--<-/.style={decoration={
  markings,
  mark=at position .35 with {\arrow[scale=0.75]{<}}},postaction={decorate}}}
\tikzset{->--/.style={decoration={
  markings,
  mark=at position .5 with {\arrow[scale=0.75]{>}}},postaction={decorate}}}
\tikzset{-<--/.style={decoration={
  markings,
  mark=at position .5 with {\arrow[scale=0.75]{<}}},postaction={decorate}}}
\ytableausetup{mathmode, boxsize=18pt, aligntableaux=center, centerboxes}

\DeclareSymbolFont{Letters}{OML}{cmm}{m}{it}
\DeclareMathSymbol{\psi}{\mathalpha}{Letters}{"20}
\newtheoremstyle{dotless}{}{}{\itshape}{}{\bfseries}{}{ }{}
\theoremstyle{dotless}

\newcommand{\Le}{\textnormal{\reflectbox{L}}}
\newcommand{\BLe}{\textnormal{\textbf{\reflectbox{L}}}}

\newcommand{\qmod}[2]{\left.\raisebox{-.1em}{$#1$}\middle\backslash\raisebox{.1em}{$#2$}\right.}

\renewcommand{\contentsname}{Table of contents}

\renewcommand{\thepart}{\Roman{part}}
\renewcommand{\thesection}{\thechapter.\arabic{section}}

\graphicspath{{figures/}{c2/}}
\newcommand{\QED}{%
  \begin{tikzpicture}[main node/.style={inner sep=0,minimum size =.09cm,circle,fill=black!80,draw},node distance = 0.75cm,scale=0.4]
    \node[main node] (a) at (0,0) {};
    \node[main node] (b) at (4,0) {};
    \node [above left=8mm of a] (a1) {};
    \node [below left=8mm of a] (a2) {};
    \node [above right=8mm of b] (b1) {};
    \node [below right=8mm of b] (b2) {};

    \draw[->-] (a) -- (a1);
    \draw[-<-] (a) -- (a2);
    \draw[-<-] (b) -- (b1);
    \draw[->-] (b) -- (b2);
    \draw[decoration=snake, decorate] (a) -- (b);
  \end{tikzpicture}
}

\newcommand{\EdgePosition}{%
    \node[main node, label=above:$x_1$] (a) at (0,0) {};
    \node[main node, label=above:$x_2$] (b) at (4,0) {};
    \node [above left=5mm of a] (a1) {};
    \node [below left=5mm of a] (a2) {};
    \node [above right=5mm of b] (b1) {};
    \node [below right=5mm of b] (b2) {};

    \draw (a) -- (a1);
    \draw (a) -- (a2);
    \draw (b) -- (b1);
    \draw (b) -- (b2);
    \draw[->-] (a) -- (b);

    \draw[->] (10,0) -- node[above]{position space} (15,0);
    \node at (20,0) {$\displaystyle\frac{1}{(x_1 - x_2)^2}$};
}

\newcommand{\EdgeMomentum}{%
    \node[main node] (a) at (0,0) {};
    \node[main node] (b) at (4,0) {};
    \node [above left=5mm of a] (a1) {};
    \node [below left=5mm of a] (a2) {};
    \node [above right=5mm of b] (b1) {};
    \node [below right=5mm of b] (b2) {};

    \draw (a) -- (a1);
    \draw (a) -- (a2);
    \draw (b) -- (b1);
    \draw (b) -- (b2);
    \draw[->-] (a) -- (b);
    \draw[->] (-0.5,1.5) to [bend right, looseness=1.1] node[above]{$p_1$} (4.5,1.5);
    \draw[->] (4.5,-1.5) to [bend right, looseness=1.1] node[below]{$p_2$} (-0.5,-1.5);

    \draw[->] (10,0) -- node[above]{momentum space} (15,0);
    \node at (22,0) {$\displaystyle\frac{1}{(p_1 - p_2 + \cdots)^2} = \frac{1}{p_e^2}$};
}

\newcommand{\Feynman}{%
  \begin{tikzpicture}[main node/.style={inner sep=0,minimum size =.09cm,circle,fill=black!80,draw},node distance = 0.75cm,scale=0.4]
    \EdgePosition
    \begin{scope}[shift={(0,-4)}]
    \EdgeMomentum
    \end{scope}
  \end{tikzpicture}
}

\newcommand{\GrphOne}{%
  \begin{tikzpicture}[main node/.style={inner sep=0,minimum size =.09cm,circle,fill=black!80,draw},node distance = 0.75cm,scale=0.4]
    \node[main node] (a) at (0,0) {};
    \node[main node] (b) at (6,0) {};
    \node [above left=5mm of a] (a1) {};
    \node [below left=5mm of a] (a2) {};
    \node [above right=5mm of b] (b1) {};
    \node [below right=5mm of b] (b2) {};

    \draw (a) -- (a1);
    \draw (a) -- (a2);
    \draw (b) -- (b1);
    \draw (b) -- (b2);
    \draw (a) to [bend left] node[above]{$\alpha_1$} (b);
    \draw (b) to [bend left] node[below]{$\alpha_2$} (a);
  \end{tikzpicture}
}

\newcommand{\CmplOne}{%
  \begin{tikzpicture}[main node/.style={inner sep=0,minimum size =.09cm,circle,fill=black!80,draw},node distance = 0.75cm,scale=0.4]
    \node[main node] (a) at (0,0) {};
    \node[main node] (b) at (6,0) {};
    \node[main node, label=above:$v$] (c) at (3,4) {};

    \draw (a) to [bend left=30] (b);
    \draw (b) to [bend left=30] (a);
    \draw (a) to [bend left=30] (c);
    \draw (c) to [bend left=30] (a);
    \draw (b) to [bend left=30] (c);
    \draw (c) to [bend left=30] (b);
  \end{tikzpicture}
}

\newcommand{\BlobBaseEmpty}{
    \node[main node] (f) at (0,0.75) {};
    \node[main node] (g) at (0,0) {};
    \node[main node] (i) at (0,-0.75) {};
    \draw (f) to [bend right=60] (g) to [bend right=60] (i);
    \draw (f) to [bend right=110, in=290, looseness=2.5] (i);
}

\newcommand{\BlobDodgson}{%
  \begin{tikzpicture}[main node/.style={inner sep=0,minimum size =.08cm,circle,fill=black!80,draw},node distance = 0.75cm,scale=0.8]
    \begin{scope}
      \BlobBaseEmpty
      \node[style={inner sep=2,circle,draw}] (v) at (f) {};
      \node[style={inner sep=2,circle,draw}] (v) at (g) {};
      \node[style={inner sep=2,circle,draw}] (v) at (i) {};
    \end{scope}
    \begin{scope}[shift={(1.5,0)}]
      \BlobBaseEmpty
      \node[style={inner sep=2,circle,draw}] (v) at (f) {};
      \node[style={inner sep=2,rectangle,draw}] (v) at (g) {};
      \node[style={inner sep=2,regular polygon,regular polygon sides=3,rotate=-90,inner sep=1.5pt,draw}] (v) at (i) {};
    \end{scope}
    \node at (2.25,0) {$=$};
    \begin{scope}[shift={(4,0)}]
      \BlobBaseEmpty
      \node[style={inner sep=2,circle,draw}] (v) at (f) {};
      \node[style={inner sep=2,circle,draw}] (v) at (g) {};
      \node[style={inner sep=2,rectangle,draw}] (v) at (i) {};
    \end{scope}
    \begin{scope}[shift={(5.5,0)}]
      \BlobBaseEmpty
      \node[style={inner sep=2,circle,draw}] (v) at (f) {};
      \node[style={inner sep=2,rectangle,draw}] (v) at (g) {};
      \node[style={inner sep=2,circle,draw}] (v) at (i) {};
    \end{scope}
    \node at (6.25,0) {$+$};
    \begin{scope}[shift={(8,0)}]
      \BlobBaseEmpty
      \node[style={inner sep=2,circle,draw}] (v) at (f) {};
      \node[style={inner sep=2,circle,draw}] (v) at (g) {};
      \node[style={inner sep=2,rectangle,draw}] (v) at (i) {};
    \end{scope}
    \begin{scope}[shift={(9.5,0)}]
      \BlobBaseEmpty
      \node[style={inner sep=2,rectangle,draw}] (v) at (f) {};
      \node[style={inner sep=2,circle,draw}] (v) at (g) {};
      \node[style={inner sep=2,circle,draw}] (v) at (i) {};
    \end{scope}
    \node at (10.25,0) {$+$};
    \begin{scope}[shift={(12,0)}]
      \BlobBaseEmpty
      \node[style={inner sep=2,circle,draw}] (v) at (f) {};
      \node[style={inner sep=2,rectangle,draw}] (v) at (g) {};
      \node[style={inner sep=2,circle,draw}] (v) at (i) {};
    \end{scope}
    \begin{scope}[shift={(13.5,0)}]
      \BlobBaseEmpty
      \node[style={inner sep=2,rectangle,draw}] (v) at (f) {};
      \node[style={inner sep=2,circle,draw}] (v) at (g) {};
      \node[style={inner sep=2,circle,draw}] (v) at (i) {};
    \end{scope}
  \end{tikzpicture}
}

\newcommand{\DTR}{%
  \begin{tikzpicture}[main node/.style={inner sep=0,minimum size =.09cm,circle,fill=black!80,draw},node distance = 0.75cm,scale=0.7]
    \begin{scope}
      \node[main node, label=above:$C$] (a) at (0,0) {};
      \node[main node, label=above right:$A$] (b) at (2,1) {};
      \node[main node, label=below right:$B$] (c) at (2,-1) {};
      \node[main node, label=above:$D$] (d) at (4,0) {};
      \node[main node] (e) at (2,2.5) {};
      \node[main node] (f) at (2,-2.5) {};
      \node[main node] (g) at (-1,1) {};
      \node[main node] (h) at (-1,-1) {};
      \node[main node] (i) at (5,1) {};
      \node[main node] (j) at (5,-1) {};

      \draw (a) -- (b) -- (d) -- (c) -- (a);
      \draw (e) -- (b) -- (c) -- (f);
      \draw (g) -- (a) -- (h);
      \draw (i) -- (d) -- (j);
    \end{scope}
    \node at (7,0) {$\rightsquigarrow$};
    \begin{scope}[shift={(10,0)}]
      \node[main node, label=above:$C$] (a) at (0,0) {};
      \node[main node, label=above right:$A$] (b) at (2,1) {};
      \node[main node, label=above:$D$] (d) at (4,0) {};
      \node[main node] (e) at (2,2.5) {};
      \node[main node] (f) at (2,-2.5) {};
      \node[main node] (g) at (-1,1) {};
      \node[main node] (h) at (-1,-1) {};
      \node[main node] (i) at (5,1) {};
      \node[main node] (j) at (5,-1) {};

      \draw (a) -- (b) -- (d) -- (a);
      \draw (e) -- (b) -- (f);
      \draw (g) -- (a) -- (h);
      \draw (i) -- (d) -- (j);
    \end{scope}
  \end{tikzpicture}
}

\newcommand{\GrphThreeEL}{%
  \begin{tikzpicture}[main node/.style={inner sep=0,minimum size =.09cm,circle,fill=black!80,draw},node distance = 0.75cm,scale=0.3]
    \node[main node] (a) at (0,0) {};
    \node[main node] (b) at (6,0) {};
    \node[main node] (c) at (6,6) {};
    \node[main node] (d) at (0,6) {};

    \draw (a) -- node[below]{$5$} (b) -- node[right]{$3$} (c) -- node[above]{$1$} (d) -- node[left]{$4$} (a);
    \draw (a) -- node[near end, above]{$2$} (c);
    \draw (b) -- node[near end, below]{$6$} (d);
  \end{tikzpicture}
}

\newcommand{\Triangle}{%
  \begin{tikzpicture}[main node/.style={inner sep=0,minimum size =.09cm,circle,fill=black!80,draw},node distance = 0.75cm,scale=0.20]
    \node[main node] (a) at (0,0) {};
    \node[main node] (b) at (6,0) {};
    \node[main node] (d) at (0,6) {};

    \draw (a) -- node[below]{$5$} (b);
    \draw (d) -- node[left]{$4$} (a);
    \draw (b) -- node[above]{$6$} (d);
  \end{tikzpicture}
}

\newcommand{\LoopOne}{%
  \begin{tikzpicture}[main node/.style={inner sep=0,minimum size =.09cm,circle,fill=black!80,draw},node distance = 0.75cm,scale=0.15]
    \node[main node] (a) at (0,0) {};
    \node[main node] (b) at (6,0) {};
    \draw (a) to [bend left] node[above]{$4$} (b);
    \draw (b) to [bend left] node[below]{$5$} (a);
    \draw (b) to [out=315,in=45,looseness=30] node[right]{$6$} (b);
  \end{tikzpicture}
}

\newcommand{\LoopTwo}{%
  \begin{tikzpicture}[main node/.style={inner sep=0,minimum size =.09cm,circle,fill=black!80,draw},node distance = 0.75cm,scale=0.15]
    \node[main node] (a) at (0,0) {};
    \node[main node] (b) at (6,0) {};
    \draw (a) to [bend left] node[above]{$6$} (b);
    \draw (b) to [bend left] node[below]{$5$} (a);
    \draw (a) to [out=135,in=225,looseness=30] node[left]{$4$} (a);
  \end{tikzpicture}
}

\newcommand{\Blob}{%
  \begin{tikzpicture}[main node/.style={inner sep=0,minimum size =.08cm,circle,fill=black!80,draw},node distance = 0.75cm,scale=0.8]
    \node[main node, label={[shift={(0,-2em)}]\scriptsize$u_2$}] (f) at (0,0) {};
    \node[main node, label=left:\scriptsize$u_1$] (g) at (-1,0) {};
    \node[main node, label=right:\scriptsize$u_3$] (i) at (1,0) {};
    \node[style={inner sep=2,circle,draw}] (v) at (0,1) {};
    \node[main node, label=above:\scriptsize$u$] (v) at (0,1) {};

    \draw (v) -- node[below left=0.01pt] {\scriptsize $2$} (f);
    \draw (v) -- node[left=0.1pt]{\scriptsize $1$} (g);
    \draw (v) -- node[right=0.1pt]{\scriptsize $3$} (i);
    \draw (g) to [bend right=60] (f) to [bend right=60] (i);
    \draw (g) to [bend right=110, in=290, looseness=2.0] (i);
  \end{tikzpicture}
}

\newcommand{\BlobBase}{%
  \begin{tikzpicture}[main node/.style={inner sep=0,minimum size =.08cm,circle,fill=black!80,draw},node distance = 0.75cm,scale=0.8]
    \node[main node, label=above:\scriptsize$u_2$] (f) at (0,0) {};
    \node[main node, label=above:\scriptsize$u_1$] (g) at (-1,0) {};
    \node[main node, label=above:\scriptsize$u_3$] (i) at (1,0) {};

    \draw (g) to [bend right=60] (f) to [bend right=60] (i);
    \draw (g) to [bend right=110, in=290, looseness=2.0] (i);
  \end{tikzpicture}
}

\newcommand{\BlobTwo}{%
  \node[main node, label={[shift={(0,-1.75em)}]\scriptsize $u_2$}] (f) at (0,-0.30) {};
  \node[main node, label=above:{\scriptsize$u=u_1 = u_3$}] (g) at (0,0.25) {};

  \draw (g) to [out=195, in=90, looseness=1] (-0.75,-0.25) to [out=270, in=225, looseness=1.5] (f);
  \draw (f) to [out=315, in=270, looseness=1.5] (0.75,-0.25) to [out=90, in=345, looseness=1] (g);
  \draw (g) to [out=188, in=180, looseness=2.25] (0,-1.25) to [out=0, in=352, looseness=2.25] (g);
}

\newcommand{\BlobThree}{%
  \node[main node, label=above:{\scriptsize$u_1$}] (i) at (-0.5,0) {};
  \node[main node, label=above:{\scriptsize$u=u_2=u_3$}] (h) at (1,0) {};

  \draw (h) to [out=225, in=320, looseness=25] (h);
  \draw (h) to [bend left=45] (i);
  \draw (h) to [out=330, in=0, looseness=2] (0.25,-1.25);
  \draw (0.25,-1.25) to [out=180, in=230, looseness=1.5] (i);
}

\newcommand{\ExampleBlob}{%
  \begin{tikzpicture}[main node/.style={inner sep=0,minimum size =.08cm,circle,fill=black!80,draw},node distance = 0.75cm,scale=0.8]
    \clip (-1.5,-1.5) rectangle (8,0.75);
    \begin{scope}
      \BlobThree
    \end{scope}
    \node at (2.75,-0.5) {$\bigcap$};
    \begin{scope}[shift={(5,0)}]
      \BlobTwo
    \end{scope}
  \end{tikzpicture}
}

\newcommand{\ShapeOne}{%
  \begin{tikzpicture}[main node/.style={inner sep=0,minimum size =.08cm,circle,fill=black!80,draw},node distance = 0.75cm,scale=0.5, baseline=0]
    \node[main node] (a) at (0,0) {};
    \node[main node] (d) at (2,0) {};
    \node[main node, label=above left:\small$x$] (g) at (-1,1) {};
    \node[main node] (h) at (-1,-1) {};
    \node[main node] (i) at (3,1) {};
    \node[main node] (j) at (3,-1) {};

    \draw (d) -- (a);
    \draw (g) -- (a) -- (h);
    \draw (i) -- (d) -- (j);
  \end{tikzpicture}
}

\newcommand{\ShapeOneCut}{%
  \begin{tikzpicture}[main node/.style={inner sep=0,minimum size =.08cm,circle,fill=black!80,draw},node distance = 0.75cm,scale=0.5, baseline=0]
    \node[main node] (a) at (0,0) {};
    \node[main node, label=above:\small$v$] (d) at (2,0) {};
    \node[style={inner sep=2,rectangle,draw}] (v) at (2,0) {};
    \node[main node, label=above left:\small$x$] (g) at (-1,1) {};
    \node[main node] (h) at (-1,-1) {};
    \node[main node] (i) at (3,1) {};
    \node[main node] (j) at (3,-1) {};

    \node (c1) at (1,1) {};
    \node (c2) at (1,-1) {};

    \draw[dashed] (c1) -- (c2);
    \draw (d) --  (a);
    \draw (g) -- (a) -- (h);
    \draw (i) -- (d) -- (j);
  \end{tikzpicture}
}

\newcommand{\ShapeTwoCut}{%
  \begin{tikzpicture}[main node/.style={inner sep=0,minimum size =.08cm,circle,fill=black!80,draw},node distance = 0.75cm,scale=0.5, baseline=0]
    \node[main node, label=above:\small$v$] (a) at (0,0) {};
    \node[style={inner sep=2,rectangle,draw}] (v) at (0,0) {};
    \node[main node] (d) at (2,0) {};
    \node[main node, label=above left:\small$x$] (g) at (-1,1) {};
    \node[main node] (h) at (-1,-1) {};
    \node[main node] (i) at (3,1) {};
    \node[main node] (j) at (3,-1) {};
    \node (c1) at (1,1) {};
    \node (c2) at (1,-1) {};

    \draw[dashed] (c1) -- (c2);

    \draw (d) -- (a);
    \draw (g) -- (a) -- (h);
    \draw (i) -- (d) -- (j);
  \end{tikzpicture}
}

\newcommand{\ShapeFive}{%
  \begin{tikzpicture}[main node/.style={inner sep=0,minimum size =.08cm,circle,fill=black!80,draw},node distance = 0.75cm,scale=0.5, baseline=0]
    \node[main node, label=above:\small$v$] (a) at (0,0) {};
    \node[main node] (g) at (-1,1) {};
    \node[main node] (h) at (-1,-1) {};
    \node[main node] (b) at (2,0) {};
    \node[main node] (c) at (2,-1.75) {};
    \node[main node, label=above:\small$w$] (d) at (4,0) {};
    \node[main node] (i) at (5,1) {};
    \node[main node] (j) at (5,-1) {};

    \draw (a) -- (b) -- (d);
    \draw (b) -- (c);
    \draw (g) -- (a) -- (h);
    \draw (i) -- (d) -- (j);
  \end{tikzpicture}
}

\newcommand{\ShapeFiveCut}{%
  \begin{tikzpicture}[main node/.style={inner sep=0,minimum size =.08cm,circle,fill=black!80,draw},node distance = 0.75cm,scale=0.5, baseline=0]
    \node[style={inner sep=2,rectangle,draw}] (v) at (0,0) {};
    \node[main node, label=above:\small$v$] (a) at (0,0) {};
    \node[main node] (g) at (-1,1) {};
    \node[main node] (h) at (-1,-1) {};
    \node[main node] (b) at (2,0) {};
    \node[main node, label=below:\small$x'$] (c) at (2,-1.75) {};
    \node[style={inner sep=2,rectangle,draw}] (w) at (4,0) {};
    \node[main node, label=above:\small$w$] (d) at (4,0) {};
    \node[main node] (i) at (5,1) {};
    \node[main node] (j) at (5,-1) {};
    \node (c1) at (1,1) {};
    \node (c2) at (1,-1) {};
    \node (c3) at (3,1) {};
    \node (c4) at (3,-1) {};

    \draw[dashed] (c1) -- (c2);
    \draw[dashed] (c3) -- (c4);

    \draw (a) -- (b) -- (d);
    \draw (b) -- (c);
    \draw (g) -- (a) -- (h);
    \draw (i) -- (d) -- (j);
  \end{tikzpicture}
}

\newcommand{\lebox}[1]{%
  \,\begin{ytableau}#1\end{ytableau}\,
}

\newcommand{\leplus}{\,+\,}

\newcommand{\leconfig}{%
  \begin{ytableau}
    \none & \none & + \\
    \none & \none & \none[\vdots] \\
    + & \none[\cdots] & 0 \\
  \end{ytableau}
}

\newcommand{\smallelbow}{%
  \begin{tikzpicture}[scale=0.13]
    \draw (0, 1) arc[start angle=-90, end angle=0, x radius=1, y radius=-1];
    \draw (2, 1) arc[start angle=90, end angle=180, x radius=1, y radius=-1];
  \end{tikzpicture}
}

\newcommand{\smallcross}{%
  \begin{tikzpicture}[scale=0.13]
    \draw (0, 1) to (2,1);
    \draw (1, 0) to (1,2);
  \end{tikzpicture}
}

\newcommand{\elbow}[2][]{%
  \begin{tikzpicture}[scale=0.32]
    \draw[#1] (0, 1) arc[start angle=-90, end angle=0, x radius=1, y radius=-1];
    \draw[#2] (2, 1) arc[start angle=90, end angle=180, x radius=1, y radius=-1];
  \end{tikzpicture}
}

\newcommand{\cross}[2][]{%
  \begin{tikzpicture}[scale=0.32]
    \draw[#1] (0, 1) to (2,1);
    \draw[#2] (1, 0) to (1,2);
  \end{tikzpicture}
}

\newcommand{\exle}{%
  \begin{tikzpicture}[baseline=(current bounding box.center)]
    \useasboundingbox (0,0) rectangle (5*18pt+5*0.04em, 3*18pt+3*0.04em);
    \draw[dotted] (0,0) rectangle (5*18pt+5*0.04em,3*18pt+3*0.04em);
    \node[inner sep=0] at (3*18pt+3*0.04em-9.19pt,2*18pt+2*0.04em-9.19pt) {%
    \begin{ytableau}
      0 & 0 & + & + & 0 \\
      + & + & + & \none & \none \\
      0 & \none & \none & \none & \none  \\
    \end{ytableau}};
    \draw[decorate,decoration={calligraphic brace,amplitude=6pt,raise=4pt}] (current bounding box.north west) -- (current bounding box.north east) node [midway, above=10pt] {\scriptsize$5=8-3$};
    \draw[decorate,decoration={calligraphic brace,mirror,amplitude=6pt,raise=4pt}] (current bounding box.north west) -- (current bounding box.south west) node [midway, left=10pt] {\scriptsize$3$};
  \end{tikzpicture}
}

\newcommand{\expipe}{%
  \begin{tikzpicture}[baseline=(current bounding box.center)]
    \useasboundingbox (0,0) rectangle (5*18pt+5*0.04em, 3*18pt+3*0.04em);
    \draw[dotted] (0,0) rectangle (5*18pt+5*0.04em,3*18pt+3*0.04em);
    \node[inner sep=0] at (3*18pt+3*0.04em-9.18pt,2*18pt+2*0.04em-9.18pt) {%
    \begin{ytableau}
      \cross{} & \cross{} & \elbow{->--} & \elbow{->--} & \cross{-->-} \\
      \elbow{->--} & \elbow{->--} & \elbow{-->-} & \none & \none \\
      \cross[--<-]{-->-} & \none & \none & \none & \none  \\
    \end{ytableau}};
    \node[inner sep=0] at (-3pt, 9pt+0.04em) {\tiny$7$};
    \node[inner sep=0] at (-3pt, 3*9pt+2*0.04em) {\tiny$4$};
    \node[inner sep=0] at (-3pt, 5*9pt+3*0.04em) {\tiny$1$};
    \node[inner sep=0] at (9pt+0.04em, 3*18pt+3*0.04em+3pt) {\tiny$8$};
    \node[inner sep=0] at (3*9pt+2*0.04em, 3*18pt+3*0.04em+3pt) {\tiny$6$};
    \node[inner sep=0] at (5*9pt+3*0.04em, 3*18pt+3*0.04em+3pt) {\tiny$5$};
    \node[inner sep=0] at (7*9pt+4*0.04em, 3*18pt+3*0.04em+3pt) {\tiny$3$};
    \node[inner sep=0] at (9*9pt+5*0.04em, 3*18pt+3*0.04em+3pt) {\tiny$2$};
    \node[inner sep=0] at (1*18pt+1*0.04em+3pt, 9pt+0.04em) {\tiny$7$};
    \node[inner sep=0] at (3*18pt+3*0.04em+3pt, 3*9pt+2*0.04em) {\tiny$4$};
    \node[inner sep=0] at (5*18pt+5*0.04em+3pt, 5*9pt+3*0.04em) {\tiny$1$};
    \node[inner sep=0] at (3*9pt+2*0.04em, 1*18pt+1*0.04em-4pt) {\tiny$6$};
    \node[inner sep=0] at (5*9pt+3*0.04em, 1*18pt+1*0.04em-4pt) {\tiny$5$};
    \node[inner sep=0] at (7*9pt+4*0.04em, 2*18pt+2*0.04em-4pt) {\tiny$3$};
    \node[inner sep=0] at (9*9pt+5*0.04em, 2*18pt+2*0.04em-4pt) {\tiny$2$};
    \node[inner sep=0] at (9pt+0.04em, -4pt) {\tiny$8$};
  \end{tikzpicture}
}

\newcommand{\explabic}{%
  \begin{tikzpicture}[main node/.style={inner sep=0,minimum size =.04cm,circle,fill=black!80,draw},node distance = 0.75cm,scale=0.1]
    \clip (-15,-35) rectangle (45,2);
    \tikzstyle{black}=[inner sep=0,minimum size=1mm,circle,fill=black!80,draw]
    \tikzstyle{white}=[inner sep=0,minimum size=1mm,circle,fill=white,draw=black]
    \begin{scope}[scale=0.25]
    \draw
    (0, 0)
     .. controls (8.9797, 1.2677) and (19.3328, 1.9016) .. (31.3762, 2.5354)
     .. controls (43.4196, 3.1693) and (57.1533, 3.8032) .. (72.366, 3.6975)
     .. controls (87.5787, 3.5919) and (104.2704, 2.7467) .. (111.1372, -4.5427)
     .. controls (118.0041, -11.8321) and (115.046, -25.5658) .. (111.0316, -32.8552)
     .. controls (107.0171, -40.1446) and (101.9462, -40.9898) .. (95.6076, -41.7293)
     .. controls (89.269, -42.4688) and (81.6626, -43.1026) .. (76.803, -43.2083)
     .. controls (71.9434, -43.3139) and (69.8305, -42.8914) .. (66.9781, -44.5817)
     .. controls (64.1257, -46.272) and (60.5339, -50.0751) .. (58.6323, -55.146)
     .. controls (56.7307, -60.2169) and (56.5194, -66.5556) .. (46.272, -69.9362)
     .. controls (36.0245, -73.3168) and (15.7409, -73.7393) .. (2.5354, -80.9231)
     .. controls (-10.67, -88.1069) and (-16.7974, -97.0519) .. (-26.5166, -95.7769)
     .. controls (-36.2358, -95.5019) and (-49.5469, -75.0071) .. (-50.4977, -56.7307)
     .. controls (-51.4485, -38.4543) and (-40.039, -22.3965) .. (-33.5947, -14.1563)
     .. controls (-27.1504, -6.916) and (-25.6714, -5.4935) .. (-21.1287, -4.3314)
     .. controls (-16.5861, -2.1693) and (-8.9797, -0.9677) .. cycle;
      \node[main node, label=below:{\tiny$8$}] (8) at (-25.5166, -95.7769) {};
      \node[main node, label=below:{\tiny$7$}] (7) at (-5.5354, -85.9231) {};
      \node[white] (7a) at (-5.5354, -65.9231) {};
      \node[main node, label=below:{\tiny$6$}] (6) at (15.7409, -75.9393) {};
      \node[white] (6a) at (12.7409, -45.9393) {};
      \node[main node, label=below:{\tiny$5$}] (5) at (38.0245, -72.0168) {};
      \node[white] (5a) at (36.0245, -42.0168) {};
      \node[main node, label=right:{\tiny$4$}] (4) at (56.7307, -60.2169) {};
      \node[black] (4a) at (50.7307, -35.2169) {};
      \node[main node, label=below:{\tiny$3$}] (3) at (76.803, -43.2083) {};
      \node[white] (3a) at (76.803, -17.2083) {};
      \node[main node, label=below:{\tiny$2$}] (2) at (95.6076, -41.7293) {};
      \node[black] (2a) at (95.6076, -30.7293) {};
      \node[main node, label=right:{\tiny$1$}] (1) at (115.046, -17.5658) {};
      \draw (7) -- (7a);
      \draw (6) -- (6a);
      \draw (5) -- (5a);
      \draw (4) to [bend left=20] (4a);
      \draw (3) -- (3a);
      \draw (2) -- (2a);
      \draw (8) arc[start angle=165.786, end angle=272.6006, x radius=34.9498, y radius=-39.9498] to (6a);
      \draw (6a) to [bend left] (5a);
      \draw (5a) to (4a);
      \draw (4a) to [bend left=60, looseness=1.2] (3a);
      \draw (3a) to (1);
      \draw[->--, gray, line width=0.1mm, shift={(-30.9166, -90.7769)}]
      (0, 0)
       .. controls (1.9016, 32.327) and (11.4095, 40.8841) .. (21.9739, 46.2191)
       .. controls (32.5383, 51.5541) and (44.1591, 53.667) .. (52.8219, 54.5122)
       .. controls (61.4847, 55.3573) and (67.1894, 54.9347) .. (70.7813, 53.2444)
       .. controls (74.3732, 51.5541) and (75.8522, 48.5961) .. (77.9651, 43.9478)
       .. controls (80.078, 39.2995) and (82.8247, 32.9608) .. (82.9981, 30.7915);
      \draw[-->-, gray, line width=0.1mm, shift={(110.046, -19.5658)}]
      (0, 0)
       .. controls (-17.1143, 0) and (-23.136, 0) .. (-26.1468, -3.0108)
       .. controls (-29.1577, -6.0217) and (-29.1577, -12.0434) .. (-29.1577, -22.319);
    \end{scope}
    \begin{scope}[shift={(35,-25)}]
      \draw (0,0) circle[radius=7];
      \node[main node, label=above right:{\tiny$2$}] (2) at (360/8-10:7) {};
      \node[black] (2a) at (3.5, 2.5) {};
      \node[main node, label=above:{\tiny$1$}] (1) at (2*360/8-10:7) {};
      \node[main node, label=above left:{\tiny$8$}] (8) at (3*360/8-10:7) {};
      \node[main node, label=left:{\tiny$7$}] (7) at (4*360/8-10:7) {};
      \node[white] (7a) at (-4.5, 1) {};
      \node[main node, label=below left:{\tiny$6$}] (6) at (5*360/8-10:7) {};
      \node[white] (6a) at (-2, 1) {};
      \node[main node, label=below:{\tiny$5$}] (5) at (6*360/8-10:7) {};
      \node[white] (5a) at (-0.5, -3) {};
      \node[main node, label=below right:{\tiny$4$}] (4) at (7*360/8-10:7) {};
      \node[black] (4a) at (2.5, -2.5) {};
      \node[main node, label=right:{\tiny$3$}] (3) at (0*360/8-10:7) {};
      \node[white] (3a) at (1.5, 1) {};

      \draw (7) to (7a);
      \draw (6) to (6a);
      \draw (5) to (5a);
      \draw (4) to (4a);
      \draw (3) to (3a);
      \draw (2) to (2a);
      \draw (1) to (3a) to (4a) to (5a) to (6a) to (8);
    \end{scope}
  \end{tikzpicture}
}

\newcommand{\exnotation}{%
  \begin{tikzpicture}[baseline=(current bounding box.center)]
    \node[inner sep=0] at (3*18pt+3*0.04em-9.19pt,2*18pt+2*0.04em-9.19pt) {%
    \begin{ytableau}
      0 & 0 & + & + & 0 \\
      + & + & + & \none & \none \\
      0 & \none & \none & \none & \none  \\
    \end{ytableau}};
    \node[inner sep=0] at (-3pt, 9pt+0.04em) {\tiny$7$};
    \node[inner sep=0] at (-3pt, 3*9pt+2*0.04em) {\tiny$4$};
    \node[inner sep=0] at (-3pt, 5*9pt+3*0.04em) {\tiny$1$};
    \node[inner sep=0] at (9pt+0.04em, 3*18pt+3*0.04em+3pt) {\tiny$8$};
    \node[inner sep=0] at (3*9pt+2*0.04em, 3*18pt+3*0.04em+3pt) {\tiny$6$};
    \node[inner sep=0] at (5*9pt+3*0.04em, 3*18pt+3*0.04em+3pt) {\tiny$5$};
    \node[inner sep=0] at (7*9pt+4*0.04em, 3*18pt+3*0.04em+3pt) {\tiny$3$};
    \node[inner sep=0] at (9*9pt+5*0.04em, 3*18pt+3*0.04em+3pt) {\tiny$2$};
    \node[inner sep=0] at (1*18pt+1*0.04em+3pt, 9pt+0.04em) {\tiny$7$};
    \node[inner sep=0] at (3*18pt+3*0.04em+3pt, 3*9pt+2*0.04em) {\tiny$4$};
    \node[inner sep=0] at (5*18pt+5*0.04em+3pt, 5*9pt+3*0.04em) {\tiny$1$};
    \node[inner sep=0] at (3*9pt+2*0.04em, 1*18pt+1*0.04em-4pt) {\tiny$6$};
    \node[inner sep=0] at (5*9pt+3*0.04em, 1*18pt+1*0.04em-4pt) {\tiny$5$};
    \node[inner sep=0] at (7*9pt+4*0.04em, 2*18pt+2*0.04em-4pt) {\tiny$3$};
    \node[inner sep=0] at (9*9pt+5*0.04em, 2*18pt+2*0.04em-4pt) {\tiny$2$};
    \node[inner sep=0] at (9pt+0.04em, -4pt) {\tiny$8$};
    \draw[->] (3,2.35) to (0.25,2.35);
    \draw[->] (-0.4,1.75) to (-0.4,0.25);
  \end{tikzpicture}
}

\newcommand{\hatD}{%
  \begin{tikzpicture}[main node/.style={inner sep=0,minimum size =.04cm,circle,fill=black!80,draw},node distance = 0.75cm,scale=0.7]
    \draw (0,0) -- ++(0,4) -- ++(6,0) -- ++(0,-1) node[midway,right] {\small $b_{1}$} -- ++(-1,0) -- ++(0,-1) node[midway,right] {\small $b_2$} -- ++(-1.5,0) -- ++(0,-1) node[midway,right] {\vdots} -- ++(-1.5,0)
      -- ++(0,-1) node[midway,right] {\small $b_{k}$} -- cycle;
    \node at (2,2.5) {$\hat{D}$};
    \draw[dotted] (0,0) rectangle (6,4);
    \draw[decorate,decoration={calligraphic brace,mirror,amplitude=6pt,raise=5pt}] (current bounding box.south west) -- (6,0) node [midway, below=12pt] {\small$n-k$};
    \draw[decorate,decoration={calligraphic brace,mirror,amplitude=6pt,raise=5pt}] (current bounding box.north west) -- (0,0) node [midway, left=12pt] {\small$k$};
  \end{tikzpicture}
}

\newcommand{\constructD}{%
  \begin{tikzpicture}[main node/.style={inner sep=0,minimum size =.04cm,circle,fill=black!80,draw},node distance = 0.75cm,scale=0.5]
    \begin{scope}
      \draw (0,0) -- ++(0,4) -- ++(6,0) -- ++(0,-1) node[midway,right] {\vdots} -- ++(-1,0) -- ++(0,-1) node[midway,right] {\tiny$b_{j-1}$} -- ++(-2,0) -- ++(0,-1) node[midway,right] {\tiny$b_{j}$} -- ++(-1,0)
      -- ++(0,-1) node[midway,right] {\vdots} -- cycle;
    \node at (2,2.5) {$\hat{D}$};
    \draw[pattern=north east lines, pattern color=gray!50] (3.85,2) rectangle (4.35,4) ;
    \node[below] at (4.1,2) {\tiny $a_n$};
    \end{scope}
    \node at (7,2) {$\longrightarrow$};
    \begin{scope}[shift={(8,0)}]
      \draw (0,0) -- ++(0,4) -- ++(5.5,0) -- ++(0,-1) node[midway,right] {\vdots} -- ++(-1,0) -- ++(0,-1) node[midway,right] {\tiny$b_{j-1}$} -- ++(-1.5,0) -- ++(0,-1) node[midway,right] {\tiny$b_{j}$} -- ++(-1,0)
      -- ++(0,-1) node[midway,right] {\vdots} -- cycle;
      \draw[dashed] (3.85,2) to (3.85,4);
    \end{scope}
    \node at (14.5,2) {$\longrightarrow$};
    \begin{scope}[shift={(15.5,-0.25)}]
      \draw (0,0) -- ++(0,4.5) -- ++(5.5,0) -- ++(0,-1) node[midway,right] {\vdots} -- ++(-1,0) -- ++(0,-1) node[midway,right] {\tiny$b_{j-1}$} -- ++(-0.75,0) -- ++(0,-0.5) node[above=2pt,right] {\tiny$a_n$} -- ++(-0.75,0) -- ++(0,-1) node[midway,right] {\tiny$b_{j}$} -- ++(-1,0)
      -- ++(0,-1) node[midway,right] {\vdots} -- cycle;
      \node at (2,3.5) {$D$};
    \draw[pattern=north east lines, pattern color=gray!50] (0,2) rectangle (3.75,2.5) ;
    \end{scope}
  \end{tikzpicture}
}

\newcommand{\shapeD}{%
  \begin{tikzpicture}[main node/.style={inner sep=0,minimum size =.04cm,circle,fill=black!80,draw},node distance = 0.75cm,scale=0.7]
    \draw (0,0) -- ++(0,5) -- ++(6,0) -- ++(0,-1) node[midway,right] {$1$} -- ++(-1,0) -- ++(0,-1) node[midway,right] {\vdots} -- ++(-1,0) -- ++(0,-1) node[midway,right] {\small$a_n$} -- ++(-1.5,0) -- ++(0,-1) node[midway,right] {\vdots} -- ++(-1,0)
    -- ++(0,-1) node[midway,right] {\small$b_k$} -- cycle;
    \node at (1.5,3.5) {$D$};
    \draw[dotted] (0,0) rectangle (6,5);
    \draw[decorate,decoration={calligraphic brace,mirror,amplitude=6pt,raise=5pt}] (current bounding box.south west) -- (6,0) node [midway, below=12pt] {\small$n-(k+1)$};
    \draw[decorate,decoration={calligraphic brace,mirror,amplitude=6pt,raise=5pt}] (current bounding box.north west) -- (0,0) node [midway, left=12pt] {\small$k+1$};
  \end{tikzpicture}
}

\newcommand{\steptwoii}{%
  \begin{tikzpicture}[main node/.style={inner sep=0,minimum size =.04cm,circle,fill=black!80,draw},node distance = 0.75cm,scale=0.7]
    \draw (-0.15,-0.15) -- ++(0,-4*0.85) -- ++(4,0) -- ++(0,0.85) -- ++(3,0) -- ++(0,0.85) -- ++(1,0) -- ++(0,0.85) -- ++(1,0) -- ++(0,0.85) node[midway,right] {\small$b_u$};
    \node at (-2,-2) {$\hat{D}$};
    \node at (0.25, -0.5) {$0$};
    \node at (0.85, -0.5) {$\rule{15pt}{0.5pt}$};
    \node at (1.45, -0.5) {$0$};
    \node at (1.95, -0.5) {$\leplus$};
    \node at (1.95, -4.0) {$L_u$};

    \node at (0.25, -1.35) {$0$};
    \node at (2.35, -1.35) {$\rule{70pt}{0.5pt}$};
    \node at (4.45, -1.35) {$0$};
    \node[draw,circle,inner sep=1pt] at (5.1, -1.35) {$\leplus$};
    \fill[pattern=north east lines, pattern color=gray!50] (3.85,-2.7) rectangle (4.65,0) ;

    \node at (0.25, -2.2) {$0$};
    \node at (3.00, -2.2) {$\rule{90pt}{0.5pt}$};
    \node at (5.75, -2.2) {$0$};
    \node at (6.25, -2.2) {$\leplus$};

    \fill[pattern=north east lines, pattern color=gray!50] (2.2,-3.55) rectangle (2.65,0) ;
    \node[draw,circle,inner sep=1pt] at (3.1, -3.05) {$\leplus$};
  \end{tikzpicture}
}

\newcommand{\beginrow}{%
  \begin{tikzpicture}[main node/.style={inner sep=0,minimum size =.04cm,circle,fill=black!80,draw},node distance = 0.75cm,scale=0.7]
    \draw (0,0) -- ++(0,1) -- ++(3,0) -- ++(0,1) node[midway,right] {\small $b_{u}$} -- ++(1,0);
    \draw[decorate,decoration={calligraphic brace,mirror,amplitude=6pt,raise=5pt}] (0.25,1) -- (3,1) node [midway, below=12pt] {\small{\text{cols. }$\ell$}};
    \node at (0.25,1.5) {$\leplus$};
    \node at (1.5, 1.45) {$\cdots\cdots$};
    \node at (2.65, 1.5) {$\leplus$};
    \begin{scope}[shift={(8,0)}]
    \draw (0,0) -- ++(0,1) -- ++(3,0) -- ++(0,1) node[midway,right] {\small $b_{u}$} -- ++(1,0);
    \draw[decorate,decoration={calligraphic brace,mirror,amplitude=6pt,raise=5pt}] (1.25,1) -- (3,1) node [midway, below=12pt] {\small{\text{cols. }$\ell$}};
    \node at (0, 1.5) {$\rule{20pt}{0.5pt}$};
    \node at (0.75, 1.5) {$0$};
    \node at (0.75, 0.6) {\small{$L_u$}};
    \node at (1.25, 1.5) {$\leplus$};
    \node at (1.95, 1.45) {$\cdots$};
    \node at (2.65, 1.5) {$\leplus$};
    \end{scope}
  \end{tikzpicture}
}

\newcommand{\lastrow}{%
  \begin{tikzpicture}[main node/.style={inner sep=0,minimum size =.04cm,circle,fill=black!80,draw},node distance = 0.75cm,scale=0.7]
    \draw (-1,1.5) -- ++(0,-1.5) -- ++(4,0) -- ++(0,1) node[midway,right] {\small $b_{k}$} -- ++(1,0);
    \node at (-0.75, 0.5) {$0$};
    \node at (0.0, 0.5) {$\rule{20pt}{0.5pt}$};
    \node at (0.75, 0.5) {$0$};
    \node at (0.75, -0.5) {\small{$L_k$}};
    \node at (1.25, 0.5) {$\leplus$};
    \node at (1.95, 0.45) {$\cdots$};
    \node at (2.65, 0.5) {$\leplus$};
    \begin{scope}[shift={(8,0)}]
    \draw (0,1.5) -- ++(0,-1.5) -- ++(3,0) -- ++(0,1) node[midway,right] {\small $a_n$} -- ++(1,0);
    \node at (0.25,0.5) {$\leplus$};
    \node at (1.5, 0.45) {$\cdots\cdots$};
    \node at (2.65, 0.5) {$\leplus$};
    \end{scope}
  \end{tikzpicture}
}

\newcommand{\anend}{%
  \begin{tikzpicture}[main node/.style={inner sep=0,minimum size =.04cm,circle,fill=black!80,draw},node distance = 0.75cm,scale=0.7]
    \draw (0,1.5) -- ++(0,-2) node[midway,left] {\small $a_{n}$};
    \draw (0,0.9) -- ++ (3,0);
    \draw (0,0.125) -- ++ (3,0);
    \node at (0.25, 0.5) {$\leplus$};
    \node at (0.95, 0.45) {$\cdots$};
    \node at (1.65, 0.5) {$\leplus$};
  \end{tikzpicture}
}

\newcommand{\rowend}{%
  \begin{tikzpicture}[main node/.style={inner sep=0,minimum size =.04cm,circle,fill=black!80,draw},node distance = 0.75cm,scale=0.7]
    \draw (0,1.5) -- ++(0,-2) node[midway,left] {\small $b_u$};
    \draw (0,0.9) -- ++ (3,0);
    \draw (0,0.125) -- ++ (3,0);
    \node at (0.25, 0.5) {$0$};
    \node at (0.95, 0.5) {$\rule{18pt}{0.5pt}$};
    \node at (1.65, 0.5) {$0$};
    \node at (2.3, -0.25) {$L_u$};
  \end{tikzpicture}
}

\newcommand{\leconfiglabels}{%
  \begin{ytableau}
    \none & \none & + & \none[\leftarrow] & \none[\text{row } $i$]\\
    \none & \none & \none[\vdots] \\
    + & \none[\cdots] & 0 & \none[\leftarrow] & \none[\text{row } $j$] \\
    \none[\uparrow] & \none & \none[\uparrow]\\
    \none[\text{col. } $m$] & \none & \none[\text{col. } \ell]
  \end{ytableau}
}

\newcommand{\lediagproof}{%
  \begin{tikzpicture}[main node/.style={inner sep=0,minimum size =.04cm,circle,fill=black!80,draw},node distance = 0.75cm,scale=0.7]
    \draw (-1,0) -- ++(0,-2.5);
    \draw (3,-2.5) -- ++(0,1) -- ++(1,0) -- ++(0,1) node[midway,right] {\small $j$} -- ++(1,0);
    \node at (-0.75, -1) {$0$};
    \node at (0.0, -1) {$\rule{20pt}{0.5pt}$};
    \node at (0.75, -1) {$0$};
    \node at (1.25, -1) {$\leplus$};
    \node at (1.25, -1.65) {$\uparrow$};
    \node at (1.25, -2.25) {\small{col. $L \geq m$}};
    \node at (6.5,-1.25) {$\hat{D}$};
    \begin{scope}[shift={(10,0)}]
    \draw (-1,0) -- ++(0,-2.5);
    \draw (3,-2.5) -- ++(0,1) node[midway,right] {\small $b$}  -- ++(1,0);
    \draw (4, -0.5) -- ++(2,0) node[midway,below] {\small $j = a_n$};
    \node at (3.5, -1) {$\vdots$};
    \node at (-0.75, -2) {$0$};
    \node at (0.0, -2) {$\rule{20pt}{0.5pt}$};
    \node at (0.75, -2) {$0$};
    \node at (1.25, -2) {$\leplus$};
    \node at (1.25, -2.65) {$\uparrow$};
    \node at (1.25, -3.25) {\small{col. $L \geq m$}};
    \end{scope}
  \end{tikzpicture}
}

\newcommand{\algexhatD}{%
  \begin{tikzpicture}[scale=0.8]
    \pgfmathsetmacro{\unit}{0.7}
      \draw (-0.25,0) -- ++(0,4*\unit) -- ++(12*\unit+0.25,0) -- ++(0,-1*\unit) node[midway,right] {\tiny{$1$}} -- ++(-1.05*\unit,0) -- ++(0,-1*\unit) node[midway,right] {\tiny{$i$}} -- ++(-3.35*\unit,0) -- ++(0,-1*\unit) node[midway,right] {\tiny$o$} -- ++(-1.8*\unit,0) -- ++(0,-1*\unit) node[midway,right] {\tiny{$q$}} -- cycle;
      \node at (-1.5, 1.5) {\small{$\hat{D}$}};
      \node at (0, -0.25*\unit) {\tiny$n$};
      \node at (0.75, -0.25*\unit) {\tiny$u$};
      \node at (1.5, -0.25*\unit) {\tiny$t$};
      \node at (2.5, -0.25*\unit) {\tiny$s$};
      \node at (3.25, -0.25*\unit) {\tiny$r$};
      \node at (4.75, 0.75*\unit) {\tiny$p$};
      \node at (6, 1.75*\unit) {\tiny$m$};
      \node at (6.75, 1.75*\unit) {\tiny$\ell$};
      \node[draw,circle,inner sep=1pt] at (0.75, 3.5*\unit) {\small$+$};
      \node at (6, 3.5*\unit) {\small$+$};
      \node[draw,circle,inner sep=1pt] at (6.75, 2.5*\unit) {\small$+$};
      \node[draw,circle,inner sep=1pt] at (2.5, 1.5*\unit) {\small$+$};
      \node at (4.75, 1.5*\unit) {\small$+$};
      \node[draw,circle,inner sep=1pt] at (1.5, 0.5*\unit) {\small$+$};
      \node at (2.5, 0.5*\unit) {\small$+$};
      \node at (3.25, 0.5*\unit) {\small$+$};
  \end{tikzpicture}
}

\newcommand{\algexDshape}{%
  \begin{tikzpicture}[scale=0.8]
    \pgfmathsetmacro{\unit}{0.7}
    \draw (-0.25,0) -- ++(0,5*\unit) -- ++(12*\unit+0.25,0) -- ++(0,-1*\unit) node[midway,right] {\tiny{$1$}} -- ++(-1.05*\unit,0) -- ++(0,-1*\unit) node[midway,right] {\tiny{$i$}} -- ++(-2.15*\unit,0) --++(0,-1*\unit) node[midway, right] {\tiny$m$} -- ++(-1.05*\unit,0) -- ++(0,-1*\unit) node[midway,right] {\tiny$o$} -- ++(-1.8*\unit,0) -- ++(0,-1*\unit) node[midway,right] {\tiny{$q$}} -- cycle;
    \node at (-1.5, 2) {\small{${D}$}};
    \node at (0, -0.25*\unit) {\tiny$n$};
    \node at (0.75, -0.25*\unit) {\tiny$u$};
    \node at (1.5, -0.25*\unit) {\tiny$t$};
    \node at (2.5, -0.25*\unit) {\tiny$s$};
    \node at (3.25, -0.25*\unit) {\tiny$r$};
    \node at (4.75, 0.75*\unit) {\tiny$p$};
    \node at (6.75, 2.75*\unit) {\tiny$\ell$};
  \end{tikzpicture}
}

\newcommand{\exrowone}{%
  \begin{tikzpicture}[baseline=(current bounding box.center),remember picture]
    \pgfmathsetmacro{\unit}{0.7}
    \draw (0,3*\unit) -- ++(0,2*\unit) -- ++(12*\unit,0) -- ++(0,-1*\unit) node[midway,right] {\small{$1$}} -- ++(-1*\unit,0) -- ++(0,-1*\unit) node[midway,right] {\small{$i$}};
    \node at (0.25, 4.5*\unit) {\small$0$};
    \node at (0.5, 4.5*\unit) {$\rule{7pt}{0.5pt}$};
    \node at (0.75, 4.5*\unit) {\small$0$};
    \node at (1, 4.5*\unit) {\small$+$};
    \node at (1, 4*\unit) {\small$\uparrow$};
    \node at (1, 3.5*\unit) {\small$u$};
    \node at (1.25, 4.5*\unit) {\small$0$};
    \node at (3.5, 4.5*\unit) {$\rule{120pt}{0.5pt}$};
    \node at (5.75, 4.5*\unit) {\small$0$};
    \node at (6, 4.5*\unit) {\small$+$};
    \node at (6, 4*\unit) {\small$\uparrow$};
    \node at (6, 3.5*\unit) {\small$m-1$};
    \node at (6.25, 4.45*\unit) {\tiny$\cdots$};
    \node at (6.5, 4.5*\unit) {\small$+$};
    \node at (6.75, 4.5*\unit) {\small$+$};
    \node at (6.75, 4*\unit) {\small$\uparrow$};
    \node at (6.75, 3.5*\unit) {\small$\ell$};
    \node at (7, 4.5*\unit) {\small$0$};
    \node at (7.25, 4.5*\unit) {$\rule{7pt}{0.5pt}$};
    \node at (7.5, 4.5*\unit) {\small$0$};
    \node at (7.75, 4.5*\unit) {\small$+$};
    \node at (8.0, 4.475*\unit) {\tiny$\cdots$};
    \node at (8.25, 4.5*\unit) {\small$+$};
    \draw[decorate,decoration={calligraphic brace,amplitude=4pt,raise=4pt}] (7.7,4.95*\unit) -- (8.3,4.95*\unit) node [midway, above=10pt] {\scriptsize$1.$};
    \draw[decorate,decoration={calligraphic brace,amplitude=6pt,raise=4pt}] (5.9,4.95*\unit) -- (7.6,4.95*\unit) node [midway, above=10pt] {\scriptsize$2.$};
    \draw[decorate,decoration={calligraphic brace,amplitude=6pt,raise=4pt}] (0.2,4.95*\unit) -- (5.8,4.95*\unit) node [midway, above=10pt] {\scriptsize$3.$};
  \end{tikzpicture}
}

\newcommand{\exrowi}{%
  \begin{tikzpicture}[baseline=(current bounding box.center),remember picture]
    \pgfmathsetmacro{\unit}{0.7}
    \draw (0,0) -- ++(0,-1.5*\unit);
    \draw (12*\unit,0) -- ++ (-1*\unit,0) -- ++(0,-1*\unit) node[midway,right] {\small{$i$}} -- ++(-2*\unit,0) -- ++(0,-0.5*\unit);
    \node at (0.25, -0.5*\unit) {\small$0$};
    \node at (3.5, -0.5*\unit) {$\rule{170pt}{0.25pt}$};
    \node at (6.75, -0.5*\unit) {\small$0$};
    \node at (7, -0.5*\unit) {\small$+$};
    \node at (7.25, -0.53*\unit) {\tiny$\cdots$};
    \node at (7.5, -0.5*\unit) {\small$+$};
    \node at (6.75, -1.35*\unit) {\small$\ell$};
    \draw[decorate,decoration={calligraphic brace,amplitude=4pt,raise=4pt}] (6.9,-0.25*\unit) -- (7.6,-0.25*\unit) node [midway, above=10pt] {\scriptsize$1.$};
    \draw[decorate,decoration={calligraphic brace,amplitude=6pt,raise=4pt}] (0.2,-0.25*\unit) -- (6.8,-0.25*\unit) node [midway, above=10pt] {\scriptsize$3.$};
  \end{tikzpicture}
}

\newcommand{\exrowm}{%
  \begin{tikzpicture}[baseline=(current bounding box.center),remember picture]
    \pgfmathsetmacro{\unit}{0.7}
    \draw (0,0) -- ++(0,-2*\unit);
    \draw (10*\unit,0) -- ++ (-0.85*\unit,0) -- ++(0,-1*\unit) node[midway,right] {\small{$m$}} -- ++(-1*\unit,0) -- ++(0,-1*\unit) node[midway,right] {\small{$o$}};
    \node at (0.25, -0.5*\unit) {\small$+$};
    \node at (0.75, -0.55*\unit) {\small$\cdots\cdots$};
    \node at (1.25, -0.5*\unit) {\small$+$};
    \node at (1.5, -0.5*\unit) {\small$+$};
    \node at (1.5, -1*\unit) {\small$\uparrow$};
    \node at (1.5, -1.5*\unit) {\small$t$};
    \node at (1.75, -0.5*\unit) {\small$0$};
    \node at (2, -0.5*\unit) {$\rule{7pt}{0.5pt}$};
    \node at (2.25, -0.5*\unit) {\small$0$};
    \node at (2.5, -0.5*\unit) {\small$+$};
    \node at (2.5, -1*\unit) {\small$\uparrow$};
    \node at (2.5, -1.5*\unit) {\small$s$};
    \node at (2.75, -0.5*\unit) {\small$0$};
    \node at (3.625, -0.5*\unit) {$\rule{40pt}{0.5pt}$};
    \node at (4.5, -0.5*\unit) {\small$0$};
    \node at (4.75, -0.5*\unit) {\small$+$};
    \node at (4.75, -1*\unit) {\small$\uparrow$};
    \node at (4.75, -1.5*\unit) {\small$p$};
    \node at (5.0, -0.5*\unit) {\small$0$};
    \node at (5.25, -0.5*\unit) {$\rule{7pt}{0.5pt}$};
    \node at (5.5, -0.5*\unit) {\small$0$};
    \node at (5.75, -0.5*\unit) {\small$+$};
    \node at (6, -0.53*\unit) {\tiny$\cdots$};
    \node at (6.25, -0.5*\unit) {\small$+$};
    \draw[decorate,decoration={calligraphic brace,amplitude=4pt,raise=4pt}] (5.7,-0.25*\unit) -- (6.3,-0.25*\unit) node [midway, above=10pt] {\scriptsize$1.$};
    \draw[decorate,decoration={calligraphic brace,amplitude=6pt,raise=4pt}] (0.2,-0.25*\unit) -- (5.6,-0.25*\unit) node [midway, above=10pt] {\scriptsize$2.$};
  \end{tikzpicture}
}

\newcommand{\exrowo}{%
  \begin{tikzpicture}[baseline=(current bounding box.center),remember picture]
    \pgfmathsetmacro{\unit}{0.7}
    \draw (0,0) -- ++(0,-2*\unit);
    \draw (9*\unit,0) -- ++ (-0.85*\unit,0) -- ++(0,-1*\unit) node[midway,right] {\small{$o$}} -- ++(-1.85*\unit,0) -- ++(0,-1*\unit) node[midway,right] {\small{$q$}};
    \node at (5, -1.35*\unit) {\small$p$};
    \node at (0.25, -0.5*\unit) {\small$0$};
    \node at (1.375, -0.5*\unit) {$\rule{55pt}{0.5pt}$};
    \node at (2.5, -0.5*\unit) {\small$0$};
    \node at (2.75, -0.5*\unit) {\small$0$};
    \node at (2.75, -1*\unit) {\small$\uparrow$};
    \node at (2.75, -1.5*\unit) {\small$s-1$};
    \node at (3, -0.5*\unit) {$\rule{7pt}{0.5pt}$};
    \node at (3.25, -0.5*\unit) {\small$0$};
    \node at (3.5, -0.5*\unit) {\small$+$};
    \node at (3.5, -1*\unit) {\small$\uparrow$};
    \node at (3.5, -1.5*\unit) {\small$r$};
    \node at (3.75, -0.5*\unit) {\small$0$};
    \node at (4, -0.5*\unit) {$\rule{7pt}{0.5pt}$};
    \node at (4.25, -0.5*\unit) {\small$0$};
    \node at (4.5, -0.5*\unit) {\small$+$};
    \node at (5, -0.55*\unit) {\small$\cdots\cdots$};
    \node at (5.5, -0.5*\unit) {\small$+$};
    \draw[decorate,decoration={calligraphic brace,amplitude=4pt,raise=4pt}] (4.4,-0.25*\unit) -- (5.6,-0.25*\unit) node [midway, above=10pt] {\scriptsize$1.$};
    \draw[decorate,decoration={calligraphic brace,amplitude=6pt,raise=4pt}] (2.7,-0.25*\unit) -- (4.3,-0.25*\unit) node [midway, above=10pt] {\scriptsize$2.$};
    \draw[decorate,decoration={calligraphic brace,amplitude=6pt,raise=4pt}] (0.2,-0.25*\unit) -- (2.55,-0.25*\unit) node [midway, above=10pt] {\scriptsize$3.$};
  \end{tikzpicture}
}

\newcommand{\exrowq}{%
  \begin{tikzpicture}[baseline=(current bounding box.center),remember picture]
    \pgfmathsetmacro{\unit}{0.7}
    \draw (0,0) -- ++(0,-1*\unit) -- ++(6*\unit,0) -- ++(0,1*\unit) node[midway,right] {\small{$q$}} -- ++(1*\unit,0);
    \node at (0.25, -0.5*\unit) {\small$0$};
    \node at (0.875, -0.5*\unit) {$\rule{20pt}{0.5pt}$};
    \node at (1.5, -0.5*\unit) {\small$0$};
    \node at (1.5, -1.35*\unit) {\small$t$};
    \node at (1.75, -0.5*\unit) {\small$+$};
    \node at (2.875, -0.55*\unit) {\small$\cdots\cdots\cdots\cdots\cdots$};
    \node at (4, -0.5*\unit) {\small$+$};
    \draw[decorate,decoration={calligraphic brace,amplitude=4pt,raise=4pt}] (1.7,-0.25*\unit) -- (4.1,-0.25*\unit) node [midway, above=10pt] {\scriptsize$1.$};
    \draw[decorate,decoration={calligraphic brace,amplitude=6pt,raise=4pt}] (0.2,-0.25*\unit) -- (1.55,-0.25*\unit) node [midway, above=10pt] {\scriptsize$3.$};
  \end{tikzpicture}
}

\newcommand{\algexD}{%
  \begin{tikzpicture}
    \pgfmathsetmacro{\unit}{0.7}
    \draw (-0.25,0) -- ++(0,5*\unit) -- ++(12*\unit+0.25,0) -- ++(0,-1*\unit) node[midway,right] {\tiny{$1$}} -- ++(-1.05*\unit,0) -- ++(0,-1*\unit) node[midway,right] {\tiny{$i$}} -- ++(-2.15*\unit,0) --++(0,-1*\unit) node[midway, right] {\tiny$m$} -- ++(-1.05*\unit,0) -- ++(0,-1*\unit) node[midway,right] {\tiny$o$} -- ++(-1.8*\unit,0) -- ++(0,-1*\unit) node[midway,right] {\tiny{$q$}} -- cycle;
    \node at (-1.5, 2) {\small{${D}$}};
    \node at (0, -0.25*\unit) {\tiny$n$};
    \node at (0.75, -0.25*\unit) {\tiny$u$};
    \node at (1.5, -0.25*\unit) {\tiny$t$};
    \node at (2.5, -0.25*\unit) {\tiny$s$};
    \node at (3.25, -0.25*\unit) {\tiny$r$};
    \node at (4.75, 0.75*\unit) {\tiny$p$};
    \node at (6.75, 2.75*\unit) {\tiny$\ell$};
    \begin{scope}[shift={(0,0)}]
      \node at (0, 4.5*\unit) {\small$0$};
      \node at (0.25, 4.5*\unit) {$\rule{7pt}{0.5pt}$};
      \node at (0.5, 4.5*\unit) {\small$0$};
      \node at (0.75, 4.5*\unit) {\small$+$};
      \node at (1, 4.5*\unit) {\small$0$};
      \node at (3.5, 4.5*\unit) {$\rule{130pt}{0.5pt}$};
      \node at (6, 4.5*\unit) {\small$0$};
      \node at (6.25, 4.5*\unit) {\small$+$};
      \node at (6.5, 4.47*\unit) {\tiny$\cdots$};
      \node at (6.75, 4.5*\unit) {\small$+$};
      \node at (7, 4.5*\unit) {\small$0$};
      \node at (7.25, 4.5*\unit) {$\rule{7pt}{0.5pt}$};
      \node at (7.5, 4.5*\unit) {\small$0$};
      \node at (7.75, 4.5*\unit) {\small$+$};
      \node at (8.0, 4.47*\unit) {\tiny$\cdots$};
      \node at (8.25, 4.5*\unit) {\small$+$};
    \end{scope}
    \begin{scope}[shift={(0,4*\unit)}]
      \node at (0, -0.5*\unit) {\small$0$};
      \node at (3.375, -0.5*\unit) {$\rule{180pt}{0.5pt}$};
      \node at (6.75, -0.5*\unit) {\small$0$};
      \node at (7, -0.5*\unit) {\small$+$};
      \node at (7.25, -0.53*\unit) {\tiny$\cdots$};
      \node at (7.5, -0.5*\unit) {\small$+$};
    \end{scope}
    \begin{scope}[shift={(0,3*\unit)}]
      \node at (0, -0.5*\unit) {\small$+$};
      \node at (0.75, -0.55*\unit) {\small$\cdot\cdots\cdots\cdot$};
      \node at (1.5, -0.5*\unit) {\small$+$};
      \node at (1.75, -0.5*\unit) {\small$0$};
      \node at (2, -0.5*\unit) {$\rule{7pt}{0.5pt}$};
      \node at (2.25, -0.5*\unit) {\small$0$};
      \node at (2.5, -0.5*\unit) {\small$+$};
      \node at (2.75, -0.5*\unit) {\small$0$};
      \node at (3.625, -0.5*\unit) {$\rule{40pt}{0.5pt}$};
      \node at (4.5, -0.5*\unit) {\small$0$};
      \node at (4.75, -0.5*\unit) {\small$+$};
      \node at (5, -0.5*\unit) {\small$0$};
      \node at (5.125, -0.5*\unit) {$\rule{2pt}{0.5pt}$};
      \node at (5.25, -0.5*\unit) {\small$0$};
      \node at (5.5, -0.5*\unit) {\small$+$};
      \node at (5.75, -0.53*\unit) {\tiny$\cdots$};
      \node at (6, -0.5*\unit) {\small$+$};
    \end{scope}
    \begin{scope}[shift={(0,2*\unit)}]
      \node at (0, -0.5*\unit) {\small$0$};
      \node at (1.5, -0.5*\unit) {$\rule{75pt}{0.5pt}$};
      \node at (3, -0.5*\unit) {\small$0$};
      \node at (3.25, -0.5*\unit) {\small$+$};
      \node at (3.5, -0.5*\unit) {\small$0$};
      \node at (3.75, -0.5*\unit) {$\rule{7pt}{0.5pt}$};
      \node at (4, -0.5*\unit) {\small$0$};
      \node at (4.25, -0.5*\unit) {\small$+$};
      \node at (4.75, -0.55*\unit) {\small$\cdots\cdots$};
      \node at (5.25, -0.5*\unit) {\small$+$};
    \end{scope}
    \begin{scope}[shift={(0,\unit)}]
      \node at (0, -0.5*\unit) {\small$0$};
      \node at (0.75, -0.5*\unit) {$\rule{30pt}{0.5pt}$};
      \node at (1.5, -0.5*\unit) {\small$0$};
      \node at (1.75, -0.5*\unit) {\small$+$};
      \node at (2.875, -0.55*\unit) {\small$\cdots\cdots\cdots\cdots\cdots$};
      \node at (4, -0.5*\unit) {\small$+$};
    \end{scope}
  \end{tikzpicture}
}

\newcommand{\algexDshapes}{%
  \begin{tikzpicture}[scale=0.85]
    \pgfmathsetmacro{\unit}{0.7}
    \draw (-0.25,0) -- ++(0,5*\unit) -- ++(12*\unit+0.25,0) -- ++(0,-1*\unit) node[midway,right] {\tiny{$1$}} -- ++(-1.05*\unit,0) -- ++(0,-1*\unit) node[midway,right] {\tiny{$i$}} -- ++(-2.15*\unit,0) --++(0,-1*\unit) node[midway, right] {\tiny$m$} -- ++(-1.05*\unit,0) -- ++(0,-1*\unit) node[midway,right] {\tiny$o$} -- ++(-1.8*\unit,0) -- ++(0,-1*\unit) node[midway,right] {\tiny{$q$}} -- cycle;
    \node at (2, -1) {\small{${D}$}};
    \node at (0, -0.25*\unit) {\tiny$n$};
    \node at (0.75, -0.25*\unit) {\tiny$u$};
    \node at (1.5, -0.25*\unit) {\tiny$t$};
    \node at (2.5, -0.25*\unit) {\tiny$s$};
    \node at (3.25, -0.25*\unit) {\tiny$r$};
    \node at (4.75, 0.75*\unit) {\tiny$p$};
    \node at (6.75, 2.75*\unit) {\tiny$\ell$};
    \draw[dotted] (0.625,5*\unit) -- ++(0.25,0) -- ++(0,-2*\unit) -- ++(0.5,0) -- ++(0,-1*\unit) -- ++ (-0.75,0) -- cycle;
    \draw[dotted] (1.375,3*\unit) -- ++(0.25,0) -- ++(0,-2*\unit) -- ++(0.75,0) -- ++(0,-1*\unit) -- ++ (-1,0) -- cycle;
    \draw[dotted] (2.375,3*\unit) -- ++(0.25,0) -- ++(0,-2*\unit) -- ++(0.5,0) -- ++(0,-1*\unit) -- ++ (-0.75,0) -- cycle;
    \draw[dotted] (3.125,2*\unit) -- ++(0.25,0) -- ++(0,-1*\unit) -- ++(0.79,0) -- ++(0,-1*\unit) -- ++ (-1.04,0) -- cycle;
    \draw[dotted] (4.625,3*\unit) -- ++(0.25,0) -- ++(0,-1*\unit) -- ++(0.54,0) -- ++(0,-1*\unit) -- ++ (-0.79,0) -- cycle;
    \draw[dotted] (6.625,5*\unit) -- ++(0.25,0) -- ++(0,-1*\unit) -- ++(0.79,0) -- ++(0,-1*\unit) -- ++ (-1.04,0) -- cycle;
    \begin{scope}[shift={(0,0)}]
      \node at (0.75, 4.5*\unit) {\small$+$};
      \node at (6.2, 4.5*\unit) {\small$+$};
      \node at (6.35, 4.47*\unit) {\tiny$\cdot$};
      \node at (6.5, 4.5*\unit) {\small$+$};
      \node at (6.75, 4.5*\unit) {\small$+$};
      \node at (7.75, 4.5*\unit) {\small$+$};
      \node at (8.0, 4.47*\unit) {\tiny$\cdots$};
      \node at (8.25, 4.5*\unit) {\small$+$};
    \end{scope}
    \begin{scope}[shift={(0,4*\unit)}]
      \node at (0.75, -0.5*\unit) {\small$0$};
      \node at (6.75, -0.5*\unit) {\small$0$};
      \node at (7, -0.5*\unit) {\small$+$};
      \node at (7.25, -0.53*\unit) {\tiny$\cdots$};
      \node at (7.5, -0.5*\unit) {\small$+$};
    \end{scope}
    \begin{scope}[shift={(0,3*\unit)}]
      \node at (0, -0.5*\unit) {\small$+$};
      \node at (0.25, -0.53*\unit) {\tiny$\cdots$};
      \node at (0.5, -0.5*\unit) {\small$+$};
      \node at (0.75, -0.5*\unit) {\small$+$};
      \node at (1, -0.53*\unit) {\tiny$\cdots$};
      \node at (1.25, -0.5*\unit) {\small$+$};
      \node at (1.5, -0.5*\unit) {\small$+$};
      \node at (2.5, -0.5*\unit) {\small$+$};
      \node at (4.75, -0.5*\unit) {\small$+$};
      \node at (5.5, -0.5*\unit) {\small$+$};
      \node at (5.75, -0.53*\unit) {\tiny$\cdots$};
      \node at (6, -0.5*\unit) {\small$+$};
    \end{scope}
    \begin{scope}[shift={(0,2*\unit)}]
      \node at (1.5, -0.5*\unit) {\small$0$};
      \node at (2.5, -0.5*\unit) {\small$0$};
      \node at (3.25, -0.5*\unit) {\small$+$};
      \node at (4.2, -0.5*\unit) {\small$+$};
      \node at (4.35, -0.53*\unit) {\tiny$\cdot$};
      \node at (4.5, -0.5*\unit) {\small$+$};
      \node at (4.75, -0.5*\unit) {\small$+$};
      \node at (5, -0.53*\unit) {\tiny$\cdots$};
      \node at (5.25, -0.5*\unit) {\small$+$};
    \end{scope}
    \begin{scope}[shift={(0,\unit)}]
      \node at (1.5, -0.5*\unit) {\small$0$};
      \node at (1.75, -0.5*\unit) {\small$+$};
      \node at (2, -0.53*\unit) {\tiny$\cdots$};
      \node at (2.25, -0.5*\unit) {\small$+$};
      \node at (2.5, -0.5*\unit) {\small$+$};
      \node at (2.75, -0.53*\unit) {\tiny$\cdots$};
      \node at (3, -0.5*\unit) {\small$+$};
      \node at (3.25, -0.5*\unit) {\small$+$};
      \node at (3.625, -0.55*\unit) {\small$\cdots$};
      \node at (4, -0.5*\unit) {\small$+$};
    \end{scope}
    \begin{scope}[shift={(9.75,0.5)}]
      \draw (-0.25,0) -- ++(0,4*\unit) -- ++(12*\unit+0.25,0) -- ++(0,-1*\unit) node[midway,right] {\tiny{$1$}} -- ++(-1.05*\unit,0) -- ++(0,-1*\unit) node[midway,right] {\tiny{$i$}} -- ++(-3.35*\unit,0) -- ++(0,-1*\unit) node[midway,right] {\tiny$o$} -- ++(-1.8*\unit,0) -- ++(0,-1*\unit) node[midway,right] {\tiny{$q$}} -- cycle;
      \node at (2, -1.5) {\small{$\hat{D}$}};
      \node at (0, -0.25*\unit) {\tiny$n$};
      \node at (0.75, -0.25*\unit) {\tiny$u$};
      \node at (1.5, -0.25*\unit) {\tiny$t$};
      \node at (2.5, -0.25*\unit) {\tiny$s$};
      \node at (3.25, -0.25*\unit) {\tiny$r$};
      \node at (4.75, 0.75*\unit) {\tiny$p$};
      \node at (6, 1.75*\unit) {\tiny$m$};
      \node at (6.75, 1.75*\unit) {\tiny$\ell$};
      \node at (0.75, 3.5*\unit) {\small$+$};
      \node at (6, 3.5*\unit) {\small$+$};
      \node at (6.75, 2.5*\unit) {\small$+$};
      \node at (2.5, 1.5*\unit) {\small$+$};
      \node at (4.75, 1.5*\unit) {\small$+$};
      \node at (1.5, 0.5*\unit) {\small$+$};
      \node at (2.5, 0.5*\unit) {\small$+$};
      \node at (3.25, 0.5*\unit) {\small$+$};
    \end{scope}
  \end{tikzpicture}
}

\newcommand*\circled[1]{%
  \tikz[baseline=(char.base)]{\node[shape=circle,draw,inner sep=1pt] (char) {$#1$};}
}

\newcommand{\lastplustoright}{%
  \begin{tikzpicture}[main node/.style={inner sep=0,minimum size =.04cm,circle,fill=black!80,draw},node distance = 0.75cm,scale=0.7]
    \draw (0,0) -- ++(0,1) -- ++(3,0) -- ++(0,1);
    \draw[decorate,decoration={calligraphic brace,mirror,amplitude=6pt,raise=5pt}] (1.5,1) -- (3,1) node [midway, below=12pt] {\small{last $\leplus$'s here}};
    \fill[pattern=north east lines, pattern color=gray!50] (1.25,1) rectangle (-1,2) ;
    \node[rotate=-45] at (-0.5,1) {$\uparrow$};
    \node at (-1,0.5) {\small all $0$'s};
    \node at (1,2.25) {$\leplus$};
    \node at (1,0.6) {\small$\ell$};
    \node at (5,1.5) {or};
    \begin{scope}[shift={(8,0)}]
      \draw (3,1) -- ++(0,1);
      \draw[decorate,decoration={calligraphic brace,mirror,amplitude=6pt,raise=5pt}] (1.5,1) -- (3,1) node [right,below=12pt] {\small{last $\leplus$'s here}};
      \fill[pattern=north east lines, pattern color=gray!50] (1.25,1) rectangle (-1,2) ;
      \node[rotate=-45] at (-0.5,1) {$\uparrow$};
      \node at (-1,0.5) {\small all $0$'s};
      \node at (1,2.25) {$\leplus$};
      \node at (1,0.6) {$+$};
      \node at (1,0.1) {$\uparrow$};
      \node at (1,-0.5) {\small$\ell$};
    \end{scope}
  \end{tikzpicture}
}

\newcommand{\lshapelabels}{%
  \begin{tikzpicture}[main node/.style={inner sep=0,minimum size =.04cm,circle,fill=black!80,draw},node distance = 0.75cm,scale=0.7]
    \draw (0,0) rectangle (0.75,2.75);
    \draw (0,0) rectangle (2.75,0.75);
    \node at (0.375, 2.375) {$+$};
    \node at (0.375, -0.375) {$\uparrow$};
    \node at (0.375, -1) {\small\text{col. }$\ell$};
    \node at (1.125, 0.375) {$+$};
    \node at (1.75, 0.34) {$\cdots$};
    \node at (2.375, 0.375) {$+$};
    \node at (2.375, -0.375) {$\uparrow$};
    \node at (2.375, -1) {\small\text{col. }$m$};
    \node at (3.25, 0.375) {$\leftarrow$};
    \node at (4.5, 0.375) {\small{row $b_B$}};
    \node at (1.25, 2.375) {$\leftarrow$};
    \node at (2.5, 2.375) {\small{row $b_T$}};
    \draw[decorate,decoration={calligraphic brace,amplitude=6pt,raise=5pt}] (0,0) -- (0,2.75) node [midway,left=12pt,text width=1.69cm,align=center] {\small{$s_{\ell}+t_{\ell}+1$ $\leplus$'s}};
  \end{tikzpicture}
}

\newcommand{\lproofone}{%
  \begin{tikzpicture}[main node/.style={inner sep=0,minimum size =.04cm,circle,fill=black!80,draw},node distance = 0.75cm,scale=0.7]
    \draw (0,0) -- ++(3.5,0);
    \fill[pattern=north east lines, pattern color=gray!50] (1.25,0) rectangle (2.25,2) ;
    \node[rotate=45] at (2,2) {$\longleftarrow$};
    \node at (2.5,2.75) {\small all $0$'s};
    \node at (1,1.5) {$\leplus$};
    \node at (0,1.5) {$\rightarrow$};
    \node at (-1,1.5) {\small row $g$};
    \node at (1,-0.4) {\small$\ell$};
    \node at (2.75,-0.4) {\small$m-1$};
    \node at (-3,1) {\small$\hat{D}$};
    \draw[decorate,decoration={calligraphic brace,amplitude=6pt,raise=5pt}] (2.75,1.75) -- (2.75,0) node [midway,right=12pt,text width=2.5cm,align=center] {\small{$+$'s here in col. $m-1$}};
  \end{tikzpicture}
}

\newcommand{\lprooftwo}{%
  \begin{tikzpicture}[main node/.style={inner sep=0,minimum size =.04cm,circle,fill=black!80,draw},node distance = 0.75cm,scale=0.7]
    \draw (0,0) -- ++(3.5,0) --++(0,1) node[midway,right] {\small$g$};
    \node at (1,0.5) {$\leplus$};
    \node at (1,-0.4) {\small$\ell$};
    \node at (1.35, 0.5) {\small$0$};
    \node at (1.75, 0.5) {$\rule{7pt}{0.5pt}$};
    \node at (2.15, 0.5) {\small$0$};
    \node at (2.5,0.5) {$\leplus$};
    \node at (2.5,-0.4) {\small$m-1$};
    \node at (2.5,-0.9) {$\uparrow$};
    \node at (3,-1.4) {\small either $+$ at $(g,m-1)$ or let $m-1=g$};
    \node at (-1,0.5) {\small$\hat{D}$};
  \end{tikzpicture}
}

\newcommand{\lproofthree}{%
  \begin{tikzpicture}[main node/.style={inner sep=0,minimum size =.04cm,circle,fill=black!80,draw},node distance = 0.75cm,scale=0.7]
    \draw (0,0) -- ++(3.25,0) -- ++(0,0.75);
    \draw (5.5,2) -- ++(0,-1) node[midway,right] {\small$g$};
    \fill[pattern=north east lines, pattern color=gray!50] (1.25,0) rectangle (2.25,2) ;
    \node at (1.75,2) {$\downarrow$};
    \node at (1.75,2.75) {\small all $0$'s};
    \node at (1,1.5) {$\leplus$};
    \node at (1,-0.4) {\small$\ell$};
    \node[rotate=45] at (3,2.25) {$\longleftarrow$};
    \node[text width=3cm,align=center] at (5,3) {\small $+$'s here or {below} in col. $m-1$};
    \node at (2.5,-0.4) {\small$m-1$};
    \node at (4.5,1.5) {$\leplus$};
    \node at (4.5,0.9) {$\uparrow$};
    \node at (4.5,0.3) {\small col. $a_n$};
    \node at (-1,1) {\small$\hat{D}$};
  \end{tikzpicture}
}

\newcommand{\lprooffour}{%
  \begin{tikzpicture}[main node/.style={inner sep=0,minimum size =.04cm,circle,fill=black!80,draw},node distance = 0.75cm,scale=0.7]
    \draw (0,0) -- ++(3.75,0) -- ++(0,0.75) node[midway,right] {\small$b$};
    \fill[pattern=north east lines, pattern color=gray!50] (0.5,0) rectangle (3,1.25) ;
    \draw (2.25,1.25) rectangle (2.75,1.75);
    \node at (0.5,0.5) {$\rightarrow$};
    \node at (-0.4,0.5) {\small all $0$'s};
    \node at (1,1.5) {$\leplus$};
    \node at (4.25,1.45) {$\leftarrow$};
    \node at (5.5,1.5) {\small row $b_B$};
    \node at (1,-0.4) {\small$L$};
    \node at (2.5,-0.4) {\small$m$};
    \node at (3.25,0.5) {$\leplus$};
    \node at (-2,1) {\small$\hat{D}$};
  \end{tikzpicture}
}

\newcommand{\stringplus}{%
  \begin{tikzpicture}[main node/.style={inner sep=0,minimum size =.04cm,circle,fill=black!80,draw},node distance = 0.75cm,scale=0.7]
    \draw (0,0) rectangle (2.75,0.75);
    \node at (0.375, 0.375) {$+$};
    \node at (1.375, 0.32) {$\cdots\cdots$};
    \node at (2.375, 0.375) {$+$};
  \end{tikzpicture}
}

\newcommand{\sectionrow}{%
  \begin{tikzpicture}[main node/.style={inner sep=0,minimum size =.04cm,circle,fill=black!80,draw},node distance = 0.75cm,scale=0.7]
        \draw (0.5,0) -- (0,0) -- ++(0,3.75) --++(5.75,0) -- ++(0,-0.5);
        \draw (1.5,0.5) -- ++(0,0.5) --++(3,0) --++(0,0.75) node[midway,right] {\small$b$} --++(0.5,0);
        \node[rotate=-45] at (1,0.375) {$\vdots$};
        \node[rotate=-45] at (5.25,2.5) {$\vdots$};
        \fill[pattern=north east lines,pattern color=gray!50] (1.5,1) rectangle (4.5,3.75);
        \draw[dotted] (1.5,1) rectangle (4.5,3.75);
        \node[rotate=-90] at (3,1) {$\longleftarrow$};
        \node at (3, 0.25) {\small section $b$};
    \end{tikzpicture}
}

\newcommand{\chain}{%
  \begin{tikzpicture}[main node/.style={inner sep=0,minimum size =.04cm,circle,fill=black!80,draw},node distance = 0.75cm,scale=0.7]
        \draw (0,3.75) --++(4.75,0);
        \draw (0,0) -- ++(0,0.75) --++(4.75,0) --++(0,0.45) node[midway,right] {\tiny$b$} --++(0.75,0);
        \draw[dotted] (0,0.75) rectangle (4.75,3.75);
        \fill[pattern=north east lines,pattern color=gray!50] (0,0.75) rectangle (4.75,3.75);
        \begin{scope}[scale=0.6,shift={(5.67,1.25)}]
          \lshapesmall
        \end{scope}
        \node at (3.07,1.85) {\tiny$\cdots$};
        \begin{scope}[scale=0.6,shift={(2.25,2.75)}]
          \lshapesmall
        \end{scope}
        \begin{scope}[scale=0.6,shift={(0,3.75)}]
          \stringplussmall
        \end{scope}
    \end{tikzpicture}
}

\newcommand{\lshapesmall}{%
  \draw[fill=white] (0,0) rectangle (0.625,2.25);
  \draw[fill=white] (0,0) rectangle (2.25,0.75);
  \draw (0.625,0) --++(0,0.75);
  \node at (0.32, 2) {\tiny$+$};
  \node at (1, 0.375) {\tiny$+$};
  \node at (1.5, 0.33) {\tiny$\cdots$};
  \node at (2, 0.375) {\tiny$+$};
}

\newcommand{\lvertsmall}{%
  \draw[fill=white] (0,0) rectangle (0.625,2.25);
  \draw (0,0.75) -- (0.625,0.75);
  \node at (0.32, 2) {\tiny$+$};
}

\newcommand{\stringplussmall}{%
  \draw[fill=white] (0,0) rectangle (2.25,0.75);
  \node at (0.25, 0.375) {\tiny$+$};
  \node at (1.125, 0.32) {\tiny $\cdots\cdots$};
  \node at (2, 0.375) {\tiny $+$};
}

\newcommand{\sectiononlystring}{%
  \begin{tikzpicture}[main node/.style={inner sep=0,minimum size =.04cm,circle,fill=black!80,draw},node distance = 0.75cm,scale=0.7]
        \draw (0,2.5) --++(2,0);
        \draw (0,0.25) -- ++(0,0.5) --++(2,0) --++(0,0.45) node[midway,right] {\tiny$b$} --++(0.75,0);
        \draw[dotted] (0,0.75) rectangle (2,2.5);
        \fill[pattern=north east lines,pattern color=gray!50] (0,0.75) rectangle (2,2.5);
        \draw[fill=white] (0,0.75) rectangle (2,1.2);
        \node at (0.25, 0.975) {\tiny$+$};
        \node at (1, 0.95) {\tiny $\cdots\cdots\cdots$};
        \node at (1.75, 0.975) {\tiny $+$};
    \end{tikzpicture}
}

\newcommand{\blocks}{%
  \begin{tikzpicture}[main node/.style={inner sep=0,minimum size =.04cm,circle,fill=black!80,draw},node distance = 0.75cm,scale=0.42]
    \draw (0,3) --++(2.25,0);
    \draw (0,-1) -- ++(2.25,0);
    \fill[pattern=north east lines,pattern color=gray!50] (0,-1) rectangle (2.25,3);
    \lshapesmall
    \node at (2,-1.5) {\tiny $m$};
    \node at (0.25,-1.5) {\tiny $\ell$};
    \begin{scope}[shift={(8,0.5)}]
    \draw (0,2.5) --++(2.25,0);
    \draw (0,-1.5) -- ++(2.25,0);
    \fill[pattern=north east lines,pattern color=gray!50] (0,-1.5) rectangle (2.25,2.5);
    \stringplussmall
    \node at (2,-2) {\tiny $m$};
    \node at (0.25,-2) {\tiny $\ell$};
    \end{scope}
  \end{tikzpicture}
}

\newcommand{\rowinbetweenl}{%
  \begin{tikzpicture}[main node/.style={inner sep=0,minimum size =.04cm,circle,fill=black!80,draw},node distance = 0.75cm,scale=0.42]
    \draw (0,-1.75) -- ++(0.625,0) -- ++(0,0.5) -- ++(0.75,0) -- ++(0,0.5) -- ++(0.875,0);
    \lshapesmall
    \node at (2,-1.25) {\tiny $m$};
    \node at (0.25,-2.25) {\tiny $\ell$};
  \end{tikzpicture}
}

\newcommand{\firstl}{%
  \begin{tikzpicture}[main node/.style={inner sep=0,minimum size =.04cm,circle,fill=black!80,draw},node distance = 0.75cm,scale=0.7]
    \draw (0,0) -- ++(0,0.5) --++(2.5,0) --++(0,0.45) node[midway,right] {\tiny$b$} --++(0.75,0);
    \begin{scope}[scale=0.6,shift={(1.92,0.833)}]
      \lshapesmall
      \node at (0.25,-0.5) {\tiny $\ell_1$};
    \end{scope}
  \end{tikzpicture}
}

\newcommand{\glueproofone}{%
  \begin{tikzpicture}[main node/.style={inner sep=0,minimum size =.04cm,circle,fill=black!80,draw},node distance = 0.75cm,scale=0.42]
    \stringplussmall
    \node at (2,-0.5) {\tiny $m$};
    \node at (2.85,0.35) {\tiny $b_B$};
    \draw[thick] (1.75,0.75) --++ (0.51,0) --++(0,-0.75) --++(-0.5,0);
    \node at (4.5,0.5) {\tiny or};
    \begin{scope}[shift={(6,-1)}]
      \lshapesmall
      \node at (2,-0.5) {\tiny $m$};
      \node at (0.25,-0.5) {\tiny $\ell$};
      \node at (2.85,0.35) {\tiny $b_B$};
      \draw[thick] (1.75,0.75) --++ (0.51,0) --++(0,-0.75) --++(-0.5,0);
    \end{scope}
    \node at (9.75,0.5) {\tiny or};
    \begin{scope}[shift={(11,-1)}]
      \lvertsmall
      \node at (0.25,-0.5) {\tiny $\ell=m$};
      \node at (1.35,0.35) {\tiny $b_B$};
      \draw[thick] (0,0.75) --++ (0.635,0) --++(0,-0.75) --++(-0.635,0);
    \end{scope}
    \node at (15,0.5) {$\Longrightarrow$};
    \node at (15,1) {\small $?$};
    \begin{scope}[scale=1.2,shift={(16,-1)}]
      \lshapesmall
      \node at (1,2) {\tiny $i$};
      \node at (0.25,-0.5) {\tiny $m-1$};
      \node at (2.85,0.35) {\tiny $j$};
      \draw[thick] (-0.5,1.625) --++ (0.5,0) --++(0,-0.625) --++(-0.5,0);
      \node at (1.325,1.3125) {\tiny$\longleftarrow$};
      \node at (2,1.3125) {\tiny$b_B$};
    \end{scope}
  \end{tikzpicture}
}

\newcommand{\glueprooftwo}{%
  \begin{tikzpicture}[main node/.style={inner sep=0,minimum size =.04cm,circle,fill=black!80,draw},node distance = 0.75cm,scale=0.6]
    \draw (0.625,2.25) --++ (0,-2.25) --++(-0.625,0) --++(0,2.25);
    \draw[fill=white] (0,0) rectangle (2.25,0.625);
    \draw (0.625,0) --++(0,0.625);
    \node at (0.32, 1.75) {\tiny$+$};
    \node at (0.32, 1.25) {\tiny$0$};
    \node at (0.32, 0.75) {\small$\vert$};
    \node at (0.32, 0.315) {\tiny$0$};
    \node at (1, 0.315) {\tiny$+$};
    \node at (1.5, 0.29) {\tiny$\cdots$};
    \node at (2, 0.315) {\tiny$+$};

    \node at (0.25,-0.5) {\tiny $m-1$};
    \node at (3,0.315) {\tiny $j \geq b_B$};
    \draw[thick] (-0.5,2.06) --++ (0.5,0) --++(0,-0.625) --++(-0.5,0);
    \node at (-0.25, 1.75) {\tiny$0$};
    \node at (1.325,1.75) {\tiny$\longleftarrow$};
    \node at (2,1.75) {\tiny$b_B$};
  \end{tikzpicture}
}

\newcommand{\glueproofthree}{%
  \begin{tikzpicture}[main node/.style={inner sep=0,minimum size =.04cm,circle,fill=black!80,draw},node distance = 0.75cm,scale=0.6]
    \draw (0.625,2.5) --++ (0,-2.5) --++(-0.625,0) --++(0,2.5);
    \draw[fill=white] (0,0) rectangle (2.25,0.625);
    \draw (0.625,0) --++(0,0.625);
    \node at (0.32, 1.75) {\tiny$+$};
    \node at (0.32, 1.25) {\tiny$0$};
    \node at (0.32, 0.75) {\small$\vert$};
    \node at (0.32, 0.315) {\tiny$0$};
    \node at (1, 0.315) {\tiny$+$};
    \node at (1.5, 0.29) {\tiny$\cdots$};
    \node at (2, 0.315) {\tiny$+$};

    \node at (0.25,-0.5) {\tiny $m-1$};
    \node at (3,0.315) {\tiny $j \geq b_B$};
    \draw (-0.5,2.5) --++ (0,-1.065) --++(0.5,0) --++(0,1.065);
    \draw[thick] (-0.5,2.06) --++ (0.5,0) --++(0,-0.625) --++(-0.5,0);
    \node at (-0.25, 1.75) {\tiny$0$};
    \node at (1.325,1.75) {\tiny$\longleftarrow$};
    \node at (2,1.75) {\tiny$b_B$};
  \end{tikzpicture}
}

\newcommand{\chainnostring}{%
  \begin{tikzpicture}[main node/.style={inner sep=0,minimum size =.04cm,circle,fill=black!80,draw},node distance = 0.75cm,scale=0.7]
        \draw (0,3.5) --++(3.5,0);
        \draw (0,0) -- ++(0,0.75) --++(3.5,0) --++(0,0.45) node[midway,right] {\tiny$b$} --++(0.75,0);
        \draw[dotted] (0,0.75) rectangle (3.5,3.5);
        \fill[pattern=north east lines,pattern color=gray!50] (0,0.75) rectangle (3.5,3.5);
        \begin{scope}[scale=0.6,shift={(3.58,1.25)}]
          \lshapesmall
        \end{scope}
        \node at (1.80,1.85) {\tiny$\cdots$};
        \begin{scope}[scale=0.6,shift={(0,2.75)}]
          \lshapesmall
        \end{scope}
    \end{tikzpicture}
}

\newcommand{\glueproof}{%
  \begin{tikzpicture}[main node/.style={inner sep=0,minimum size =.04cm,circle,fill=black!80,draw},node distance = 0.75cm,scale=0.6]
    \node at (0,0) {\chain};
    \node at (5.5,0.5) {\tiny or};
    \node at (10,2.5) {\sectiononlystring};
    \node at (10.5,-2.5) {\chainnostring};
  \end{tikzpicture}
}

\newcommand{\doublel}{%
  \begin{tikzpicture}[main node/.style={inner sep=0,minimum size =.04cm,circle,fill=black!80,draw},node distance = 0.75cm,scale=0.6]
    \draw (0,0) rectangle (0.625,3.5);
    \draw (0,1.5) rectangle (3,2.25);
    \draw (0,0) rectangle (3,0.75);
    \draw (1.25,0) -- (1.25,0.75);
    \node[rotate=-90] at (0.92, -0.1) {\small$\longleftarrow$};
    \node at (0.93, -0.75) {\small could be $\leplus$ or $0$};
    \node at (0.32, 3.125) {\small$+$};
    \node at (1.55, 0.375) {\small$+$};
    \node at (2.15, 0.32) {\small$\cdots$};
    \node at (2.75, 0.375) {\small$+$};
    \node at (0.9, 1.875) {\small$+$};
    \node at (1.825, 1.82) {\small$\cdots\cdots$};
    \node at (2.75, 1.875) {\small$+$};
  \end{tikzpicture}
}

\newcommand{\permproofn}{%
  \begin{tikzpicture}[main node/.style={inner sep=0,minimum size =.04cm,circle,fill=black!80,draw},node distance = 0.75cm,scale=0.6]
    \draw (0,0) --++(0,0.5) node[midway,right] {\tiny $n$} --++(1.25,0);
    \draw (0,0) --++(0,1.25);
    \node at (-0.25,0.25) {\tiny $n$};
  \end{tikzpicture}
}

\newcommand{\permproofncolzero}{%
  \begin{tikzpicture}[main node/.style={inner sep=0,minimum size =.04cm,circle,fill=black!80,draw},node distance = 0.75cm,scale=0.6]
    \draw (1,2) --++(-1,0) --++(0,-2) --++(1,0);
    \draw[pattern=north east lines, pattern color=gray!50] (0,0) rectangle (0.5,2);
    \draw[fill=white] (0.75,1.25) --++(-0.75,0) --++(0,-0.5) node[midway,left] {\tiny $a_n$} --++(0.75,0);
    \draw (0.5,1.25) -- (0.5,0.75);
    \node at (0.25,-0.25) {\tiny $n$};
    \node at (0.25, 1) {\tiny $+$};
  \end{tikzpicture}
}

\newcommand{\permproofncol}{%
  \begin{tikzpicture}[main node/.style={inner sep=0,minimum size =.04cm,circle,fill=black!80,draw},node distance = 0.75cm,scale=0.6]
    \draw (2,3) --++(-2,0) --++(0,-3) --++(2,0);
    \draw[pattern=north east lines, pattern color=gray!50] (0,0) rectangle (0.5,3);
    \begin{scope}[scale=0.8,shift={(0,0.75)}]
      \lshapesmall
      \node at (0.3125, 0.375) {\tiny$+$};
      \node at (-0.5,0.375) {\tiny $a_n$};
    \end{scope}
    \node at (0.25,-0.25) {\tiny $n$};
  \end{tikzpicture}
}

\newcommand{\permprooffixedrow}{%
  \begin{tikzpicture}[main node/.style={inner sep=0,minimum size =.04cm,circle,fill=black!80,draw},node distance = 0.75cm,scale=0.6]
    \draw (0,1.75) --++(2.25,0);
    \draw (0,0) --++(1.5,0) --++(0,0.5) node[midway,right] {\tiny $i$} --++(0.75,0);
    \draw[pattern=north east lines, pattern color=gray!50] (1,0.5) rectangle (1.5,1.75);
    \draw (1,0) --++ (0,0.5);
    \node at (1.25,0.25) {\tiny $+$};
    \node at (1.25,-0.25) {\tiny $i+1$};
    \node at (1.25,2) {\tiny $i+1$};
  \end{tikzpicture}
}

\newcommand{\permprooffixedcolzero}{%
  \begin{tikzpicture}[main node/.style={inner sep=0,minimum size =.04cm,circle,fill=black!80,draw},node distance = 0.75cm,scale=0.6]
    \draw (0,1.75) --++(1.5,0);
    \draw (0,0) --++(1.5,0);
    \draw[pattern=north east lines, pattern color=gray!50] (0.25,0) rectangle (0.75,1.75);
    \draw[pattern=north east lines, pattern color=gray!50] (0.75,0) rectangle (1.25,1.75);
    \draw[fill=white] (0.25,0.625) rectangle (0.75,1.125);
    \draw[fill=white] (0.75,0.625) rectangle (1.25,1.125);
    \node at (0.5,0.875) {\tiny $+$};
    \node at (1,0.875) {\tiny $+$};
    \node at (0.5,-0.25) {\tiny $i+1$};
    \node at (1.05,-0.25) {\tiny $i$};
    \node at (0.5,2) {\tiny $i+1$};
    \node at (1.05,2) {\tiny $i$};
  \end{tikzpicture}
}

\newcommand{\permprooffixedcol}{%
  \begin{tikzpicture}[main node/.style={inner sep=0,minimum size =.04cm,circle,fill=black!80,draw},node distance = 0.75cm,scale=0.6]
    \draw (0,3) --++(1.5,0);
    \draw (0,-0.25) --++(1.5,0);
    \draw[pattern=north east lines, pattern color=gray!50] (0.25,-0.25) rectangle (0.75,3);
    \draw[pattern=north east lines, pattern color=gray!50] (0.75,-0.25) rectangle (1.25,3);
    \draw[fill=white] (0.25,1.4) rectangle (0.75,1.9) node[midway,left=4pt] {\tiny$b_L$};
    \node at (0.5,1.65) {\tiny $+$};
    \begin{scope}[scale=0.8,shift={(0.9375,0.4)}]
      \draw[fill=white] (0,0) rectangle (0.625,2.65);
      \draw[fill=white] (0,0) rectangle (2.25,0.75);
      \draw (0.625,0) --++(0,0.75);
      \node at (1, 0.375) {\tiny$+$};
      \node at (1.5, 0.33) {\tiny$\cdots$};
      \node at (2, 0.375) {\tiny$+$};
      \node at (0.3125, 0.375) {\tiny$0$};
      \node at (0.3125, 1) {\tiny$0$};
      \node at (0.3125, 1.66) {\tiny$+$};
    \end{scope}
    \draw (0.75,1.4) rectangle (1.25,1.9);
    \node at (0.5,-0.5) {\tiny $i+1$};
    \node at (1.05,-0.5) {\tiny $i$};
    \node at (0.5,3.25) {\tiny $i+1$};
    \node at (1.05,3.25) {\tiny $i$};
  \end{tikzpicture}
}

\newcommand{\permproofexpath}{%
  \begin{tikzpicture}[main node/.style={inner sep=0,minimum size =.04cm,circle,fill=black!80,draw},node distance = 0.75cm,scale=0.7]
    \node at (-1,0.75) {\tiny $D$};
    \draw (0,1) --++(0,1) node[midway,left] {\tiny$j$}--++(2.5,0);
    \draw (1.5,0) --++(1,0) node[midway,below] {\tiny $i$};
    \node at (2,0.75) {\tiny$+$};
    \node at (1,0.75) {\tiny$+$};
    \node at (1,1.5) {\tiny$+$};
    \draw[->,decorate,decoration=snake] (4,1) -- node[above] {\tiny pipe dream} node[below=2pt] {\tiny $P(D)$} (6,1);
    \begin{scope}[shift={(8,0)}]
      \draw (0,1) --++(0,1) node[midway,left] {\tiny$j$}--++(2.5,0);
      \draw (1.5,0) --++(1,0) node[midway,below] {\tiny $i$};
      \node at (2,0.325) {\tiny $\uparrow$};
      \node at (2,0.75) {\smallelbow};
      \node at (1.5,0.7) {\tiny $\leftarrow$};
      \node at (1,0.75) {\smallelbow};
      \node at (1,1.125) {\tiny $\uparrow$};
      \node at (1,1.5) {\smallelbow};
      \node at (0.5,1.45) {\tiny $\leftarrow$};
    \end{scope}
    \draw[->,decorate,decoration=snake] (12,1) -- node[above] {\tiny path/pipe $i$} (14,1);
    \begin{scope}[shift={(16,0)}]
      \node[main node, label=below:{\tiny$i$}] (i) at (2,0) {};
      \node[main node, label=left:{\tiny$j$}] (j) at (0,1.75) {};
      \draw[rounded corners,postaction={decorate,decoration={markings,mark=at position 0.1 with {\arrow[scale=0.75]{>}},mark=at position 0.9 with {\arrow[scale=0.75]{>}}}}] (i) --++ (0,0.75) --++(-1,0) --++(0,1) -- (j);
    \end{scope}
  \end{tikzpicture}
}

\newcommand{\permproofgeneralpath}{%
  \begin{tikzpicture}[main node/.style={inner sep=0,minimum size =.04cm,circle,fill=black!80,draw},node distance = 0.75cm,scale=0.7]
    \node at (-1,2) {\small $\hat{D}$};
    \draw[thick] (0,4.5) --++(1,0) node[midway,above] {\tiny$j=\ell_{m+1}$};
    \draw[thick] (4,0) --++(1,0) node[midway,below] {\tiny $i+1=\ell_1$};
    \fill[pattern=north east lines, pattern color=gray!50] (0.375,3.75) rectangle (0.625,4.5);
    \fill[pattern=north east lines, pattern color=gray!50] (0.75,3.625) rectangle (1.25,3.375);
    \node at (0.5,3.5) {\tiny$+$};
    \node at (1.5,3.5) {\tiny$+$};
    \fill[pattern=north east lines, pattern color=gray!50] (1.375,3.25) rectangle (1.625,3);
    \fill[pattern=north east lines, pattern color=gray!50] (1.75,2.875) rectangle (2.25,2.625);
    \node at (1.5,2.75) {\tiny$+$};
    \node at (2.5,2.25) {\tiny$\cdots$};
    \node at (3.5,1.75) {\tiny$+$};
    \fill[pattern=north east lines, pattern color=gray!50] (2.75,1.625) rectangle (3.25,1.875);
    \fill[pattern=north east lines, pattern color=gray!50] (3.375,1.5) rectangle (3.625,1.25);
    \node at (3.5,1) {\tiny$+$};
    \node at (4.5,1) {\tiny$+$};
    \fill[pattern=north east lines, pattern color=gray!50] (3.75,0.875) rectangle (4.25,1.125);
    \fill[pattern=north east lines, pattern color=gray!50] (4.375,0.75) rectangle (4.625,0);
    \node at (3.5,-0.30) {\tiny$\ell_2$};
    \node at (2.5,-0.30) {\tiny$\cdots$};
    \node at (1.5,-0.30) {\tiny$\ell_m$};
    \node at (6.0,1) {\tiny$b_1$};
    \node at (6.0,1.75) {\tiny$b_2$};
    \node[rotate=90] at (6.05,2.25) {\tiny$\cdots$};
    \node at (6.0,2.75) {\tiny$b_{m-1}$};
    \node at (6.0,3.5) {\tiny$b_{m}$};
    \draw[->,decorate,decoration=snake] (8,2) -- node[above] {\tiny path $\hat\pi(i+1)$} (10,2);
    \begin{scope}[shift={(11,0)}]
      \node[main node, label=below:{\tiny$i+1$}] (i) at (3.5,0) {};
      \node[main node, label=above:{\tiny$j$}] (j) at (0,4.5) {};
      \draw[rounded corners,postaction={decorate,decoration={markings,mark=at position 0.1 with {\arrow[scale=0.75]{>}}}}] (i) |- (2.5,1) |- (2,1.75);
      \node[rotate=-45] at (1.75,2.2) {\tiny$\cdots$};
      \draw[rounded corners,postaction={decorate,decoration={markings,mark=at position 0.9 with {\arrow[scale=0.75]{>}}}}] (1.5,2.75) -| (1,3.5) -| (j);
    \end{scope}
  \end{tikzpicture}
}

\newcommand{\permproofpathzeros}{%
  \begin{tikzpicture}[main node/.style={inner sep=0,minimum size =.04cm,circle,fill=black!80,draw},node distance = 0.75cm,scale=0.7]
    \draw[fill=white] (0,1.25) rectangle (0.5,1.75);
    \node at (0.25,1.5) {\tiny$+$};
    \fill[pattern=north east lines, pattern color=gray!50] (0.5,0.5) rectangle (-1,1.25);
    \fill[pattern=north east lines, pattern color=gray!50] (0.5,0.5) rectangle (1.5,2.5);
    \node at (3.5,1.25) {\small$\rightarrow$};
    \begin{scope}[shift={(7,0.5)}]
      \draw[fill=white] (0,0) rectangle (0.5,0.5);
      \draw[fill=white] (0,1.25) rectangle (0.5,1.75);
      \node at (0.25,1.5) {\tiny$+$};
      \node at (0.25,0.25) {\tiny$+$};
      \fill[pattern=north east lines, pattern color=gray!50] (0.5,0.5) rectangle (-1,1.25);
    \end{scope}
    \begin{scope}[shift={(10,1)}]
      \draw[fill=white] (0,0) rectangle (0.5,0.5);
      \draw[fill=white] (1.25,0) rectangle (1.75,0.5);
      \node at (1.5,0.25) {\tiny$+$};
      \node at (0.25,0.25) {\tiny$+$};
      \fill[pattern=north east lines, pattern color=gray!50] (0.5,0) rectangle (1.25,1.5);
    \end{scope}
  \end{tikzpicture}
}

\newcommand{\permproofihat}{%
  \begin{tikzpicture}[main node/.style={inner sep=0,minimum size =.04cm,circle,fill=black!80,draw},node distance = 0.75cm,scale=0.7]
    \node at (-2,2) {\small $\hat{D}$};
    \draw[dashed] (-1,0) node[below=2pt] {\tiny $a_n$} --++ (0,4.5);
    \draw[thick] (0,4.5) --++(1,0) node[midway,above] {\tiny$j=\ell_{m+1}$};
    \draw[thick] (4,0) --++(1,0) node[midway,below] {\tiny $i+1=\ell_1$};
    \node at (0.5,3.5) {\tiny$+$};
    \node[draw,circle,inner sep=1pt] at (1.5,3.5) {\tiny$+$};
    \node at (1.5,2.75) {\tiny$+$};
    \node at (2.5,2.25) {\tiny$\cdots$};
    \node[draw,circle,inner sep=1pt] at (3.5,1.75) {\tiny$+$};
    \node at (3.5,1) {\tiny$+$};
    \node[draw,circle,inner sep=1pt] at (4.5,1) {\tiny$+$};
    \node at (3.5,-0.30) {\tiny$\ell_2$};
    \node at (2.5,-0.30) {\tiny$\cdots$};
    \node at (1.5,-0.30) {\tiny$\ell_m$};
    \node at (6.0,1) {\tiny$b_1$};
    \node at (6.0,1.75) {\tiny$b_2$};
    \node[rotate=90] at (6.05,2.25) {\tiny$\cdots$};
    \node at (6.0,2.75) {\tiny$b_{m-1}$};
    \node at (6.0,3.5) {\tiny$b_{m}$};
    \node at (8,2) {\small or};
    \begin{scope}[shift={(11,0)}]
      \draw[dashed] (-1,0) node[below=2pt] {\tiny $a_n$} --++ (0,4.5);
      \draw[thick] (0,4.5) --++(1,0) node[midway,above] {\tiny$j=\ell_{m+1}$};
      \draw[thick] (5.5,-0.125) --++(0,0.75) node[midway,right] {\tiny $i+1$};
      \node at (0.5,3.5) {\tiny$+$};
      \node[draw,circle,inner sep=1pt] at (1.5,3.5) {\tiny$+$};
      \node at (1.5,2.75) {\tiny$+$};
      \node at (2.5,2.25) {\tiny$\cdots$};
      \node[draw,circle,inner sep=1pt] at (3.5,1.75) {\tiny$+$};
      \node at (3.5,1) {\tiny$+$};
      \node[draw,circle,inner sep=1pt] at (4.5,1) {\tiny$+$};
      \node at (4.5,0.25) {\tiny$+$};
      \node at (4.5,-0.30) {\tiny$\ell_1$};
      \node at (3.5,-0.30) {\tiny$\ell_2$};
      \node at (2.5,-0.30) {\tiny$\cdots$};
      \node at (1.5,-0.30) {\tiny$\ell_m$};
      \node at (6.0,1) {\tiny$b_1$};
      \node at (6.0,1.75) {\tiny$b_2$};
      \node[rotate=90] at (6.05,2.25) {\tiny$\cdots$};
      \node at (6.0,2.75) {\tiny$b_{m-1}$};
      \node at (6.0,3.5) {\tiny$b_{m}$};
    \end{scope}
  \end{tikzpicture}
}

\newcommand{\permprooficols}{%
  \begin{tikzpicture}[main node/.style={inner sep=0,minimum size =.04cm,circle,fill=black!80,draw},node distance = 0.75cm,scale=0.7]
    \node at (2,1) {\small $\hat{D}$};
    \node[draw,circle,inner sep=1pt] at (3.5,2) {\tiny$+$};
    \node at (3.5,1) {\tiny$+$};
    \node[draw,circle,inner sep=1pt] at (4.5,1) {\tiny$+$};
    \node at (4.5,0.25) {\tiny$+$};
    \node at (4.5,-0.75) {\tiny$\ell_s$};
    \node at (3.5,-0.75) {\tiny$\ell_{s+1}$};
    \node at (6.0,0.25) {\tiny$b_{s-1}$};
    \node at (6.0,1) {\tiny$b_s$};
    \node at (6.0,2) {\tiny$b_{s+1}$};
    \node at (8.25,1) {$\longmapsto$};
    \begin{scope}[shift={(8,0)}]
      \node at (2,1) {\small ${D}$};
      \draw (3.25,1) -- (3.25,2.5);
      \draw (3.75,1) -- (3.75,2.5);
      \draw (4.25,0.5) -- (4.25,2.5);
      \draw (4.75,0.5) -- (4.75,2.5);
      \fill[pattern=north east lines, pattern color=gray!50] (3.75,1.75) rectangle (4.25,2.25);
      \fill[pattern=north east lines, pattern color=gray!50] (4.25,1.75) rectangle (4.75,1.25);
      \node[draw,circle,inner sep=1pt] at (3.5,2) {\tiny$+$};
      \node at (4.5,2) {\tiny$+$};
      \node[draw,circle,inner sep=1pt] at (4.5,1) {\tiny$+$};
      \node at (4.5,-0.5) {\tiny$\ell_s$};
      \node at (3.5,-0.5) {\tiny$\ell_{s+1}$};
      \node at (6.0,1) {\tiny$b_s$};
      \node at (6.0,2) {\tiny$b_{s+1}$};
    \end{scope}
  \end{tikzpicture}
}

\newcommand{\permproofi}{%
  \begin{tikzpicture}[main node/.style={inner sep=0,minimum size =.04cm,circle,fill=black!80,draw},node distance = 0.75cm,scale=0.7]
    \node at (-1,2) {\small ${D}$};
    \draw[thick] (0,4.5) --++(1,0) node[midway,above] {\tiny$j=\ell_{m+1}$};
    \draw (0.25,3) --++(0,0.75) --++(0.5,0) --++ (0,-0.75);
    \fill[pattern=north east lines, pattern color=gray!50] (0.25,3.75) rectangle (0.75,4.5);
    \fill[pattern=north east lines, pattern color=gray!50] (0.75,3.75) rectangle (1.25,3.25);
    \node at (0.5,3.5) {\tiny$+$};
    \draw (1.25,4) -- (1.25,2.25);
    \draw (1.75,4) -- (1.75,2.25);
    \fill[pattern=north east lines, pattern color=gray!50] (1.25,3.25) rectangle (1.75,3);
    \fill[pattern=north east lines, pattern color=gray!50] (1.75,3) rectangle (2.25,2.5);
    \node at (1.5,3.5) {\tiny$+$};
    \node[draw,circle,inner sep=1pt] at (1.5,2.75) {\tiny$+$};
    \node at (2.5,2.25) {\tiny$\cdots$};
    \draw (3.25,2.25) -- (3.25,0.5);
    \draw (3.75,2.25) -- (3.75,0.5);
    \fill[pattern=north east lines, pattern color=gray!50] (2.75,1.5) rectangle (3.25,2);
    \fill[pattern=north east lines, pattern color=gray!50] (3.25,1.5) rectangle (3.75,1.25);
    \node at (3.5,1.75) {\tiny$+$};
    \node[draw,circle,inner sep=1pt] at (3.5,1) {\tiny$+$};
    \draw (4.25,1.5) -- (4.25,-0.25);
    \draw (4.75,1.5) -- (4.75,-0.25);
    \fill[pattern=north east lines, pattern color=gray!50] (3.75,0.75) rectangle (4.25,1.25);
    \fill[pattern=north east lines, pattern color=gray!50] (4.25,0.75) rectangle (4.75,0.5);
    \node at (4.5,1) {\tiny$+$};
    \node[draw,circle,inner sep=1pt] at (4.5,0.25) {\tiny$+$};
    \node at (4.5,-0.50) {\tiny$\ell_1$};
    \node at (3.5,-0.50) {\tiny$\ell_2$};
    \node at (2.5,-0.50) {\tiny$\cdots$};
    \node at (1.5,-0.50) {\tiny$\ell_m$};
    \node at (6.0,0.25) {\tiny$b_1$};
    \node at (6.0,1) {\tiny$b_2$};
    \node at (6.0,1.75) {\tiny$b_3$};
    \node[rotate=90] at (6.05,2.25) {\tiny$\cdots$};
    \node at (6.0,2.75) {\tiny$b_{m}$};
    \node at (6.0,3.5) {\tiny$b$};
    \draw[->,decorate,decoration=snake] (8,2) -- node[above] {\tiny portion of path $\pi(i)$} (10,2);
    \begin{scope}[shift={(11,0)}]
      \node[main node, label=below:{\tiny$\ell_1$}] (i) at (3.5,0.25) {};
      \node[circle,inner sep=1pt,draw] at (3.5,0.25) {};
      \node[right=2pt of i] {\tiny$b_1$};
      \node[main node, label=above:{\tiny$j$}] (j) at (0,4.5) {};
      \draw[rounded corners,postaction={decorate,decoration={markings,mark=at position 0.1 with {\arrow[scale=0.75]{>}}}}] (i) |- (2.5,1) |- (2,1.75);
      \node[rotate=-45] at (1.75,2.2) {\tiny$\cdots$};
      \draw[rounded corners,postaction={decorate,decoration={markings,mark=at position 0.9 with {\arrow[scale=0.75]{>}}}}] (1.5,2.75) -| (1,3.5) -| (j);
    \end{scope}
  \end{tikzpicture}
}

\newcommand{\permproofibothcols}{%
  \begin{tikzpicture}[main node/.style={inner sep=0,minimum size =.04cm,circle,fill=black!80,draw},node distance = 0.75cm,scale=0.7]
    \draw (0,0) --++(1.5,0);
    \draw (0.25,1.5) --++(0,-0.875) --++(0.5,0) --++(0,0.875);
    \fill[pattern=north east lines, pattern color=gray!50] (0.75,0) rectangle (1.25,0.625);
    \node[circle,inner sep=1pt,draw] at (0.5,0.875) {\tiny $+$};
    \node at (1,0.875) {\tiny $+$};
    \node at (0,-0.5) {\tiny $\ell_1 = i+1$};
    \node at (1,-0.5) {\tiny $i$};
    \node at (-0.25,0.875) {\tiny $b_1$};
    \draw[->,decorate,decoration={snake,amplitude=2pt}] (2.25,0.5) -- (2.75,0.5);
    \begin{scope}[shift={(3.25,0)}]
      \node[main node, label=below:{\tiny$i$}] (i) at (1,0) {};
      \node[circle,inner sep=1pt,draw] at (0.5,0.875) {};
      \node[main node, label=above:{\tiny$\ell_1$}] (j) at (0.5,0.875) {};
      \draw[rounded corners,postaction={decorate,decoration={markings,mark=at position 0.25 with {\arrow[scale=0.75]{>}}}}] (i) |- (j);
    \end{scope}
  \end{tikzpicture}
}

\newcommand{\permprooficolrow}{%
  \begin{tikzpicture}[main node/.style={inner sep=0,minimum size =.04cm,circle,fill=black!80,draw},node distance = 0.75cm,scale=0.7]
    \draw (0,0) --++(0.75,0) --++(0,0.75) node[midway,right] {\tiny $i = b_1$};
    \node[circle,inner sep=1pt,draw] at (0.375,0.375) {\tiny $+$};
    \node at (0,-0.25) {\tiny $\ell_1 = i+1$};
    \draw[->,decorate,decoration={snake,amplitude=2pt}] (2.25,0.375) -- (2.75,0.375);
    \begin{scope}[shift={(3.25,0)}]
      \node[main node, label=right:{\tiny$i$}] (i) at (1.25,0.375) {};
      \node[circle,inner sep=1pt,draw] at (0.5,0.375) {};
      \node[main node, label=above:{\tiny$\ell_1$}] (j) at (0.5,0.375) {};
      \draw[postaction={decorate,decoration={markings,mark=at position 0.5 with {\arrow[scale=0.75]{>}}}}] (i) -- (j);
    \end{scope}
  \end{tikzpicture}
}

\newcommand{\permproofibothrows}{%
  \begin{tikzpicture}[main node/.style={inner sep=0,minimum size =.04cm,circle,fill=black!80,draw},node distance = 0.75cm,scale=0.7]
    \draw (1.5,0) --++(0,1.5);
    \draw (0.125,1.5) -- (0.125,0.25);
    \draw (0.625,1.5) -- (0.625,0.25);
    \fill[pattern=north east lines, pattern color=gray!50] (0.625,0.625) rectangle (1.5,1.125);
    \node[circle,inner sep=1pt,draw] at (0.375,0.875) {\tiny $+$};
    \node at (2,0.25) {\tiny $i+1$};
    \node at (2,0.875) {\tiny $i = b_1$};
    \node at (0.375,-0.25) {\tiny $\ell_1$};
    \draw[->,decorate,decoration={snake,amplitude=2pt}] (3,0.875) -- (3.5,0.875);
    \begin{scope}[shift={(4,0.5)}]
      \node[main node, label=right:{\tiny$i$}] (i) at (1.25,0.375) {};
      \node[circle,inner sep=1pt,draw] at (0.5,0.375) {};
      \node[main node, label=above:{\tiny$\ell_1$}] (j) at (0.5,0.375) {};
      \draw[postaction={decorate,decoration={markings,mark=at position 0.5 with {\arrow[scale=0.75]{>}}}}] (i) -- (j);
    \end{scope}
  \end{tikzpicture}
}

\newcommand{\permproofirowcol}{%
  \begin{tikzpicture}[main node/.style={inner sep=0,minimum size =.04cm,circle,fill=black!80,draw},node distance = 0.75cm,scale=0.7]
    \draw (1.25,-1.25) --++(0,1) --++(1,0);
    \draw (0.125,1.5) -- (0.125,0.25);
    \draw (0.625,1.5) -- (0.625,0.25);
    \fill[pattern=north east lines, pattern color=gray!50] (0.625,0.625) rectangle (1.5,1.125);
    \fill[pattern=north east lines, pattern color=gray!50] (1.5,0.625) rectangle (2,-0.25);
    \node[circle,inner sep=1pt,draw] at (0.375,0.875) {\tiny $+$};
    \node at (1.75,0.875) {\tiny $+$};
    \node at (1.75,-0.5) {\tiny $i$};
    \node at (1.75,-1) {\tiny $i+1$};
    \node at (0.375,2) {\tiny $\ell_1$};
    \node at (-0.5,0.875) {\tiny $b_1$};
    \draw[->,decorate,decoration={snake,amplitude=2pt}] (3,0.25) -- (3.5,0.25);
    \begin{scope}[shift={(4,0)}]
      \node[main node, label=below:{\tiny$i$}] (i) at (1.75,-0.25) {};
      \node[circle,inner sep=1pt,draw] at (0.375,0.875) {};
      \node[main node, label=above:{\tiny$\ell_1$}] (j) at (0.375,0.875) {};
      \draw[rounded corners,postaction={decorate,decoration={markings,mark=at position 0.25 with {\arrow[scale=0.75]{>}}}}] (i) |- (j);
    \end{scope}
  \end{tikzpicture}
}

\newcommand{\permproofiv}{%
  \begin{tikzpicture}[main node/.style={inner sep=0,minimum size =.04cm,circle,fill=black!80,draw},node distance = 0.75cm,scale=0.7]
    \node at (-2,2.5) {\small $\hat{D}$};
    \draw[dashed] (0,5.25) --++ (5,0)  node[right=4pt] {\tiny $a_n$};
    \draw[thick] (-0.5,3.75) --++(0,1) node[midway,left] {\tiny$j$};
    \node at (0.5,4.25) {\tiny$+$};
    \node at (0.5,3.5) {\tiny$+$};
    \node[draw,circle,inner sep=1pt] at (1.5,3.5) {\tiny$+$};
    \node at (1.5,2.75) {\tiny$+$};
    \node at (2.5,2.25) {\tiny$\cdots$};
    \node[draw,circle,inner sep=1pt] at (3.5,1.75) {\tiny$+$};
    \node at (3.5,1) {\tiny$+$};
    \node[draw,circle,inner sep=1pt] at (4.5,1) {\tiny$+$};
    \node at (4.5,-0.30) {\tiny$\ell_1$};
    \node at (3.5,-0.30) {\tiny$\ell_2$};
    \node at (2.5,-0.30) {\tiny$\cdots$};
    \node at (1.5,-0.30) {\tiny$\ell_m$};
    \node at (0.5,-0.30) {\tiny$\ell_{m+1}$};
    \node at (6.0,1) {\tiny$b_1$};
    \node at (6.0,1.75) {\tiny$b_2$};
    \node[rotate=90] at (6.05,2.25) {\tiny$\cdots$};
    \node at (6.0,2.75) {\tiny$b_{m-1}$};
    \node at (6.0,3.5) {\tiny$b_{m}$};
    \node at (8,2.4) {$\longmapsto$};
    \begin{scope}[shift={(11,0.75)}]
      \node at (-1.5,1.75) {\small ${D}$};
      \draw[dashed] (0,4.5) --++ (5.5,0)  node[right=4pt] {\tiny $a_n$};
      \draw[thick] (-0.5,3) --++(0,1);
      \node at (-0.75,3.5) {\tiny$j$};
      \draw (0.25,3) -- (0.25,4);
      \draw (0.75,3) -- (0.75,4);
      \fill[pattern=north east lines, pattern color=gray!50] (0.25,3.75) rectangle (-0.5,3.25);
      \fill[pattern=north east lines, pattern color=gray!50] (0.75,3.75) rectangle (1.25,3.25);
      \node at (0.5,3.5) {\tiny$0$};
      \draw (1.25,4) -- (1.25,2.25);
      \draw (1.75,4) -- (1.75,2.25);
      \fill[pattern=north east lines, pattern color=gray!50] (1.25,3.25) rectangle (1.75,3);
      \fill[pattern=north east lines, pattern color=gray!50] (1.75,3) rectangle (2.25,2.5);
      \node at (1.5,3.5) {\tiny$+$};
      \node[draw,circle,inner sep=1pt] at (1.5,2.75) {\tiny$+$};
      \node at (2.5,2.25) {\tiny$\cdots$};
      \draw (3.25,2.25) -- (3.25,0.5);
      \draw (3.75,2.25) -- (3.75,0.5);
      \fill[pattern=north east lines, pattern color=gray!50] (2.75,1.5) rectangle (3.25,2);
      \fill[pattern=north east lines, pattern color=gray!50] (3.25,1.5) rectangle (3.75,1.25);
      \node at (3.5,1.75) {\tiny$+$};
      \node[draw,circle,inner sep=1pt] at (3.5,1) {\tiny$+$};
      \draw (4.25,1.5) -- (4.25,-0.25);
      \draw (4.75,1.5) -- (4.75,-0.25);
      \fill[pattern=north east lines, pattern color=gray!50] (3.75,0.75) rectangle (4.25,1.25);
      \fill[pattern=north east lines, pattern color=gray!50] (4.25,0.75) rectangle (4.75,0.5);
      \node at (4.5,1) {\tiny$+$};
      \node[draw,circle,inner sep=1pt] at (4.5,0.25) {\tiny$+$};
      \node at (4.5,-1) {\tiny$\ell_1$};
      \node at (3.5,-1) {\tiny$\ell_2$};
      \node at (2.5,-1) {\tiny$\cdots$};
      \node at (1.5,-1) {\tiny$\ell_m$};
      \node at (0.5,-1) {\tiny$\ell_{m+1}$};
      \node at (6.0,0.25) {\tiny$b_1$};
      \node at (6.0,1) {\tiny$b_2$};
      \node at (6.0,1.75) {\tiny$b_3$};
      \node[rotate=90] at (6.05,2.25) {\tiny$\cdots$};
      \node at (6.0,2.75) {\tiny$b_{m}$};
      \node at (6.0,3.5) {\tiny$b=j$};
    \end{scope}
  \end{tikzpicture}
}

\newcommand{\permproofiihat}{%
  \begin{tikzpicture}[main node/.style={inner sep=0,minimum size =.04cm,circle,fill=black!80,draw},node distance = 0.75cm,scale=0.7]
    \node at (-2,2.5) {\small $\hat{D}$};
    \draw[thick] (0,6.5) --++(1,0) node[midway,above] {\tiny$j=\ell_{m+1}$};
    \node at (0.5,5.25) {\tiny$+$};
    \node at (1.5,5.25) {\tiny$+$};
    \node at (1.5,4.5) {\tiny$+$};
    \node at (2.25,4) {\tiny$\cdots$};
    \node at (3,3.5) {\tiny$+$};
    \node at (3,2.75) {\tiny$+$};
    \node at (3.75,2.25) {\tiny$\cdots$};
    \node[draw,circle,inner sep=1pt] at (4.5,1.75) {\tiny$+$};
    \node at (4.5,1) {\tiny$+$};
    \node[draw,circle,inner sep=1pt] at (5.5,1) {\tiny$+$};
    \node at (5.5,-0.30) {\tiny$\ell_1$};
    \node at (4.5,-0.30) {\tiny$\ell_2$};
    \node at (3.85,-0.30) {\tiny$\cdots$};
    \node at (3,-0.30) {\tiny$a_n = \ell_s$};
    \draw (2.75,0.5) rectangle (3.25,6.5);
    \node[main node] at (3,6.25) {};
    \node[main node] at (3,5.25) {};
    \node[main node] at (3,4.5) {};
    \node at (2.15,-0.30) {\tiny$\cdots$};
    \node at (1.5,-0.30) {\tiny$\ell_m$};
    \node at (7.0,1) {\tiny$b_1$};
    \node at (7.0,1.75) {\tiny$b_2$};
    \node[rotate=90] at (7.05,2.25) {\tiny$\cdots$};
    \node at (7.0,2.75) {\tiny$b_{s-1}$};
    \node at (7.0,3.5) {\tiny$b_{s}$};
    \node[rotate=90] at (7.05,4) {\tiny$\cdots$};
    \node at (7.0,4.5) {\tiny$b_{m-1}$};
    \node at (7.0,5.25) {\tiny$b_{m}$};
    \node[rotate=90] at (7.05,5.75) {\tiny$\cdots$};
    \node at (7.0,6.25) {\tiny$1$};
    \draw[decorate,decoration={calligraphic brace,amplitude=6pt,raise=5pt}] (7.25,6) -- (7.25,2.75) node [midway, right=12pt] {\tiny rows of type (II)};
    \node at (9,2.5) {\small or};
    \begin{scope}[shift={(12,0)}]
      \draw[thick] (-0.5,5.5) --++(0,1) node[midway,left] {\tiny$j$};
      \node at (0.5,6) {\tiny$+$};
      \node at (0.5,5.25) {\tiny$+$};
      \node at (1.5,5.25) {\tiny$+$};
      \node at (1.5,4.5) {\tiny$+$};
      \node at (2.25,4) {\tiny$\cdots$};
      \node at (3,3.5) {\tiny$+$};
      \node at (3,2.75) {\tiny$+$};
      \node at (3.75,2.25) {\tiny$\cdots$};
      \node[draw,circle,inner sep=1pt] at (4.5,1.75) {\tiny$+$};
      \node at (4.5,1) {\tiny$+$};
      \node[draw,circle,inner sep=1pt] at (5.5,1) {\tiny$+$};
      \node at (5.5,-0.30) {\tiny$\ell_1$};
      \node at (4.5,-0.30) {\tiny$\ell_2$};
      \node at (3.85,-0.30) {\tiny$\cdots$};
      \node at (3,-0.30) {\tiny$a_n = \ell_s$};
      \draw (2.75,0.5) rectangle (3.25,6.5);
      \node[main node] at (3,6) {};
      \node[main node] at (3,5.25) {};
      \node[main node] at (3,4.5) {};
      \node at (2.15,-0.30) {\tiny$\cdots$};
      \node at (1.5,-0.30) {\tiny$\ell_m$};
      \node at (0.5,-0.30) {\tiny$\ell_{m+1}$};
      \node at (7.0,1) {\tiny$b_1$};
      \node at (7.0,1.75) {\tiny$b_2$};
      \node[rotate=90] at (7.05,2.25) {\tiny$\cdots$};
      \node at (7.0,2.75) {\tiny$b_{s-1}$};
      \node at (7.0,3.5) {\tiny$b_{s}$};
      \node[rotate=90] at (7.05,4) {\tiny$\cdots$};
      \node at (7.0,4.5) {\tiny$b_{m-1}$};
      \node at (7.0,5.25) {\tiny$b_{m}$};
      \draw[decorate,decoration={calligraphic brace,mirror,amplitude=6pt,raise=4pt}] (-1,6) -- (-1,2.75);
    \end{scope}
  \end{tikzpicture}
}

\newcommand{\permproofii}{%
  \begin{tikzpicture}[main node/.style={inner sep=0,minimum size =.04cm,circle,fill=black!80,draw},node distance = 0.75cm,scale=0.7]
    \node at (-1.5,1) {\small ${D}$};
    \draw[thick] (0,4.5) --++(1,0) node[midway,above] {\tiny$j=\ell_{m+1}$};
    \draw (0.25,3) --++(0,0.75) --++(0.5,0) --++ (0,-0.75);
    \fill[pattern=north east lines, pattern color=gray!50] (0.25,3.75) rectangle (0.75,4.5);
    \fill[pattern=north east lines, pattern color=gray!50] (0.75,3.75) rectangle (1.25,3.25);
    \node at (0.5,3.5) {\tiny$+$};
    \node at (1.25,3) {\tiny$\cdots$};
    \draw (1.75,3) -- (1.75,1.25);
    \draw (2.25,3) -- (2.25,1.25);
    \fill[pattern=north east lines, pattern color=gray!50] (1.75,2.25) rectangle (1.25,2.75);
    \fill[pattern=north east lines, pattern color=gray!50] (1.75,2.25) rectangle (2.25,2);
    \fill[pattern=north east lines, pattern color=gray!50] (2.25,2) rectangle (2.75,1.5);
    \node at (2,2.5) {\tiny$+$};
    \node at (2,1.75) {\tiny$+$};
    \draw (2.75,2.25) -- (2.75,0.5);
    \draw (3.25,2.25) -- (3.25,0.5);
    \fill[pattern=north east lines, pattern color=gray!50] (2.75,1.5) rectangle (3.25,1.25);
    \fill[pattern=north east lines, pattern color=gray!50] (3.25,0.75) rectangle (3.75,1.25);
    \node at (3,1.75) {\tiny$+$};
    \node[draw,circle,inner sep=1pt] at (3,1) {\tiny$+$};
    \node at (3.75,0.5) {\tiny$\cdots$};
    \draw (4.25,0.5) -- (4.25,-1.25);
    \draw (4.75,0.5) -- (4.75,-1.25);
    \fill[pattern=north east lines, pattern color=gray!50] (3.75,-0.25) rectangle (4.25,0.25);
    \fill[pattern=north east lines, pattern color=gray!50] (4.25,-0.25) rectangle (4.75,-0.5);
    \node at (4.5,0) {\tiny$+$};
    \node[draw,circle,inner sep=1pt] at (4.5,-0.75) {\tiny$+$};
    \node at (4.5,-1.750) {\tiny$\ell_1$};
    \node at (3.75,-1.750) {\tiny$\cdots$};
    \node at (3,-1.750) {\tiny$\ell_{s-1}$};
    \node at (2,-1.750) {\tiny$\ell_{s+1}$};
    \node at (1.25,-1.750) {\tiny$\cdots$};
    \node at (0.5,-1.750) {\tiny$\ell_{m+1}$};
    \draw[decorate,decoration={calligraphic brace,mirror,amplitude=6pt,raise=3pt}] (0.25,-2) -- (2.25,-2) node[midway,below=12pt] {\tiny original path};
    \draw[decorate,decoration={calligraphic brace,mirror,amplitude=6pt,raise=3pt}] (2.75,-2) -- (4.75,-2) node[midway,below=12pt] {\tiny same as (i), (iv)};
    \node at (6.0,-0.75) {\tiny$b_1$};
    \node at (6.0,0) {\tiny$b_2$};
    \node[rotate=90] at (6.05,0.5) {\tiny$\cdots$};
    \node at (6.0,1) {\tiny$b_{s-1}$};
    \node at (6.0,1.75) {\tiny$b_{s}$};
    \node at (6.0,2.5) {\tiny$b_{s+1}$};
    \node[rotate=90] at (6.05,3) {\tiny$\cdots$};
    \node at (6.0,3.5) {\tiny$b_m$};
    \node at (8.5,1) {\small or};
    \begin{scope}[shift={(11,0)}]
      \draw[thick] (-0.5,3.75) --++(0,1);
      \node at (-0.75,4.25) {\tiny$j$};
      \fill[pattern=north east lines, pattern color=gray!50] (0.25,4.5) rectangle (-0.5,4.0);
      \node at (0.5,4.25) {\tiny$+$};
      \node at (0.5,3.5) {\tiny$+$};
      \draw (0.25,4.75) -- (0.25,2.75);
      \draw (0.75,4.75) -- (0.75,2.75);
      \fill[pattern=north east lines, pattern color=gray!50] (0.25,4.0) rectangle (0.75,3.75);
      \fill[pattern=north east lines, pattern color=gray!50] (0.75,3.75) rectangle (1.25,3.25);
      \node at (1.25,3) {\tiny$\cdots$};
      \draw (1.75,3) -- (1.75,1.25);
      \draw (2.25,3) -- (2.25,1.25);
      \fill[pattern=north east lines, pattern color=gray!50] (1.75,2.25) rectangle (1.25,2.75);
      \fill[pattern=north east lines, pattern color=gray!50] (1.75,2.25) rectangle (2.25,2);
      \fill[pattern=north east lines, pattern color=gray!50] (2.25,2) rectangle (2.75,1.5);
      \node at (2,2.5) {\tiny$+$};
      \node at (2,1.75) {\tiny$+$};
      \draw (2.75,2.25) -- (2.75,0.5);
      \draw (3.25,2.25) -- (3.25,0.5);
      \fill[pattern=north east lines, pattern color=gray!50] (2.75,1.5) rectangle (3.25,1.25);
      \fill[pattern=north east lines, pattern color=gray!50] (3.25,0.75) rectangle (3.75,1.25);
      \node at (3,1.75) {\tiny$+$};
      \node[draw,circle,inner sep=1pt] at (3,1) {\tiny$+$};
      \node at (3.75,0.5) {\tiny$\cdots$};
      \draw (4.25,0.5) -- (4.25,-1.25);
      \draw (4.75,0.5) -- (4.75,-1.25);
      \fill[pattern=north east lines, pattern color=gray!50] (3.75,-0.25) rectangle (4.25,0.25);
      \fill[pattern=north east lines, pattern color=gray!50] (4.25,-0.25) rectangle (4.75,-0.5);
      \node at (4.5,0) {\tiny$+$};
      \node[draw,circle,inner sep=1pt] at (4.5,-0.75) {\tiny$+$};
      \node at (4.5,-1.750) {\tiny$\ell_1$};
      \node at (3.75,-1.750) {\tiny$\cdots$};
      \node at (3,-1.750) {\tiny$\ell_{s-1}$};
      \node at (2,-1.750) {\tiny$\ell_{s+1}$};
      \node at (1.25,-1.750) {\tiny$\cdots$};
      \node at (0.5,-1.750) {\tiny$\ell_{m+1}$};
      \draw[decorate,decoration={calligraphic brace,mirror,amplitude=6pt,raise=3pt}] (0.25,-2) -- (2.25,-2) node[midway,below=12pt] {\tiny original path};
      \draw[decorate,decoration={calligraphic brace,mirror,amplitude=6pt,raise=3pt}] (2.75,-2) -- (4.75,-2) node[midway,below=12pt] {\tiny same as (i), (iv)};
      \node at (6.0,-0.75) {\tiny$b_1$};
      \node at (6.0,0) {\tiny$b_2$};
      \node[rotate=90] at (6.05,0.5) {\tiny$\cdots$};
      \node at (6.0,1) {\tiny$b_{s-1}$};
      \node at (6.0,1.75) {\tiny$b_{s}$};
      \node at (6.0,2.5) {\tiny$b_{s+1}$};
      \node[rotate=90] at (6.05,3) {\tiny$\cdots$};
      \node at (6.0,3.5) {\tiny$b_m$};
    \end{scope}
  \end{tikzpicture}
}

\newcommand{\permproofiiihat}{%
  \begin{tikzpicture}[main node/.style={inner sep=0,minimum size =.04cm,circle,fill=black!80,draw},node distance = 0.75cm,scale=0.7]
    \node at (-2,3.125) {\small $\hat{D}$};
    \draw[thick] (0,6.5) --++(1,0) node[midway,above] {\tiny$j=\ell_{m+1}$};
    \node at (0.5,5.25) {\tiny$+$};
    \node at (1.5,5.25) {\tiny$+$};
    \node at (1.5,4.5) {\tiny$+$};
    \node at (2.25,4.062) {\tiny$\cdots$};
    \node at (3,3.625) {\tiny$+$};
    \draw[dashed] (0,3.125) --++ (6.5,0)  node[right=3pt] {\tiny $a_n$};
    \node at (3,2.625) {\tiny$+$};
    \node at (3.75,2.187) {\tiny$\cdots$};
    \node[draw,circle,inner sep=1pt] at (4.5,1.75) {\tiny$+$};
    \node at (4.5,1) {\tiny$+$};
    \node[draw,circle,inner sep=1pt] at (5.5,1) {\tiny$+$};
    \node at (5.5,-0.30) {\tiny$\ell_1$};
    \node at (4.5,-0.30) {\tiny$\ell_2$};
    \node at (3.85,-0.30) {\tiny$\cdots$};
    \node at (3,-0.30) {\tiny$\ell_s$};
    \node at (2.25,-0.30) {\tiny$\cdots$};
    \node at (1.5,-0.30) {\tiny$\ell_m$};
    \node at (7.0,1) {\tiny$b_1$};
    \node at (7.0,1.75) {\tiny$b_2$};
    \node[rotate=90] at (7.05,2.187) {\tiny$\cdots$};
    \node at (7.0,2.625) {\tiny$b_{s-1}$};
    \node at (7.0,3.625) {\tiny$b_{s}$};
    \node[rotate=90] at (7.05,4.062) {\tiny$\cdots$};
    \node at (7.0,4.5) {\tiny$b_{m-1}$};
    \node at (7.0,5.25) {\tiny$b_{m}$};
    \draw[decorate,decoration={calligraphic brace,amplitude=6pt,raise=5pt}] (7.25,6) -- (7.25,3.5) node [midway, right=12pt] {\tiny rows of type (II)};
    \node at (9.25,3.125) {\small or};
    \begin{scope}[shift={(12,0)}]
      \draw[thick] (-0.5,5.5) --++(0,1) node[midway,left] {\tiny$j$};
      \node at (0.5,6) {\tiny$+$};
      \node at (0.5,5.25) {\tiny$+$};
      \node at (1.5,5.25) {\tiny$+$};
      \node at (1.5,4.5) {\tiny$+$};
      \node at (2.25,4.062) {\tiny$\cdots$};
      \node at (3,3.625) {\tiny$+$};
      \draw[dashed] (0,3.125) --++ (6.5,0)  node[right=3pt] {\tiny $a_n$};
      \node at (3,2.625) {\tiny$+$};
      \node at (3.75,2.187) {\tiny$\cdots$};
      \node[draw,circle,inner sep=1pt] at (4.5,1.75) {\tiny$+$};
      \node at (4.5,1) {\tiny$+$};
      \node[draw,circle,inner sep=1pt] at (5.5,1) {\tiny$+$};
      \node at (5.5,-0.30) {\tiny$\ell_1$};
      \node at (4.5,-0.30) {\tiny$\ell_2$};
      \node at (3.85,-0.30) {\tiny$\cdots$};
      \node at (3,-0.30) {\tiny$\ell_s$};
      \node at (2.25,-0.30) {\tiny$\cdots$};
      \node at (1.5,-0.30) {\tiny$\ell_m$};
      \node at (0.5,-0.30) {\tiny$\ell_{m+1}$};
      \node at (7.0,1) {\tiny$b_1$};
      \node at (7.0,1.75) {\tiny$b_2$};
      \node[rotate=90] at (7.05,2.187) {\tiny$\cdots$};
      \node at (7.0,2.625) {\tiny$b_{s-1}$};
      \node at (7.0,3.625) {\tiny$b_{s}$};
      \node[rotate=90] at (7.05,4.062) {\tiny$\cdots$};
      \node at (7.0,4.5) {\tiny$b_{m-1}$};
      \node at (7.0,5.25) {\tiny$b_{m}$};
      \draw[decorate,decoration={calligraphic brace,mirror,amplitude=6pt,raise=4pt}] (-1,6) -- (-1,3.5);
    \end{scope}
  \end{tikzpicture}
}

\newcommand{\permproofiii}{%
  \begin{tikzpicture}[main node/.style={inner sep=0,minimum size =.04cm,circle,fill=black!80,draw},node distance = 0.75cm,scale=0.7]
    \node at (-2.5,1.75) {\small ${D}$};
    \draw[thick] (-1,5.25) --++(1,0) node[midway,above] {\tiny$j=\ell_{m+1}$};
    \draw (-0.75,3.75) --++(0,0.75) --++(0.5,0) --++ (0,-0.75);
    \fill[pattern=north east lines, pattern color=gray!50] (-0.75,4.5) rectangle (-0.25,5.25);
    \fill[pattern=north east lines, pattern color=gray!50] (-0.25,4.5) rectangle (0.25,4);
    \node at (-0.5,4.25) {\tiny$+$};
    \node at (0.25,3.75) {\tiny$\cdots$};
    \draw (0.75,3.75) -- (0.75,2);
    \draw (1.25,3.75) -- (1.25,2);
    \fill[pattern=north east lines, pattern color=gray!50] (0.75,3) rectangle (0.25,3.5);
    \fill[pattern=north east lines, pattern color=gray!50] (0.75,3) rectangle (1.25,2.75);
    \fill[pattern=north east lines, pattern color=gray!50] (1.25,2.75) rectangle (1.75,2.25);
    \node at (1,3.25) {\tiny$+$};
    \node at (1,2.5) {\tiny$+$};
    \draw (1.75,3) -- (1.75,1.25);
    \draw (2.25,3) -- (2.25,1.25);
    \fill[pattern=north east lines, pattern color=gray!50] (1.75,2.25) rectangle (2.25,2);
    \fill[pattern=north east lines, pattern color=gray!50] (2.25,2) rectangle (2.75,1.5);
    \node at (2,2.5) {\tiny$+$};
    \node at (2,1.75) {\tiny$+$};
    \draw (2.75,2.25) -- (2.75,0.5);
    \draw (3.25,2.25) -- (3.25,0.5);
    \fill[pattern=north east lines, pattern color=gray!50] (2.75,1.5) rectangle (3.25,1.25);
    \fill[pattern=north east lines, pattern color=gray!50] (3.25,0.75) rectangle (3.75,1.25);
    \node at (3,1.75) {\tiny$+$};
    \node[draw,circle,inner sep=1pt] at (3,1) {\tiny$+$};
    \node at (3.75,0.5) {\tiny$\cdots$};
    \draw (4.25,0.5) -- (4.25,-1.25);
    \draw (4.75,0.5) -- (4.75,-1.25);
    \fill[pattern=north east lines, pattern color=gray!50] (3.75,-0.25) rectangle (4.25,0.25);
    \fill[pattern=north east lines, pattern color=gray!50] (4.25,-0.25) rectangle (4.75,-0.5);
    \node at (4.5,0) {\tiny$+$};
    \node[draw,circle,inner sep=1pt] at (4.5,-0.75) {\tiny$+$};
    \node at (4.5,-1.750) {\tiny$\ell_1$};
    \node at (3.75,-1.750) {\tiny$\cdots$};
    \node at (3,-1.750) {\tiny$\ell_{s-1}$};
    \node at (2,-1.750) {\tiny$\ell_{s}$};
    \node at (1,-1.750) {\tiny$\ell_{s+1}$};
    \node at (0.25,-1.750) {\tiny$\cdots$};
    \node at (-0.5,-1.750) {\tiny$\ell_{m+1}$};
    \draw[decorate,decoration={calligraphic brace,amplitude=6pt,raise=3pt}] (6.25,5) -- (6.25,2.5) node[midway,right=12pt] {\tiny original path};
    \draw[decorate,decoration={calligraphic brace,amplitude=6pt,raise=3pt}] (6.25,1) -- (6.25,-0.75) node[midway,right=12pt] {\tiny same as (i), (iv)};
    \node at (6.0,-0.75) {\tiny$b_1$};
    \node at (6.0,0) {\tiny$b_2$};
    \node[rotate=90] at (6.05,0.5) {\tiny$\cdots$};
    \node at (6.0,1) {\tiny$b_{s-1}$};
    \node at (6.0,1.75) {\tiny$a_{n}$};
    \node at (6.0,2.5) {\tiny$b_{s}$};
    \node at (6.0,3.25) {\tiny$b_{s+1}$};
    \node[rotate=90] at (6.05,3.75) {\tiny$\cdots$};
    \node at (6.0,4.25) {\tiny$b_m$};
    \node at (8,1.75) {\small or};
    \begin{scope}[shift={(12,0)}]
      \draw[thick] (-1.5,4.5) --++(0,1);
      \node at (-1.75,5) {\tiny$j$};
      \fill[pattern=north east lines, pattern color=gray!50] (-0.75,5.25) rectangle (-1.5,4.75);
      \node at (-0.5,5) {\tiny$+$};
      \node at (-0.5,4.25) {\tiny$+$};
      \draw (-0.75,5.5) -- (-0.75,3.5);
      \draw (-0.25,5.5) -- (-0.25,3.5);
      \fill[pattern=north east lines, pattern color=gray!50] (-0.75,4.5) rectangle (-0.25,4.75);
      \fill[pattern=north east lines, pattern color=gray!50] (-0.25,4.5) rectangle (0.25,4);
      \node at (0.25,3.75) {\tiny$\cdots$};
      \draw (0.75,3.75) -- (0.75,2);
      \draw (1.25,3.75) -- (1.25,2);
      \fill[pattern=north east lines, pattern color=gray!50] (0.75,3) rectangle (0.25,3.5);
      \fill[pattern=north east lines, pattern color=gray!50] (0.75,3) rectangle (1.25,2.75);
      \fill[pattern=north east lines, pattern color=gray!50] (1.25,2.75) rectangle (1.75,2.25);
      \node at (1,3.25) {\tiny$+$};
      \node at (1,2.5) {\tiny$+$};
      \draw (1.75,3) -- (1.75,1.25);
      \draw (2.25,3) -- (2.25,1.25);
      \fill[pattern=north east lines, pattern color=gray!50] (1.75,2.25) rectangle (2.25,2);
      \fill[pattern=north east lines, pattern color=gray!50] (2.25,2) rectangle (2.75,1.5);
      \node at (2,2.5) {\tiny$+$};
      \node at (2,1.75) {\tiny$+$};
      \draw (2.75,2.25) -- (2.75,0.5);
      \draw (3.25,2.25) -- (3.25,0.5);
      \fill[pattern=north east lines, pattern color=gray!50] (2.75,1.5) rectangle (3.25,1.25);
      \fill[pattern=north east lines, pattern color=gray!50] (3.25,0.75) rectangle (3.75,1.25);
      \node at (3,1.75) {\tiny$+$};
      \node[draw,circle,inner sep=1pt] at (3,1) {\tiny$+$};
      \node at (3.75,0.5) {\tiny$\cdots$};
      \draw (4.25,0.5) -- (4.25,-1.25);
      \draw (4.75,0.5) -- (4.75,-1.25);
      \fill[pattern=north east lines, pattern color=gray!50] (3.75,-0.25) rectangle (4.25,0.25);
      \fill[pattern=north east lines, pattern color=gray!50] (4.25,-0.25) rectangle (4.75,-0.5);
      \node at (4.5,0) {\tiny$+$};
      \node[draw,circle,inner sep=1pt] at (4.5,-0.75) {\tiny$+$};
      \node at (4.5,-1.750) {\tiny$\ell_1$};
      \node at (3.75,-1.750) {\tiny$\cdots$};
      \node at (3,-1.750) {\tiny$\ell_{s-1}$};
      \node at (2,-1.750) {\tiny$\ell_{s}$};
      \node at (1,-1.750) {\tiny$\ell_{s+1}$};
      \node at (0.25,-1.750) {\tiny$\cdots$};
      \node at (-0.5,-1.750) {\tiny$\ell_{m+1}$};
      \draw[decorate,decoration={calligraphic brace,mirror,amplitude=6pt,raise=3pt}] (-2,5) -- (-2,2.5);
      \draw[decorate,decoration={calligraphic brace,mirror,amplitude=6pt,raise=3pt}] (-2,1) -- (-2,-0.75);
      \node at (-2.25,1.75) {\tiny$a_{n}$};
      \node at (6.0,-0.75) {\tiny$b_1$};
      \node at (6.0,0) {\tiny$b_2$};
      \node[rotate=90] at (6.05,0.5) {\tiny$\cdots$};
      \node at (6.0,1) {\tiny$b_{s-1}$};
      \node at (6.0,1.75) {\tiny$a_{n}$};
      \node at (6.0,2.5) {\tiny$b_{s}$};
      \node at (6.0,3.25) {\tiny$b_{s+1}$};
      \node[rotate=90] at (6.05,3.75) {\tiny$\cdots$};
      \node at (6.0,4.25) {\tiny$b_m$};
    \end{scope}
  \end{tikzpicture}
}

\newcommand{\permproofiiispecialone}{%
  \begin{tikzpicture}[main node/.style={inner sep=0,minimum size =.04cm,circle,fill=black!80,draw},node distance = 0.75cm,scale=0.7]
    \node at (2,1) {\small $\hat{D}$};
    \draw[thick] (3,2.5) --++(1,0) node[midway,above] {\tiny$j=\ell_{m+1}$};
    \draw[dashed] (3,1.75) --++ (2.5,0) node[right=3pt] {\tiny $a_n$};
    \node at (3.5,1) {\tiny$+$};
    \node[draw,circle,inner sep=1pt] at (4.5,1) {\tiny$+$};
    \node at (4.5,0.25) {\tiny$+$};
    \node at (4.5,-0.75) {\tiny$\ell_m$};
    \node at (6.0,0.25) {\tiny$b_{m-1}$};
    \node at (6.0,1) {\tiny$b_m$};
    \node at (8.25,1) {$\longmapsto$};
    \begin{scope}[shift={(8,0)}]
      \node at (2,1) {\small ${D}$};
      \draw[thick] (3,2.5) --++(1,0) node[midway,above] {\tiny$j=\ell_{m+1}$};
      \draw (3.25,1) -- (3.25,2) --++ (0.5,0) -- (3.75,1);
      \draw (4.25,0.5) -- (4.25,2.15);
      \draw (4.75,0.5) -- (4.75,2.15);
      \fill[pattern=north east lines, pattern color=gray!50] (3.25,2) rectangle (3.75,2.5);
      \fill[pattern=north east lines, pattern color=gray!50] (3.75,1.5) rectangle (4.25,2);
      \fill[pattern=north east lines, pattern color=gray!50] (4.25,1.5) rectangle (4.75,1.25);
      \node at (3.5,1.75) {\tiny$+$};
      \node at (4.5,1.75) {\tiny$+$};
      \node[draw,circle,inner sep=1pt] at (4.5,1) {\tiny$+$};
      \node at (4.5,-0.5) {\tiny$\ell_m$};
      \node at (6.0,1) {\tiny$b_m$};
      \node at (6.0,1.75) {\tiny$a_{n}$};
    \end{scope}
  \end{tikzpicture}
}

\newcommand{\permproofiiispecialtwo}{%
  \begin{tikzpicture}[main node/.style={inner sep=0,minimum size =.04cm,circle,fill=black!80,draw},node distance = 0.75cm,scale=0.7]
    \node at (2,0) {\small $\hat{D}$};
    \draw[thick] (4,-1.25) --++(1,0) node[midway,below] {\tiny$i+1=\ell_{1}$};
    \draw[dashed] (3,-0.5) --++ (2.5,0) node[right=3pt] {\tiny $a_n$};
    \node at (3.5,1) {\tiny$+$};
    \node at (3.5,0.25) {\tiny$+$};
    \node at (4.5,0.25) {\tiny$+$};
    \node at (3.5,-1.58) {\tiny$\ell_2$};
    \node at (6.0,0.25) {\tiny$b_{1}$};
    \node at (6.0,1) {\tiny$b_2$};
    \node at (8.25,0) {$\longmapsto$};
    \begin{scope}[shift={(8,0)}]
      \node at (2,0) {\small ${D}$};
      \draw[thick] (4,-1.25) --++(1,0) node[midway,below] {\tiny$i+1=\ell_{1}$};
      \node at (6,-0.5){\tiny $a_n$};
      \node at (3.5,1) {\tiny$+$};
      \node at (3.5,0.25) {\tiny$+$};
      \node at (4.5,0.25) {\tiny$+$};
      \node at (4.5,-0.5) {\tiny$+$};
      \draw (3.25,1.5) -- (3.25,-0.25);
      \draw (3.75,1.5) -- (3.75,-0.25);
      \draw (4.25,0.75) -- (4.25,-0.75) --++(0.5,0) -- (4.75,0.75);
      \fill[pattern=north east lines, pattern color=gray!50] (3.25,0.5) rectangle (3.75,0.75);
      \fill[pattern=north east lines, pattern color=gray!50] (3.75,0) rectangle (4.25,0.5);
      \fill[pattern=north east lines, pattern color=gray!50] (4.25,0) rectangle (4.75,-0.25);
      \node at (3.5,-1.58) {\tiny$\ell_2$};
      \node at (6.0,0.25) {\tiny$b_{1}$};
      \node at (6.0,1) {\tiny$b_2$};
    \end{scope}
  \end{tikzpicture}
}

\newcommand{\permproofpathex}{%
  \begin{tikzpicture}[main node/.style={inner sep=0,minimum size =.04cm,circle,fill=black!80,draw},node distance = 0.75cm,scale=0.7]
    \pgfmathsetmacro{\unit}{0.75}
    \draw (0,5*\unit) --++ (7*\unit,0) --++(0,-1*\unit) node[main node,midway,label=right:{\tiny (i)}] (i) {} --++(-1*\unit,0) --++(0,-1*\unit) --++(-1*\unit,0) node[main node,midway,label=below:{\tiny (i)}] (i2) {} --++(0,-1*\unit) --++(-1*\unit,0) node[main node,midway, label=below:{\tiny (ii)}] (ii) {} --++(0,-1*\unit) node[main node,midway] (ii2) {} --++(-3*\unit,0) --++(0,-1*\unit);
    \node[main node] (ib) at (4.875,3.75) {};
    \node[main node] (i2b) at (2.625,3.75) {};
    \node[main node] (iib) at (-0.25,2.75) {};
    \node[main node] (ii2b) at (1.125,3.75) {};
    \draw (1.625,0.5) rectangle (2.125,4);
    \node at (1.875, 0.15) {\tiny $a_n$};
    \draw[rounded corners,postaction={decorate,decoration={markings,mark=at position 0.8 with {\arrow[scale=0.75]{>}}}}] (i) -| (ib);
    \draw[rounded corners,postaction={decorate,decoration={markings,mark=at position 0.25 with {\arrow[scale=0.75]{>}},mark=at position 0.95 with {\arrow[scale=0.75]{>}}}}] (i2) |- (3.325,2.625) |- (2.625,3.325) -- (i2b);
    \draw[rounded corners,postaction={decorate,decoration={markings,mark=at position 0.15 with {\arrow[scale=0.75]{>}},mark=at position 0.9 with {\arrow[scale=0.75]{>}}}}] (ii) |- (2.625,1.975) |- (1.875,2.375) |- (iib);
    \draw[rounded corners,postaction={decorate,decoration={markings,mark=at position 0.2 with {\arrow[scale=0.75]{>}},mark=at position 0.9 with {\arrow[scale=0.75]{>}}}}] (ii2) -| (2.425,1.775) -| (1.875,2.175) -| (ii2b);
    \node at (7,2) {\small$\hat{D}$};
    \begin{scope}[shift={(10,0)}]
      \draw (0,5*\unit) --++(0,-5*\unit) --++ (3*\unit,0) --++(0,1*\unit) node[main node,midway,label=right:{\tiny (iv)}] (iv) {}  --++(2*\unit,0) --++(0,1*\unit) --++(1*\unit,0) --++(0,2*\unit) --++(1*\unit,0);
      \node[main node,label=right:{\tiny (iii)}] (iii) at (4.5,1.75) {};
      \node[main node,label=below:{\tiny (iii)}] (iii2) at (3.375,0.75) {};
      \node[main node,label=below:{\tiny (iv)}] (iv2) at (0.625,0) {};
      \node[main node] (iiib) at (0,3.375) {};
      \node[main node] (iii2b) at (1.75,4) {};
      \node[main node] (ivb) at (0,1.75) {};
      \node[main node] (iv2b) at (0,0.375) {};
      \draw (-0.25, 2) rectangle (4.75,2.5);
      \node at (5.1,2.25) {\tiny $a_n$};
      \draw[rounded corners,postaction={decorate,decoration={markings,mark=at position 0.8 with {\arrow[scale=0.75]{>}}}}] (iv2) |- (iv2b);
      \draw[rounded corners,postaction={decorate,decoration={markings,mark=at position 0.15 with {\arrow[scale=0.75]{>}},mark=at position 0.95 with {\arrow[scale=0.75]{>}}}}] (iv) -| (1.25,1.125) -| (0.625,1.75) -- (ivb);
      \draw[rounded corners,postaction={decorate,decoration={markings,mark=at position 0.15 with {\arrow[scale=0.75]{>}},mark=at position 0.95 with {\arrow[scale=0.75]{>}}}}] (iii2) |- (2.625,1.125) |- (1.75,1.75) -- (iii2b);
      \draw[rounded corners,postaction={decorate,decoration={markings,mark=at position 0.03 with {\arrow[scale=0.75]{>}},mark=at position 0.9 with {\arrow[scale=0.75]{>}}}}] (iii) -| (4.125,2.75) -| (3.375,3.375) -- (iiib);
    \end{scope}
  \end{tikzpicture}
}

\bibliography{le-diagrams/lediagram,c2/c2}

\title{A Combinatorial Tale of Two Scattering Amplitudes: \\ See Two Bijections}
\author[S. Simone Hu]{Simeng Simone Hu}

\begin{document}
\maketitle
\pagenumbering{roman}

\begin{center}
\vfill
A thesis \\
presented to the University of Waterloo \\
in fulfillment of the \\
thesis requirement for the degree of \\
Master of Mathematics \\
in \\
Combinatorics and Optimization \\

\vfill
Waterloo, Ontario, Canada, 2021 \\
\vspace*{1.0cm}
\copyright{ Simeng Simone Hu 2021}
\vfill
\end{center}

\newpage
\chapter*{Author's declaration}

I hereby declare that I am the sole author of this thesis.
This is a true copy of the thesis, including any required final revisions, as accepted by my examiners.

I understand that my thesis may be made electronically available to the public.

\newpage
\chapter*{Abstract}

In this thesis, we take a journey through two different but not dissimilar stories with an underlying theme of combinatorics emerging from scattering amplitudes in quantum field theories.

The first part tells the tale of the $c_2$-invariant, an arithmetic invariant related to the Feynman integral in $\phi^4$-theory, which studies the zeros of the Kirchoff polynomial and related graph polynomials.
Through reformulating the $c_2$-invariant as a purely combinatorial problem, we show how enumerating certain edge bipartitions through fixed-point free involutions can complete a special case of the long sought after $c_2$ completion conjecture.

The second part tells the tale of the positive Grassmannian and a combinatorial T-duality map on its cells, as related to scattering amplitudes in planar $\mathcal{N}=4$ SYM theory.
In particular, T-duality is a bridge between triangulations of the hypersimplex and triangulations of the amplituhedron, two objects that appear as images of the positive Grassmannian.
We give an algorithm for viewing T-duality as a map on Le diagrams and characterize a nice structure to the Le diagrams (which can then be used in lieu of the algorithm).
Through this Le diagram perspective on T-duality, we show how the dimensional relationship between the positroid cells on either side of the map can be directly explained.

\newpage
\chapter*{Acknowledgements}

I would first and foremost like to thank my supervisor Karen Yeats for introducing me to this wonderfully  interesting field and for all the support, encouragement and advice throughout the years (especially this one), I will be forever grateful!

A big thank you to my thesis readers, Logan Crew and Olya Mandelshtam, for all their helpful comments and valuable feedback.

I would also like to thank my friends and family for camping with me through this whole process.
Thank you to the you who has always believed in me.
(:

\newpage
\setcounter{tocdepth}{1}
\pdfbookmark[0]{\contentsname}{toc}
\tableofcontents
\setcounter{tocdepth}{2}

\newpage
\newlength\oldparskip
\setlength\oldparskip\parskip
\setlength{\parskip}{0pt}
\listoffigures

\newpage
\setlength{\parskip}{\oldparskip}
\makeatletter
\def\@evenhead{%
  \setTrue{runhead}%
  \normalfont\scriptsize
  \hfil
  \def\thanks{\protect\thanks@warning}%
  \leftmark{}{}\hfil}%
\def\@oddhead{%
  \setTrue{runhead}%
  \normalfont\scriptsize \hfil
  \def\thanks{\protect\thanks@warning}%
  \rightmark{}{}\hfil
}%
\let\@mkboth\markboth
\def\partmark{\@secmark\markboth\partrunhead\partname}%
\def\chaptermark{%
  \@secmark\markright\chapterrunhead{}}%
\def\sectionmark{%
  \@secmark\markright\sectionrunhead\sectionname}%
\let\sectionmark\@gobble
\markboth{\MakeUppercase\shortauthors}{\MakeUppercase\shorttitle}%
\makeatother

\pagenumbering{arabic}
\newrefsegment
\chapter{Introduction}\label{C:intro}

\section[Combinatorics in disguise in the world of physics]{\texorpdfstring{Combinatorics in disguise in the world of physics\except{toc}{: \\ The story of two types of scattering amplitudes}}{Combinatorics in disguise in the world of physics}}

One of the fundamental problems in high energy physics is to understand how particles interact, how they scatter when they collide.
For example, as illustrated in Figure~\ref{fig:higgs}, we can think of experiments done in particle accelerators like the Large Hadron Collider (LHC) in which some number of particles are brought together, and the particles resulting from that interaction are detected.
\textbf{Quantum field theory} (QFT) is one such framework to understand these processes, and one which unifies quantum mechanics and special relativity.
A natural question that arises is: what is the probability of a given configuration of particle scatterings?
In quantum field theory, this probability is encoded in the \textbf{scattering amplitude}~\cite{iz}.

However, we don't know what actually happens when the particles interact and what paths the particles are taking, the "history" of the collision.
Perturbative quantum field theory thus expresses each of these scatterings as a formal power series, a perturbative expansion, in the \textbf{coupling constant} which is measuring the strength of the interactions.
In this expansion, each term corresponds to a particular "history" that the particles trace out when interacting.
In other words, in the spirit of quantum physics, we are taking a weighted sum of all possible interactions that get us to some particular scattering configuration.
The weights in this sum are called \textbf{Feynman integrals} and the scattering amplitude of some configuration is then the sum of all the Feynman integrals in the expansion.
However, in general these integrals are notoriously hard to calculate and the number of terms in each sum rapidly expands as the number of particles grows.

The fascinating part about these particular problems is the emergence of rather beautiful connections to combinatorics, and the study of discrete structures.
Here begins our tale of two different connections to combinatorics as related to scattering amplitudes in quantum field theory. \\

\begin{minipage}{\textwidth}
  \addtocontents{lof}{\protect\contentsline{chapter}%
    {\protect\tocchapter{\chaptername}{\thechapter}{\thischaptertitle}}{}{}}
  \global\togglefalse{noFigs}
  \centering
  \includegraphics[scale=0.14]{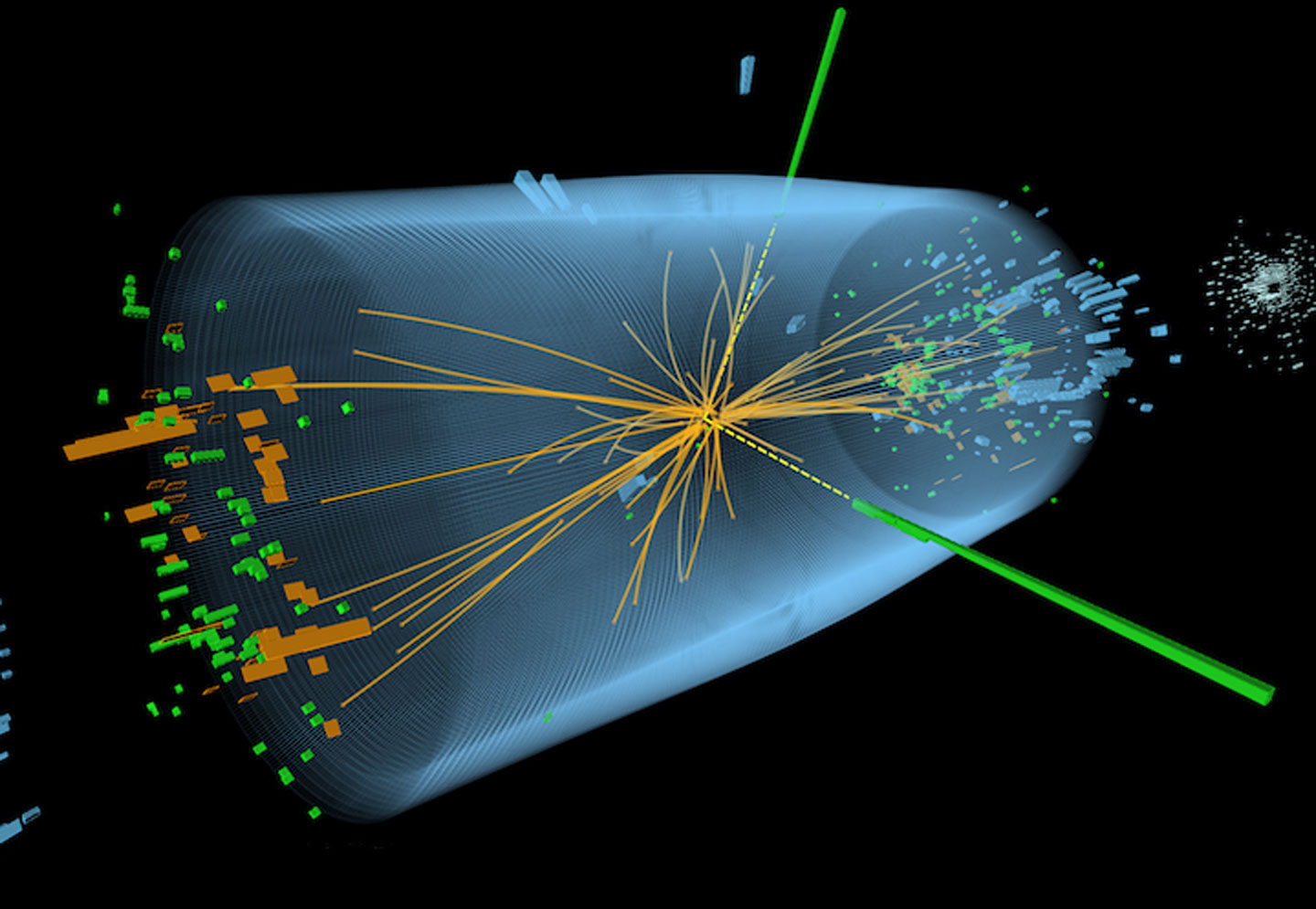}
  \captionof{figure}[A candidate Higgs boson event.]{A candidate Higgs boson event from collisions between protons in the LHC where the particle is decaying into two protons indicated by the yellow dashed-lines with green towers (Image: CMS/CERN)}
  \label{fig:higgs}
\end{minipage}

\bpoint{Feynman diagrams}

One way to view these particle scatterings is through the underlying graphs of the interactions, called Feynman diagrams or \textbf{Feynman graphs}, formed by the paths that the particles take.
The edges of the graphs correspond to the propagations of particles, while the vertices correspond to the collisions or interactions.
Incoming and outgoing particles are then represented by external half-edges (legs).
As such, the complexity of a problem involving Feynman graphs can be measured in two ways, via the number of cycles, which physicists call \textbf{loops}, or via the structure of the external legs.
For a particular field theory, we can further encode the different fundamental particles through having different edge types, and the permissible interactions between them can be encoded through allowable vertex types (as defined by the number of each edge type adjacent to a vertex).

\begin{figure}[ht]
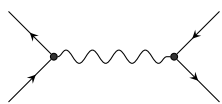

  \centering
  \QED
  \caption[A Feynman diagram in quantum electrodynamics.]{A Feynman diagram in quantum electrodynamics (QED), a quantum field theory of the electromagnetic force. The squiggly edge represents the propagation of a proton and the directed half-edges represent incoming and outgoing fermions.}
  \label{fig:qed}
\end{figure}

Each graph can then be associated to a Feynman integral though the \textbf{Feynman rules}.
That is the rules dictate how to build the integral based on the structure of the graph.
In the standard approach to quantum field theory, the Feynman rules can be seen as arising from the Lagrangian in the path integral formulation for the scattering amplitude.
However, by taking a graphs-first approach to particle interactions, this opens the door to using tools from the well established and classical field of graph theory to tackle problems in perturbative quantum field theory~\cite{perspective}.

In practice, as many of the Feynman integrals are often divergent, physicists introduced a technique called \textbf{renormalization} in order to extract meaningful quantities from these integrals.
As an example, one such scheme for renormalizing Feynman integrals is called \textbf{BPHZ renormalization}, which we can think of as recursively subtracting off divergent parts of the integral corresponding to divergent subgraphs, "subdivergences".
In particular, renormalization itself has an interesting connection to combinatorics via the combinatorial Hopf algebras that appear as the underlying algebraic structure to some of these schemes.
For BPHZ renormalization, this is the Connes-Kreimer Hopf algebra of rooted trees.
Renormalization, Hopf algebras, and recursive structures on Feynman diagrams then also open the door to Dyson-Schwinger equations, Green's functions, and their related combinatorics such as chord diagrams~\cite{perspective}.
With renormalization, we can much more easily study graphs with high loop orders and a low number of external legs.

Going back to Feynman integrals, other than their relation to the scattering amplitude, one of the other reasons to study them lies in the many interesting numbers that they can evaluate to.
As an example, there are Euler sums, some of which evaluate to zeta values, and polylogarithms.
However as we mentioned before, these integrals are quite complicated due to the need for renormalization and very hard to compute.
Thus in order to study the numbers that appear as Feynman integrals, we instead restrict to a particular residue of the integral called the \textbf{Feynman period}~\cite{periods} which is renormalization scheme independent and captures much of the number-theoretic content for massless scalar field theories in four-dimensional space-time.
We can think of this period as taking only the "primitive divergent" part of the leading contribution to the integral, that is the part corresponding to subdivergence free Feynman graphs.

We will focus on periods in massless Euclidean \textbf{$\mathbf{\phi^4}$-theory}, which is a scalar field theory with one particle type and whose interactions are all quartic, that is, there is only one type of edge and all vertices are 4-valent.
In graph theory language, the Feynman graphs are just 4-regular graphs with some number of external half-edges.
Historically, $\phi^4$-theory has been an important toy model for more complicated field theories.
As an motivating example, some of the numerical structure seen in the anomalous magnetic moment of the electron in perturbative QED is also manifest in the beta-function of $\phi^4$-theory via its coefficients, which are periods.

\begin{figure}[t]
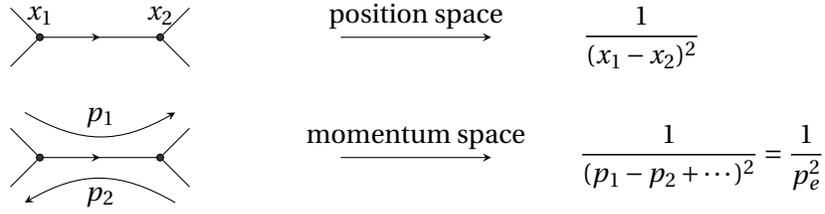

  \centering
  \Feynman
  \caption[Feynman rules for massless scalar field theories.]{Feynman rules for massless scalar field theories. For momentum space, $p_e$ refers to the signed sum of the cycles running through edge $e$.}
  \label{fig:feynman-rules}
\end{figure}

Now, to build the period of a given primitive divergent graph $G$, we can use the Feynman rules which give equivalent representations of the period in the different spaces: position, momentum, parametric, and dual parametric.
The first two of these corresponds to assigning 4-dimensional vectors to vertices and cycles (in a cycle basis) respectively, with each edge of the graph contributing a factor to the integrand following the rules in Figure~\ref{fig:feynman-rules}.
Note that we are arbitrarily orienting the edges and cycles of the graph and choosing a cycle basis in momentum space.
For a graph with loop order $\ell$, i.e. the number of cycles or first Betti number, the period of $G$ in momentum space is
\[ P_G = \pi^{-2(\ell-1)}\bigintsss \frac{\;\mathrm{d}p_2 \cdots \,\mathrm{d}p_{\ell}\;}{\displaystyle\prod_e p_e^2 \,\bigg|_{p_1=\mathbb{1}}}, \]
where we are integrating over $\mathbb{R}^{4(\ell - 1)}$ and we take $p_e^2$ to mean taking the norm squared.
Note that we can set one momentum vector, say $p_1$, to a fixed choice of unit $\mathbb{1}$ as we can always normalize the momentum variables.
In momentum space, we can then think of the period as the Feynman integral where we are setting all external properties, masses and momenta, to zero (and setting a variable to 1).

We could similarly represent the period in position space using the rules or through a Fourier transform of the momentum space period.
Using the Schwinger trick, we can transform the period to parametric space, which will be our main focus,
\[ P_G = \bigintsss_{\; 0}^{\infty} \frac{\;\mathrm{d}\alpha_1 \cdots \,\mathrm{d}\alpha_{\abs{E(G)}-1}\;}{\Psi_G^2 \,\bigg|_{\alpha_{\abs{E(G)}}=1}}. \]
Here $\Psi_G$ is the \textbf{Kirchhoff polynomial}.
Continuing with the graphs-first approach to study the period, here starts one story of combinatorics emerging in scattering amplitudes.
We will continue with this story in Part~\ref{P:c2}.  \\

\bpoint{On-shell diagrams}

One of the limitations of the classical formulation of perturbative quantum field theory is the explosion in the number of Feynman diagrams that can appear for any particular scattering configuration, even in the simple case of two particles interacting.
This is a consequence of taking an all possible "histories" expansion, which introduces a large amount of redundancy adding to the complexity of computing scattering amplitudes.
These redundancies can be seen as coming from the fact that we don't actually know the paths that the particles are taking and the internal edges in the Feynman diagrams are actually representing \textbf{virtual}, or \textbf{off-shell particles}.
They are called virtual particles as they do not satisfy the energy-momentum relation and thus are non-physical particles whose momenta cannot exist.

As an answer to this problem, a small group of researchers~\cite{grassampl} (as part of an ongoing program) developed a new way of thinking about scattering amplitudes in quantum field theory, one using only physical, or \textbf{on-shell particles}.
This reformulation illuminated the many beautiful mathematical structures underlying on-shell processes and the remarkable simplicity of their amplitudes.

At the heart of it all, this theory was built upon the extensive exploration of scattering amplitudes in planar \textbf{$\mathcal{N}=4$ supersymmetric Yang-Mills (SYM) theory}.
To try to simply explain what this theory is, we take a brief detour using the language of differential geometry.
Instead of working directly on space-time like scalar field theories, \emph{gauge theories} are models that are defined on fibre bundles over space-time in which the Lagrangian is invariant under certain symmetries, \emph{gauge transformations}.
To get back to space-time, one would then need to pick a local section, a \emph{gauge}.
In particular, \emph{Yang-Mills theory} is one such important gauge theory when trying to understand the Standard model.
To study this rather complicated theory, researchers turn to $\mathcal{N}=4$ SYM, which is a maximally supersymmetric theory often considered as a toy model for four-dimensional Yang-Mills theory, and thus for quantum chromodynamics (QCD), which is a quantum field theory (and a particular Yang-Mills theory) describing the strong interactions between quarks and gluons.
The $\mathcal{N}=4$ here refers to the number of supersymmetries.
While $\mathcal{N}=4$ SYM incorporates many complicated particles into its theory, its highly symmetric nature lends itself to having many special properties and dualities, and is what makes this theory attractive.

Unlike the Feynman integral story where the main study is on diagrams with high loop orders but low numbers of external legs, much of the work in these gauge theory amplitudes started with the opposite complexity; looking at low loop orders but high numbers of external legs.
In particular, at tree-level (no loops) one of the breakthroughs has been in the discovery of the \textbf{BCFW recurrence relations}~\cite{BCFW} for amplitudes in Yang-Mills theory, which could decompose the amplitude in various ways depending on the choices at each recursive step.
Furthermore, each step only ever involved on-shell particles.

Building upon these relations and related work at the time, a surprising connection between scattering amplitudes in $\mathcal{N}=4$ SYM theory and the \textbf{positive Grassmannian}, which is a certain "non-negative" subset of the Grassmannian, started to emerge.
In the flurry of work which culminated in~\cite{grassampl}, \textbf{on-shell diagrams} were introduced to describe on-shell processes like those that decompose the scattering amplitude and which arise from glueing together three-particle amplitudes.

\begin{figure}[ht]
  \centering
  \includegraphics[scale=0.47]{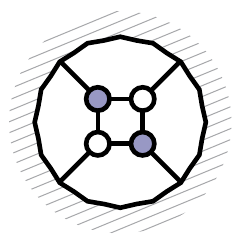}
  \caption[An on-shell diagram.]{An on-shell diagram which represents a four-particle tree-level amplitude. Note that "tree" here refers to the amplitude in the context of the Feynman diagram expansion. (Image:~\cite{grassampl})}
  \label{fig:on-shell}
\end{figure}

In graph theory language, these diagrams are bi-coloured trivalent graphs which are planar when in the context of planar $\mathcal{N}=4$ SYM.
In this story about scattering amplitudes, the on-shell diagram plays the role of the Feynman diagram.
These on-shell diagrams, which satisfied certain equivalences, turned out to be \textbf{plabic graphs}, which Postnikov~\cite{tpgrass} introduced as one of the many combinatorial objects that index cells of the positive Grassmannian!

Inspired by these connections and the question of where the on-shell diagram decompositions of the amplitude comes from, Arkani-Hamed and Trnka~\cite{ampl} introduced a new mathematical object called the \textbf{amplituhedron}.
They proposed that the $m=4$ amplituhedron was a geometric representation of scattering amplitudes in planar $\mathcal{N}=4$ SYM, which would be encoded through its "volume", and one in which the BCFW recursion relations could be manifested through its "triangulations".
Through the positive Grassmannian, the amplituhedron and triangulations, herein lies our second combinatorial story, this time in a different type of scattering amplitude.
We will pursue this story in Part~\ref{P:le}.  \\

\section{Overview and organization}

This thesis is composed of two self-contained parts, consisting of two stories with the underlying theme of combinatorics related to scattering amplitudes in quantum field theories.

In the first part, Part~\ref{P:c2}, we delve deeper into the combinatorics of an arithmetic invariant, called the \textbf{$c_2$-invariant}, related to a particular residue of the Feynman integral called the Feynman period. This invariant studies the zeros of the denominators of these periods that arise after several integrations.

We start off in Chapter~\ref{C:c2background} with an overview of the main object of interest in this story, the $c_2$-invariant.
Section~\ref{S:phi4} provides some background on $\phi^4$-theory and Feynman periods, setting up the motivation for defining the $c_2$-invariant which we present in Section~\ref{S:c2-def}.
In Section~\ref{S:graph-polys} we develop the framework of graph polynomials, specifically Dodgsons and spanning forest polynomials, which then allows us to start interpreting the $c_2$ combinatorially via spanning trees and spanning forests.
Going back to the period and using this framework, in Section~\ref{S:denominator} we present the main algorithms used to compute the $c_2$ which are based on the denominators appearing in periods after successive integrations.
In particular in Section~\ref{SS:denom-count} we see how we can reformulate the $c_2$ into a purely combinatorial problem, one about taking particular coefficients of graph polynomials!
We end this introductory chapter with Section~\ref{S:families}, giving some of the symmetries of and computational results on the $c_2$-invariant.

In Chapter~\ref{C:completion} we answer a special case of an over 10-year old conjecture on the $c_2$-invariant using combinatorial methods involving enumerating certain edge bipartitions, completing the argument as first started by Yeats~\cite{specialc2}.
We first set-up the problem in Section~\ref{S:conjecture} and then give a complete, self-contained proof of each of the three cases in the following sections, Sections~\ref{S:T-case},~\ref{S:S-case} and ~\ref{S:R-case}.
In particular in Sections~\ref{SS:S-new} and~\ref{SS:R-new} we prove two new involutions that together with results from~\cite{specialc2} completes the proof of the conjecture in their respective cases.
We assemble everything in Section~\ref{S:completing-conjecture} and give the main result of Part~\ref{P:c2} in Theorem~\ref{p2-completed}.

The second part, Part~\ref{P:le}, then moves to exploring the combinatorics of a geometric object called the \textbf{positive Grassmannian} and a particular map on cells of the positive Grassmannian called \textbf{T-duality}.
These cells on one side of the map are then related to the on-shell diagrams from $\mathcal{N}=4$ SYM theory.

We open this part with an introductory chapter, Chapter~\ref{C:lebackground}, on the many objects behind T-duality.
In Section~\ref{S:grass} we motivate the study of the positive Grassmannian through one particular decomposition of the Grassmannian called its matroid stratification.
Section~\ref{S:positroids} then introduces the main combinatorial objects of interest that are in bijection with cells of the positive Grassmannian, which we call positroids cells as they are indexed by positroids.
We end this chapter with Section~\ref{S:t-duality}, where we look at the geometric objects arising from images of the positive Grassmannian under two different maps and how they are related.
In particular in Section~\ref{SS:t-duality}, we define T-duality as a combinatorial map which arose from looking at triangulations of these geometric objects, and set-up the motivation for the central question of Chapter~\ref{C:lebijection}.

In Chapter~\ref{C:lebijection} we investigate what T-duality looks as a map on Le diagrams, one of the combinatorial objects in bijection with positroids, doing so in two different ways.
Starting with Section~\ref{S:algorithm}, we define an explicit algorithm in Section~\ref{SS:algorithm} for determining the Le diagram resulting from the T-duality map.
At the end of the section, in Section~\ref{ex:algorithm}, we work through an example of this intricate algorithm.
In Section~\ref{S:visual} we notice that there is actually a nice structure to the Le diagrams created from the above algorithm, giving a more visual perspective of how T-duality affects Le diagrams, the main result being Theorem~\ref{glueing}.
Finally in Section~\ref{S:proof}, we prove that the algorithm does indeed give the correct Le diagram, resulting in Theorem~\ref{letheorem}.
In particular, here we see in Theorem~\ref{dimension} how viewing T-duality on Le diagrams directly explains the dimensional relationship between the positroid cells on either side of the map.

We conclude this thesis in Part~\ref{P:final} with Chapter~\ref{C:conclusion}, giving some further research directions for both stories in Section~\ref{S:further} and tying everything together with some final thoughts in Section~\ref{S:final}.

\newpage
\newrefsegment
\part{The \texorpdfstring{$c_2$}{c2}-invariant}\label{P:c2}

\chapter{On the \texorpdfstring{$c_2$}{c2}-invariant}\label{C:c2background}

In the first part of this thesis, we present one story on the combinatorics arising in perturbative quantum field theory.
Among the many combinatorial connections, this is specifically a story about an arithmetic invariant called the $c_2$-invariant that is related to scattering amplitudes in scalar $\phi^4$ theory and the Feynman integrals in their expansion.

We begin with an introductory chapter on the $c_2$-invariant, providing the background needed for the specific problem that we will focus on in Chapter~\ref{C:completion}, which aims to answer one of the conjectured symmetries of the $c_2$-invariant via combinatorial techniques.
As we will see in Section~\ref{SS:denom-count}, a combinatorial picture emerges where the $c_2$-invariant can be thought of as an enumeration problem on graphs, counting certain edge partitions. \\

\section{Motivation from \texorpdfstring{$\phi^4$}{phi-4-}-theory}\label{S:phi4}

\bpoint{4-point graphs}

\begin{figure}[th]
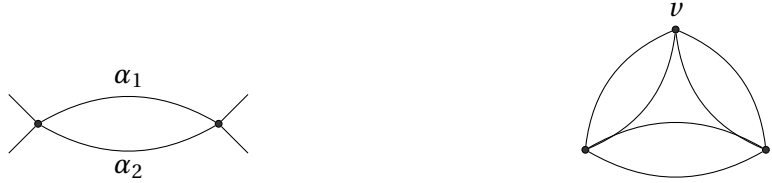

  \centering
  \captionsetup{justification=centering}
  \begin{subfigure}[b]{.45\linewidth}
    \centering
    \GrphOne
    \caption{A primitive one-loop Feynman graph.}
    \label{fig:G1}
  \end{subfigure}
  \begin{subfigure}[b]{.45\linewidth}
    \centering
    \CmplOne
    \caption{The unique completion of (a).}
    \label{fig:Comp1}
  \end{subfigure}
  \caption{The primitive graph of one loop and its completion.}
  \label{fig:l1}
\end{figure}

The Feynman graphs we are interested in are those arising from 4-point Feynman integrals in four-dimensional ($D=4$) massless Euclidean $\phi^4$-theory.
The "4-point" here corresponds to the number of particles in the scattering configuration that we are interested in and the "massless" corresponds to setting masses to 0.
Viewing this scalar field theory as a \textbf{combinatorial physical theory} (see \textsection 5.2 of~\cite{perspective}), the Feynman graphs in question correspond to graphs with 4 external half-edges, or "legs", where the (internal) edges are thought of as consisting of two half-edges, and where every vertex is of degree 4 (with external edges contributing to the degree).
As we will be considering the Feynman period in the massless case, these external legs can usually be disregarded outside of playing a role in vertex degree counting.
An example of a 4-point graph in $\phi^4$-theory is given in Figure~\ref{fig:G1}.

Let $G$ be a 4-point graph in $\phi^4$-theory, that is a 4-regular graph, not necessarily simple, with 4 external legs.
We will also assume throughout that $G$ is connected.
Here $\abs{E(G)}$ is the number of internal edges of $G$ and $\ell(G)$ is the loop order of $G$.

As we want to talk about Feynman integrals in this combinatorial context, recalling that renormalization was needed to circumvent the divergent nature of these integrals, we need to translate these notions into our setting.
To measure how badly a Feynman integral diverges as the energies get large, we can use power counting (counting the exponents) on the integration variables in comparison with the number of such variables, which can then be directly distilled via the Feynman rules into a property of the underlying graph.
More precisely, the \textbf{superficial degree of divergence} (sdd) of a graph is defined as
\[\text{sdd}(G) = D\ell(G) -\sum_e w(e) - \sum_v w(v), \]
where $D$ is the dimension of space-time, $w$ is the power counting weight function associated to the specific combinatorial physical theory in question and the sums run over all internal edges $e$ and all vertices $v$, respectively.
In scalar $\phi^4$-theory, which has only one type of edge and vertex, edges have weight $2$ ($w(e) = 2$) and vertices have weight $0$ ($w(v) = 0$).
We can see this directly through the Feynman rules in momentum space (see Figure~\ref{fig:feynman-rules}), where each edge contributes a power of 2 of its associated variable in the denominator of the integrand.
Generally, there may be different weights for different edge types and for different vertex types.

We call a graph \textbf{divergent} when its sdd is non-negative and \textbf{logarithimcally divergent} when it is exactly $0$.
We then say that a combinatorial physical theory is \textbf{renormalizable} if the sdd of a graph in the theory only depends on the structure of the external legs.

In our case, $\phi^4$-theory is renormalizable in $D = 4$ and thus the superficial degree of divergence of $G$ is
\[ \text{sdd}(G) = 4\ell(G) - 2\abs{E(G)} = 4 - q, \]
where $q$ is the number of external legs.
Here we use that there are ${2\abs{E(G)} + q = 4\abs{V(G)}}$ half-edges and Euler's formula for connected graphs,
\[ \abs{V(G)} - \abs{E(G)} + \ell(G) = 1, \]
which gives
\[ 4\ell(G) - 2\abs{E(G)} = 4(1 + \abs{E(G)} - \abs{V(G)}) - 2\abs{E(G)} = 4 + 2\abs{E(G)} - 4\abs{V(G)} = 4 - q. \]

As $q=4$, every 4-point graph is logarithmically divergent (because $\text{sdd}(G) = 0$), giving the equality $\abs{E(G)} = 2\ell(G)$.
Using Euler's formula again gives us that $\abs{V(G)} = \ell(G) + 1$. \\

\bpoint{Feynman periods}

Instead of looking at the full Feynman integral, we will look at a particular residue of it called the \textbf{Feynman period}.
In parametric space, the Feynman period (also referred to just as the period) of a graph $G$ is defined as
\begin{equation}\label{eq:period}
  P_G \coloneqq \bigintsss_{\; 0}^{\infty} \frac{\;\mathrm{d}\alpha_1 \cdots \,\mathrm{d}\alpha_{\abs{E(G)}-1}\;}{\Psi_G^2 \,\bigg|_{\alpha_{\abs{E(G)}}=1}},
\end{equation}
where we associate a Schwinger parameter $\alpha_e$ to each edge $e$ in $G$ and
\[ \Psi_G = \sum\limits_{\substack{T \\ \text{spanning tree}}} \prod_{e \notin T} \alpha_e \]
is the \textbf{graph polynomial} or \textbf{Kirchhoff polynomial} of G.
Note that this integral is independent of the choice of edge $e=\abs{E(G)}$ to set $\alpha_e = 1$.
There are many equivalent formulations of this integral, in the different spaces (position, momentum etc.) as well projective versions of each.
We refer the reader to ~\cite{census} and \textsection 3.1, 3.2 of~\cite{periods} for further details.

Note that the term "period" comes from algebraic geometry:
looking at $P_G$ in its parametric form $\Psi_G$ (if it exists) is simply a polynomial in the variables $\alpha_e$ with integer coefficients.
Thus $\Psi_G^{-2}$ is a rational function with $P_G$ the evaluation of its integral over $\alpha_e \geq 0$.
That is, $P_G$ is a period as defined by Kontsevich and Zagier ~\cite{konzag} and in the same sense as how multiple zeta values are periods.

It turns out that $P_G$ is well-defined when the graph $G$ is primitive and logarithmically divergent.
In fact, these are necessary and sufficient conditions for the convergence of $P_G$, see Proposition 5.2 of~\cite{motives}.
In this case $P_G$ defines a positive real number which we call the \textbf{period of $\mathbf{G}$}.
In physics language, $G$ is a 4-point graph with no 1PI divergent subgraphs and algebraically primitive here means primitive for the co-product of the renormalization Hopf algebra on Feynman graphs.
In terms of graphs we can define primitivity with logarithmic divergence included as follows.

\tpointn{Definition}\label{def:primitive}
\statement{
  A graph $G$ is \textbf{primitive} (or \textbf{primitive divergent}) if
  \begin{itemize}
    \item $\abs{E(G)} = 2\ell$; where $\abs{E(G)}$ is the number of edges in $G$ and $\ell = \ell(G)$ is the loop number of $G$, and
    \item every non-empty proper subgraph $\gamma \subset G$ has $\abs{E(\gamma)} > 2\ell(\gamma)$.
  \end{itemize}
}

We note that the second condition tells us our primitive graphs have no divergent subgraphs (subdivergences), as
\[ \text{sdd}(\gamma) = 4\ell(\gamma) - 2\abs{E(\gamma)} \geq 0 \quad\Longleftrightarrow\quad 2\ell(\gamma) \geq \abs{E(\gamma)}. \]

We also have that any primitive graph with at least three vertices must be simple, namely without multiple edges and self-loops.
The only non-simple primitive graph is given in Figure~\ref{fig:G1}, which is also the only graph of loop order $\ell = 1$.
This graph also (most likely) corresponds to the only rational $\phi^4$ period and has by far the easiest period to calculate.

\tpointn{Example} (A period calculation)
\statement{
Consider the (unique) primitive graph with one loop $G$ in Figure~\ref{fig:G1}.
Using Equation~\eqref{eq:period} we get
\begin{align*}
  P_G
  &= \bigintsss_{\; 0}^{\infty} \frac{\;\mathrm{d}\alpha_1}{\Psi_G^2 \,\bigg|_{\alpha_2=1}}
  = \bigintsss_{\; 0}^{\infty} \frac{\;\mathrm{d}\alpha_1}{(\alpha_1 + \alpha_2)^2 \,\bigg|_{\alpha_2=1}}
  = \bigintsss_{\; 0}^{\infty} \frac{\;\mathrm{d}\alpha_1}{(\alpha_1 + 1)^2}
  = - \left. \frac{1}{1+\alpha_1} \right\rvert_0^{\infty} = 1.
\end{align*}
}

As the loop order $\ell$ increases, so does the number of non-isomorphic Feynman graphs and the difficulty of the integrations.
Starting from $\ell=3$, we already see the appearance of special numbers, namely multiple zeta values (see Section~\ref{SS:mzv}), as periods.
For example, the complete graph on four vertices $K_4$ can be viewed as a primitive $\phi^4$ graph and has a period of $6\zeta(3)$.
\\

\bpoint{Period symmetries}\label{SS:symmetries}

One way to tackle the explosion of graphs is to find symmetries of the period, allowing us to create equivalence classes of graphs with the same period.
One particularly important symmetry is called \textbf{completion}, which reduces the problem of calculating periods for primitive divergent graphs to one on 4-regular graphs.
Note that 4-regular graphs can also be thought of as $\phi^4$ graphs with no external legs and thus no external momenta, sometimes called \textbf{vacuum graphs}.

As we are looking at graphs with 4 external legs where each vertex is of degree 4, notice that we can uniquely "complete" such graphs by adding a new vertex connected to all the external half-edges creating a 4-regular graph which we say is the \textbf{completion} of the primitive graph.
Conversely, given a 4-regular graph and then "decompleting" it by removing a vertex, we are not guaranteed to always get a primitive divergent graph, nor will we always get the same resulting graph when picking a different vertex to remove.
Thus first we need a notion of primitivity for these 4-regular graphs such that removing any vertex indeed gives a primitive divergent graph which we then call a \textbf{decompletion} of the 4-regular graph.

\tpointn{Definition} (Definition $2.4$ of~\cite{census})\label{def:completedprim}
\statement{
  A 4-regular graph $G$ with $\geq 3$ vertices is called \textbf{completed primitive} if the only way to split $G$ into multiple connected components with 4 edge cuts is to separate off a vertex, in other words, there are only trivial 4-edge cuts.
  In other words, $G$ is \textbf{internally 6-edge connected}.
  In this case, we say that $G$ has \textbf{loop order} $\ell$ if $G - v$ has loop order $\ell$ for any vertex $v$.
}

As an example, Figure~\ref{fig:l1} gives a primitive graph and its unique completion, which is completed primitive.
The complete graph on 5 vertices $K_5$ is another example of a completed primitive graph.
Notice that any decompletion of $K_5$ gives the primitive graph $K_4$.
In Section~\ref{SS:trivialc2}, we briefly discuss why this completed primitive condition is necessary by relating it back to having no subdivergences in primitive graphs.

\tpointn{Proposition} (Proposition $2.6$ of~\cite{census})
\statement{
  Let $G$ be a 4-regular graph and $v$ any vertex in $G$.
  Then $G$ is completed primitive if and only if $G - v$ is primitive.
}

With this notion of primitivity for completed graphs, Schnetz~\cite{census} then proved that the period is completion invariant: any two decompletions of the same 4-regular completed primitive graph have the same period.

\tpointn{Theorem} (Theorem $2.7$ of~\cite{census})\label{period-comp}
\statement{
  Let $G$ be a 4-regular completed primitive graph of loop order $\ell$.
  The period of $G - v$ for any vertex $v$ is the same for all choices of $v$.
}

This completion symmetry tells us that rather than looking at periods of 4-point graphs, we can instead focus on 4-regular graphs.
Note that completion considerably reduces the number of relevant graphs at each loop order $\geq 5$.

Other than completion, there are currently four other known period symmetries which we briefly describe below.
For each of the symmetries, we are taking our graphs $G$ to be completed primitive, that is 4-regular and internally 6-edge connected.
The last three symmetries tells us that under their respective transforms, the period stays invariant.

\begin{itemize}[itemsep=5pt]
  \item \textbf{Product identity} for 3-vertex joins: When there is a 3-vertex cut, in which case we call $G$ \textbf{reducible}, the period decomposes as a product of two periods corresponding to the two sides of the cut where to each side a triangle is added to the three vertices in the cut.
    In other words, this occurs when $G$ is the 3-vertex join of two connected completed primitive graphs.
    For primitive divergent graphs, this product identity holds for 2-vertex joins.
  \item \textbf{Planar duality}, also called the \textbf{Fourier identity}: One natural period identity arises from reinterpreting the Fourier transform as a graph transform.
    Graphically, this is taking the planar dual of a decompletion $G - v$, and then completing it (if possible).
  \item \textbf{Twist identity}: This identity due to Schnetz arises from the twist transform which uses a 4-vertex cut that separates $G$ into two connected subgraphs and reattaches one side of this separation in a particular manner.
  \item \textbf{Fourier split identity} -- This is a relatively new period symmetry which combines the ideas of the Fourier and twist transforms.
    Instead of just reattaching the subgraph in the twist transform, the Fourier split transform first takes a dual of this subgraph and then reattaches it in a particular manner.
\end{itemize}

For more details on the first three symmetries, see \textsection 2.5-2.7 of~\cite{census} and for the Fourier split see \textsection 2.3 of~\cite{further}.
We note that in the context of Feynman integrals and primitive divergent graphs, the product identity for 2-vertex joins is a well-known result and the Fourier identity has been used as early as ~\cite{knots}.

While perhaps not as hard as the Feynman integral to compute, the period still poses quite the computational challenge even after utilizing these symmetries.
As such, related invariants were introduced to study properties of the period such as the $c_2$-invariant, which is our main focus here, the extended graph permanent~\cite{permanent} and the Hepp bound~\cite{hepp}.
It is either conjectured or proven that these invariants have all the above symmetries of the period. \\

\bpoint{Multiple zeta values}\label{SS:mzv}

The number-theoretic content of the period is itself an interesting problem, in particular due to its relationship to the underlying geometries of the period, see for instance~\cite{motives,periods,k3}.
By "geometries", we refer to types of varieties that appear as defined by the vanishing of the graph polynomial.
In some sense these varieties control the period, giving an algebro-geometric feel to the period.
In particular, many of the early computable periods~\cite{knots,census} were all found to be rational linear combinations of multiple zeta values.

For $s_1,\ldots,s_k$ positive integers and $s_1 > 1$, the \textbf{multiple zeta value} (MZV) is defined to be
\[ \zeta(s_1,\ldots,s_k) = \sum_{n_1 > \cdots > n_k > 0} \frac{1}{n_1^{s_1}\cdots n_k^{s_k}}, \]
and we call the sum of the exponents $n = s_1 + \cdots + s_k$ the \textbf{weight} of the MZV.

Viewing an MZV through its iterated integral representation, we can then also think of its weight in the sense of Kontsevich-Zagier periods, as the minimum number of nested integrals needed to write it as such a period.
As such, we refer to the \textbf{weight} of a Feynman period in this same way, as the mininum number of integrals needed to write it as an integral of an algebraic function over an algebraic domain.
From the definition of the Feynman period, we can already say that an upper bound for its weight is $2\abs{E(G)} - 1 = 2\ell(G) - 1$.
We will discuss more about weight in Section~\ref{S:denominator}.

We end this section on $\phi^4$-theory and periods with a beautiful result~\cite{zigzag} that gives a formula for the period of an infinite family of graphs called \emph{zig-zag graphs}, denoted by $Z_{\ell}$ where $\ell$ is the loop order of the graph.
The formula was first conjectured by Broadhurst and Kreimer in~\cite{knots}.
This special family includes the complete graph $Z_3 = K_4$ and their completions form the circulant graphs $C_{\ell+2}(1,2)$.

\tpointn{Theorem} (Theorem 1.1 of~\cite{zigzag})\label{zigzag}
\statement[eq]{
  The period of the zig-zag graph $Z_{\ell}$ is given by
  \[ P_{Z_{\ell}} = \frac{4(2\ell - 2)!}{\ell!(\ell - 1)!}\left(1 - \frac{1 - (-1)^{\ell}}{2^{2\ell - 3}}\right) \zeta(2\ell - 3).\]
} \\

\section{Defining the \texorpdfstring{$c_2$}{c2}-invariant}\label{S:c2-def}

Looking at the period in its parametric form~\eqref{eq:period}, notice we are integrating over the denominator $\Psi_G^2$, which is just a polynomial in $\abs{E(G)}$ variables.
In particular, in order to understand and characterize properties of the period, we need to understand the structure of $\Psi_G$.
This motivates the study of the zeros of $\Psi_G$ and the polynomials (denominators) that arise after several steps of integration.

In 1997, inspired by the work of Broadhurst and Kreimer~\cite{knots,association}, Kontsevich informally conjectured that for all graphs, the function taking prime powers $q$ to the number of zeros (also known as point counts) of $\Psi_G$ over the finite finite $\mathbb{F}_q$ is a polynomial in $q$.
This stemmed from the appearance of multiple zeta values in certain $\phi^4$ periods, which also led to the conjecture that all $\phi^4$ periods were linear combinations of products of multiple zeta values.
As such, it was thought that there must have been some reason these special numbers were showing up in the period, namely due to specific structures~\cite{motives} of the algebraic varieties defined by $\Psi_G$, which in turn would imply Kontsevich's conjecture.
However, the conjecture turned out to be false in general.
Belkale and Brosnan~\cite{BB} first proved that Kontsevich's conjecture was generically false in 2000 and later Doryn~\cite{doryn} and Schnetz~\cite{fq} independently found counterexamples specifically in the class of $\phi^4$ graphs.
Brown and Schnetz~\cite{k3} then further showed that even under certain "physical" constraints on the graphs, Kontsevich's conjecture still could not hold.

Nevertheless, this point count function, which we denote by $[\Psi_G]_q$, still has an interesting connection to the period.
When Kontsevich's conjecture is indeed true for some graph $G$, then the coefficient of $q^2$ in $[\Psi_G]_q$, which is now a constant with respect to $q$, was found to give information on the weight of the multiple zeta values that appear in the period of this graph.
In~\cite{fq}, Schnetz introduced the following arithmetic invariant.

\tpointn{Definition} (Theorem $2.9$ of~\cite{fq})\label{def:c2}
\statement{
  Let $q = p^n$ be a prime power and $\mathbb{F}_q$ the finite field with $q$ elements.
  Let $G$ be a connected graph with at least $3$ vertices.
  Then the \textbf{$\mathbf{c_2}$-invariant} of $G$ at $q$ is
  \begin{equation}\label{eq:originalc2}
    c_2^{(q)}(G) \equiv \frac{\left[\Psi_G\right]_q}{q^2} \mod q,
  \end{equation}
  where $\left[\Psi_G\right]_q$ is the number of zeros of $\,\Psi_G$ in $\mathbb{F}_q^{\abs{E(G)}}$. \\
  Denote by $\mathbf{c_2(G)}$ the sequence of $c_2^{(q)}(G)$ for all prime powers $q$.
}

To show that this invariant is well-defined, Schnetz proved that $[\Psi_G]_q$ was indeed divisible by $q^2$.
In particular, the $c_2^{(q)}(G)$ would then exactly be the quadratic coefficient of $[\Psi_G]_q$ if it was indeed a polynomial in $q$, hence the name "$c_2$".
When specifically referring to the $c_2$ at primes $q=p$, we will often use $p$ instead of $q$.

In a sense, the $c_2$-invariant is measuring how badly Kontsevich's conjecture fails for a graph.
More importantly, other than its combinatorial flavour which we will see shortly, the interest in the $c_2$-invariant lies in how it seems to detect the types of numbers that appear in the period and thus telling us something about the geometries underlying them~\cite{k3,modular,geometries}.
In Section~\ref{S:families}, we give a brief overview of some of the results in this direction.

This connection between the $c_2$-invariant and the period is further strengthened by the following conjecture, which says that the $c_2$ is a period invariant.

\tpointn{Conjecture} (Remark 2.11 (2) of~\cite{fq}, Conjecture $2$ of ~\cite{modular})\label{c2-period-conj} \\
\statement[eq]{
  Let $G_1$ and $G_2$ be primitive divergent graphs.
  If their periods coincide, that is, $P_{G_1} = P_{G_2}$, then
  \[ c_2(G_1) = c_2(G_2). \]
}

Currently, this conjecture holds for all known examples.
However, whether all the symmetries of the period also hold for the $c_2$-invariant is still unknown, which we would need for the above conjecture to be true.
There has been some progress along this front: first shown in 2013, Doryn~\cite{isinv} proved that the $c_2$-invariant is indeed an invariant across the four representations of the period under a condition for the graph called \emph{duality admissibility}.
This condition encompasses all planar graphs, and thus the $c_2$ remains invariant under planar duality, as does the period.
For the completion symmetry, which for the $c_2$ we refer to as the \textbf{completion conjecture}, in 2018 Yeats~\cite{specialc2} first made headway on this in some special cases, and in Chapter~\ref{C:completion} we complete the argument for when $q=2$.

Finally, there wouldn't be much interest in the $c_2$-invariant if computationally it was as hard as the period.
In the following two sections, we will see how to transform the $c_2$ into a more tractable form, along the way illuminating its combinatorial nature.
In this transformed form, the $c_2$-invariant becomes relatively easy to calculate compared to the period, at least for low loop orders or small primes $q=p$, while still encapsulating many of the properties of the period and its underlying graph. \\

\section{Graph polynomials}\label{S:graph-polys}

Before we can work on massaging Equation~\eqref{eq:originalc2}, we need to familiarize ourselves with the function in the numerator, the graph polynomial.
Let $G$ be a connected graph.
Recall that the graph polynomial of $G$ is defined as
\[ \Psi_{G} = \sum_{\substack{T \\\text{spanning tree}}} \prod_{e \not\in T} \alpha_e, \] \\
where the $\alpha_e$ are parameters associated to each edge $e$ and the sum runs over all spanning trees of $G$.
Notice that $\Psi_G$ is a homogeneous polynomial of degree $|E(G)| - |V(G)| + 1$ in $|E(G)|$ variables.
There is also a nice deletion-contraction relation, where for any edge $e$ in $G$ we have \\
\begin{equation}\label{eq:del-cont-psi}
  \Psi_G = \alpha_e \Psi_{G \setminus e} + \Psi_{G / e},
\end{equation}\\
since we can partition the sum based on whether a spanning tree of $G$ contains $e$ or not.
To see this, notice that the set of spanning trees of $G$ not containing $e$ is equal to the set of spanning trees on $G \setminus e$ and thus $\alpha_e$ can be factored out of the monomials.
The set of spanning trees of $G$ containing $e$ is then equal to the set of spanning trees of $G / e$.

As the graphs get larger, these polynomials explode as well, and thus the first step in transforming the $c_2$ is to find smaller polynomials in which we can take point counts of.
To do this, we need the theory of some related graph polynomials, which following Brown we call \textbf{Dodgson polynomials}, and which stem from the deletion-contraction relation of the graph polynomial and Kirchhoff's matrix-tree theorem.

As first introduced in ~\cite{periods} in the context of calculating periods and the theory of \textbf{denominator reduction}, which we will cover in the next section, Dodgson polynomials were originally defined up to overall sign.
Up until very recently, in the subsequent work on the $c_2$-invariant this overall sign did not matter.
In 2019, Schnetz~\cite{geometries} re-derived much of the work on Dodgson polynomials to account for the correct signs, which we present here.
Note that in the $q=2$ case, which is our main focus in Chapter~\ref{C:completion}, these signs won't be needed. \\

\bpoint{Determinantal framework}

The underlying framework of this theory on graph polynomials relies on the fact that the graph polynomial can be seen as the determinant of a particular matrix.
To represent the graph polynomial as such, we first define the following.

\tpointn{Definition}\label{def:L}
\statement{
  Given a connected graph $G$, choose an arbitrary orientation on the edges and an order $\iota$ on the edges and vertices of $G$, where the edges are ordered before the vertices,
  \[ \iota : E(G) \cup V(G) \quad\longrightarrow\quad \{1, \ldots, \abs{E(G)}+\abs{V(G)} \}.\]
  We denote by $\iota_e$ or $\ \iota_v$, to mean the order of edge $e$ or vertex $v$ under $\iota$, respectively. \\\\
  Let $\cE_G$ be the $\abs{V(G)} \times \abs{E(G)}$ signed incidence matrix, with any one row (corresponding to a vertex) removed.
  That is, for all vertices $v$ except one, and all edges $e$
  \[ [\cE_G]_{v,e} = \begin{cases} 1 & \text{ if }\, e = (v,u), u \in V(G), \\ -1 & \text{ if }\, e = (u,v), u \in V(G), \\ 0 & \text{ otherwise. }\end{cases} \]
  Let $\Lambda$ be the diagonal matrix of indeterminates $\alpha_e$ for $e$ in $E(G)$ in the chosen order, that is
  \[ [\Lambda]_{e_i,e_j} = \begin{cases} \alpha_{e_i} & \text{ if }\, i=j, \\ 0 & \text{ otherwise. }\end{cases} \]
  Then we define the \textbf{expanded Laplacian} of $G$ to be
  \[ L_G = \left[
    \renewcommand\arraystretch{1.5}
    \begin{array}{c|c}
      \Lambda & {\cE_G}^T \\
      \hline
      \cE_G & 0 \\
    \end{array}
    \right],
  \]
  which is an $(\abs{E(G)}+ \abs{V(G)} - 1) \times (\abs{E(G)} + \abs{V(G)} - 1)$ matrix, with rows and columns ordered by $\iota$.
}

While this matrix is not well-defined, as it depends on the choice of row removed in $\cE_G$ as well as the choice of orderings and orientation, we have, for any such choice\\
\begin{equation}\label{eq:det-psi}
  \Psi_G = (-1)^{|V(G)|-1}\det\left(L_G\right).
\end{equation}

To see why this holds, following the proof in \textsection 2.2 of~\cite{periods}, we need a lemma due to Kirchhoff which we can think of as essentially the matrix-tree theorem.

\tpointn{Lemma} (\cite{kirchhoff})\label{matrixtree}
\statement[eq]{
  Let $U$ be a subset of edges of $G$ such that
  \[{\abs{E(G \setminus U)} = \ell(G) = \abs{E(G)} - \abs{V(G)} + 1}. \]
  Let $\cE_G(G\setminus U)$ denote the square $(\abs{V(G)} - 1) \times (\abs{V(G)} - 1)$ matrix obtained from $\cE_G$ by removing the columns indexed by the edges of $G \setminus U$ (recall that $\cE_G$ already has one row removed).
  Then
  \[ \det \cE_G(G\setminus U) =
  \begin{cases}
    \pm 1 & \text{if } U \text{ is a spanning tree of } G, \\
    0 & \text{otherwise.}
  \end{cases}\]
}

We also need the Leibniz formula for determinants of square $n\times n$ matrices $A$ with entries $a_{i,j}$
\[ \det(A) = \sum_{\sigma \in S_{n}} \sgn(\sigma) \prod_{i=1}^n a_{i, \sigma(i)}, \]
where $S_{n}$ is the symmetric group on $n$ elements and $\sgn(\sigma)$ is the sign of the permutation $\sigma$
\[ \sgn(\sigma) = (-1)^{N(\sigma)}. \]
Here $N(\sigma)$ is the number of inversions in $\sigma$, i.e. the number of pairs $(i,j)$ such that $i < j$, but $\sigma(i) > \sigma(j)$.

Then, expanding the determinant of $L_G$ via the Leibniz formula, we have that \\
\[ \det(L_G) = \sum_{U \subset E(G)} \det\left[
    \begin{array}{c|c}
      0 & {\cE_G(G\setminus U)}^T \\
      \hline
      \cE_G (G \setminus U) & 0 \\
    \end{array}
    \right] \prod_{e \not\in U} \alpha_e\
  = \hspace{-3.5pt}\sum_{\substack{U \subset E(G) \\ \abs{E(G \setminus U)} = \ell(G)}} \hspace{-3.5pt} (-1)^{\abs{V(G)}-1}\det\left(\cE_G(G \setminus U)\right)^2 \prod_{e \not\in U} \alpha_e. \]

The first equality comes from choosing $E(G \setminus U)$ to be the fixed points of permutation $\sigma$ and passing to the minor with the corresponding rows and columns removed, setting $\alpha_e=0$ for $e \in U$.
Now, the terms corresponding to any permutation with fixed points in the last $\abs{V(G)}-1$ columns vanish, and thus the only fixed points for non-zero terms in the sum are those in the first $\abs{E(G)}$ columns.
Additionally, the sign of a permutation $\sigma$ is the same as the sign of the permutation $\pi$ created by removing the fixed points of $\sigma$.
The second equality comes from a determinantal identity on block matrices and the restriction in the sum occurs as all the terms with $\abs{E(G \setminus U)} \neq \ell(G)$ vanish since the rank of $\cE_G$ is $\abs{V(G)} - 1$ (recall we are removing the columns corresponding to $E(G \setminus U)$).
Then, applying Lemma~\ref{matrixtree} and rearranging gives Equation~\eqref{eq:det-psi} since $\det(\cE_G(G \setminus U))^2 = 1$ if $U$ is a spanning tree of $G$. \\

\bpoint{Dodgson polynomials}\label{SS:dodgsons}

Under this determinantal framework and looking at the deletion-contraction relation~(\ref{eq:del-cont-psi}) for $\Psi_G$, we see that removing an edge $e$ to get $\Psi_{G\setminus e}$ corresponds to the minor where the row and column indexed by $e$ is removed.
This comes from noticing that $\Psi_{G\setminus e}$ is the coefficient of $\alpha_e$ in $\det(L_G)$.
Contracting an edge $e$ to get $\Psi_{G / e}$ then corresponds to setting $\alpha_e$ to zero.
In either case, some consideration is needed for the signs.
Motivated by this observation, we can extend the definition to minors of $L_G$ which we call \textbf{Dodgson polynomials}, or just Dodgsons.

Before defining Dodgsons, we need to take care of the signs.
As originally defined in \textsection 2.3 of~\cite{periods}, the edges corresponding to the rows and columns being removed from $L_G$ were taken as sets and thus the original Dodgsons were only defined up to overall sign, depending on the minor taken.
We refer to these as \textbf{unsigned Dodgsons}.

The main idea here is that to get the correct signs, we need to keep track of the order in which the rows and columns are being removed, as this affects the sign in the coefficient of the corresponding minor in the cofactor expansion of $L_G$ following this order.
In particular, this sign is determined by the indices of the rows and columns being removed and the orders of these removals.

To get the correct signs, we lift the sets of edges being removed to words on $E(G)$, where the letters are the edges ordered with respect to the chosen ordering $\iota$ from Definition~\ref{def:L}.
For convenience, we interchangeably use the terms edges and letters, as well as refer to edges $e$ by their numerical position $\iota_e$ as given by $\iota$.
We also use set notation with these words when we are referring to the set of edges in the words.
For such a word $w=w_1\cdots w_n$ in the letters $E(G)$, let $\iota_w$ be the sum of the orders $\iota_{w_j}$ of the letters $w_j$ of $w$ under $\iota$.
We also define
\[ \sgn(w) = \begin{cases}
      0 & \text{ if $w$ has repeated letters,} \\
\sgn(\pi) & \text{ otherwise, where $\pi(i) = j$ if $w_j$ is the $i$th smallest letter in $w$ with respect to $\iota$.} \end{cases} \]
As an example, the word $w=31$, where the edges are ordered naturally, has $\sgn(31) = -1$ as the permutation corresponding to $w$ is $\pi = 21$.
To add a letter $x$ to $w$, we have
\[ \sgn(wx) = (-1)^{N(\sigma)} = (-1)^{N(\pi) + \abs{\{j \;:\; i \;<\; j \;,\; \sigma(i) = n+1 \;>\; \sigma(j) \}}} = (-1)^{\abs{\{w_j \;\in\; w \;:\; x \;<\; w_j\}}}\sgn(w),  \]
where $\sigma$ is the permutation corresponding to the word $wx$.
The equalities hold as the inversions in $\sigma$ coming from $x$ are exactly those involving $\sigma^{-1}(n+1)$, which then correspond to the letters of $w$ larger than $x$ (with respect to $\iota$).
In other words, the sign of a word is the parity of how many transpositions it takes to get the edges back in their order under $\iota$.
Now we are ready to define our Dodgson polynomials.

\tpointn{Definition}\label{def:dodgson} (Definition 7 of~\cite{geometries})
\statement{
  Let $I$ and $J$ be words in the edges such that $\abs{I} = \abs{J}$ and let $K$ be a subset of edges of $G$.
  Denote $L_G{(I, J)}_K$ the matrix obtained from $L_G$ by removing rows indexed by $I$ and columns indexed by $J$, and setting $\alpha_e = 0$ for $e \in K$.
  Then the \textbf{Dodgson polynomial} is defined to be
  \[ \Psi_{G, K}^{I, J} = (-1)^{\abs{V(G)} + \iota_I + \iota_J -1} \sgn(I) \sgn(J) \det L_G(I, J)_K, \]
  We define $\Psi_{G,K}^{I,J} = 0$ if $\abs{I} \neq \abs{J}$ and define the empty determinant to be 1.
}

While there is a choice of ordering $\iota$ and vertex removed in $\cE_G$ for the matrix $L_G$, Dodgsons do not depend on this choice (see Lemma 9 of~\cite{geometries}).
When the graph $G$ is clear from the context, we will drop the subscript of $G$.
We also omit empty indices.
When $I, J, K$ are all empty, we recover $\Psi_G$.

Most, if not all, of the properties of unsigned Dodgsons (see \textsection 2.2, 2.3 of~\cite{periods}) also hold for Dodgsons, and in particular Dodgsons also satisfy a deletion-contraction relation.
The following proposition captures some of these properties, see \textsection 2 of~\cite{geometries}.

\tpointn{Proposition}\label{dodgson-props}
\statement[eq]{
  Let $I$ and $J$ be words in the edges such that $\abs{I} = \abs{J}$, and $K$ a subset of edges of $G$.
  \begin{itemize}
    \setlength\itemsep{0.7em}
    \item $\Psi^{I,J}_{G,K} = \Psi^{J,I}_{G,K}$.
    \item For any edge $e \in E(G) \setminus (I \cup J \cup K)$, there is a deletion-contraction relation
      \[ \Psi_{G, K}^{I, J} = \alpha_e \Psi_{G, K}^{Ie, Je} + \Psi_{G, K \cup e}^{I, J} = \alpha_e \Psi_{G \setminus e, K}^{I, J} + \Psi_{G/e, K}^{I,J}. \]
    \item It follows from the relation that we can always pass to a minor of $G$
      \[ \Psi_{G,K}^{I,J} = \Psi_{G'}^{I', J'}, \]
      where $G' = G \setminus (I \cap J) / (K \setminus (I \cup J))$, $I'=I\setminus (I \cap J)$, and $J'=J\setminus (I \cap J)$. \\
      By doing so, we can assume that $I \cap J = K = \emptyset$.
    \item If $\,\Psi_{G,K}^{I,J} \neq 0$, the degree of $\,\Psi_{G,K}^{I,J}$ is $\ell(G) - \abs{I}$. \\
  \end{itemize}
}

Lastly, like how the monomials of $\Psi_G$ correspond to spanning trees of $G$, there is a nice combinatorial interpretation of the monomials appearing in Dodgsons.
Once again, like for Equation~\eqref{eq:det-psi}, expanding the determinant of a minor of $L_G$ and using Lemma~\ref{matrixtree} gives the following theorem.

\tpointn{Theorem}\label{dodgsonsum} (Proposition 23 of~\cite{periods})
\statement{
  Suppose $I \cap J = K = \emptyset$.
  Then we have
  \[ \Psi_{G}^{I, J} = (-1)^{\iota_I + \iota_J} \sgn(I) \sgn(J) \sum_{U \subset G \setminus (I \cup J)}  \det\left( \cE_G( G \setminus (U \cup I)) \right) \det\left( \cE_G( G \setminus (U \cup J)) \right)\prod_{e \not\in U} \alpha_e, \]
  where the sum runs over all subgraphs $U$ such that $U \cup I$ and $U \cup J$ are both spanning trees in $G$.
}

Ignoring the signs, what this theorem is saying is that we can view Dodgson polynomials through the possible shapes of the underlying graph after some deletions and contractions.
For $U \cup I$ to be a spanning tree in $G$, we must have that $U$ is a spanning tree in $G \setminus J / I$ since $U$ are edges not in $I$ and $J$.
Similarly for $U \cup J$ to be a spanning tree, $U$ must be a spanning tree in $G \setminus I / J$.
The polynomial can then be thought of as an "intersection" of graphs $\left(G \setminus J / I\right) \cap \left(G \setminus I / J\right)$ where the graphs represent the spanning trees and $\cap$ is taken to mean the resulting polynomial of common terms which are spanning trees in each minor.
When $K \neq \emptyset$, this corresponds to contracting the edges in $K$ (that are not in $I, J$) in both graphs.
When $I \cap J \neq \emptyset$, this corresponds to removing the edges in $I \cap J$ in both graphs.\\

\bpoint{Dodgson identities}\label{SS:dodgson-ids}

In addition to having these nice properties for general Dodgsons, there is also a multitude of other relations falling into two types:
\begin{itemize}
  \item graph-specific identities,
  \item determinantal identities.
\end{itemize}

The first type stems from the combinatorial nature of Dodgsons, which are enumerating particular spanning trees of minors.
Thus we can get identities, sometimes also using determinantal identities, based on specific structures within a graph.
Some structures include cuts, cycles, vertex-joins, 3-valent vertices, and triangles, see ~\cite{periods,k3,properties,geometries}.

As a simple example, consider a graph with a 3-valent vertex $v$.
We immediately see that if we label the three adjacent edges of $v$ as $1,2,3$, then we have $\Psi^{123,123} = 0$, as removing all three edges from the graph disconnects $v$ and thus there are no possible spanning trees.
Another identity occurs by noticing that removing any two of the edges and contracting the third edge adjacent to $v$ results in the same underlying graph.
Thus we have the identity $\Psi^{12,12}_{3} = \Psi^{13,13}_{2} = \Psi^{23,23}_{1}$.

The second type stems from the determinantal framework of Dodgsons which allows us to use known determinantal identities and formulas in the context of Dodgsons.
Two determinantal identities that we will highlight here are what are called Pl{\"u}cker identities and Dodgson identities.

Based on the Pl{\"u}cker relations, we have the following formula (\textsection 2.4 of~\cite{periods}) which we call the \textbf{Pl{\"u}cker identities}
\[ \sum_{k=n}^{2n} (-1)^k \Psi_G^{(i_1\cdots i_{n-1}i_k), (i_n\cdots\hat{i}_k\cdots i_{2n})} = 0, \] \\
where $i_1 < \cdots < i_{2n}$ is an increasing sequence of edges and $\hat{i}_k$ means removing letter $i_k$.
For example, in the case when $n = 2$, we have
\[ \Psi_G^{12,34} - \Psi_G^{13,24} + \Psi_G^{14,23} = 0. \]

An important family of identities are what we call \textbf{Dodgson identities}, which are based on Jacobi's determinant formula and Dodgson condensation (\textsection 2.5 of~\cite{periods}, \textsection 2 of~\cite{geometries}).
One special case is due to the classical Dodgson identity
\[ \det(M)\det(M^{ij,ij}) = \det(M^{i,i})\det(M^{j,j}) - \det(M^{i,j})\det(M^{j,i}), \]
where $M$ is a square matrix, $i \neq j$ are indices and $M^{I,J}$ is the matrix obtained from $M$ by removing the rows indexed by $I$ and the columns indexed by $J$.
Translating this to our Dodgson polynomials gives the quadratic identity
\[ \Psi_{K}^{I, J} \Psi_K^{Ief,Jef} = \Psi_K^{Ie,Je} \Psi_K^{If,Jf} - \Psi_K^{Ie,Jf} \Psi_K^{If,Je}, \]\\
where $I,J$, are words of equal length in the edges, $K = K' \cup \{e,f\}$ for some subset of edges $K'$ and any two edges $e,f \not\in I \cup J \cup K'$.
In particular, when $I = J = K' = \emptyset$
\begin{equation}\label{eq:dodgson-id}
  \Psi_{12}\Psi^{12,12} = \Psi_2^{1,1}\Psi_1^{2,2} - \left(\Psi^{1,2}\right)^2,
\end{equation}
where we are using the property of Dodgsons that we can always pass to $K \setminus (I \cup J)$ in the subscript since the entries $\alpha_e$ for $e \in I \cup J$ are removed from the matrix via the removal of rows $I$ and columns $J$.
We can also think of Equation~\eqref{eq:dodgson-id} as taking a particular coefficient of the more general quadratic identity above.\\

\bpoint{Spanning forest polynomials}\label{SS:forest-polys}

Following the discussion in Section~\ref{SS:dodgsons}, we saw that we could interpret Dodgson polynomials via their underlying graphs.
However, we are still dealing with the common spanning trees between two possibly different graphs.
In an effort to better understand these polynomials combinatorially, Brown and Yeats introduced spanning forest polynomials~\cite{forest}.

\tpointn{Definition} (Definition 9 of~\cite{forest})\label{def:span-forest}
\statement{
  Let $P = P_1 \cup \cdots \cup P_k$ be a set partition of a subset of the vertices of $G$.
  Let $F$ be a spanning forest that partitions the vertices of $P$ exactly into $P_1 \cup \cdots \cup P_k$.
  More precisely, if $F = T_1 \cup \cdots \cup T_k$ then each tree $T_i$ of $F$ contains all the vertices in $P_i$ and no other vertices of $P$, and possibly other vertices in $V(G) \setminus P$.
  Then we say that the spanning forest $F$ is \textbf{compatible} with the vertex partition $P$.
  Note that the vertices not in $P$ can belong to any tree of $F$. \\\\
  We define a \textbf{spanning forest polynomial} of $G$ to be
  \[ \Phi_G^P = \sum_F \prod_{e \not\in F} \alpha_e, \]
  where the sum runs over all spanning forests $F$ that are compatible with $P$.
}

To represent spanning forest polynomials graphically, we can associate a shape or colour to each part of $P$ and
draw the vertices in $P$ accordingly on the graph $G$.

Like Dodgsons, spanning forest polynomials also satisfy deletion-contraction relations.

\tpointn{Proposition} (Proposition 10 of~\cite{forest})
\statement{
  Let $e$ be an edge in $G$ and let $P$ be a partition of a subset of vertices of $G$.
  Then we have
  \[ \Phi_G^P = \begin{cases}
    \alpha_e \Phi_{G \setminus e}^P & \text{if the ends of $e$ are in different parts of $P$,}\\
    \alpha_e \Phi_{G \setminus e}^P + \Phi_{G / e}^{P / e} & \text{otherwise,} \end{cases}\]
  where $P / e$ is the partition created from $P$ by identifying the ends of $e$ if they appear in $P$.
}

To relate spanning forest polynomials back to Dodgsons, instead of viewing Dodgsons through the matrix-tree theorem like in Theorem~\ref{dodgsonsum}, we can interpret them via the all-minors matrix-tree theorem~\cite{minorsmatrixtree}.
As a consequence, Dodgsons can be viewed as sums of spanning forest polynomials.

\tpointn{Theorem} (Proposition 12 of~\cite{forest})\label{tree-to-forest}
\statement{
  Let $I$, $J$ and $K$ be sets of edges of $G$ where $\abs{I} = \abs{J}$.
  Then we can write
  \[ \Psi_{G,K}^{I,J} = \sum f_P \Phi_{G \setminus (I \cup J \cup K)}^{P}, \]
  where $f_P \in \{-1,1\}$, and the sum runs over all partitions $P$ of the vertices in $(I \cup J \cup K) \setminus (I \cap J)$ such that all forests compatible with $P$ are spanning trees in both $G \setminus I / (J \cup K)$ and $G \setminus J / (I \cup K)$.
}

Following the proof in~\cite{forest}, comparing with Theorem~\ref{dodgsonsum} we can determine the sign $f_P$ via the appropriate determinants.
Passing to the case when $I \cap J = K = \emptyset$, notice from that theorem the subgraphs $U$ are exactly the forests in the spanning forest polynomials above since $U \cup I$ and $U \cup J$ must be spanning trees in $G$.
Equivalently, for a fixed forest $U$, this condition means that $I$ and $J$ must be spanning trees in $G \setminus \bar{U} / U$ where $\bar{U}$ are the edges in $G \setminus (I \cup J \cup U)$.

Now in $G \setminus \bar{U} / U$, each tree of $U$ has been contracted into a vertex and since the only edges left are those in $I \cup J$, each contracted vertex also tells us which vertices of $I \cup J$ are in the same tree of $U$.
That is, the vertices of $G \setminus \bar{U} / U$ partition the vertices in $I \cup J$ exactly into $P$, which is the partition that $U$ is compatible with.
In particular, any forest that is compatible with $P$ gives a graph isomorphic to $G \setminus \bar{U} / U$.
Thus, in $\Psi_G^{I,J}$ the monomials of $\Phi_G^P$ all have the same coefficient $f_P$ which is the following product of determinants for any forest $U$ compatible with $P$
\[\det\left( \cE_{G}( G \setminus (U \cup I)) \right) \det\left( \cE_{G}( G \setminus (U \cup J)) \right), \]
up to an overall sign.

Similarly to Dodgsons, there are many identities and relations that spanning forest polynomials satisfy, see \textsection 2, 3 of~\cite{forest}.
In particular, based on the Dodgson identity~(\ref{eq:dodgson-id}) we have the following identity which we represent graphically:
\[ \BlobDodgson \]
Here each blob represents a spanning forest polynomial, all with the same underlying graph, and the shapes (circle, square, and triangle) indicate the partitions of the three distinguished vertices for each polynomial.
Each term in the identity is a product of two spanning forest polynomials. \\

\section{Denominator reduction}\label{S:denominator}

With the theory of graph polynomials in place, we return to our main goal of transforming the definition of the $c_2$-invariant~(\ref{eq:originalc2}) into a more tractable form.
The main idea behind this is an algorithm called \textbf{denominator reduction} due to Brown~\cite{periods,massless}, which goes back to the period of a primitive divergent $\phi^4$ graph $G$ and looks at its denominators after several integrations.

Order the edges according to $\iota$ and let $N = \abs{E(G)}$.
Recall that the period is defined as
\[ P_G \coloneqq \int_0^{\infty} \frac{\;\mathrm{d}\alpha_1 \cdots \,\mathrm{d}\alpha_{N-1}\;}{\left(\Psi_G\right)^2 \,\bigg|_{\alpha_N=1}}. \]
Using the deletion-contraction relation~(\ref{eq:del-cont-psi}), Dodgson identities, and the theory of hyperlogarithms (see \textsection 10.2 of~\cite{periods}, \textsection 5 of~\cite{massless} and lecture notes on iterated integrals in this context~\cite{iteratedQFT}), we can explicitly compute the first few integrations, where we are integrating with respect to the first five edges.
The first two integrations are relatively simple, starting with integrating with respect to $\alpha_1$
\begin{align*}
  I_1
    &= \int_{\alpha_N=1} \frac{\;\mathrm{d}\alpha_1 \cdots \,\mathrm{d}\alpha_{N-1}\;}{(\alpha_1\Psi^{1,1} + \Psi_1)^2}
    = \int_{\alpha_N=1} -\frac{\;\mathrm{d}\alpha_2 \cdots \,\mathrm{d}\alpha_{N-1}\;}{(\alpha_1\Psi^{1,1} + \Psi_1)\Psi^{1,1}} \,\Bigg|_{\alpha_1 = 0}^{\infty}
    = \int_{\alpha_N=1} \frac{\;\mathrm{d}\alpha_2 \cdots \,\mathrm{d}\alpha_{N-1}\;}{\Psi^{1,1} \Psi_1}. \\
\end{align*}
For the second integration with respect to $\alpha_2$, we use the Dodgson identity~(\ref{eq:dodgson-id}) in the denominator
\begin{align*}
  \renewcommand\arraystretch{1.25}
  I_2
    &= \int_{\alpha_N=1} \frac{\;\mathrm{d}\alpha_2 \cdots \,\mathrm{d}\alpha_{N-1}\;}{(\alpha_2 \Psi^{12,12} + \Psi_{2}^{1,1}) (\alpha_2\Psi_1^{2,2} + \Psi_{12})}
    = \int_{\alpha_N=1} \frac{\log\left(\frac{\alpha_2\Psi^{12,12} + \Psi^{1,1}_2}{\alpha_2\Psi^{2,2}_1 + \Psi_{12}}\right)\,\bigg|_{\alpha_2 = 0}^{\infty}}{\Psi^{12,12}\Psi_{12} - \Psi_2^{1,1}\Psi_1^{2,2}}\;\mathrm{d}\alpha_3 \cdots \,\mathrm{d}\alpha_{N-1} \\
    &= \int_{\alpha_N=1} \frac{\log(\Psi^{1,1}_2) + \log(\Psi^{2,2}_1) - \log(\Psi^{12,12}) - \log(\Psi_{12})}{\left( \Psi^{1,2}\right)^2}\;\mathrm{d}\alpha_3 \cdots \,\mathrm{d}\alpha_{N-1}.
\end{align*}
Starting from here is where the integrations get more complicated, though reducing with respect to $\alpha_3$ is still doable.
By integration by parts, we first note that
\[ \int^{\infty}_0 \frac{\log(cx+d)}{(ax+b)^2}\;\mathrm{d}x = \frac{\log(d)}{ab} + \frac{c\left(\log(ab) - \log(cd)\right)}{a(ad-bc)}. \]
Applying this to $I_2$ after expanding the Dodgsons with respect to edge 3 gives the next integral
\begin{align*}
  I_3 &= \int_{\alpha_N=1} \frac{\log(\Psi^{1,1}_{23})+\log(\Psi^{2,2}_{13}) - \log(\Psi^{12,12}_3) - \log(\Psi_{123})}{\Psi^{13,23}\Psi^{1,2}_3} \;\mathrm{d}\alpha_4 \cdots \,\mathrm{d}\alpha_{N-1} + \cdots,
\end{align*}
where we've combined the first terms of integration of each of the four $\log$s in $I_2$ and omit the rest of the terms in the sums.
Another form of $I_3$ as given in \textsection 10.2 of~\cite{periods} is
\begin{align*}
  I_3 = \int_{\alpha_N=1} \frac{\Psi^{123,123}\log(\Psi^{123,123})}{\Psi^{12,13}\Psi^{12,23}\Psi^{13,23}} &- \frac{\Psi_{123}\log(\Psi_{123})}{\Psi_1^{2,3}\Psi_2^{1,3}\Psi_3^{1,2}} \\
    &+ \sum_{\{i,j,k\}}\frac{\Psi^{i,i}_{jk}\log(\Psi^{i,i}_{jk})}{\Psi^{ij,ik}\Psi^{i,k}_j\Psi^{i,j}_k} - \frac{\Psi^{ij,ij}_k\log(\Psi_{k}^{ij,ij})}{\Psi^{ij,ik}\Psi^{ij,jk}\Psi^{i,j}_k}\;\mathrm{d}\alpha_4 \cdots \,\mathrm{d}\alpha_{N-1},
\end{align*}
where the sum runs over the 6 permutations of $\{1,2,3\}$.
For the last two integrations, we only present the results, also from the same section of~\cite{periods}
\begin{align*}
  I_4 &= \int_{\alpha_N=1} \frac{A}{\Psi^{12,34}\Psi^{13,24}} + \frac{B}{\Psi^{13,24}\Psi^{14,23}} + \frac{C}{\Psi^{12,34}\Psi^{14,23}}\;\mathrm{d}\alpha_5 \cdots \,\mathrm{d}\alpha_{N-1}, \\
  I_5 &= \int_{\alpha_N=1} \frac{F}{^5\Psi(1,2,3,4,5)} \;\mathrm{d}\alpha_6 \cdots \,\mathrm{d}\alpha_{N-1}, \\
\end{align*}
where $A,B,C$ are di-logarithms, and $F$ is a hyperlogarithm of weight 3.
The denominator in $I_5$ is an important combination of Dodgson polynomials called the 5\nobreakdash-invariant, which is defined as follows.

\tpointn{Definition}\label{def:5-inv}
\statement[eq]{Given edges $1,\,\dots\,,5$ for a graph $G$, define the \textbf{5-invariant} of G, $\ ^5\Psi_G(1,\,\dots\, ,5)$, as
\[
  ^5\Psi_G (1,\,\dots\, ,5) = \pm \det \left(
  \renewcommand\arraystretch{1.25}
  \begin{array}{cc}
    \Psi_5^{12,34} & \Psi_5^{13,24} \\
    \Psi^{125,345} & \Psi^{135,245}\\
  \end{array}
  \right) = \pm \left(\Psi^{12,34}_5\Psi^{135,245} - \Psi_5^{13,24}\Psi^{125,345} \right).
\]
}

The 5-invariant is defined up to overall sign and holds for any five distinct edges.
Furthermore, permuting the order of the edges only changes the sign of $\ ^5\Psi_G(1,\,\dots\, ,5)$, see Lemma 87 in ~\cite{periods}.

The things to notice in these integrations are that:
\begin{itemize}[itemsep=4pt]
  \item The denominators are specific products of Dodgsons.
    Specifically the next denominator is the resultant of the two factors of the previous denominator (the determinant of their associated Sylvester matrix), assuming it is not a square.
  \item The weight of the numerators, in the sense of MZVs and hyperlogarithms, increases by exactly 1 when the previous denominator is not a square, i.e. $1 \to 1 \to \log \to \log \to \text{di-log} \to \text{weight 3}$.
  \item After the fifth integration, the denominator is in a special form.
\end{itemize}

These observations form the basis of denominator reduction, which in the context of the $c_2$-invariant allows us to reduce $\Psi_G$ with respect to a sequence of edges.  \\

\bpoint{Reducing variables}

As we noticed above, after the fifth integration there was a special denominator which is the 5-invariant.
Thus using the theory of hyperlogarithms, when possible, we can then completely determine the denominators of successive integrations, see Corollary 126 of~\cite{periods}.
This gives the following algorithm called \textbf{denominator reduction}, see \textsection 10 of~\cite{periods}, which gives higher invariants.

\tpointn{Definition}\label{def:denom-red} (Proposition 130 of~\cite{periods})
\statement{
  Given a graph $G$ and a sequence of edges $1, \ldots, \abs{E(G)}$, we define
  \[ D^5_G = {^5\Psi}_G(1,2,3,4,5). \]
  For $n > 5$, we recursively define the \textbf{n-invariant} $D^n_G(1,\ldots,n)$, also referred to as $\ ^n\Psi_G(1,\ldots,n)$, as follows:
  Suppose $D^n_G$ for $n \geq 5$ factors into the following product of linear factors in $\alpha_{n+1}$
  \[ D^n_G(1,\ldots,n) = (A\alpha_{n+1} + B)(C\alpha_{n+1} + D). \]
  Then we define
  \[ D^{n+1}_G(1,\ldots,n+1) = \pm (AD - BC), \]
  which is the resultant of the two factors of $D_n$.
  Otherwise, if $D^{n+1}_G = 0$ or $D^n_G$ cannot be factored in such a form (in which case the $n+1$-invariant and higher do not exist), we say that \textbf{denominator reduction} ends.
  If $D^{n+1}_G = 0$ (and thus also for higher invariants) for some sequence of edges,  we say that $G$ has \textbf{weight drop}.
  Note that in particular, when $D^n_G$ is a perfect square, then $G$ has weight drop.
}

Like the 5-invariant, these higher invariants (if they exist) are defined up to overall sign and do not depend on the order of edges up to that point.
However, a different sequence of edges may give a sequence of denominators of a different length.
As these $D^n_G$ are indeed the denominators in the period after integrating out the $n$ edges variables, we can say we are "reducing" with respect to the $n$ edges.
Note that the factors in the denominators are usually Dodgsons or combinations of Dodgsons.

The notion of weight drop for a graph $G$ directly relates back to the weight of its Feynman period $P_G$, which is the minimum number of nested integrals needed to write $P_G$ as a period.
As the number of edges in $\phi^4$ graphs is $\abs{E(G)} = 2\ell$ and we have set one edge variable to 1, we start with $2\ell-1$ integrals.
From the initial integrations above, we notice that the weight of the numerators does not increase in two of the integrations.
Thus the maximum weight of a $\phi^4$ period is $2\ell - 3$.
Note that this maximum is frequently achieved like in the case of the zig-zag graphs, see Theorem~\ref{zigzag}.
Now, when $G$ has a weight drop, this means that the weight of its period $P_G$ will be less than this maximum weight of $2\ell - 3$.
Hence denominator reduction, and thus the $c_2$-invariant, can tell us about the weight of a graph's period.
In Section~\ref{SS:trivialc2}, we discuss weight drop in the context of the $c_2$-invariant.

We can also extend the definition of the $n$-invariant to $n=3,4$ at the expense of uniqueness, where now the denominator depends on the order of the edges chosen.
For $n=3$ we have three choices for $D^3_G$ and taking the edges $1,2,3$ in that order we can define
\[ D^3_G(1,2,3) = \pm \Psi_G^{13,23}\Psi^{1,2}_{G,3}. \]
Similarly for $n=4$ we also have three choices for $D^4_G$, as seen in $I_4$, but notice we can apply denominator reduction to the $D^3_G$ above and use a Dodgson identity to get one such choice
\[ D^4_G(1,2,3,4) = \pm \Psi_G^{13,24}\Psi_G^{14,23}. \]
As we can denominator reduce any choice of $D^3_G$ or $D^4_G$ to the 5-invariant, these choices are sufficient.

This brings us to the first equivalent definition of the $c_2$-invariant where simultaneously we can get rid of dividing by $q^2$ and reduce the number of variables in the polynomial we need take point counts of.

\tpointn{Theorem} (Theorem $29$ of ~\cite{k3} with Corollary 28 of ~\cite{k3} for $n<5$)\label{denomred}
\statement{
  Let $G$ be a connected graph with $2\ell(G) \leq \abs{E(G)}$ and $\abs{E(G)} \geq 3$.
  Suppose that $D^n_G(e_1,\,\dots\, ,\, e_n)$ is the result of the denominator reduction after $3 \leq n < \abs{E(G)}$ steps.
  Then
  \[ c_2^{(q)}(G) \equiv (-1)^n \left[ D^n_G(e_1,\,\dots\, ,\, e_n)\right]_q \mod q.\]
  If $\ 2\ell(G) < \abs{E(G)} \geq 4$, then $c_2^{(q)} \equiv 0 \mod q$.
}

In particular, for any three distinct edges labelled $1,2,3$ of $G$ we have
\[ c_2^{(q)}(G) \equiv - \left[ \Psi_G^{13,23} \Psi_{G,3}^{1,2} \right]_q \mod q.\]

This theorem becomes a powerful tool in computing the $c_2$-invariant, at least for small $q$, as we can usually reduce via denominator reduction for many steps before needing to take point counts.
One strategy then becomes to find "good" sequences of edges in which we can denominator reduce for as long as possible.
In~\cite{modular}, Brown and Schnetz was able to exhaustively compute the $c_2^{(p)}$ for all $\phi^4$ graphs up to loop order 10 and for small primes $p$ using these methods. \\

\bpoint{Becoming a counting problem}\label{SS:denom-count}

While Theorem~\ref{denomred} brings us to a much more tractable form of the $c_2$-invariant, we are still left with the problem of point counting, which itself is a hard problem as the denominators can still have many variables.
To tackle this we turn to a consequence of a theorem from number theory, see \textsection 2 of~\cite{ax}, which will allow us to determine coefficients instead of counting zeros.

\tpointn{Theorem} (Corollary of Chevalley-Warning Theorem)\label{chevalley}
\statement{
  Let $F$ be a polynomial of degree $N$ in $N$ variables, $x_1,\,\dots\, ,x_N$, with integer coefficients. Then we have
  \[ \text{coefficient of } x_1^{q-1}\,\cdots\, x_N^{q-1} \text{ in } F^{q-1} \equiv (-1)^{N-1}[F]_q \mod p, \]
  where $q$ is some power of the prime $p$.
}

Note that this theorem holds when reducing modulo prime $p$ and not when reducing modulo $q$.
An example of this is given in \textsection 5 of~\cite{geometries}.
Thus, to use this theorem in the context of the $c_2$-invariant, we will restrict to when $q=p$ is a prime (see also Conjecture~\ref{prime-powers}).

To apply this theorem to $D^n_G$, we just need to make sure that the degree of the polynomial matches up with the number of variables.
We verify this recursively, recalling that $\abs{E(G)} = 2\ell(G)$ for our graphs of interest (primitive divergent graphs):
\begin{itemize}[itemsep=4pt]
  \item Starting with $\Psi_G^2$, as the monomials in $\Psi_G$ correspond to the complement edges of spanning trees, $\Psi_G^2$ is a homogeneous polynomial of degree $2(\abs{E(G)} - \abs{V(G)} + 1) = 2\ell(G) = \abs{E(G)}$.
  \item As $\Psi_G = \alpha_1\Psi^{1,1} + \Psi_1$, we have that $\Psi^{1,1}$ must be a homogeneous polynomial of degree $\ell(G) - 1$ and $\Psi^{1,1}$ is of degree $\ell(G)$.
    This means that $\Psi^{1,1}\Psi_1$ has degree $2\ell(G) - 1 = \abs{E(G)} - 1$ which is exactly the number of variables left after reducing with respect to edge 1.
  \item As $\Psi^{1,1} = \alpha_2 \Psi^{12,12} + \Psi^{1,1}_2$ and $\Psi_1 = \alpha_2 \Psi^{2,2}_1 + \Psi_{12}$ we have that $\Psi^{12,12}\Psi_{12}$ and $\Psi^{1,1}_2\Psi_1^{2,2}$ are both polynomials of degree $\ell(G) - 1 + \ell(G) - 1 = \abs{E(G)} - 2$.
    Thus the second denominator is of degree $\abs{E(G)} - 2$.
  \item Following a similar process for $D^3_G$, we get a degree of $\abs{E(G)} - 3$.
  \item For $n \geq 3$, suppose $D^n_G$ is a homogeneous polynomial of degree $\abs{E(G)} - n$ in $\abs{E(G)} - n$ variables.
    If $D^{n}_G$ factors as $(A\alpha_{n+1} + B)(C\alpha_{n+1} + D)$ then we must have that $AD$ and $BC$ are both homogeneous polynomials of degree $\abs{E(G)} - n - 1$.
    Thus $D^{n+1}_G = \pm (AD - BC)$ must be homogeneous of degree $\abs{E(G)} - (n+1)$ and in $\abs{E(G)} - (n+1)$ variables.
\end{itemize}

Thus, $D^n_G(e_1,\,\dots\, ,e_n)$ satisfies the criterion for Theorem~\ref{chevalley} for $n \geq 3$, where $N = \abs{E(G)} - n$.
Applying this to Theorem~\ref{denomred} and returning the notation $[\cdot]$ back to meaning taking coefficients, we arrive at the second transformation for the $c_2$-invariant.

\tpointn{Corollary}
\statement[eq]{
  Let $G$ be a connected graph with $2\ell(G) \leq \abs{E(G)}$ and $\abs{E(G)} \geq 3$.
  Suppose that $D^n_G(e_1,\,\dots\, ,\, e_n)$ is the result of the denominator reduction after $3 \leq n < \abs{E(G)}$ steps.
  Then
  \[ c_2^{(p)}(G) \equiv (-1)^{\abs{E(G)}-1} [ \alpha_{n+1}^{p-1} \cdots \alpha_{\abs{E(G)}}^{p-1} ]\; \left(D^n_G(e_1,\,\dots\, ,\, e_n)\right)^{p-1} \mod p.\]
}

In particular, for any three distinct edges labelled $1,2,3$ of $G$ and when $\abs{E(G)} = 2\ell(G)$ we have
\begin{equation}\label{eq:c2-version}
  c_2^{(p)}(G) \equiv -[\alpha_4^{p-1}\cdots\alpha_{|E(G)|}^{p-1}]\; \left(\Psi_G^{13,23} \Psi_{G,3}^{1,2} \right)^{p-1} \mod p.
\end{equation} \\
An example of using this equation to calculate the $c_2$-invariant is given in Example~\ref{ex:c2-calc}.

What's nice about this version of the $c_2$-invariant is that we have reformulated what was a problem on counting zeros of polynomials into a problem about finding the coefficient of a particular monomial in some power of said polynomials.
Furthermore, these polynomials inherently have a combinatorial nature as they arise from spanning trees or spanning forests of particular graph minors.
Thus, we can interpret this coefficient as counting the number of ways to distribute edges into some of these spanning trees or spanning forests, which is a purely combinatorial problem!
We will expand on this interpretation in Chapter~\ref{C:completion}.

Using this method with the theory of spanning forest polynomials and combinatorial techniques, Yeats independently and also with Chorney were able to compute the $c_2^{(p)}$ for entire families of graphs (those which are called \emph{recursively constructible}) and fixed small primes $p$~\cite{circulants,prefixes,families}.
Some examples of these families include particular circulants and toroidal grids.

Combining the two methods of finding "good" sequences of edges for denominator reduction and taking coefficients, Shaw et al.~\cite{further} were able to extend Brown and Schnetz's computations of $c_2^{(p)}$ to all $\phi^4$ graphs up to loop order 11 and $p=7$, with many graphs up to $p=13$. \\\\

\bpoint{Quadratic denominator reduction}\label{SS:qdr}

Up until now, most of the results on the $c_2$-invariant did not require us to know the specific signs of Dodgsons.
It was not until ~\cite{further} that the question of signs became needed, where point counts of sums of Dodgsons were being taken.
Furthermore, up until now the main tool in computing the $c_2$-invariant has been denominator reduction.
As we saw, denominator reduction only guarantees up to the first five integrations.
After the fifth integration, the use of denominator reduction requires the previous denominator to be factored into two different linear factors.
If this does not hold or we have a perfect square, we cannot reduce any further.
Often for high loop ordered graphs, even after denominator reduction, point counting would be too computationally intensive.
In~\cite{geometries}, Schnetz addresses both of these issues; fixing a sign convention for Dodgsons as we presented in Section~\ref{SS:dodgsons} and introducing an improved reduction called \textbf{quadratic denominator reduction}.

The main idea of quadratic denominator reduction is that instead of using point counts directly, we can use sums of \textbf{Legendre symbols}, which enable us to move to quadratic factors instead of linear factors.
We can then define a similar reduction algorithm for particular forms of these quadratic factors.
This idea stemmed from the observation that if the denominator $D^n_G$ does not factor into linear factors and is a general quadratic in variable $\alpha_{n+1}$, integrating with respect to $\alpha_{n+1}$ produces a square root of its discriminant in the next denominator.
Notice that in the denominator reduction case, this fact is still true, except now the discriminant is a perfect square (hence the ambiguity of the sign)!
Thus quadratic denominator reduction is generalizing the standard denominator reduction.

The Legendre symbol for any $a \in \mathbb{F}_q$ where $q$ is an odd prime power is defined as
\[ \left(\frac{a}{q}\right) = \abs{\{x \in \mathbb{F}_q \,:\, x^2 = a \}} - 1  \ \in \{-1,0,1\}. \]
For any polynomial $F \in \mathbb{Z}[\alpha_1, \ldots, \alpha_N]$, we can define
\[ (F)_q = \sum_{\alpha \in \mathbb{F}_q^N} \left(\frac{F(\alpha)}{q} \right) \ \in \mathbb{Z}. \]
Then the point count $[F]_q$ can be expressed as
\[ [F]_q = q^N - (F^2)_q  \equiv -(F^2)_q \mod q. \]
Note that we can also extend these definitions and results to when $q=2$ (Remark 3.7 of~\cite{geometries}). For more details see \textsection 5 of~\cite{geometries}.

Relating this back to the $c_2$, we are now reducing with respect to $\left(\Psi_G\right)^4$.
Based on the results of an additional integration after denominator reduction ends, Schnetz found two cases of denominators which behave nicely, with respect to the geometry of the varieties, under such an integration, see \textsection 6 of~\cite{geometries}.

\tpointn{Definition}\label{def:qdr} (Definition 34 of~\cite{geometries})
\statement{
  Given a graph $G$ with at least three edges and a sequence of edges $1, \ldots, \abs{E(G)}$, we define
  \[ ^3\Psi^2_G(1,2,3) = \left(\pm {^3}\Psi_G(1,2,3)\right)^2 = \left( \Psi^{13,23}_G\Psi^{1,2}_{G,3}\right)^2. \]
  For $n \geq 3 $, we recursively define the \textbf{quadratic n-invariant} $\ ^n\Psi^2_G(1,\ldots,n)$, as follows. \\
  If $\ {^n}\Psi^2_G$ for $n \geq 3$ is of the form
  \[ ^n\Psi^2_G(1,\ldots,n) = (A\alpha^2_{n+1} + B\alpha_{n+1} + C)^2, \]
  we define
  \[ ^{n+1}\Psi_G(1,\ldots,n+1) = B^2 - 4AC. \]
  If $\ {^n}\Psi^2_G$ for $n \geq 3$ is of the form
  \[ ^n\Psi^2_G(1,\ldots,n) = (D\alpha^2_{n+1} + E\alpha_{n+1} + F)(H\alpha_{n+1} + J)^2, \]
  we define
  \[ ^{n+1}\Psi_G(1,\ldots,n+1) = DJ^2 - EHJ + FH^2. \]
  Otherwise, if $\ ^{n+1}\Psi^2_G = 0$ or $\ ^n\Psi^2_G$ cannot be factored in such a form (in which case the quadratic $n+1$-invariant and higher do not exist), we say that \textbf{quadratic denominator reduction} ends at step $n$.
  If $\ ^{n}\Psi^2_G = 0$ (and thus also for higher invariants), we say that $G$ has \textbf{weight drop}.
}

Note that for quadratic denominator reduction, if $\ ^n\Psi^2_G$ factors into squares of linear factors, the signs are completely determined, and both cases reduce to the standard denominator reduction.
That is, when $D^n_G$ exists, we have that $D^n_G = \pm ({^n}\Psi_G^2)^{\frac{1}{2}}$.

Using quadratic invariants, we can once again reformulate the $c_2$-invariant.

\tpointn{Theorem} (Theorem $36$ of ~\cite{geometries})\label{qdr}
\statement{
  Let $q$ be an odd prime power.
  Suppose $G$ is a connected graph with $2\ell(G) \leq \abs{E(G)}$ and $\abs{E(G)} \geq 3$.
  Then
  \[ c_2^{(q)}(G) \equiv (-1)^{n-1} \left( ^n\Psi^2_G(e_1,\,\dots\, ,\, e_n)\right)_q \mod q,\]
  whenever $\ ^n\Psi^2_G$ exists.
}

The power of quadratic denominator reduction allowed Schnetz to compute and extend previous results on initial prime sequences for the $c_2^{(p)}$ for all $\phi^4$ graphs up to loop order 11 and partially to loop orders 12 and 13, see \textsection 10 of~\cite{geometries}.
This comes from the fact that we can always do a minimum of nine initial reductions (and in numerous graphs we can do many more) using quadratic denominator reduction as compared to the previous seven for standard denominator reduction, see \textsection 8 of~\cite{geometries}.
Furthermore, using this new technique, Schnetz and Yeats were able to determine a formula for a special family of infinite graphs~\cite{hourglass}, which we briefly describe in Section~\ref{SS:hourglass}. \\

\section{Families of graphs and numerical results}\label{S:families}

As an avatar for the period, continued interest in the $c_2$-invariant mostly falls into two categories which are deeply intertwined.
The first, based on its relative ease of computability, lies in the types of sequences that can arise as a $c_2$ and what that means for the geometries of the varieties underlying the corresponding periods.
This direction then motivates us to come up with better algorithms for computing the $c_2$-invariant that can handle higher values of $q$ and larger graphs.
In Section~\ref{S:denominator}, we introduced two algorithms, denominator reduction and quadratic denominator reduction, which currently form the basis for most of the computational results.

The second category is based on finding symmetries and properties of the $c_2$-invariant.
For example, resolving the period symmetries from Section~\ref{SS:symmetries} in the context of the $c_2$ or finding graph transforms that preserve the $c_2$, both of which gives families of graphs with the same $c_2$.
On the one hand, these aid computational results and the discovery of better algorithms as they allow us to either reduce the number of graphs we need to calculate the $c_2$ for or extend computations by transforming the problem to a simpler form (with respect to computability).
On the other hand, many symmetries and properties of the $c_2$-invariant have an inherently combinatorial flavour and thus are interesting combinatorial problems in their own right.
These combinatorial problems then also use known algorithms to solve some aspect of the problem or even to initially get the problem into a combinatorial form.
In Section~\ref{S:graph-polys}, we introduced the notion of Dodgsons and spanning forest polynomials which allows us to interpret the $c_2$ through spanning trees and spanning forests.
Chapter~\ref{C:completion} uses this interpretation to solve the $c_2$ completion conjecture when $q=2$.

In this final section, we present some of the results in both of these categories: the types of sequences that appear as a $c_2$, and graph transforms that give infinite families of graphs with the same $c_2$.

First, we note that we will restrict the $c_2$ sequence to all primes $q=p$ instead of prime powers $q=p^n$.
This is based on the following conjecture by Schnetz which says that prime powers do not give any extra information.

\tpointn{Conjecture} (Conjecture 2 of~\cite{geometries})\label{prime-powers}
\statement[eq]{
  Let $G_1$ and $G_2$ be two graphs with equivalent $c_2^{(p)}$ for all primes $p$, that is $c_2^{(p)}(G_1) \equiv c_2^{(p)}(G_2) \mod p$.
  Then for all prime powers $q = p^n$
  \[ c_2^{(q)}(G_1) \equiv c_2^{(q)}(G_2) \mod q. \]
}

Part of the motivation for this conjecture (other than for computational reasons) came from the use of only the primes to identify certain $c_2$ sequences, see Section~\ref{SS:numerics}.
Some evidence for this conjecture was also given in~\cite{fq} for a few small prime powers.
Very recently, Esipova and Yeats ~\cite{maria} were able to prove a version of this conjecture where instead of working modulo $q$ they worked modulo $p$.
Specifically, they proved that $c_2^{(q)}(G_1) \equiv c_2^{(q)}(G_2) \mod p$.

Secondly, we note that in the context of numerical results on $\phi^4$ graphs, we normally assume the completion conjecture for the $c_2$-invariant which says that for any completed primitive graph $G$, decompleting at any vertex $v$ gives the same $c_2$ (see Section~\ref{SS:symmetries}).
Finally, when the $c_2$ is a constant sequence consisting of all $c$'s, we may say that $c_2 = c$.

For the rest of this section, $G$ will be a connected graph with $2\ell(G) \leq \abs{E(G)}$ and at least three edges. \\

\bpoint{Weight drop, subdivergences and the trivial \texorpdfstring{$c_2$}{c2}}\label{SS:trivialc2}

The first families of graphs that we will consider is when the $c_2$ is trivial.
These are graphs $G$ with $c_2^{(p)}(G) \equiv 0 \mod p$ for all primes $p$.

One type of family is based on the denominator reduction algorithm (Definition~\ref{def:denom-red}) which tells us that when $G$ has weight drop (i.e. when the denominator is 0 after reducing with respect to some sequence of edges), by Theorem~\ref{denomred} we have that
\[ c_2^{(p)}(G) \equiv 0 \mod p. \]
Thus, any graph $G$ that has weight drop will have a trivial $c_2$.
Intuitively, recall that for a graph the weight of the (possibly transcendental) number that is its period can be thought of as the minimum number of nested integrals needed to write the period.
From denominator reduction, we saw that the maximal weight of a period is $2\ell(G) - 3$, and having weight drop then means that the weight is actually less than this maximal value.

Two graph properties that give rise to weight drop are:
\begin{itemize}
  \item doubled edges in graphs with at least five edges (Lemma 92 of~\cite{periods}), and
  \item 2-vertex reducible (Proposition 36 of~\cite{forest}): These are graphs which have a 2-vertex cut, in other words, they have vertex connectivity 2. For completed primitive graphs, the corresponding property is what we call reducible, having 3-vertex cuts.
\end{itemize}
Note that being (2- or 3-vertex) reducible and having $c_2 = 0$ is the $c_2$-invariant equivalent of the product identity (see Section~\ref{SS:symmetries}) for the period!

The second type of graph we will look at that has trivial $c_2$ are those from $\phi^4$-theory which are not primitive but still logarithmic divergent.
Recall from Definition~\ref{def:primitive} that these graphs $G$ contain non-trivial divergent subgraphs $\gamma$ which are subgraphs such that $\abs{E(\gamma)} \leq 2\ell(\gamma)$.
To translate this condition into an equivalent graph property, we think of $\gamma$ as being 4-regular with $k$ external legs giving
\[ \abs{E(\gamma)} = \frac{4\abs{V(\gamma)} - k}{2} \leq 2\ell(\gamma) \quad\implies\quad \abs{V(\gamma)} \leq \ell(\gamma) + \frac{k}{4}. \]
Then, using Euler's formula with the divergent subgraph condition, we have
\[ \abs{E(\gamma)} = \abs{V(\gamma)} + \ell(\gamma) - 1 \leq 2\ell(\gamma) + \frac{k}{4} - 1, \]
which tells us that $k \leq 4$.
Now, thinking of these $k$ external legs in the context of the full graph $G$, this means that $G$ has a cut with $k$ edges where one side of the cut is the subgraph $\gamma$.
Thus, having subdivergences corresponds to having non-trivial cuts of size at most 4.
In \textsection 5 of~\cite{properties}, Brown, Schnetz and Yeats were able to prove that if a graph has a non-trivial cut of size at most 4, then its $c_2$ is trivial!

Note that after completing $G$ (see Section~\ref{SS:symmetries}), the same argument still holds for any subdivergence $\gamma$.
This shows why we need the completed primitive condition (being internally 6-edge connected) for 4-regular graphs to ensure we get primitive divergent decompletions. \\

\bpoint{Double triangle reduction}\label{SS:DTR}

Now that we know when graphs have $c_2 = 0$, we switch our focus to a graph transformation called the \textbf{double triangle reduction} (DTR).
While double triangle reduction does not preserve the period, it does preserve the $c_2$-invariant.
Thus, having this transformation allows us to reduce our problems to smaller graphs and so we only need to look at the $c_2$'s of those without double triangles.

\begin{figure}[ht]
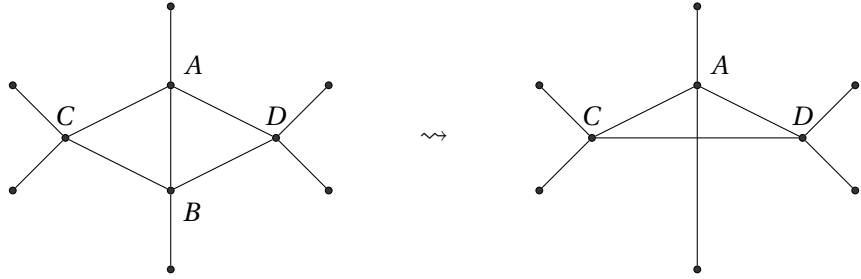

  \centering
  \DTR
  \caption[Double triangle reduction.]{Double triangle reduction from left to right. \\
  Note: $A, B, C$ and $D$ do not need to be 4-valent.}
  \label{fig:DTR}
\end{figure}

Suppose a graph $G$ has an edge that is shared by exactly two triangles.
Call this edge $(A,B)$ with triangles $(A,B,C)$ and $(A,B,D)$.
The double triangle reduced graph of $G$ is $G$ with one of the vertices of $(A,B)$, say $B$, replaced with the edge $(C,D)$.
If $B$ has a neighbour not in the triangles, then it is now adjacent to the remaining vertex $A$.
This is illustrated in Figure~\ref{fig:DTR}.
We note that in fact $A, B, C$ and $D$ do not need to be 4-valent.

\tpointn{Theorem} (Theorem 3.5 of~\cite{further} using Corollary $34$ in ~\cite{k3}, Theorem $35$ in ~\cite{forest})\label{DTRinvar} \\
\statement[eq]{
  Let $G$ be a connected 4-regular graph and $G'$ be the double triangle reduction of $G$.
  Then, decompleting at any (same) vertex $v$, we get
  \[ c_2(G - v) = c_2(G' - v). \] \\
  In particular, when $G_1$ is a primitive divergent graph in $\phi^4$ and $G'_1$ is its double triangle reduction, then
  \[ c_2(G_1) = c_2(G'_1). \]
} \\

\bpoint{\texorpdfstring{$K_5$}{K5} and its descendants}

Now, assuming the completion conjecture for the $c_2$-invariant allows us to only look at completed primitive graphs.
The results in Section~\ref{SS:trivialc2} says that we can then disregard those graphs with vertex connectivity $< 4$ and with non-trivial 4-edge cuts since they yield trivial $c_2$'s.
Finally, double triangle reduction further reduces our graphs to those without double triangles.
Combining all three reductions, we are led to the following type of graphs.

\tpointn{Definition} (Definition 22 of~\cite{geometries})
\statement{
  A graph with $\geq 5$ vertices is a \textbf{prime ancestor} if all of the following properties hold:
  \begin{enumerate}
    \item it is 4-regular,
    \item it is internally 6-edge connected,
    \item it has vertex connectivity 4, and
    \item it has no edge which is shared by exactly two triangles.
  \end{enumerate}
}

Thus, we only need to investigate the $c_2$-invariant for decompletions of prime ancestors as each prime ancestor defines an equivalence class of graphs all with the same $c_2$.
In each such class, the graphs are related by double triangle reductions and thus we call these graphs \textbf{descendants} of the ancestor.
There is one prime ancestor of particular interest which is conjectured to exactly correspond to the class of graphs with $c_2 = -1$.

\tpointn{Conjecture} (Conjecture 25 of~\cite{modular})
\statement{
  The only $\phi^4$ prime ancestor with $c_2^{(q)} \equiv -1 \mod q$ for all $q$ is the complete graph $K_5$.
}

In~\cite{double} Laradji, Mishna and Yeats studied the structure of $K_5$ descendants and gave some partial results on properties of these graphs.
One interesting result is that the minimum number of triangles in a $K_5$ descendant is four.

We verify here that decompletions of $K_5$ indeed have a $c_2$ of $-1$, giving an example of a $c_2$ calculation.

\tpointn{Example} (A $c_2$ calculation)\label{ex:c2-calc}
\statement{
  Consider the complete graph $K_4$ on four vertices whose completion is $K_5$, with the edges labelled as follows:
  \[ \GrphThreeEL \]
  To compute it's $c_2$-invariant, we use Equation~\eqref{eq:c2-version}.
  First we compute the two Dodgsons in the equation with respect to $K_4$ using Theorem~\ref{dodgsonsum} and the discussion thereafter.
  For $\Psi^{12,13}$ notice that $K_4 \setminus 12 / 3$ is isomorphic to $K_4 \setminus 13 / 2$ which is the triangle with edges $4,5,6$ and thus
  \[ \Psi^{12,13} = \Psi\left( \raisebox{-9mm}[10mm][10mm]{\Triangle} \right) = \alpha_4 + \alpha_5 + \alpha_6. \]
  For $\Psi^{2,3}_1$ we compute the common spanning trees of $K_4 \setminus 2 / 13$ and $K_4 \setminus 3 / 12$, which is just the edge $5$,
  \[ \Psi^{2,3}_1 = \Psi\left( \raisebox{-7mm}[7mm][7mm]{\LoopOne} \cap \raisebox{-7mm}[7mm][7mm]{\LoopTwo} \right) = \alpha_4\alpha_6. \]
  Putting everything together we get
  \begin{align*}
    c_2^{(q)}(K_4)
      &\equiv - [\alpha_{4}^{q-1}\alpha_5^{q-1}\alpha_6^{q-1}]\ (\Psi^{12,13}\Psi^{2,3}_1)^{q-1} \mod q \\[5pt]
      &\equiv - [\alpha_{4}^{q-1}\alpha_5^{q-1}\alpha_6^{q-1}]\ (\alpha_4\alpha_6(\alpha_4 + \alpha_5 + \alpha_6))^{q-1} \mod q \\[5pt]
      &\equiv - [\alpha_5^{q-1}]\ (\alpha_4 + \alpha_5 + \alpha_6)^{q-1} \mod q \\[5pt]
      &\equiv -1 \mod q.
  \end{align*}
  So the $c_2$ of $K_4$ is the constant sequence $-1$!
}

In fact all zig-zag graphs $Z_{\ell}$ with $\ell \geq 3$, which includes $K_4 = Z_3$, have $c_2(Z_{\ell}) = -1$ (see Corollary 54 of~\cite{k3}).
This makes zig-zag graphs particularly special as they are an infinite family of graphs for which both the period (Theorem~\ref{zigzag}) and the $c_2$ are completely known! \\

\bpoint{Hourglass chains}\label{SS:hourglass}

After double triangle reduction, up until recently, there hasn't been much progress in the direction of new graph transforms giving families of graphs with the same $c_2$.
A breakthrough came after Schnetz~\cite{geometries} introduced quadratic denominator reduction (see Section~\ref{SS:qdr}) as a new technique to compute the $c_2$-invariant.
The potential of this new technique was illustrated in~\cite{hourglass} where Schnetz and Yeats were able to determine a formula (using both cases in Definition~\ref{def:qdr}) for the $c_2$ of special infinite families of graphs.

The interesting part is that these graphs, which are formed from attaching \textbf{hourglass chains} to a \textbf{kernel} graph, were able to be reduced to their kernels (which then may be further reduced).
That is, the $c_2$ of such graphs only depends on said kernel, thus creating infinite families of graphs (of all loop orders) with the same $c_2$ and for which the $c_2$ is completely determined.
Furthermore, for suitable kernels these hourglass chain graphs actually form infinite families of prime ancestors!

This is exciting as we have previously only seen families of graphs for which the $c_2$ is trivial or $c_2 = -1$ (like zig-zag graphs), or we have double-triangle families, some of which the $c_2$ is unknown.
With hourglass chains, we now have families with a variety of interesting $c_2$s, including sequences that have never been found before! \\

\bpoint{Legendre symbols and modular forms}\label{SS:numerics}

Finally we very briefly mention some of the (other) sequences that appear as $c_2$-invariants of $\phi^4$ graphs.
Of the identified sequences, these $c_2$s currently fall into two types.
The first type corresponds to constant and quasi-constant $c_2$s, more generally those that are Legendre symbols.
The second type are what we call modular, where the $c_2$ arises as the Fourier coefficients of certain modular forms.
We refer the reader to ~\cite{modular,geometries} and the references therein for complete details. \\

\chapter{Completing the \texorpdfstring{$c_2$}{c2} completion conjecture for \texorpdfstring{$p=2$}{p=2}}\label{C:completion}

One of the motivating reasons to study the $c_2$-invariant (Definition~\ref{def:c2}) is its close relation to the Feynman period of a graph, as given in Equation~\eqref{eq:period}, where these are primitive divergent graphs in $\phi^4$-theory.
Recall from Section~\ref{S:phi4} that these particular graphs can be described as connected 4-regular graphs with one vertex removed and with the property that $\abs{E} = 2\ell$, where $\ell$ is the loop order, and $\abs{E(\lambda)} > 2\ell(\lambda)$ for any non-empty proper subgraph $\lambda$.

Not only is the $c_2$-invariant essentially capturing the behaviour of the zeros of the denominators of the period, but the connection is further strengthened by Conjecture~\ref{c2-period-conj}, which says that the $c_2$ should be a period invariant.
For this conjecture to hold, we would need the $c_2$-invariant to have all the symmetries of the period as described in Section~\ref{SS:symmetries}.
In particular, it should be invariant under the period symmetry called completion.

Schnetz proved in~\cite{census} that the period could be defined for the graphs created by adding back the vertex that was "removed" from primitive divergent graphs, and thus "completing" them to become 4-regular.
We call this 4-regular graph $G$ the \textbf{completion} and call the original primitive divergent graph $G-v$ for some vertex $v$ a \textbf{decompletion} of $G$.
Note that decompletions of 4-regular graphs could be non-isomorphic.
To utilize this, Schnetz further proved that if two primitive divergent graphs have the same completion, then their periods are also the same (see Theorem~\ref{period-comp}).
Thus the power of completion is that it allows us to lift the problem of calculating periods on primitive divergent graphs to one on 4-regular graphs, in the process cutting down the number of graphs to compute periods for.

For the $c_2$-invariant, whether or not the completion symmetry holds is still an open conjecture.
\tpointn{The $c_2$ Completion Conjecture} (Conjecture $4,35$ of ~\cite{k3})\label{c2-comp-conj} \\
\statement[eq]{
  Let $G$ be a connected 4-regular graph, and let $v$ and $w$ be vertices of $G$.
  Then,
  \[ c_2(G - v) = c_2(G - w). \]
}

While the $c_2$-invariant is defined for all prime powers $q$ and more general connected graphs, we will restrict to when $q = p$ is a prime (see also Conjecture~\ref{prime-powers}) and to decompletions of connected {4-regular} graphs $G$.
Note that for any vertex $v$, by Euler's formula, the decompletion $G - v$ indeed satisfies $\abs{E(G-v)} = 2\ell$ where $\ell$ is the loop order of $G-v$.

Yeats in~\cite{specialc2} first made progress in the special case of $p=2$ by reducing the conjecture to a combinatorial counting problem; one involving enumerating certain edge bipartitions.
Using two different constructions, Yeats was able to partially prove the $p=2$ case.
However, the obstruction to a full proof lay in the second construction, which required restricting the argument to completed graphs with an odd number of vertices.

In this chapter, inspired by the first construction in~\cite{specialc2}, we present this argument and show how by using the same ideas we can get rid of the mysterious parity condition on the vertices, thereby completing the $c_2$ completion conjecture for this special case.

\section{The conjecture and the \texorpdfstring{$p=2$}{p=2} case}\label{S:conjecture}

Let $G$ be a connected 4-regular graph, and let $v$ and $w$ be vertices of $G$.
Before even approaching the conjecture, we notice that it suffices to prove the equation holds when $v$ and $w$ are adjacent.
This is because as $G$ is connected, there is a path between any two vertices, and successively applying the equation to each pair of neighbours along the path gives the more general result.

\begin{figure}[th]
  \centering
  \includegraphics[scale=0.75]{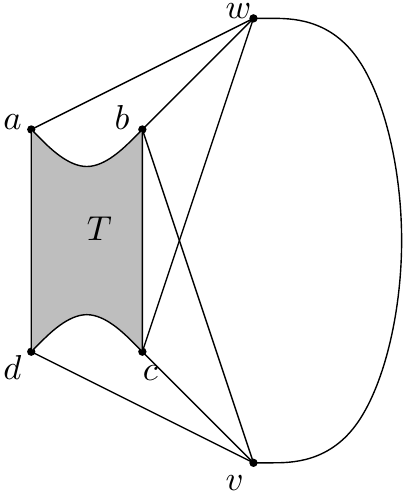}
  \quad\quad\quad \includegraphics[scale=0.75]{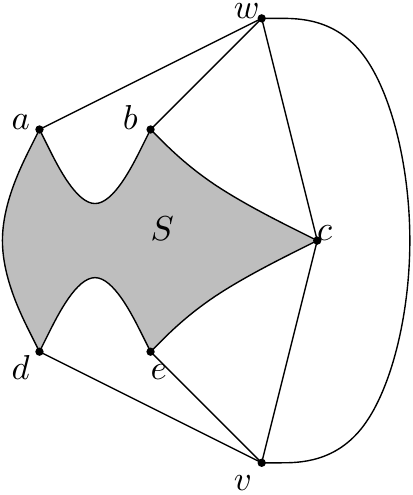}
  \quad\quad\quad \raisebox{4mm}[0pt][0pt]{\includegraphics[scale=0.75]{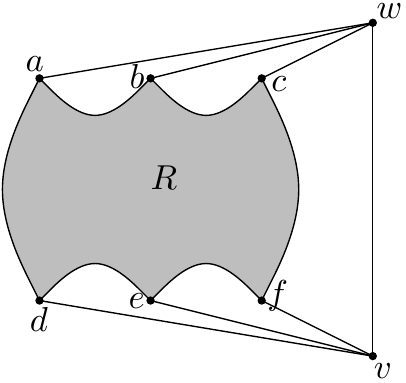}}
  \caption{The graph $G$, and the three cases $T$, $S$, and $R$.}
  \label{fig:TSR}
\end{figure}

As we will see shortly, the neighbours of $v$ and $w$ play an important role, and thus we split the conjecture into cases based on the number of neighbours that $v$ and $w$ have in common.
Since $G$ is a 4-regular graph and $v$ and $w$ are adjacent, there are four possibilities: they share all three neighbours, or exactly two or one neighbour, or they have no neighbours in common.
In the first case, when $v$ and $w$ share all neighbours, $G-v$ and $G-w$ are isomorphic, and thus the conjecture holds trivially.
For the others, we will deal with each case separately and label them with $T$, $S$, and $R$, respectively.
The $T$-case is when $v$ and $w$ share exactly two neighbours, the $S$-case is when they share only one neighbour, and the $R$-case is when they do not share any neighbours.
The three cases are depicted in Figure~\ref{fig:TSR} where the grey blobs are the graphs $G - \{v,w\}$ in each case, and will be addressed in Sections~\ref{S:T-case},~\ref{S:S-case}, and~\ref{S:R-case}, respectively.

With the set-up complete, the first step in tackling the $c_2$ completion conjecture is the question of how to reduce it into an approachable combinatorial counting problem.
Recall that in Section~\ref{SS:denom-count}, we were able to transform the $c_2$-invariant for primes $p$ into the following~\eqref{eq:c2-version} for a graph $G - v$ \\
\[ c_2^{(p)}(G-v) \equiv -[\alpha_4^{p-1}\cdots\alpha_{|E|}^{p-1}]\; \left(\Psi^{13,23} \Psi_{3}^{1,2} \right)^{p-1} \mod p, \] \\
where $1,2,3$ are the labels of any three distinct edges of $G - v$, and the Dodgson polynomials in the equation are with respect to $G - v$.

Now suppose $u$ is a 3-valent vertex in $G - v$.
Let the edges $1,2,3$ be the three edges incident to $u$, and let $u_1, u_2, u_3$ be the respective vertices adjacent to these edges as below:
\[ \Blob \]

Looking at the term $\Psi^{1,2}_3$, as discussed in Section~\ref{SS:dodgsons} we can view this Dodgson polynomial as coming from the following intersection of spanning trees
\[ \ExampleBlob \]
where the first blob depicts spanning trees on $(G -v) \setminus 1 / \{2,3\}$, and the second blob is on $(G - v) \setminus 2 / \{1,3\}$.

To translate this into spanning forests on $(G - v) \setminus \{1,2,3\}$ using Theorem~\ref{tree-to-forest}, we just need to determine the partitions $P$ of $\{u, u_1, u_2, u_3\}$ (as these are the vertices of edges $\{1,2,3\}$) such that all the spanning forests compatible with $P$ are spanning trees in the above intersection.
First, spanning trees on the first blob correspond to spanning forests such that $u$, $u_2$, and $u_3$ are in different trees, and $u_1$ can be in any of the three trees.
Similarily, spanning trees on the second blob correspond to spanning forests such that $u$, $u_1$, and $u_3$ are in different trees, and there is no restriction for which tree $u_2$ is in.
As these spanning forests have exactly three parts, the intersection of these blobs are spanning forests such that $u$ is in one tree, $u_3$ is in another tree, and $u_1$ and $u_2$ are in the last tree.
Furthermore, this is the only possible partition of $\{u, u_1, u_2, u_3\}$ that satisfies the required condition.

In particular, there is only one spanning forest polynomial in the sum, which is $\Phi^{\{u\},\{u_3\},\{u_1,u_2\}}_{(G - v) \setminus \{1,2,3\}}$.
However, notice that $u$ is an isolated vertex in $(G - v) \setminus \{1,2,3\}$ since edges $1,2,3$ are exactly the edges incident to $u$.
This means we can remove the vertex $u$, giving the equation
\[ \Psi^{1,2}_{3} = \Phi^{\{u_3\},\{u_1,u_2\}}_{G - \{u,v\}}. \]
Thus we are left with spanning forests with exactly two parts, which we call \textbf{spanning 2-forests}, such that they are compatible with the vertex partition $\{u_3\},\{u_1,u_2\}$.

Following the same process for $\Psi^{13,23}$, in either case we want spanning trees on the graph $G - \{u,v\}$
\[ \BlobBase \]
where for $(G - v) \setminus \{1,3\} / 2$, we have $u = u_2$, and for $(G - v) \setminus \{2,3\} / 1$, we have $u = u_1$.
Once again, we can translate these to spanning forests on $(G - v) \setminus \{1,2,3\}$, this time with exactly two parts, giving the polynomial $\Phi^{\{u\},\{u_1,u_2,u_3\}}_{(G - v) \setminus \{1,2,3\}}$.
Noticing that $u$ is an isolated vertex in this graph, removing vertex $u$ gives
\[ \Psi^{13,23} = \Phi^{\{u_1,u_2,u_3\}}_{G - \{u,v\}} = \Psi_{G - \{u,v\}}. \]

Putting everything together, we have reformulated the $c_2$-invariant to \\
\[ c_2^{(p)}(G-v) \equiv -[\alpha_4^{p-1}\cdots\alpha_{|E|}^{p-1}]\; \left(\Phi_{G-\{u,v\}}^{\{u_3\},\{u_1,u_2\}} \Psi_{G-\{u,v\}} \right)^{p-1} \mod p, \]

which we call \textbf{reducing with respect to vertex $u$}.
Note that the ordering of the vertices $u_1$, $u_2$, and $u_3$ was completely arbitrary, and this formulation of the $c_2$ does not depend on this choice.

Finally, by working with the $c_2$ in this form, there's a particularly nice interpretation of what this equation is counting.
From the coefficient extraction, the monomial we are looking for can be viewed as distributing $p-1$ copies of each edge in $G - \{u,v\}$ across $2p-2$ polynomials, where now $G - \{u,v\}$ is the underlying graph for the polynomials.
Of these $2p-2$  polynomials, $p-1$ of them arise from spanning 2-forests compatible with the vertex partition $\{u_3\},\{u_1,u_2\}$ and the other $p-1$ come from spanning trees.
Then, since all the monomials in $\Psi$ and $\Phi$ appear with a coefficient of $1$, we are enumerating all the ways to do such a distribution.

Now, recalling that the terms in $\Psi$ and $\Phi$ correspond to the edges \emph{not} in the respective trees or forests, we need to shift our view slightly to make sense of what this coefficient is counting.
Given $p-1$ spanning trees of $G - \{u,v\}$, the corresponding monomial in $\left(\Psi_{G - \{u,v\}}\right)^{p-1}$ comes from combining all the edges of the graph not in each spanning tree.
Thus the edges left over to be distributed to the $p-1$ $\Phi$ polynomials are exactly the edges of the $p-1$ spanning trees.
Furthermore, these edges must be grouped in a such a way that removing them gives $p-1$ spanning 2-forests.

Similarily, given $p-1$ spanning 2-forests the corresponding monomial in $\left(\Phi_{G-\{u,v\}}^{\{u_3\},\{u_1,u_2\}}\right)^{p-1}$ combines all the edges of the graph not in each forest.
Notice this monomial must then match the edges left over from the monomial given by the $p-1$ spanning trees above, which are exactly the edges of the spanning trees.
By the same logic, the monomial arising from the $p-1$ spanning trees must then match the edges of the $p-1$ spanning 2-forests.
Thus, we can swap the roles of the spanning trees and forests, and equivalently express the enumeration as counting the number of ways to partition $p-1$ copies of each edge into $p-1$ spanning trees and $p-1$ spanning 2-forests compatible with $\{u_3\},\{u_1,u_2\}$!

When $p=2$, this interpretation is much simpler as the equation simplifies to \\
\begin{equation}\label{eq:c2-counting}
  c_2^{(2)}(G-v) \equiv [\alpha_4\cdots\alpha_{|E|}]\; \Phi_{G-\{u,v\}}^{\{u_3\},\{u_1,u_2\}} \Psi_{G-\{u,v\}} \mod 2.
\end{equation}

Now we only have one copy of each edge, and we are counting the number of edge bipartitions where one part is a spanning tree and the other is a spanning 2-forest compatible with $\{u_3\},\{u_1,u_2\}$.
We will usually denote these edge partitions as $(\psi, \phi)$.
Additionally, as we are counting modulo $2$, we really only care about the parity of this count.
Thus one strategy to tackling the $c_2$ completion conjecture when $p=2$ is to find fixed-point free involutions on the appropriate sets of edge bipartitions! \\

\section{The \texorpdfstring{$T$}{T}-case}\label{S:T-case}

The first case we will look at is the $T$-case, which is also special in and of itself.
It is the starting point for many of the arguments needed later on.
Let $G$ be a connected 4-regular graph, and let $v$ and $w$ be two adjacent vertices of $G$ such that they share two common neighbours.
Let $T = G - \{v,w\}$ be the graph obtained by removing vertices $v$ and $w$ from $G$, the first grey blob in Figure~\ref{fig:TSR}, with the neighbours of $v$ and $w$ labelled $\{a,b,c,d\}$ as in the figure.
Here $w$ has neighbours $\{a,b,c,v\}$, and $v$ has neighbours $\{b,c,d,w\}$.

What's special about this case is when $v$ and $w$ share two common neighbours, these four vertices form a double triangle as defined in Section~\ref{SS:DTR}.
From the results in Section~\ref{SS:DTR}, we know that double triangles are important since the $c_2$ is invariant under double triangle reductions.
Now adding in decompleting at a double triangle vertex, as proved in ~\cite{further} (see also Theorem~\ref{DTRinvar}), we can resolve the $c_2$ completion conjecture in the $T$-case for all values of $p$ as a consequence of double triangle reduction.
This proof uses a different argument, namely via proving an equivalence of certain graph polynomials.
Even though we do have this alternative proof, it is still worthwhile to present the enumerative argument from~\cite{specialc2} in this section, as it showcases the proof techniques and overall framework that the $S$ and $R$ cases build upon.

Furthermore, another reason for the specialness of the $T$-case is that by using the same counting methods as for the $p=2$ case, these results can be extended to all values of $p$~\cite{higherp}.
While the $S$ and $R$ cases are not as simple, this extension gives some hope for generalizing the arguments needed in those cases for higher values of $p$. \\

\bpoint{Set-up}

To start off, we define the particular sets of edge partitions that we are counting using the vertex bipartition that the spanning 2-forests are compatible with (see Definition~\ref{def:span-forest}).
We will use similar definitions for the $S$ and $R$ cases.

\tpointn{Definition} (Definition 3.1 of~\cite{specialc2})\label{T-defs}
\statement{
  Suppose $P$ is a bipartition of $\{a,b,c,d\}$. \\
  Let $\mathcal{T}_P$ be the set of bipartitions $(\psi, \phi)$ of the edges of $T$ such that $\psi$ is a spanning tree and $\phi$ is a spanning 2-forest compatible with $P$.
  Let $t_P = \abs{\mathcal{T}_P}$.
}

Now, since we can choose the order of the edges adjacent to the 3-valent vertex we are reducing by in Equation~\eqref{eq:c2-counting}, and thus choose the partition of its neighbours into a 1-part and 2-part for the 2-forest, in the $T$-case we can do so in such a way that exploits the property that $v$ and $w$ share two common neighbours.

\tpointn{Proposition} (Proposition 3.2, 3.3 of~\cite{specialc2})\label{T-eqs}
\statement[eq]{
  When $v$ and $w$ have common neighbours $b, c$,
  \[ c_2^{(2)}(G - v) = t_{\{a\},\{b,c\}} \mod 2, \]
  \[ c_2^{(2)}(G - w) = t_{\{d\},\{b,c\}} \mod 2, \]
  and thus we have
  \[ c_2^{(2)}(G - v) - c_2^{(2)}(G - w) = {t}_{\{a\},\{b,c,d\}} + {t}_{\{d\},\{a,b,c\}} \mod 2. \]
}
\begin{proof}
  In the $G-v$ case, from Equation~\eqref{eq:c2-counting} we can reduce with respect to vertex $w$, since now $w$ is a 3-valent vertex, and choose the order of its adjacent edges such that we get vertex partition $\{a\}, \{b,c\}$ for the spanning 2-forest.
  Then, by the discussion in Section~\ref{S:conjecture}, this is equivalent to counting the number of ways to partition the edges of $G - \{v,w\} = T$ into a spanning tree and a spanning 2-forest compatible with $\{a\},\{b,c\}$, which gives
  \[ c_2^{(2)}(G - v) \equiv [\alpha_4\cdots\alpha_{|E(G-v)|}]\; \Phi_{G-\{v,w\}}^{\{a\},\{b,c\}} \Psi_{G-\{v,w\}} \equiv t_{\{a\},\{b,c\}} \mod 2. \]
  Similarily, in the $G-w$ case, $v$ is now a 3-valent vertex, and we choose vertex partition $\{d\},\{b,c\}$ to get
  \[ c_2^{(2)}(G - w) \equiv t_{\{d\},\{b,c\}} \mod 2. \]

  Finally, enumerating over all possibilities for the vertex set $\{a,b,c,d\}$ in the spanning 2-forests
  \[t_{\{a\},\{b,c\}} = t_{\{a\},\{b,c,d\}} + t_{\{a,d\},\{b,c\}}, \]
  \[t_{\{d\},\{b,c\}} = t_{\{d\},\{a,b,c\}} + t_{\{a,d\},\{b,c\}}, \]
  and thus since $t_{\{a,d\},\{b,c\}}$ appears in both equations and we are working modulo 2
  \[ c_2^{(2)}(G - v) - c_2^{(2)}(G - w) = t_{\{a\},\{b,c\}} - t_{\{d\},\{b,c\}} \equiv {t}_{\{a\},\{b,c,d\}} + {t}_{\{d\},\{a,b,c\}} \mod 2. \]
\end{proof}

\bpoint{Swapping around a two-valent vertex}

The main property of the $T$-case, that makes it much simplier to deal with, is that since $v$ and $w$ share two common neighbours, after removing both $v$ and $w$ we are left with two 2-valent vertices.
As we are dealing with edge bipartitions involving spanning trees and forests of a particular form, these 2-valent vertices cannot be isolated in said forests, and thus we know exactly how the two edges of a 2-valent vertex must be distributed across the edge partition.
Furthermore, we can swap these edges to obtain another valid edge partition.

\tpointn{Lemma} (Lemma 4.1 of~\cite{specialc2})\label{swap-two}
\statement{
  Let $V$ be a set of marked vertices of a connected graph $G$, and suppose $c \in V$ is a 2-valent vertex.
  Let $(\psi,\phi)$ be an edge partition of $G$ where $\psi$ is a spanning tree and $\phi$ is a spanning 2-forest that is compatible with a bipartition of $V$ where each part (and thus forest) contains at least one of the vertices in $V \setminus \{c\}$. \\\\
  Then, of the the two edges incident to $c$, exactly one edge is in each part of $(\psi, \phi)$, and swapping which edge is in which part yields a new edge partition $(\psi', \phi')$ with the same properties, except possibly the bipartition of $V$ that $\phi'$ is compatible with.
}
\begin{proof}
  Looking at the edge partition $(\psi, \phi)$, first, since $\psi$ is a spanning tree of $G$ there must be at least one edge in $\psi$ that is adjacent to the vertex $c$.
  Now since $\phi$ is a spanning 2-forest in which each forest contains at least one of the vertices in $V \setminus \{c\}$, we have that in particular, $c$ is not an isolated vertex in $\phi$ and must be connected to the rest of its forest via at least one edge.
  As $c$ is a 2-valent vertex adjacent to at least one edge in part of the edge partition, exactly one edge incident to $c$  must be in each of $(\psi, \phi)$.

  Because $c$ is a leaf in both $\psi$ and $\phi$, removing the edges incident to $c$ disconnects and furthermore, isolates $c$ from the rest of the tree and the 2-forest, respectively.
  Since both $\psi$ and $\phi$ are spanning, they must contain both neighbours of $c$, and thus we can reconnect $c$ via the opposite edges incident to $c$ as in $\psi$ and $\phi$.
  This gives a new edge partition $(\psi', \phi')$ where $\psi'$ is a spanning tree and $\phi'$ is a spanning 2-forest, the only difference between the two edge partitions being the swapping of the edges incident to $c$.
  However, in $\phi'$, $c$ could be reconnected to a different tree than in $\phi$ thus changing the bipartition of $V$ that $\phi'$ is compatible with.
  Specifically, either the bipartition of $V$ remains the same as for $\phi$ or $c$ swaps between the parts of the bipartition.
\end{proof}

With this observation that we can swap the edges around any 2-valent vertex to get two different edge partitions $(\psi, \phi)$ and $(\psi', \phi')$, we can obtain fixed-point free involutions on $\mathcal{T}_P$ for some particular bipartitions $P$ of $\{a,b,c,d\}$.

\tpointn{Theorem} (from \textsection 6 of ~\cite{specialc2})\label{T-bijection}
\statement[eq]{
  There is a fixed-point free involution on
  \[ \mathcal{T}_{\{a\},\{b,c,d\}} \cup \mathcal{T}_{\{a,b\},\{c,d\}}, \]
  and on
  \[ \mathcal{T}_{\{d\},\{a,b,c\}} \cup \mathcal{T}_{\{a,b\},\{c,d\}}. \]
  Thus we have
  \[ {t}_{\{a\},\{b,c,d\}} + {t}_{\{a,b\},\{c,d\}} \equiv 0 \mod 2,\]
  \[ {t}_{\{d\},\{a,b,c\}} + {t}_{\{a,b\},\{c,d\}} \equiv 0 \mod 2,\]
  which together gives
  \[ {t}_{\{a\},\{b,c,d\}} + {t}_{\{d\},\{a,b,c\}} \equiv 0 \mod 2.\]
}
\begin{proof}
  Consider the edge partitions in $\mathcal{T}_{\{a\},\{b,c,d\}} \cup \mathcal{T}_{\{a,b\},\{c,d\}}$.
  Since $b$ is a 2-valent vertex in $T$, and both $\{a\}, \{b,c,d\}$ and $\{a,b\},\{c,d\}$ are bipartitions of $\{a,b,c,d\}$ in which each part contains at least one of $\{a,c,d\}$, these edge partitions satisfy the conditions of Lemma~\ref{swap-two}.

  Take any edge partition $(\psi, \phi) \in \mathcal{T}_{\{a\},\{b,c,d\}}$ and swap the edges incident to $b$ between the two parts to get a new edge partition $(\psi', \phi')$.
  Looking at the new spanning 2-forest $\phi'$, either both neighbours of $b$ were in the same tree in $\phi$, and thus $b$ was reconnected to that same tree corresponding to the part $\{b,c,d\}$ in the vertex bipartition.
  Or the neighbours of $b$ were in different trees in $\phi$, and thus $b$ was reconnected but to the tree coresponding to the part $\{a\}$.
  In the first case, $(\psi', \phi')$ is once again an edge partition in $\mathcal{T}_{\{a\},\{b,c,d\}}$.
  In the second case, $(\psi', \phi')$ in now an edge partition in $\mathcal{T}_{\{a,b\},\{c,d\}}$.
  Similarily, for edge partitions in $\mathcal{T}_{\{a,b\},\{c,d\}}$, swapping around $b$ either gives a new edge partition in the same set (first case) or in $\mathcal{T}_{\{a\},\{b,c,d\}}$ (second case).

  In either case, since we are swapping the edges incident to $b$ in each part, $(\psi', \phi')$ must be different from $(\psi, \phi)$, and swapping around $b$ once again returns us to $(\psi, \phi)$.
  Thus, this action of swapping around $b$ gives a fixed-point free involution on $\mathcal{T}_{\{a\},\{b,c,d\}} \cup \mathcal{T}_{\{a,b\},\{c,d\}}$, and so the size of this set must be even, giving the equation
  \[ {t}_{\{a\},\{b,c,d\}} + {t}_{\{a,b\},\{c,d\}} \equiv 0 \mod 2. \]

  Now since $c$ is also 2-valent in $T$, using the exact same argument on $\mathcal{T}_{\{d\},\{a,b,c\}} \cup \mathcal{T}_{\{a,b\},\{c,d\}}$ with $c$ as the vertex we are swapping around, we obtain a fixed-point free involution on these edge partitions, and thus
  \[ {t}_{\{d\},\{a,b,c\}} + {t}_{\{a,b\},\{c,d\}} \equiv 0 \mod 2. \]

  Finally, notice that ${t}_{\{a,b\},\{c,d\}}$ appears in both equations, and since we are working modulo 2, adding the two equations gives
  \[ {t}_{\{a\},\{b,c,d\}} + 2{t}_{\{a,b\},\{c,d\}} + {t}_{\{d\},\{a,b,c\}} \equiv {t}_{\{a\},\{b,c,d\}} + {t}_{\{d\},\{a,b,c\}} \equiv 0 \mod 2. \]
\end{proof}

From Proposition~\ref{T-eqs}, the final equation in Theorem~\ref{T-bijection} is exactly what we need to prove the $c_2$ completion conjecture for $p=2$ in the $T$-case.

\tpointn{Corollary}\label{T-case}
\statement[eq]{
  Let $G$ be a connected 4-regular graph.
  Let $v$ and $w$ be adjacent vertices of $G$ such that they share two common neighbours.
  Then,
  \[ c_2^{(2)}(G - v) = c_2^{(2)}(G - w). \]
} \\

\section{The \texorpdfstring{$S$}{S}-case}\label{S:S-case}

With the $T$-case complete, we move to the slightly more complicated $S$-case.
Let $G$ be a connected 4-regular graph, and let $v$ and $w$ be two adjacent vertices of $G$ such that they share one neighbour.
Let $S = G - \{v,w\}$ be the graph obtained by removing vertices $v$ and $w$ from $G$, the middle grey blob in Figure~\ref{fig:TSR}, with the neighbours of $v$ and $w$ labelled $\{a,b,c,d,e\}$ as in the figure.
Here $w$ has neighbours $\{a,b,c,v\}$, and $v$ has neighbours $\{c,d,e,w\}$.

While the $S$-case is more complicated, it does use ideas that stem from the $T$-case.
Firstly, as $v$ and $w$ still share a neighbour, after removing both vertices $S$ will have a 2-valent vertex, and thus we can use the same swapping around a 2-valent vertex argument as in the $T$-case.
In fact, originally in ~\cite{specialc2} this argument was first used for the $S$-case!
However, it was not enough to complete the proof of the conjecture for $p=2$.
A different type of argument involving compatible cycles (see \textsection 5 of~\cite{specialc2}) was given for the rest of the sets not covered by swapping around a 2-valent vertex, but under the condition that $G$ had an odd number of vertices.

Secondly, this idea of swapping around a 2-valent vertex can be generalized to swapping around a particular vertex called the control vertex, under certain conditions.
While first used in ~\cite{specialc2} for the $R$-case, we prove that this argument can be modified to apply in the $S$-case as well.
This modification was able to cover all the sets from the compatible cycle argument without the odd vertex condition on $G$, and thus complete the proof of the $c_2$ completion conjecture for $p=2$ in the $S$-case.

Being not as simple as the $T$-case but also perhaps not as hard as the $R$-case due to the presence of a 2-valent vertex, the $S$-case seems to be the perfect playground and stepping stone for trying to further extend these ideas for higher values of $p$.

To give a complete, self-contained proof of the $S$-case, we begin by presenting the set-up and some of the  results from~\cite{specialc2} with proof, and in Section~\ref{SS:S-new} we introduce our new results which allows us to finish the $S$-case for $p=2$! \\

\bpoint{Set-up}

As with the $T$-case, we start with defining the particular sets of edge partitions that we are counting.

\tpointn{Definition} (Definition 3.1 of~\cite{specialc2})\label{S-defs}
\statement{
  Suppose $P$ is a bipartition of $\{a,b,c,d,e\}$. \\
  Let $\mathcal{S}_P$ be the set of bipartitions $(\psi, \phi)$ of the edges of $S$ such that $\psi$ is a spanning tree and $\phi$ is a spanning 2-forest compatible with $P$.
  Let $s_P = \abs{\mathcal{S}_P}$.
}

Once again, we can pick the vertex bipartitions for the spanning 2-forests to exploit the property that $v$ and $w$ share a neighbour.

\tpointn{Proposition} (Proposition 3.2, 3.3 of~\cite{specialc2})\label{S-eqs}
\statement[eq]{
  When $v$ and $w$ have one common neighbour $c$
  \[ c_2^{(2)}(G - v) = s_{\{c\},\{a,b\}} \mod 2, \]
  \[ c_2^{(2)}(G - w) = s_{\{c\},\{d,e\}} \mod 2, \]
  and thus we have
  \begin{align*}\refstepcounter{equation}\tag{\theequation}\label{eq:S-counts}
    c_2^{(2)}(G - v) - c_2^{(2)}(G - w)
      =\ &{s}_{\{c,d\},\{a,b,e\}} + {s}_{\{c,e\},\{a,b,d\}} + {s}_{\{a,b\},\{c,d,e\}} \\
      &+ {s}_{\{a,c\},\{b,d,e\}} + {s}_{\{b,c\},\{a,d,e\}} + {s}_{\{d,e\},\{a,b,c\}} \mod 2.
  \end{align*}
}
\begin{proof}
  Using the same argument as in the proof of Proposition~\ref{T-eqs}, we can pick the vertex partition $\{c\},\{a,b\}$ for the $G-v$ case and $\{c\},\{d,e\}$ for the $G-w$ case to equivalently determine each $c_2$ via counting edge partitions into a spanning tree and a spanning 2-forest compatible with the respective vertex partition.

  Then, enumerating over all possibilities for the vertex set $\{a,b,c,d,e\}$ in the spanning 2-forests, and noticing that $s_{\{c\},\{a,b,d,e\}}$ appears for both $G-v$ and $G-w$, we obtain the final equation by subtracting the expanded counts for $G-v$ and $G-w$, and simplifying modulo $2$.
\end{proof}

Here is a nice visual representation of Equation~\eqref{eq:S-counts} from ~\cite{specialc2} where we depict the vertex bipartition for the spanning 2-forest using different shapes:

\[ c_2^{(2)}(G - v) - c_2^{(2)}(G - w) = \raisebox{-3mm}[0pt][0pt]{\includegraphics[scale=1]{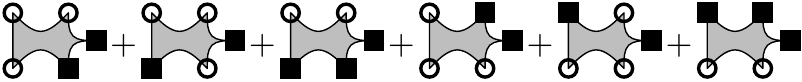}} \mod 2 \] \\

\bpoint{Swapping around \texorpdfstring{$c$}{c}}

The main property of the $S$-case is that $v$ and $w$ share a neighbour, and in particular, after removing both $v$ and $w$ from $G$, we are left with a 2-valent vertex.
Since we already know how to deal with 2-valent vertices, we immediately obtain the following theorem.

\tpointn{Theorem} (Lemma 4.2 of~\cite{specialc2})\label{S-swapc}
\statement[eq]{
  There is a fixed-point free involution on
  \[ \mathcal{S}_{\{a,b,c\},\{d,e\}} \cup \mathcal{S}_{\{a,b\},\{c,d,e\}}, \]
  and thus we have
  \[ {s}_{\{a,b,c\},\{d,e\}} + {s}_{\{a,b\},\{c,d,e\}} \equiv 0 \mod 2. \]
}
\begin{proof}
  Consider the edge partitions in $\mathcal{S}_{\{a,b,c\},\{d,e\}} \cup \mathcal{S}_{\{a,b\},\{c,d,e\}}$.
  Since $c$ is a 2-valent vertex in $S$, and both $\{a,b,c\},\{d,e\}$ and $\{a,b\}, \{c,d,e\}$ are bipartitions of $\{a,b,c,d,e\}$ in which each part contains at least one of $\{a,b,d,e\}$, these edge partitions satisfy the conditions of Lemma~\ref{swap-two}.

  Then, using the same argument as in the proof of Theorem~\ref{T-bijection}, swapping the edges incident to $c$ between the two parts of the edge partitions of interest gives a new edge partition with a spanning 2-forest either compatible with the same vertex bipartition as we started with or the vertex bipartition where $c$ swaps between the parts of the original vertex bipartition.
  Once again we have a fixed-point free involution, this time on $\mathcal{S}_{\{a,b,c\},\{d,e\}} \cup \mathcal{S}_{\{a,b\},\{c,d,e\}}$, and so the size of this set must also be even
  \[ {s}_{\{a,b,c\},\{d,e\}} + {s}_{\{a,b\},\{c,d,e\}} \equiv 0 \mod 2. \]
\end{proof}

Swapping around $c$ takes care of two specific vertex bipartitions on the right-hand side of Equation~\eqref{eq:S-counts}, and what we are left with are sets of the form $\mathcal{S}_{\{c,*\}, \{*,*,*\}}$.
However, when we try to swap around $c$ for these sets we obtain vertex bipartitions not in the equation nor ones that coincide with each other, like in the $T$-case with $t_{\{a,b\},\{c,d\}}$.
For example, for $\{a,c\},\{b,d,e\}$, we would get the equation $s_{\{a,c\},\{b,d,e\}} + s_{\{a\},\{b,c,d,e\}} \equiv 0 \mod 2$ and similarly for the rest.
To tackle this problem, we generalize the idea of swapping around a 2-valent vertex! \\

\bpoint{Swapping around a control vertex}\label{SS:S-new}

For the remaining sets in Equation~\eqref{eq:S-counts}, we generalize the swapping argument to one using a vertex that is not 2-valent in $S$.
While this idea of using control vertices is not new, we adapt it to be able to apply it to the $S$-case.
In the process we develop a new swapping method involving multiple vertices, starting with using a control vertex in conjunction with a 2-valent vertex.
Generalizing once again, the notion of using multiple control vertices will then be used for the $R$-case.

Recall that in the 2-valent vertex case, we are swapping the edges incident to this vertex between the two parts of the edge partition $(\psi, \phi)$.
We can do so uniquely since under some requirements for $(\psi, \phi)$, this vertex is a leaf in both parts.

As we no longer have an appropriate 2-valent vertex, we need to define a different vertex to swap edges around. We do so under certain conditions such that the same vertex can be recovered after swapping an incident edge in both $\psi$ and $\phi$, and such that we can get a new 2-forest compatible with a vertex bipartition that we care about!

\begin{figure}[t]
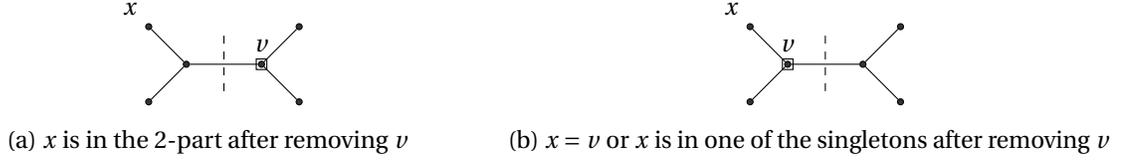

  \centering
  \begin{subfigure}[b]{0.5\linewidth}
    \centering
    \ShapeOneCut%
    \caption{$x$ is in the 2-part after removing $v$}
    \label{fig:S-control}
  \end{subfigure}%
  \begin{subfigure}[b]{0.5\linewidth}
    \centering
    \ShapeTwoCut%
    \caption{$x=v$ or $x$ is in one of the singletons after removing $v$}
    \label{fig:R-control}
  \end{subfigure}%
  \caption[The two possible control vertices for Lemma~\ref{control-vertex}.]{The two possible control vertices $v$ depending on the specification for special vertex $x$ for Lemma~\ref{control-vertex}. The three unlabelled leaves are the vertices in $p \setminus \{x\}$. The dashed lines indicate the edges we will be picking for one of the involutions in the $S$-case and in the $R$-case, respectively. Some of the leaf edges are possibly contracted, with the exception that two leaf edges incident to the same vertex cannot simultaneously both be contracted.}
  \label{fig:control}
\end{figure}

\tpointn{Lemma} (generalizing Lemma 4.3 of~\cite{specialc2})\label{control-vertex}
\statement{
  Let $G$ be a connected 4-regular graph with two adjacent vertices removed, and let $V$ be the set of neighbours of those vertices.
  Fix a vertex $x \in V$, which we call the \textbf{special vertex}.
  Suppose we have a bipartition $P$ of $V$ such that each part contains at least one of the vertices in $V$ and one part is of the form $p = \{x, *, *, *\}$. \\\\
  Let $(\psi, \phi)$ be an edge partition of $G$ where $\psi$ is a spanning tree and $\phi$ is a spanning 2-forest compatible with $P$.
  Let $t$ be the tree corresponding to $p$ in the 2-forest $\phi$.
  Then, there are two distinct vertices $v$ such that removing $v$ from $t$ partitions $\{x,*,*,*\}$ into a 2-part, a singleton, and possibly a third part which is also a singleton. \\\\
  If we further specify that $x$ must satisfy one of the following properties, either
  \begin{enumerate}
    \item $x$ is in the 2-part, or
    \item $x=v$ or $x$ is in one of the singleton components,
  \end{enumerate}
  then there is a unique vertex $v$, which we call the \textbf{control vertex}, that satisfies all the required properties.
  Figure~\ref{fig:control} depicts the control vertex in each case.
  Furthermore, in $t$ this control vertex $v$ will either be 2-valent with $v \in p$ or 3-valent and $v \not\in p$.
  In both cases $v$ is a leaf in the spanning tree $\psi$.
}
\begin{proof}
  Looking at the tree $t$ and the subtree created by taking the union of all the paths in $t$ between any two vertices in $p$, since $p$ has four vertices this subtree will be of the form
  \[ \ShapeOne \]\\
  in which the leaves are vertices in $p$, the edges are paths in $t$, and some paths are possibly contracted.

  First, we argue that this is indeed the correct form for the subtree.
  Since $\psi$ is a spanning tree and we are partitioning the edges of $G$, the vertices in $\phi$, and thus the non-leaves in the subtree, are at most 3-valent.
  Furthermore, in $\phi$ the vertices in $P$, and thus in $p$, are at most 2-valent.
  This immediately implies that we have the subtree form above, where some of the paths may be contracted.
  The valency constraints also mean that the path between the two non-leaf vertices cannot be contracted.
  Note trivially, the paths between two leaves adjacent to the same closest non-leaf vertex cannot both be contracted simultaneously.
  All possible configurations are drawn in Figure 6 of~\cite{specialc2}.

  Now, because of the restrictions on the subtree configurations, removing any one of the two non-leaf vertices in the subtree above from $t$ creates components that partitions $p$ into a 2-part, a singleton, and possibly a third singleton.
  By the valency of the vertices in $\phi$, these two vertices can only be 2-valent if they are in $p$ or otherwise 3-valent.
  In either case, by the properties of $G$, they will be leaves in $\psi$.

  Let $v$ be one of the two non-leaf vertices in the subtree above.
  Finally, to get uniqueness of $v$ we need to fix where one vertex $x \in p$ must be after removing $v$ from $t$.
  In the first case, fixing $x$ to be in the 2-part component, we get the control vertex $v$ as pictured in Figure~\ref{fig:S-control}.
  Notice here $x$ could never be the control vertex $v$.
  In the second case, if we want $x$ to be in a singleton component we choose the other non-leaf vertex to be $v$ as in Figure~\ref{fig:R-control}.
  However, with this choice of $v$, we could also have that $x=v$ (the path between $x$ and $v$ is contracted in $t$).
  Thus in this case, the property we are fixing is either $x=v$ or $x$ is in one of the singleton components.
\end{proof}

Now we are ready to tackle the remaining sets in the Equation~\eqref{eq:S-counts}, which all are of the form $\mathcal{S}_{\{c,*\},\{*,*,*\}}$.
Instead of dealing with these sets directly, since we know there is a fixed-point free involution via swapping around $c$ (as it is a 2-valent vertex), we can equivalently show that there is a fixed-point free inovlution on sets of the form $\mathcal{S}_{\{*\},\{c,*,*,*\}}$ (the "swapped" version).

The involution uses a two-phase process, where we utilize that we can swap around $c$ for certain edge bipartitions.
For the others, we modify the involution used in Theorem~\ref{R-control-bij}, which uses the same idea of swapping edges incident to the control vertex between the two parts of $(\psi, \phi)$.

\tpointn{Theorem}\label{S-bijection}
\statement[eq]{
  There is a fixed-point free involution on
  \[ \mathcal{S}_{\{a\},\{b,c,d,e\}} \cup \mathcal{S}_{\{b\},\{a,c,d,e\}} \cup \mathcal{S}_{\{d\},\{a,b,c,e\}} \cup \mathcal{S}_{\{e\},\{a,b,c,d\}}, \]
  and thus we have
  \[ {s}_{\{a\},\{b,c,d,e\}} + {s}_{\{b\},\{a,c,d,e\}} + {s}_{\{d\},\{a,b,c,e\}} + {s}_{\{e\},\{a,b,c,d\}} \equiv 0 \mod 2. \]
}
\begin{proof}
  Consider the edge partitions in the union of sets of the form $\mathcal{S}_{\{*\},\{c,*,*,*\}}$ where the $*$'s are the vertices $\{a,b,d,e\}$ in any order.
  Since partitioning $\{a,b,c,d,e\}$ (the neighbours of the two adjacent vertices removed to obtain the graph $S$) into the form $\{*\},\{c,*,*,*\}$ leaves each part containing at least one of $\{a,b,c,d,e\}$ and a 4-part containing $c$, these edge partitions satisfy the conditions of both Lemma~\ref{swap-two} and Lemma~\ref{control-vertex}.
  For Lemma~\ref{control-vertex}, we are taking $c$ to be the special vertex and specifying that $c$ be in the 2-part after removing the control vertex from the approriate tree in the spanning 2-forest (the first property).

  Let $(\psi, \phi)$ be an edge partition in any set of the form $\mathcal{S}_{\{x\},\{c,*,*,*\}}$, and let $v$ be its control vertex.
  Let $t$ be the tree in $\phi$ corresponding to part $\{c,*,*,*\}$ which also contains $v$.
  Since $v$ is a leaf in the spanning tree $\psi$, there is exactly one edge incident to $v$ in $\psi$, which we will call $\eta_v$, and let $n_v$ be the neighbour of $v$ in $\psi$.
  Since $c$ is a 2-valent vertex, and thus also a leaf in $\psi$, let $\eta_c$ be the edge incident to $c$ in $\psi$ and $n_c$ its neighbour.

  Now to describe the involution giving a new edge partition:
  \begin{enumerate}[label=(\arabic*)]
    \item\label{S-step-1} \emph{Swapping $c$ stays in} -- If $n_c \in t$, then swap the edges incident to $c$ between $\psi$ and $\phi$.
    \item\label{S-step-2} \emph{Swapping $c$ goes out} -- Otherwise:
      \begin{enumerate}[label=(\roman*), ref=(2\roman*)]
        \item\label{S-step-in} \emph{Swapping $v$ stays in} -- If $n_v \in t$, so the control vertex and its neighbour in $\psi$ are in the same tree of $\phi$, let $\eta$ be the edge incident to $v$ in $t$ in the path to $n_v$.\\
          Then, swap the edges $\eta_v$ and $\eta$ between $\psi$ and $\phi$.
        \item\label{S-step-out} \emph{Swapping $v$ goes out} -- If $n_v \not\in t$, so the control vertex and its neighbour in $\psi$ are in different trees of $\phi$, let $\eta$ be the edge incident to $v$ in $t$ in the path to $c$; in Figure~\ref{fig:S-control} this is the edge towards the 2-part indicated by the dashed line. \\
          Then swap the edges $\eta_v$ and $\eta$ between $\psi$ and $\phi$, and finally swap the edges incident to $c$.
      \end{enumerate}
  \end{enumerate}

  First, we show that we get a valid edge partition in a set of the form $\mathcal{S}_{\{*\},\{c,*,*,*\}}$.
  For ~\ref{S-step-1}, this is just a specific case of the swapping around a 2-valent vertex argument.
  As in Lemma~\ref{swap-two}, since $c$ is a 2-valent vertex and is a leaf in both $\psi$ and $\phi$, we can swap the edges incident to $c$ between the two parts to get $(\psi', \phi')$ where $\psi'$ is still a spanning tree and $\phi'$ is still a spanning 2-forest.
  Looking at $\phi'$, as $n_c \in t$ and thus both neighbours of $c$ were in the same tree in $\phi$, $c$ gets reconnected to $t$ after the swap.
  Thus $\phi'$ is compatible with the same vertex bipartition as $\phi$, and $(\psi', \phi')$ is in the same set $\mathcal{S}_{\{x\},\{c,*,*,*\}}$ as before.

  For ~\ref{S-step-2}, we can use a similar argument to swapping around a 2-valent vertex, except replacing the 2-valent vertex with the control vertex $v$.
  From Lemma~\ref{control-vertex}, we know that $v$ is the unique vertex such that removing $v$ from $t$ partitions the vertices in part $\{c,*,*,*\}$ into a 2-part containing $c$, a singleton, and possibly a third singleton.
  Unlike the 2-valent argument, as $v$ is no longer a leaf in $t$ and we now have a choice of which edge incident to $v$ in $t$ we are swapping, we need to be careful of creating cycles when reconnecting $v$ to the rest of $t$ and make sure that we can recover the edge swapped and the same control vertex to indeed get an involution.

  To deal with the possibility of creating cycles, we notice that the only way to get a cycle after swapping edges incident to $v$ between $\psi$ and $\phi$ is if the neighbour of $v$ in $\psi$ was in the same tree as $v$ in $\phi$ and the edge picked to be swapped in $\phi$ was not on the cycle created by adding the edge incident to $v$ in $\psi$ to $\phi$.
  Thus for ~\ref{S-step-in}, if $n_v \in t$, then we pick the edge $\eta$ to be the edge incident to $v$ on the path to $n_v$ in $t$.
  In this case all the neighbours of $v$ are in the same tree in $\phi$.
  Removing $\eta$ from $t$ breaks the tree $t$ into two components, one with $v$ and one containing $n_v$, splitting $\phi$ into three components.
  Adding $\eta_v$ reconnects the two components from $t$ via $v$ and $n_v$ to create a new spanning 2-forest $\phi'$.
  As $n_v$ and $v$ were both in $t$, $\phi'$ is compatible with the same vertex bipartition as $\phi$.
  Since $v$ is a leaf in $\psi$, removing $\eta_v$ from $\psi$ and adding $\eta$ to create $\psi'$ maintains the spanning tree structure.
  Thus swapping the edges $\eta_v$ and $\eta$ between $\psi$ and $\phi$ creates a new edge partition $(\psi', \phi')$ in the same set $\mathcal{S}_{\{x\},\{c,*,*,*\}}$ as we started with.

  For ~\ref{S-step-out} when $n_v \not\in t$, we pick $\eta$ to be the edge incident to $v$ on the path to $c$ in $t$.
  Removing $\eta$ from $t$ splits the tree into a component with $v$ and a component which contains $c$ and another vertex in part $\{c,*,*,*\}$, call it $y$, breaking $\phi$ into three components.
  Since $n_v \not\in t$, adding $\eta_v$ connects the component with $v$ to the other tree in $\phi$ (corresponding to part $\{x\}$) creating a new spanning 2-forest $\psi'$.
  Once again as $v$ is a leaf in $\psi$, removing $\eta_v$ from $\psi$ and adding $\eta$ to create $\psi'$ maintains the spanning tree structure.
  Thus swapping the edges $\eta_v$ and $\eta$ between $\psi$ and $\phi$ creates a new edge partition $(\psi', \phi')$,
  However, now $\psi'$ is compatible with the vertex bipartition $\{c,y\},\{x,*,*\}$.
  To get back to a vertex bipartition of the form $\{*\},\{c,*,*,*\}$, we notice that we are in the case where $n_c \not\in t$, and in particular, $n_c$ was in the other tree of $\phi$ which is now connected to $x$.
  Since $c$ can never be the control vertex, swapping the edges incident to $c$ between $\psi'$ and $\phi'$ gives a new edge partition $(\psi'', \phi'')$ now in the set $\mathcal{S}_{\{y\},\{c,x,*,*\}}$, which is of the correct form.
  Note that $y \neq x$ could be any of $\{a,b,d,e\}$.

  In all three cases, as we are always changing the edges incident to $c$ and/or $v$ in $\psi$ and $\phi$, the new edge partitions can never be identical to $(\psi, \phi)$.
  Thus this transformation is fixed-point free.

  Lastly, we prove that this transformation is indeed an involution by showing that the control vertex and the two edges being swapped remains the same in the new edge partition.
  Let $(\psi', \phi')$ be the transformed edge partition.
  For~\ref{S-step-1}, since $c$ is a 2-valent vertex, there is no choice of which edge is to be swapped and as both $n_c$ and $c$ are in the same tree in $\phi'$, we are in the same case of the transformation.
  Thus applying the transformation again returns us to $(\psi, \phi)$.

  For~\ref{S-step-2}, it will be easier to first formulate the control vertex in another way.
  We can define the control vertex $v$ to be the first common vertex in $t$ in the paths from exactly two of the vertices in part $\{c, *, *, *\} \setminus \{c\}$ to $c$.
  Under this formulation, we immediately see that any reconnecting of $v$ via $\eta_v$ does not change this property, and thus the control vertex remains the same after the swapping of edges incident to $v$.
  Additionally for ~\ref{S-step-out}, the swapping of edges incident to $c$ afterwards also does not affect this property of $v$ since $n_c$ was in the tree connected to $x$, and thus once again the control vertex remains the same in $\psi'$.

  For the edges, notice that~\ref{S-step-in} occurs when all neighbours of $v$ are in the same tree of $\phi$, and thus by the choice of $\eta$, all the neighbours of $v$ are still in the same tree of $\phi'$.
  Additionally in $\psi'$, $\eta_v$ is exactly the edge incident to $v$ on the cycle created by adding $\eta$, which is now the edge incident to $v$ in $\psi'$.
  For~\ref{S-step-out}, since we are reconnecting $v$ to the other tree of $\phi$ via $\eta_v$, the new neighbour of $v$ in $\psi'$ is now in a different tree than $v$ in $\phi'$.
  Additionally in $\psi'$, because $n_c$ was in the tree connected to $x$ in $\phi$, after the two swaps $\eta$ is exactly the edge incident to $v$ on the path to $c$.
  In both cases $(\psi', \phi')$ remain in the same case of ~\ref{S-step-2} as before with $\eta'_v = \eta$, as $v$ is a leaf in $\psi'$ and so there is only one choice of edge, and $\eta' = \eta_v$ as shown above.
  Thus, applying the transformation again returns us to $(\psi, \phi)$.

  Putting everything together, we indeed have a fixed-point free involution in all cases.
  Finally, this holds for any edge partition in sets of the form $\mathcal{S}_{\{x\},\{c,*,*,*\}}$ where we are either staying in the same set $\mathcal{S}_{\{x\},\{c,*,*,*\}}$ or swappping to a different set $\mathcal{S}_{\{y\},\{c,*,*,*\}}$ where $y \neq x$ could be any of $\{a,b,d,e\}$.
  Thus our involution is on the union of sets of the form $\mathcal{S}_{\{*\},\{c,*,*,*\}}$,
  \[ \mathcal{S}_{\{a\},\{b,c,d,e\}} \cup \mathcal{S}_{\{b\},\{a,c,d,e\}} \cup \mathcal{S}_{\{d\},\{a,b,c,e\}} \cup \mathcal{S}_{\{e\},\{a,b,c,d\}}, \]
  and therefore, the size of this set must be even, giving the equation
  \[ {s}_{\{a\},\{b,c,d,e\}} + {s}_{\{b\},\{a,c,d,e\}} + {s}_{\{d\},\{a,b,c,e\}} + {s}_{\{e\},\{a,b,c,d\}} \equiv 0 \mod 2. \]
\end{proof}

As in our discussion before Theorem~\ref{S-bijection}, since swapping around $c$ gives fixed-point free involutions $\mathcal{S}_{\{x\},\{c,*,*,*\}} \cup \mathcal{S}_{\{c,x\},\{*,*,*\}}$ for $x \in \{a,b,d,e\}$, we immediately get
\[ {s}_{\{a,c\},\{b,d,e\}} + {s}_{\{b,c\},\{a,d,e\}} + {s}_{\{c,d\},\{a,b,e\}} + {s}_{\{c,e\},\{a,b,d\}} \equiv 0 \mod 2 \]
by adding all the equations together and using the result from Theorem~\ref{S-bijection}.

Alternatively, notice that we could use the same fixed-point free involution on the union of sets of the form $\mathcal{S}_{\{c,*\},\{*,*,*\}}$ except with a slight modification.
In~\ref{S-step-2} where the two neighbours of $c$ are in different trees in the spanning 2-forest, we first swap around the edges of $c$, giving a new spanning 2-forest compatible with $\{*\},\{c,*,*,*\}$, and proceed as before with finding the control vertex and edges to swap.
Instead of swapping around $c$ again in ~\ref{S-step-out}, we do so in ~\ref{S-step-in} to ensure we get edge partitions in sets of the form $\mathcal{S}_{\{c,*\},\{*,*,*\}}$.
As everything else remains the same, we obtain the following result.

\tpointn{Corollary}\label{S-bij-cor}
\statement[eq]{
  There is a fixed-point free involution on
  \[ \mathcal{S}_{\{a,c\},\{b,d,e\}} \cup \mathcal{S}_{\{b,c\},\{a,d,e\}} \cup \mathcal{S}_{\{c,d\},\{a,b,e\}} \cup \mathcal{S}_{\{c,e\},\{a,b,d\}}, \]
  and thus we have
  \[ {s}_{\{a,c\},\{b,d,e\}} + {s}_{\{b,c\},\{a,d,e\}} + {s}_{\{c,d\},\{a,b,e\}} + {s}_{\{c,e\},\{a,b,d\}} \equiv 0 \mod 2. \]
}

\bpoint{Completing the \texorpdfstring{$S$}{S}-case}

With the two swapping arguments, first with swapping around $c$ and then with swapping around a control vertex, we have covered all the $\mathcal{S}$ sets in Equation~\eqref{eq:S-counts}!
In fact they are really the same argument as we can think of the 2-valent vertex $c$ as acting as a special type of control vertex.
Thus we have everything we need to prove the $c_2$ completion conjecture for $p=2$ in the $S$-case.

\tpointn{Theorem}\label{S-case}
\statement[eq]{
  Let $G$ be a connected 4-regular graph.
  Let $v$ and $w$ be adjacent vertices of $G$ such that they share one common neighbour.
  Then,
  \[ c_2^{(2)}(G - v) = c_2^{(2)}(G - w). \]
}
\begin{proof}
  Combining Proposition~\ref{S-eqs}, Theorem~\ref{S-swapc} and Corollary~\ref{S-bij-cor} gives the result.
\end{proof}
\vspace{\baselineskip}

\section{The \texorpdfstring{$R$}{R}-case}\label{S:R-case}

The last case we need to deal with is the $R$-case.
Let $G$ be a connected 4-regular graph, and let $v$ and $w$ be two adjacent vertices of $G$ such that they do not share any neighbours.
Let $R = G - \{v,w\}$ be the graph obtained by removing vertices $v$ and $w$ from $G$, the last grey blob in Figure~\ref{fig:TSR}, with the neighbours of $v$ and $w$ labelled $\{a,b,c,d,e,f\}$ as in the figure.
Here $w$ has neighbours $\{a,b,c,v\}$, and $v$ has neighbours $\{d,e,f,w\}$.

The difficulty of the $R$-case stems from the fact that $v$ and $w$ no longer share any neighbours, and thus there is no distinguishing feature of $R$ that we can readily exploit.
However, the lack of specialness of the vertices $\{a,b,c,d,e,f\}$ is what lends itself to the symmetric nature of the $R$-case.
To prove the conjecture for $p=2$, we needed to use this symmetry to our advantage.

As we no longer have any 2-valent vertex to swap around, we look towards the more general control-vertex argument.
Notice in the swapping around a control vertex argument for the $S$-case there was a symmetric flavour in how the new edge partitions that we ended up with could have been from any $\mathcal{S}$ set of a similar form to the original, but we didn't need to know exactly which one.
Using a simplified version of the involution in the $S$-case, which was actually first discovered for the $R$-case, we use this symmetry to partially deal with the $R$-case.

However, like in the $S$-case, this particular control vertex argument was not enough.
Once again an argument involving compatible cycles (see \textsection 5 of~\cite{specialc2}) was used for the remaining sets when $G$ had an odd number of vertices.
To get rid of this parity requirement, a more complicated control vertex argument is necessary this time involving multiple control vertices.

Once again, to give a complete, self-contained proof of the $R$-case, we begin with the same set-up and partial results as from~\cite{specialc2}, and in Section~\ref{SS:R-new} we present our new results, with allows us to finish the $R$-case for $p=2$! \\

\bpoint{Set-up}

As with the previous cases, we start with defining the sets of edge partitions that we are counting.

\tpointn{Definition} (Definition 3.1 of~\cite{specialc2})\label{R-defs}
\statement{
  Suppose $P$ is a bipartition of $\{a,b,c,d,e,f\}$. \\
  Let $\mathcal{R}_P$ be the set of bipartitions $(\psi, \phi)$ of the edges of $R$ such that $\psi$ is a spanning tree and $\phi$ is a spanning 2-forest compatible with $P$.
  Let $r_P = \abs{\mathcal{R}_P}$.
}

This time because of the symmetry of the neighbours of $v$ and $w$, and since we are counting modulo $2$, instead of picking specific vertex bipartitions for the spanning 2-forests we can add them all together.

\tpointn{Proposition} (Proposition 3.2, 3.3 of~\cite{specialc2})\label{R-eqs}
\statement[eq]{
  When $v$ and $w$ have no common neighbours
  \[ c_2^{(2)}(G - v) = r_{\{a\},\{b,c\}} + r_{\{b\},\{a,c\}} + r_{\{c\},\{a,b\}} \mod 2, \]
  \[ c_2^{(2)}(G - w) = r_{\{d\},\{e,f\}} + r_{\{e\},\{d,f\}} + r_{\{f\},\{d,e\}} \mod 2, \]
  and thus we have
  \begin{align*}\refstepcounter{equation}\tag{\theequation}\label{eq:R-counts}
    c_2^{(2)}(G - v) - c_2^{(2)}(G - w)
      =\ &{r}_{\{a\},\{b,c,d,e,f\}} + {r}_{\{b,c\},\{a,d,e,f\}} \\
      &+ {r}_{\{b\},\{a,c,d,e,f\}} + {r}_{\{a,c\},\{b,d,e,f\}} \\
      &+ {r}_{\{c\},\{a,b,d,e,f\}} + {r}_{\{a,b\},\{c,d,e,f\}} \\
      &+ {r}_{\{d\},\{a,b,c,e,f\}} + {r}_{\{e,f\},\{a,b,c,d\}} \\
      &+ {r}_{\{e\},\{a,b,c,d,f\}} + {r}_{\{d,f\},\{a,b,c,e\}} \\
      &+ {r}_{\{f\},\{a,b,c,d,e\}} + {r}_{\{d,e\},\{a,b,c,f\}} \mod 2.
  \end{align*}
}
\begin{proof}
  Like in the proof of Proposition~\ref{T-eqs}, we can pick specific vertex partitions for the $G-v$ and $G-w$ cases to equivalently determine each $c_2$ via counting edge partitions into a spanning tree and a spanning 2-forest compatible with the respective vertex partition.
  However, since the same argument holds for all vertex partitions of $\{a,b,c\}$ and $\{d,e,f\}$ into the form $\{*\},\{*,*\}$, we have that
  \begin{align*}
    c_2^{(2)}(G - v)
      &\equiv r_{\{a\},\{b,c\}}
      \equiv r_{\{b\},\{a,c\}}
      \equiv r_{\{c\},\{a,b\}} \mod 2, \\
    c_2^{(2)}(G - w)
      &\equiv r_{\{d\},\{e,f\}}
      \equiv r_{\{e\},\{d,f\}}
      \equiv r_{\{f\},\{d,e\}} \mod 2.
  \end{align*}
  Since we are working modulo $2$, we can add each of the three choices together
  \begin{align*}
    c_2^{(2)}(G - v) &\equiv r_{\{a\},\{b,c\}} + r_{\{b\},\{a,c\}} + r_{\{c\},\{a,b\}} \mod 2, \\
    c_2^{(2)}(G - w) &\equiv r_{\{d\},\{e,f\}} + r_{\{e\},\{d,f\}} + r_{\{f\},\{d,e\}} \mod 2.
  \end{align*}

  Then, enumerating over all possibilities for the vertex set $\{a,b,c,d,e,f\}$ in the spanning 2-forests, we obtain the final equation by subtracting the expanded counts for $G-v$ and $G-w$ and simplifying modulo $2$.
  To illustrate this, expanding the counts for $r_{\{a\},\{b,c\}}$ gives
  \begin{align*}
    r_{\{a\}, \{b,c\}}
      =\ &r_{\{a\}, \{b,c,d,e,f\}} \\
      &+ r_{\{a,d\}, \{b,c,e,f\}} + r_{\{a,e\}, \{b,c,d,f\}} + r_{\{a,f\}, \{b,c,d,e\}} \\
      &+ r_{\{a,d,e\}, \{b,c,f\}} + r_{\{a,d,f\}, \{b,c,e\}} + r_{\{a,e,f\}, \{b,c,d\}} \\
      &+ r_{\{a,d,e,f\}, \{b,c\}}.
  \end{align*}
  Now notice that the terms in the middle two lines appear in the expanded counts for the terms in the $G-w$ case, for example, $r_{\{a,d\},\{b,c,e,f\}}$ and $r_{\{a,e,f\},\{b,c,d\}}$ both appear in the expansion of $r_{\{d\},\{e,f\}}$ and similarily for the other four terms.
  Thus we would be left with $r_{\{a\},\{b,c,d,e,f\}}$ and $r_{\{b,c\},\{a,d,e,f\}}$ after the subtraction.

  As this holds for all $r_{\{*\},\{*,*\}}$, where the $*$'s are partitioning either $\{a,b,c\}$ or $\{d,e,f\}$, after subtracting the $G-v$ and $G-w$ counts we are left with terms of the form $r_{\{*\},\{*,*,*,*,*\}}$ and $r_{\{y,z\},\{x,*,*,*\}}$, where $\{x,y,z\}$ is either $\{a,b,c\}$ or $\{d,e,f\}$.
\end{proof}

\bpoint{Swapping around a control vertex}

Immediately from Equation~\eqref{eq:R-counts} we see that there are 6 terms corresponding to sets of the form $\mathcal{R}_{\{y,z\},\{x,*,*,*\}}$ where $\{x,y,z\}$ is either $\{a,b,c\}$ or $\{d,e,f\}$.
The thing to note is that these vertex bipartitions for the 2-forests are exactly in the form where we can apply Lemma~\ref{control-vertex}, this time using the other specification for the special vertex $x$.
As we already have an idea of how to deal with a control vertex from this lemma, we can use a simplified version of the involution in Theorem~\ref{S-bijection}.

\tpointn{Theorem} (Lemma 4.4 of~\cite{specialc2})\label{R-control-bij}
\statement[eq]{
  There is a fixed-point free involution on
  \begin{align*}
    \mathcal{R}_{\{a,b\},\{c,d,e,f\}} &\cup \mathcal{R}_{\{a,c\},\{b,d,e,f\}} \cup \mathcal{R}_{\{b,c\},\{b,d,e,f\}} \\
    &\cup \mathcal{R}_{\{d,e\},\{a,b,c,f\}} \cup \mathcal{R}_{\{d,f\},\{a,b,c,e\}} \cup \mathcal{R}_{\{e,f\},\{a,b,c,d\}},
  \end{align*}
  and thus we have
  \begin{align*}
    {r}_{\{a,b\},\{c,d,e,f\}} &+ {r}_{\{a,c\},\{b,d,e,f\}} + {r}_{\{b,c\},\{b,d,e,f\}} \\
    &+ {r}_{\{d,e\},\{a,b,c,f\}} + {r}_{\{d,f\},\{a,b,c,e\}} + {r}_{\{e,f\},\{a,b,c,d\}} \equiv 0 \mod 2.
  \end{align*}
}
\begin{proof}
  Consider the edge partitions in the union of sets of the form $\mathcal{R}_{\{y,z\},\{x,*,*,*\}}$ where $\{x,y,z\}$ is either the trio $\{a,b,c\}$ or $\{d,e,f\}$, and the $*$'s are the other three vertices in $\{a,b,c,d,e,f\} \setminus \{x,y,z\}$.
  Since we can always distinguish the vertex $x$ that is alone from its trio, we can choose $x$ uniquely to be the special vertex.
  Thus these edge partitions satisfy the conditions of Lemma~\ref{control-vertex}.
  This time because of the special form of the vertex bipartitions, we will specify for $x$ the second property from the lemma, that $x$ is either the control vertex or $x$ is in one of the singleton components after removing the control vertex from its tree in the spanning 2-forest.

  Let $(\psi, \phi)$ be an edge partition in any set of the form $\mathcal{R}_{\{y,z\},\{x,*,*,*\}}$, and let $v$ be its control vertex.
  Let $t$ be the tree in $\phi$ corresponding to $\{x,*,*,*\}$, which also contains $v$.
  Since $v$ is a leaf in $\psi$, let $\eta_v$ be the edge incident to $v$ and $n_v$ the neighbour of $v$ in $\psi$.

  Now to describe the involution which gives a new edge partition $(\psi', \phi')$:
  \begin{enumerate}
    \item\label{nv-in} \emph{Swapping stays in} -- If $n_v \in t$, so the control vertex and its neighbour in $\psi$ are in the same tree of $\phi$, let $\eta$ be the edge incident to $v$ in $t$ in the path to $n_v$.\\
      Then, swap the edges $\eta_v$ and $\eta$ between $\psi$ and $\phi$.
    \item\label{nv-out} \emph{Swapping goes out} -- Otherwise if $n_v \not\in t$, so the control vertex and its neighbour in $\psi$ are in different trees of $\phi$, let $\eta$ be the edge incident to $v$ in $t$ in the path to the component of $t - v$ that contains exactly two vertices of $\{x,*,*,*\}$; in Figure~\ref{fig:R-control} this is the edge towards the 2-part indicated by the dashed line. \\
      Then, swap the edges $\eta_v$ and $\eta$ between $\psi$ and $\phi$.
  \end{enumerate}

  Notice that this is essentially phase~\ref{S-step-2} in Theorem~\ref{S-bijection}, just without the extra swapping of edges in~\ref{S-step-out}, and thus most of the proof still holds for this simplified involution.
  Firstly, as we are always changing which edge is incident to $v$ in $\psi$ and $\phi$, $(\psi', \phi') \neq (\psi, \phi)$, and so this transformation is fixed-point free.

  Now for~\ref{nv-in}, we can argue exactly as we did for~\ref{S-step-in} in Theorem~\ref{S-bijection} to show that $(\psi', \phi')$ stays a valid edge partition in the same set $\mathcal{R}_{\{y,z\},\{x,*,*,*\}}$ as $(\psi, \phi)$.
  Furthermore, it was shown that if $v$ is still the control vertex after the transformation, which we will prove shortly, then $(\psi',\phi')$ remains in case~\ref{nv-in} and applying the transformation again brings us back to $(\psi, \phi)$.

  For~\ref{nv-out}, a little more work is needed to show that $(\psi', \phi')$ is in a valid $\mathcal{R}$ set but the general framework is the same.
  By the choice of $\eta$, removing $\eta$ from $t$ splits the tree into two components that partitions $\{x,*,*,*\}$ into the form $\{x,x'\}, \{*,*\}$ for some $x' \in \{a,b,c,d,e,f\} \setminus \{x,y,z\}$.
  Since $n_v \not\in t$, adding the edge $\eta_v$ connects the component from $t$ containing $v$ to the other tree of $\phi$ (corresponding to part $\{y,z\}$), to create a new spanning 2-forest $\phi'$.
  In particular, as $v$ is in the component corresponding to $\{x, x'\}$, $\phi'$ is now compatible with the vertex bipartition $\{*,*\},\{x',x,y,z\}$.
  As $v$ is a leaf in $\psi$, removing $\eta_v$ from $\psi$ and adding $\eta$ to create $\psi'$ maintains the spanning tree structure.
  Thus we have that $(\psi', \phi')$ is a valid edge partition in the set $\mathcal{R}_{\{*,*\},\{x',x,y,z\}}$.
  Notice that as $\{x,y,z\}$ is either $\{a,b,c\}$ or $\{d,e,f\}$, the other trio must be $\{x',*,*\}$ in any order, and so $\{*,*\}, \{x',x,y,z\}$ is in the form we want.

  Once again if $v$ is remains the control vertex, $v$ and its neighbour in $\psi'$ are in different trees of $\phi'$, and thus we stay in case~\ref{nv-out}.
  As $x'$ is now the special vertex, $\eta_v$ is exactly the edge incident to $v$ connecting the component with $\{x',x\}$ and the component with $\{y,z\}$ in $\phi'$, that is $\eta' = \eta_v$.
  Thus applying the transformation again returns us to $(\psi, \phi)$.

  Finally what's left to prove is that the control vertex $v$ stays the same across the transformation.
  As in the proof of Theorem~\ref{S-bijection}, we can formulate the control vertex in an alternate way.
  This time $v$ is the first common vertex in $t$ in all three paths from $x$ to the other vertices in part $\{x,*,*,*\}$.
  Under this formulation, we immediately see that this property does not change after the transformation, where now $v$ is the first common vertex in the tree corresponding to part $\{x',x,y,z\}$ in $\phi'$ from $x'$ to the other three vertices.
  Thus $v$ remains the control vertex in $(\psi', \phi')$.

  Putting everything together, we indeed have a fixed-point free involution on the union of sets of the form $\mathcal{R}_{\{y,z\},\{x,*,*,*\}}$ where $\{x,y,z\}$ is either $\{a,b,c\}$ or $\{d,e,f\}$,
  \begin{align*}
    \mathcal{R}_{\{a,b\},\{c,d,e,f\}} &\cup \mathcal{R}_{\{a,c\},\{b,d,e,f\}} \cup \mathcal{R}_{\{b,c\},\{b,d,e,f\}} \\
    &\cup \mathcal{R}_{\{d,e\},\{a,b,c,f\}} \cup \mathcal{R}_{\{d,f\},\{a,b,c,e\}} \cup \mathcal{R}_{\{e,f\},\{a,b,c,d\}},
  \end{align*}
  and therefore, the size of this set must be even, giving the equation
  \begin{align*}
    {r}_{\{a,b\},\{c,d,e,f\}} &+ {r}_{\{a,c\},\{b,d,e,f\}} + {r}_{\{b,c\},\{b,d,e,f\}} \\
    &+ {r}_{\{d,e\},\{a,b,c,f\}} + {r}_{\{d,f\},\{a,b,c,e\}} + {r}_{\{e,f\},\{a,b,c,d\}} \equiv 0 \mod 2.
  \end{align*}
\end{proof}

\bpoint{Two control vertices}\label{SS:R-new}

For the remaining 6 sets in Equation~\eqref{eq:R-counts}, we are no longer in a situation where we can apply Lemma~\ref{control-vertex} since these sets are of the form $\mathcal{R}_{\{x\},\{*,*,*,*,*\}}$ for $x \in \{a,b,c,d,e,f\}$.
What we need is a new notion of "control vertex", this time for a tree with five marked vertices instead of four.
However, as we no longer have any distinguishable vertices in this 5-part, we could not appropriately define one unique vertex as the control vertex.
Rather two control vertices are necessary to get an involution on the correct sets.

\begin{figure}[t]
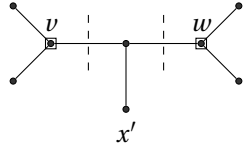

  \centering
  \ShapeFiveCut
  \caption[The two control vertices for Lemma~\ref{R-control}.]{The two control vertices $v$ and $w$ for Lemma~\ref{R-control}. The unlabelled leaves and $x'$ are all vertices in $p$. The dashed lines indicate the edges we will be picking for one of the involutions in the $R$-case. Some of the leaf edges are possibly contracted, with the exception that the edges of two leaves adjacent to the same vertex cannot simultaneously both be contracted.}
  \label{fig:two-control}
\end{figure}

\tpointn{Lemma}\label{R-control}
\statement{
  Suppose we have a bipartition of $\{a,b,c,d,e,f\}$ of the form $\{x\}, \{*,*,*, *, * \}$.
  Let $(\psi,\phi)$ be an edge partition in $\mathcal{R}_{\{x\},\{*,*,*,*,*\}}$ where $\psi$ is a spanning tree and $\phi$ is a spanning 2-forest compatible with $\{x\},\{*,*,*,*,*\}$.
  Let $t$ be the tree corresponding to part $p = \{*,*, *, *, *\}$ in the 2-forest $\phi$. \\\\
  Then, there is a unique pair of non-adjacent vertices $v$ and $w$ called the \textbf{control vertices}, such that removing $v$ and $w$ from $t$ partitions $p$ into all singletons.
  Figure~\ref{fig:two-control} depicts the two control vertices.\\
  Furthermore, in $t$ each control vertex will either be
  \begin{itemize}
    \item 2-valent and in $\{a,b,c,d,e,f\} \setminus \{x\}$, or
    \item 3-valent and not a vertex in $\{a,b,c,d,e,f\}$.
  \end{itemize}
  In all cases, the control vertices are leaves in the spanning tree $\psi$.
}
\begin{proof}
  In a similar fashion to the proof of Lemma~\ref{control-vertex}, we start by looking at the tree $t$ and the subtree created by taking the union of all the paths in $t$ between any two vertices in $p$.
  Since $p$ has five vertices this subtree will be of the form
  \[ \ShapeFive \]

  where the leaves are vertices in $p$, the edges are paths in $t$, and some of the paths may be contracted.

  Since $\psi$ is a spanning tree, the vertices in $\phi$, and thus $t$, are at most 3-valent, and the vertices in $\{a,b,c,d,e,f\}$, and thus $p$, are at most 2-valent.
  In particular, we have that the above is the only possible form for the subtree, where some of the paths may be contracted.
  Note trivially, the paths between two leaves adjacent to the same non-leaf vertex cannot be contracted simultaneously.
  Now, the valency restrictions also mean that in the subtree form above, the paths from $v$ and $w$ to the middle non-leaf vertex cannot be contracted.
  Thus $v \neq w$, and furthermore, they are not adjacent in $\phi$ i.e. there must be at least two edges in the path from $v$ to $w$.
  By the valency of the vertices in $\phi$, $v$ and $w$ can only be 2-valent if they are in the set $\{a,b,c,d,e,f\}$ or otherwise 3-valent, making them both leaves in the spanning tree $\psi$.
  As $R$ has more than two vertices, $v$ and $w$ are also not adjacent in $\psi$, and thus in $R$ the control vertices are non-adjacent.

  Finally, removing $v$ and $w$ from $t$ clearly partitions $p$ into all singletons and these are the unique pair of vertices which do so.
  The number of components of $t - \{v,w\}$ will depend on whether $v$ or $w$ are in the set $\{a,b,c,d,e,f\}$.
\end{proof}

Now that we've established the notion of control vertices, we can use a swapping argument once more for one final fixed-point free involution to cover the remaining sets in Equation~\eqref{eq:R-counts}.
As there are two vertices to swap around, like in Theorem~\ref{S-bijection}, we need a more complicated involution than previously.

\tpointn{Theorem}\label{R-bijection}
\statement[eq]{
  There is a fixed-point free involution on
  \begin{align*}
    \mathcal{R}_{\{a\},\{b,c,d,e,f\}} &\cup \mathcal{R}_{\{b\},\{a,c,d,e,f\}} \cup \mathcal{R}_{\{c\},\{a,b,d,e,f\}} \\
    &\cup \mathcal{R}_{\{d\},\{a,b,c,e,f\}} \cup \mathcal{R}_{\{e\},\{a,b,c,d,f\}} \cup \mathcal{R}_{\{f\},\{a,b,c,d,e\}},
  \end{align*}
  and thus we have
  \begin{align*}
    {r}_{\{a\},\{b,c,d,e,f\}} &+ {r}_{\{b\},\{a,c,d,e,f\}} + {r}_{\{c\},\{a,b,d,e,f\}} \\
    &+ {r}_{\{d\},\{a,b,c,e,f\}} + {r}_{\{e\},\{a,b,c,d,f\}} + {r}_{\{f\},\{a,b,c,d,e\}} \equiv 0 \mod 2.
  \end{align*}
}
\begin{proof}
  By construction the edge partitions in the union of sets of the form $\mathcal{R}_{\{x\},\{*,*,*,*,*\}}$, where $x \in \{a,b,c,d,e,f\}$, and the $*$'s are the rest of the vertices $\{a,b,c,d,e,f\} \setminus \{x\}$, satisfy the conditions of Lemma~\ref{R-control}.
  Let $(\psi, \phi)$ be an edge partition in any set of the form $\mathcal{R}_{\{x\},\{*,*,*,*,*\}}$, and let $v,w$ be its control vertices.
  Let $t$ be the tree in $\phi$ corresponding to part $\{*,*,*,*,*\}$ which contains $v$ and $w$.
  Since the control vertices are leaves in $\psi$, let $\eta_v$ and $\eta_w$ be the edges incident to $v$ and $w$ in $\psi$, respectively, and let $n_v$ and $n_w$ be their neighbours.

  Consider the following involution giving a new edge partition $(\psi',\phi')$:
  \begin{enumerate}[label=(\arabic*)]
    \item\label{R-step-1} \emph{Swapping either stays in} -- If $n_v \in t$ or $n_w \in t$:
      \begin{enumerate}[label=(\roman*), ref=(1\roman*)]
        \item\label{R-step-vin} \emph{$v$ stays in} -- If $n_v \in t$, so control vertex $v$ and its neighbour in $\psi$ are in the same tree of $\phi$, let $\eta_1$ be the edge incident to $v$ in $t$ in the path to $n_v$.\\
          Then, swap the edges $\eta_v$ and $\eta_1$ between $\psi$ and $\phi$.
        \item\label{R-step-win} \emph{$w$ stays in} -- If $n_w \in t$, so control vertex $w$ and its neighbour in $\psi$ are in the same tree of $\phi$, let $\eta_2$ be the edge incident to $w$ in $t$ in the path to $n_w$.\\
          Then, swap the edges $\eta_w$ and $\eta_2$ between $\psi$ and $\phi$.
      \end{enumerate}
      Here we are using non-exclusive or, so if $n_v, n_w \in t$, then swap around both $v$ and $w$ as above.
    \item\label{R-step-2} \emph{Swapping both goes out} -- Otherwise we have $n_v, n_w \not\in t$, where the control vertices and their neighbours in $\psi$ are in different trees of $\phi$.
      Let $\eta_1$ be the edge incident to $v$ in $t$ in the path to $w$, and let $\eta_2$ be the edge incident to $w$ in $t$ in the path to $v$; in Figure~\ref{fig:two-control} these are the edges indicated by the dashed lines. \\
      Then, swap edges $\eta_v$ and $\eta_w$ with $\eta_1$ and $\eta_2$ between $\psi$ and $\phi$.
  \end{enumerate}

  Firstly, as we are always changing which edges are incident to $v$ and/or $w$, $(\psi', \phi')$ must be different from $(\psi,\phi)$.
  Thus this transformation is fixed-point free.

  For~\ref{R-step-1}, both cases are symmetric so without loss of generality assume $n_v \in t$.
  Then, removing edge $\eta_1$ from $\phi$ breaks $t$ into two components and because of how $\eta_1$ was chosen, $n_v$ is now in a different component than $v$.
  As $v$ and its neighbour in $\psi$, $n_v$, are in the same tree of $\phi$, adding edge $\eta_v$ reconnects the two components from $t$.
  By the choice of $\eta_1$, we cannot create any cycles by adding $\eta_v$, and thus $\phi'$ is a spanning 2-forest which is compatible with the same vertex bipartition as $\phi$.
  Since $v$ is a leaf in $\psi$, removing edge $\eta_v$ and adding edge $\eta_1$ to create $\psi'$ maintains the spanning tree structure.
  Thus $(\psi', \phi')$ is a valid edge partition in the same $\mathcal{R}$ set as $(\psi, \phi)$.
  We also have if the control vertices stay the same, $(\psi', \phi')$ will still have the property that $v$ and its neighbour in $\psi'$ are in the same tree of $\phi'$, and thus applying the transformation again returns to $(\psi, \phi)$.

  When both $n_v, n_w \in t$, since $v\neq w$ and they are not adjacent, the four edges in the swap are all distinct; $\eta_1 \neq \eta_2$ and $\eta_v \neq \eta_w$.
  Thus we can simultaneously swap edges $\eta_v$ and $\eta_w$ with edges $\eta_1$ and $\eta_2$, respectively.
  Note that we need to swap around both control vertices as there is no way to distinguish $v$ from $w$, and thus no way to ensure we get an involution if we arbitrarily picked one to swap edges around.
  Applying the above swapping argument twice, we have that $(\psi', \phi')$ is a valid edge partition in the same $\mathcal{R}$ set as before and applying the transformation again returns to $(\psi, \phi)$.

  For~\ref{R-step-2} when $n_v, n_w \not\in t$, once again we know that the four edges in the swap are all distinct by the same reasoning as when $n_v, n_w \in t$.
  Then, removing edges $\eta_1$ and $\eta_2$ breaks $t$ into three subtrees (thus breaking $\phi$ into four components); one with $v$, one with $w$, and one with the vertex $x' \in \{a,b,c,d,e,f\} \setminus \{x\}$ as in Figure~\ref{fig:two-control}.
  Now since $v$ and $n_v$ are in different trees of $\phi$, adding $\eta_v$ connects the subtree from $t - \{\eta_1, \eta_2\}$ with $v$ to the tree in $\phi$ corresponding to part $\{x\}$.
  Similarily, as $w$ and $n_w$ are in different trees of $\phi$, adding $\eta_w$ connects the subtree with $w$ also to the tree with $x$ in $\phi$.

  As the last subtree, which contains $x'$, never gets connected to any other component, this transformation creates a new spanning 2-forest $\phi'$ which is compatible with the vertex bipartition $\{x'\}, \{x, *,*,*,*\}$.
  Since $v,w$ are leaves in $\psi$, removing edges $\eta_v$ and $\eta_w$ from $\psi$ disconnects the control vertices and adding edges $\eta_1$ and $\eta_2$ reconnects them, creating $\psi'$ which is once again a spanning tree.
  Thus $(\psi', \phi')$ is a valid edge partition in $\mathcal{R}_{\{x'\},\{x,*,*,*,*\}}$, which is in the correct form.
  Once again, if the control vertices stay the same, then their neighbours in $\psi$ will be in a different tree of $\phi'$ landing us back in ~\ref{R-step-2}.
  Thus applying the transformation again brings us back to $(\psi, \phi)$.

  What's left of the proof is to make sure $v$ and $w$ are still the control vertices in $(\psi', \phi')$.
  Notice we can formulate the control vertices as the two vertices in $t$ such that for a specific $x'$ in part $\{*,*,*,*,*\}$, each vertex is the last common vertex on the paths from $x'$ to the exactly two of the other vertices in the part.
  The vertex $x'$ is as depicted in Figure~\ref{fig:two-control}.
  Under this formulation, we immediately see that in~\ref{R-step-2} as $\phi'$ is created by connecting subtrees with $v$ and $w$ to the tree corresponding to part $\{x\}$ in $\phi$, the vertex $x$ now acts as the $x'$ for $\phi'$.
  Thus $v$ and $w$ remain the control vertices for $(\psi',\phi')$.

  For~\ref{R-step-1}, depending on where $n_v$ or $n_w$ were in the tree $t$, the vertex acting as $x'$ in $\phi'$ may no longer be the same $x'$ from $\phi$.
  However, one can check that in all possible situations, the property for the control vertices is still satisfied by $v$ and $w$, and thus $v$ and $w$ are still the control vertices.

  Finally putting everything together, we indeed get a fixed-point free involution on the union of sets of the form $\mathcal{R}_{\{x\},\{*,*,*,*,*\}}$,
  \begin{align*}
    \mathcal{R}_{\{a\},\{b,c,d,e,f\}} &\cup \mathcal{R}_{\{b\},\{a,c,d,e,f\}} \cup \mathcal{R}_{\{c\},\{a,b,d,e,f\}} \\
    &\cup \mathcal{R}_{\{d\},\{a,b,c,e,f\}} \cup \mathcal{R}_{\{e\},\{a,b,c,d,f\}} \cup \mathcal{R}_{\{f\},\{a,b,c,d,e\}},
  \end{align*}
  and therefore, the size of this must be even
  \begin{align*}
    {r}_{\{a\},\{b,c,d,e,f\}} &+ {r}_{\{b\},\{a,c,d,e,f\}} + {r}_{\{c\},\{a,b,d,e,f\}} \\
    &+ {r}_{\{d\},\{a,b,c,e,f\}} + {r}_{\{e\},\{a,b,c,d,f\}} + {r}_{\{f\},\{a,b,c,d,e\}} \equiv 0 \mod 2.
  \end{align*}
\end{proof}

As a remark on why two control vertices were necessary, looking at Figure~\ref{fig:two-control} one would first be tempted to pick the unlabelled non-leaf vertex as the control vertex, call it $u$, and cut the edge towards $x'$.
However, in some edge partitions, this path may be contracted with $u$ being $x'$ itself.
In this case there is no obvious way to choose an edge incident to $u$ to remove.

To mitigate this, one would then decide to cut the two edges incident to $u$ that is not towards $x'$.
However, in this case by valency arguments, $u$ is a leaf in the corresponding spanning tree, and thus there is only one edge available to swap two edges with, which is not possible.
Thus, we cannot pick $u$ to be the control vertex.

Looking to $v$ or $w$, as there are no special vertices in this 5-part of the vertex bipartition, there would be no way to distinguish between $v$ and $w$ to pick just one control vertex.
Even if there were, only picking one of $v$ or $w$ to swap edges around would not swap to a spanning 2-forest compatible with a vertex bipartition that we care about.
Therefore, we really did need to pick both vertices $v$ and $w$ to be control vertices!

Notice that we can think of this involution involving two control vertices as a generalized version of the involution used in Theorem~\ref{S-bijection}.
For that involution, the 2-valent vertex $c$ acts as a second "control vertex". \\

\bpoint{Completing the \texorpdfstring{$R$}{R}-case}

While there were no special vertices in $\{a,b,c,d,e,f\}$ that we could readily use to swap between specific vertex bipartitions, we were able to take advantage of the symmetric nature of the $R$-case and use control vertices to find involutions where we only cared about the form of the vertex bipartition we may be swapping our edge partitions to.
 With the two swapping arguments, one involving one control vertex and the other involving two control vertices, we have everything we need to prove the $c_2$ completion conjecture for $p=2$ in the $R$-case.

\tpointn{Theorem}\label{R-case}
\statement[eq]{
  Let $G$ be a connected 4-regular graph.
  Let $v$ and $w$ be adjacent vertices of $G$ such that they do not share any common neighbours.
  Then,
  \[ c_2^{(2)}(G - v) = c_2^{(2)}(G - w). \]
}
\begin{proof}
  Together Proposition~\ref{R-eqs}, Theorem~\ref{R-control-bij} and Theorem~\ref{R-bijection} gives the result.
\end{proof}
\vspace{\baselineskip}

\section{Completing the \texorpdfstring{$p=2$}{p=2} case}\label{S:completing-conjecture}

To prove the $c_2$ completion conjecture for $p=2$, in all three cases the underlying idea remained the same: finding a vertex, or vertices, and picking particular edges incident to them to swap between the two parts of an edge partition, where this action gives a fixed-point free involution on specific sets of edge partitions.

In the $T$-case, we found that we could use the two 2-valent vertices to swap edges around.
In the $S$-case, we reused the 2-valent vertex swapping argument and generalized it to give the notion of swapping around a particular vertex called the control vertex.
We found that for a second involution, we needed a two-phase swapping process where we used both the control vertex and a 2-valent vertex.
Finally in the $R$-case, we reused the control vertex argument and found that we needed to extend it to include two control vertices for a second involution, this extension itself being a generalized version of the second involution in the $S$-case.

We note here that the two new bijections, giving Theorems~\ref{S-bijection} and~\ref{R-bijection} in the $S$ and $R$ cases, respectively, were implemented in \texttt{Sage}~\cite{sagemath} to verify our results on some small graphs.

Assembling everything together, we finally complete the $c_2$ completion conjecture for $p=2$!

\tpointn{Theorem}\label{p2-completed}
\statement[eq]{
  Let $G$ be a connected 4-regular graph, and let $v$ and $w$ be vertices of $G$.
  Then,
  \[ c_2^{(2)}(G - v) = c_2^{(2)}(G - w). \]
}
\begin{proof}
  As in the discussion in Section~\ref{S:conjecture}, since $G$ is connected there is path between any two vertices of $G$, and thus it suffices to prove the conjecture when $v$ and $w$ are adjacent vertices.
  For any two non-adjacent vertices, we could then just follow any path between them, getting a chain of equivalences to obtain the required result.

  Then, as $G$ is a 4-regular graph, the four possibilities for the neighbours of $v$ and $w$ are:
  \begin{itemize}
    \item $v$ and $w$ share all neighbours, or
    \item $v$ and $w$ have exactly two common neighbours -- $T$-case, or
    \item $v$ and $w$ have only one common neighbour -- $S$-case, or
    \item $v$ and $w$ do not share any neighbours -- $R$-case.
  \end{itemize}

  When $v$ and $w$ share all neighbours, notice that $G - v$ and $G - w$ are isomorphic.
  Thus, trivially, their $c_2$'s are equivalent.
  When $v$ and $w$ share some neighbours or none at all, Corollary~\ref{T-case}, Theorem ~\ref{S-case}, and Theorem~\ref{R-case} give the result in the $T$, $S$, and $R$ cases, respectively.
  This completes the proof.
\end{proof}

\newpage
\newrefsegment
\part{Le diagrams}\label{P:le}

\chapter{On the positive Grassmannian and related objects}\label{C:lebackground}

In the second part of this thesis we move to a different story but one which is still very much in the same flavour as Part~\ref{P:c2}, on the emergence of combinatorics and using combinatorial techniques to study aspects of scattering amplitudes.
This time, instead of looking at an arithmetic invariant related to Feynman diagrams through enumerating certain edge bipartitions, we are looking at a geometric story of triangulations through a simple combinatorial map where one side of the map is related to on-shell diagrams.

We start with introducing the main objects of interest, the positive Grassmannian and some of the combinatorial objects in bijection with the cells that decompose it.
Finishing off this introductory chapter, in Section~\ref{SS:t-duality} we set-up the motivation for the question we investigate in Chapter~\ref{C:lebijection}, which looks at the map mentioned above through the lens of a different combinatorial object.
The exposition from this chapter follows mainly from~\cite{tpgrass, LPW}. \\

\section{Grassmannia}\label{S:grass}

The Grassmannian is one of those classical geometric objects in mathematics that has been extensively studied due to its nice structure and its many connections to different areas of mathematics, including (of course) combinatorics, algebraic and differential geometry, and representation theory.
For one, the Grassmannian is at the heart of \textbf{Schubert calculus}, which aims to solve problems in enumerative geometry via the Grassmannian, Schubert varieties, and flag varieties, and which lies at the intersection of combinatorics and algebraic geometry.
An example of one of the early problems studied via Schubert calculus is~\cite{schubert}: "How many lines in 3-space, in general, intersect four given lines?" (of which the answer is generically 2).
As an added benefit, the Grassmannian can be described very simply.

\tpointn{Definition}\label{def:grass}
\statement{
  For $0 \leq k \leq n$, the (real) \textbf{Grassmannian} $\:\Gr_{k,n}(\mathbb{R})$ is the space of all $k$-dimensional (linear) subspaces of $\:\mathbb{R}^n$.
}

We can think of the elements $V$ of the Grassmannian as represented by real $k \times n$ matrices $A$ of rank $k$ modulo invertible row operations, where $V$ is the subspace spanned by the row vectors of $A$.
That is,
\[ \Gr_{k,n}(\mathbb{R}) \cong \qmod{\GL(k, \mathbb{R})}{\Mat^*(k,n)}, \]
where $\Mat^*(k,n)$ is the set of all real $k \times n$ matrices of rank $k$.
The dimension of $\:\Gr_{k,n}$ is $k(n - k)$.
We often drop the reference to the field $\mathbb{R}$ and just refer to the real Grassmannian as $\Gr_{k,n}$.
As a geometric object, we can view the Grassmannian through a couple of different lenses (which we do not go through here): as a compact smooth manifold, as a projective variety, as a scheme, and as a CW-complex.

\tpointn{Example}
\statement{$\Gr_{1,n}$ consists of all lines through the origin in $\:\mathbb{R}^n$, or viewing through matrices all $1 \times n$ matrices modulo scaling, which is equivalent to projective space $\:\mathbb{P}^{n-1}$.
}

For algebraic combinatorialists, the interest in the Grassmannian lies in the various beautiful connections to combinatorics arising from its decompositions.
Through the classic Schubert decomposition of the Grassmannian into \textbf{Schubert cells}, which can be indexed by partitions, we are led to some familiar combinatorial machinery such as Young tableaux, Schur functions, and Schubert polynomials.
A standard reference for many of these topics is~\cite{Fulton}.

From a different decomposition, as first described by Gelfand, Goresky, MacPherson, and Serganova in~\cite{GGMS}, and which we will present here, we are led to matroids and matroid polytopes, opening up a different connection to combinatorics.

First, let's fix some notation.
Let $[n]$ denote the set $\{1,\ldots, n\}$ and $\binom{[n]}{k}$ be the set of all $k$-element subsets of $[n]$.
Given a $k \times n$ matrix $A$ of rank $k$ and an $I \in \binom{[n]}{k}$, let $p_I(A) \coloneqq \det(A_I)$, where $A_I$ is the $k \times k$ matrix of $A$ restricted to the columns indexed by $I$.
We call this a \textbf{maximal minor} of $A$.
For $V \in \Gr_{k,n}$ represented by $A$, the \textbf{Pl\"ucker coordinates} of $V$ are then defined to be $p_I(V) \coloneqq p_I(A)$ for all $I \in \binom{[n]}{k}$, where we note that $p_I(V)$ does not depend on the choice of representative matrix $A$.

Now, thinking about $V \in \Gr_{k,n}$ as given by the $n$ column vectors of $A$ which together span $\mathbb{R}^k$, since $A$ is of rank $k$, we can then look at all the $k$-subsets of these $n$ column vectors such that they form a basis for $\mathbb{R}^k$.
In particular, each of these subsets $I \in \binom{[n]}{k}$ must have the property that $p_I(A) \neq 0$.
Viewing through the combinatorial lens, in other words we can construct a representable matroid $\mathcal{M}$ over the ground set $[n]$ which has bases
\[ M(A) \coloneqq \left\{ I \in \binom{[n]}{k} \ :\ p_I(A) \neq 0 \right\}, \] \\
and where $A$ represents $\mathcal{M}$ over $\mathbb{R}$.
For a standard reference on matroids, we refer the reader to~\cite{matroid}.

Going back to the Grassmannian, we can then divide the elements of the Grassmannian based on which Pl\"ucker coordinates are non-zero, into \textbf{matroid strata}, also known as the Gelfand-Serganova strata.
That is, for any $M \subset \binom{[n]}{k}$ we can define
\[ \mathscr{S}_M \coloneqq \left\{V \in \Gr_{k,n} \ : \ p_I(V) \neq 0 \text{ if and only if } I \in M \right\}. \]
Equivalently, if $A$ is a $k \times n$ matrix representing $V$, then the condition reduces to $M = M(A)$ where $M(A)$ is as defined above.
Notice that $M$ defines a representable matroid of rank $k$, where here we are taking the convention of defining a matroid by its bases.
Thus we get the following (disjoint) decomposition, called the \textbf{matroid stratification}, also known as the Gelfand-Serganova stratification, of the Grassmannian, \\
\[ \Gr_{k,n} = \bigsqcup_{\text{representable matroids } M} \mathscr{S}_M. \]

However, unfortunately Mn\"ev's Universality Theorem~\cite{mnev} tells us that the structure of these $\mathscr{S}_M$ can be very complicated, in fact as complicated as any algebraic variety! \\

\bpoint{The positive Grassmannian}

Instead of looking at the full Grassmannian, Postnikov in~\cite{tpgrass} initiated the study of a certain subset of the Grassmannian called the positive Grassmannian, through giving a combinatorial description of its cells, called positroids, which turned out to have a much nicer geometric structure.
This opened the door to the extensive study of the positive Grassmannian, both combinatorially and through the multitude of emerging connections with other branches of mathematics and physics.
Here we just define the positive Grassmannian and in the next section, we will delve into some aspects of its decomposition.

As an analogue to totally positive and totally non-negative matrices, we can extend this notion of positivity to the matrices underlying the Grassmannian.

\tpointn{Definition} (Definition 3.1 of~\cite{tpgrass})\label{def:TNN}
\statement{
  We say that a $k \times n$ matrix $A$ of rank $k$ is \textbf{totally non-negative}, resp. \textbf{totally positive}, if all of its maximal minors are non-negative i.e. $p_I(A) \geq 0$, resp. all positive $p_I(A) > 0$, for all $I \in \binom{[n]}{k}$.
  We then say that $V \in \Gr_{k,n}(\mathbb{R})$ is \textbf{totally non-negative}, resp. \textbf{totally positive}, if $V$ is represented by a $k \times n$ matrix $A$ that is totally non-negative i.e. $p_I(V) \geq 0$, resp. totally positive $p_I(V) > 0$, for all $i \in \binom{[n]}{k}$. \\\\
  The \textbf{totally non-negative (TNN) Grassmannian}, also referred to as the \textbf{positive Grassmannian}, $\Gr_{k,n}^{\geq 0}$ is the set of all totally non-negative $V \in \Gr_{k,n}(\mathbb{R})$.
  The \textbf{totally positive Grassmannian} $\Gr_{k,n}^{> 0}$ is the set of all totally positive $V \in \Gr_{k,n}(\mathbb{R})$.
}

Note that the positive Grassmannian $\Gr_{k,n}^{\geq 0}$ is a closed subset of the (full) Grassmannian $\Gr_{k,n}$ of dimension $k(n-k)$, which is also the dimension of $\:\Gr_{k,n}$.
For the rest of this chapter, we will focus on the positive Grassmannian.

We end this section with an example of an element in the positive Grassmannian.
\tpointn{Example}\label{ex:pg-elt}
\statement{
  Consider the following $2 \times 4$ matrix of rank $2$, which represents an element $V$ in $\:\Gr_{2,4}$:
  \[ A = \left[\begin{array}{cccc} 1 & 0 & -1 & -2 \\ 0 & 1 & 2 & 3 \end{array}\right]. \]
  Then looking at the Pl\"ucker coordinates, for the following order on the $2$-subsets of $[4]$:
  \[ \binom{[4]}{2} = \{12, 13, 14, 23, 24, 34\}, \]
  we have the sequence
  \[ p_I(A) = 1, 2, 3, 1, 2, 1, \]
  where as a sample calculation
  \[ p_{34}(A) = \det\left[ \begin{array}{cc} -1 & -2 \\ 2 & 3 \end{array}\right] = 1. \]
  In particular, $p_I(A) \geq 0$ for all $I \in \binom{[4]}{2}$, and thus $V$ is an element in the positive Grassmannian $\Gr_{2,4}^{\geq 0}$. \\
  As we actually have $p_I(A) > 0$ for all $I$, $V$ is also in the totally positive Grassmannian $\Gr_{2,4}^{> 0}$.\\
  Here we have that $M(A) = \binom{[4]}{2}$, which is a uniform matroid!
} \\

\section{Positroids and related combinatorial objects}\label{S:positroids}

Like how the Grassmannian could be subdivided into its matroid strata, analogously we can do the same with the positive Grassmannian, this time based on which maximal minors are positive.

\tpointn{Definition} (Definition 3.2 of ~\cite{tpgrass})\label{def:positroid}
  \statement{
    Let $M$ be a subset of $\ \binom{[n]}{k}$, and define
    \[ S_M \coloneqq \mathscr{S}_{M} \cap \Gr_{k,n}^{\geq 0}  = \left\{V \in \Gr_{k,n}^{\geq 0} \ : \ p_I(V) > 0 \text{ for all } I \in M \text{ and } p_I(V) = 0 \text{ for all } I \in \binom{[n]}{k} \setminus M \right\}. \]
    If $S_M \neq \emptyset$, we call $M$ a \textbf{positroid}, and $S_M$ is its \textbf{positroid cell}.
  }

Note that $M$ is a rank $k$ representable matroid, represented by a $k \times n$ matrix with non-negative maximal minors, with the term "positroid" arising from abbreviating "positive matroid".
However, the positroid is not invariant under matroid isomorphism but is invariant under cyclic shifts of elements.
Additionally, notice that the totally positive Grassmannian corresponds to the uniform matroid $M = \binom{[n]}{k}$ and is the unique top-dimensional cell (of dimension $k(n-k)$) of the positive Grassmannian.
From Example~\ref{ex:pg-elt}, we have that the $M(A) = \binom{[4]}{2}$ is actually a positroid!

In other words, we are partitioning the positive Grassmannian based on the Pl\"ucker coordinates that are strictly positive.
We call this (disjoint) decomposition the \textbf{positroid stratification} of the positive Grassmannian
\[ \Gr^{\geq 0}_{k,n} = \bigsqcup_{M \text{ is a positroid }} S_M. \]

What's nice about the geometric structure of these positroid cells is that they are indeed cells.

\tpointn{Theorem} (Theorem 3.5 of~\cite{tpgrass}, Theorem 5.4 of~\cite{PSW})
\statement{
  Each positroid cell $S_M$ is homeomorphic to an open ball.
  The decomposition of $\:\Gr_{k,n}^{\geq 0}$ into the union of positroid cells is a $CW$-complex.
}

For us the interest in positroid cells lies in its surprisingly rich combinatorial structure, arising from the numerous families of combinatorial objects that Postnikov in ~\cite{tpgrass} showed index these cells.
Other than positroids, these objects include decorated permutations, Le diagrams, Grassmann necklaces, and equivalence classes of reduced plabic graphs.
We present these objects and one of the bijections between them, giving us a way to label positroid cells by their associated combinatorial objects. \\

\bpoint{Decorated permutations}

We start with the simplest family of combinatorial objects called decorated permutations.

\tpointn{Definition}\label{def:decperm}
\statement{
  A \textbf{decorated permutation} of $[n]$ is a permutation $\pi$ on $[n]$ where every fixed-point is designated (coloured) as a \textbf{loop} (black), denoted by $\pi(i) = \underline{i}$, or a \textbf{co-loop} (white), denoted by $\pi(i) = \overline{i}$. \\
  An \textbf{anti-excedance} of $\pi$ is an element $i$ such that either $\pi^{-1}(i) > i$ or $i$ is a co-loop.
  We refer to $\pi^{-1}(i)$ as the \textbf{position} of $i$.
}

\tpointn{Example} (from~\cite{LPW})\label{ex:perm}
\statement{
  The following is a decorated permutation of $[8]$ in two-line notation:
  \[ \pi = \left(\begin{array}{cccccccc} 1 & 2 & 3 & 4 & 5 & 6 & 7 & 8 \\ 3 & \underline{2} & 5 & 1 & 6 & 8 & \overline{7} & 4 \end{array}\right), \]
  which has 3 anti-excedances, $1, 4, 7$, as these elements when viewed on the bottom line, do not exceed their corresponding element in the top line and are not loops.
  Here $\pi$ has one loop, $2$, and one co-loop, $7$.
}

To show that a specific type of decorated permutation is in bijection with the cells of the positive Grassmannian, Postnikov utilized another object called the \textbf{Grassmann necklace} $\mathcal{I}$ which is directly in bijection with positroids and can be read off from the bases of the positroid.
Decorated permutations can then be seen as expressing the information from Grassmann necklaces in a more compact way.
We refer the reader to \textsection 16 of~\cite{tpgrass} for more details.

\tpointn{Theorem} (Lemma 16.2 of~\cite{tpgrass})
\statement{
  There is a bijection between decorated permutations $\pi$ of $[n]$ with $k$ anti-excedances and Grassmann necklaces $\mathcal{I}$ of type $(k,n)$.
  Therefore, decorated permutations of $[n]$ with $k$ anti-excedances index the cells of $\:\Gr_{k,n}^{\geq 0}$, and we denote by $S_{\pi}$ the positroid cell indexed by $\pi$.
}

Going back to Example~\ref{ex:perm}, this means that $\pi$ indexes the cell $S_{\pi}$ of $\:\Gr_{3,8}^{\geq 0}$. \\

\bpoint{Le diagrams and pipe dreams}\label{SS:le-def}

While Grassmann necklaces can directly see the bases of the positroid, and decorated permutations are simple and succinct objects that encode positroids, one property that both do not easily see is the dimension of the associated positroid cell.
For this, we need the next family of combinatorial objects called Le diagrams.

\tpointn{Definition} (Definition 6.1 of~\cite{tpgrass})\label{def:le}
\statement{
  A \textbf{\BLe-diagram} (Le diagram), is a filling $D$ of a Young diagram of shape $\lambda$ with $0$'s and $\leplus$'s such that $D$ avoids the \textbf{\BLe-configuration}:
  \[ \leconfig \] \\
  That is, no $0$ has both a $\leplus$ above and to the left of it, which we refer to as the \textbf{\BLe-condition}.\\
  For $0 \leq k \leq n$ we say that the \Le-diagram $D$ is of \textbf{type (k,n)} if the shape $\lambda$ fits inside a $k \times (n-k)$ rectangle.
}

An example of a \Le-diagram $D$ is given in Figure~\ref{fig:ex-le}, which we will see shortly is the \Le-diagram associated to the decorated permutation in Example~\ref{ex:perm}.
$D$ is indeed a \Le-diagram as it avoids the \Le-configuration, and it is of type $(3,8)$ since the shape of $D$ fits inside a $3 \times 5 = 3 \times (8 - 3)$ rectangle.

In \textsection 20 of~\cite{tpgrass}, Postnikov gave two bijections between \Le-diagrams and decorated permutations, the first through associating a series of other objects (hook diagrams, networks, and plabic graphs) to \Le-diagrams.
The second, which we describe here, uses an algorithm from \textsection 19 of~\cite{tpgrass} going through a slightly different object called pipe dreams.

Given a \Le-diagram $D$ of type $(k,n)$, we associate a decorated permutation $\pi_D$ on $[n]$ as follows:
\begin{enumerate}[itemsep=3pt]
  \item In $D$, we replace each $0$ with a cross $\,\smallcross\,$ and each $\leplus$ with an elbow joint $\,\smallelbow\;$:
    \begin{align*}
      \lebox{0}\quad&\longmapsto\quad \lebox{\cross{}} \\[4pt]
      \lebox{+} \quad&\longmapsto\quad \lebox{\elbow{}}
    \end{align*}
  \item Notice that we can view the south-east (SE) border of $D$ as a lattice path with $n$ steps, and thus label the edges with $1, \ldots, n$ from the top-right corner to the bottom-left corner.
  \item Then label the edges of the north-west (NW) border of $D$ so that, viewing $D$ as a grid, the rows and columns have the same labels (the opposite horizontal/vertical edges are labelled the same).
    We call this diagram $P$ the \textbf{pipe dream} associated to $D$.
  \item To get the decorated permutation $\pi$ associated to $P$ and $D$, we follow the "pipes" of $P$ from the SE border to the NW border. That is, $\pi(i) = j$ if the pipe starts at $i$ and ends at $j$.
    A horizontal pipe starting and ending at $i$, where $i$ labels vertical edges, is denoted a co-loop $\pi(i) = \overline{i}$.
    A vertical pipe starting and ending at $i$, where $i$ labels horizontal edges, is denoted a loop $\pi(i) = \underline{i}$.
\end{enumerate}

To illustrate this algorithm, we use an example.
\tpointn{Example}\label{ex:le}
\statement{
  Given the \Le-diagram $D$ of type $(3,8)$ in Figure~\ref{fig:ex-le}, we label the SE and NW borders with $1, \ldots, 8$ as above, and replace the $0$'s and $\leplus$'s with crosses and elbow-joints.
  This gives the pipe dream $P$ in Figure~\ref{fig:ex-pipe}, where arrows are added to show the direction to follow the pipes for the permutation. \\
\begin{figure}[ht]
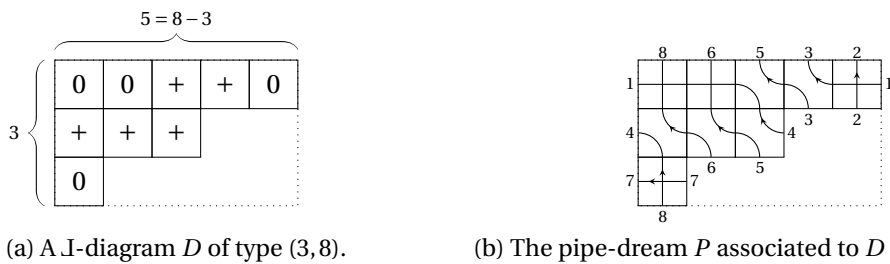

  \centering
  \begin{subfigure}[b]{.48\linewidth}
    \centering
    \[ \exle \]
    \caption{A \Le-diagram $D$ of type $(3,8)$.}
    \label{fig:ex-le}
  \end{subfigure}
  \begin{subfigure}[b]{.48\linewidth}
    \centering
    \[ \expipe \]
    \caption{The pipe-dream $P$ associated to $D$ in Figure~\ref{fig:ex-le}.}
    \label{fig:ex-pipe}
  \end{subfigure}
  \caption[The \Le-diagram and pipe dream associated to Example~\ref{ex:perm}.]{The \Le-diagram and pipe dream associated to Example~\ref{ex:perm}. The arrows on the pipes indicate the direction to follow to read off the decorated permutation.}
  \label{fig:ex-le-pipe}
\end{figure} \\
  To read off the decorated permutation, we follow the pipes in the direction of the arrows to get
  \[ \pi_D = \left(\begin{array}{cccccccc} 1 & 2 & 3 & 4 & 5 & 6 & 7 & 8 \\ 3 & \underline{2} & 5 & 1 & 6 & 8 & \overline{7} & 4 \end{array}\right), \]
  which is the decorated permutation from Example~\ref{ex:perm}!
  Note that $\pi_D$ is a decorated permutation on $[8]$ with 3 anti-excedances.
}

From the example, we make a few observations about this map $D \mapsto \pi_D$, where $P$ is the associated pipe dream to $D$:
\begin{itemize}[itemsep=3pt]
  \item Rows of all $0$'s in $D$ correspond to horizontal pipes in $P$, and thus these are the co-loops in $\pi_D$.
    In Example~\ref{ex:le} this would be the row labelled $7$ in Figure~\ref{fig:ex-le-pipe}.
  \item Similarly, columns of all $0$'s in $D$ correspond to vertical pipes in $P$, and thus these are the loops in $\pi_D$.
    In Example~\ref{ex:le} this would be the column labelled $2$ in Figure~\ref{fig:ex-le-pipe}.
  \item Because of the way the elbow joints are drawn, notice that pipes in $P$ can only ever go to the left and up from its starting point.
    Thus, pipes starting from $i$ that end on vertical edges $j$ must have $i \geq j$, with $i = j$ being a horizontal pipe.
    Whereas pipes starting from $i$ and ending on horizontal edges $j$ have $i < j$.
    That is, the labels of the rows (vertical edges) of $P$ exactly correspond to the anti-excedances of $\pi_D$.
    In Example~\ref{ex:le}, this would be the rows labelled $1,4,7$ in Figure~\ref{fig:ex-le-pipe}, which indeed are the anti-excedances of $\pi_D$.
  \item Finally, avoiding the \Le-configuration in $D$ corresponds to pipes in $P$ crossing at most once, and furthermore once two pipes cross, they will also no longer touch (via an elbow joint).
\end{itemize}

In turns out that this map is indeed a bijection, and furthermore one of the reasons to look at \Le-diagrams is because of how easily they display the dimension of the positroid cell that they index.
In particular, the dimension of the corresponding positroid cell is determined by counting the number of $\leplus$'s in the \Le-diagram.

\tpointn{Theorem} (Theorem 6.5, Corollary 20.1, Theorem 20.3 of~\cite{tpgrass})
\statement{
  The map $D \mapsto \pi_D$ is a bijection from \Le-diagrams of type $(k,n)$ to decorated permutations of $[n]$ with $k$ anti-excedances.
  Therefore, \Le-diagrams of type $(k,n)$ index the cells of $\:\Gr_{k,n}^{\geq 0}$.
  Furthermore, let $S_D$ be the positroid cell indexed by the \Le-diagram $D$.
  Then the dimension of the cell $S_D$ is the number of $\:\leplus$'s in $D$.
}

Going back to Example~\ref{ex:le}, this means that $D$ indexes the cell $S_{D}$ of $\:\Gr_{3,8}^{\geq 0}$, which corresponds to the cell $S_{\pi}$ from Example~\ref{ex:perm}, and the dimension of $S_D$ is $5$ since there are 5 $\leplus$'s in $D$ . \\

\bpoint{Plabic graphs}\label{SS:plabic}

The final family of combinatorial objects that we will define (without going into detail) are plabic graphs, where the term "plabic" comes from abbreviating "planar bicoloured".
These graphs play an important role in the foundations for the original correspondences of the above objects to positroid cells as well as the basis of many of the applications and connections from positroids to other fields.

\tpointn{Definition} (Definition 11.5 of~\cite{tpgrass})
\statement{
  A \textbf{plabic graph} is an undirected planar graph $G$ drawn inside a disk with the following properties:
  \begin{itemize}
    \item $n$ boundary vertices on the boundary of the disk labelled clockwise with $1, \ldots n$,
    \item some number of internal vertices, which are coloured either black of white, and
    \item each boundary vertex is incident to a single edge; if this edge is incident to a leaf, the leaf is called a \textbf{lollipop}.
  \end{itemize}
}

From associating \Le-diagrams to plabic graphs, Postinikov showed that what we really want to be looking at are equivalence classes of \textbf{reduced plabic graphs}, which are specific types of plabic graphs that can be obtained via certain local moves (transformations) on plabic graphs.
That is, these equivalence classes also index the cells of the positive Grassmannian!
In particular, like \Le-diagrams, reduced plabic graphs also see the dimension of the cell, via the number of faces of the graph minus one.
We refer the reader to \textsection 11--13, 20 of~\cite{tpgrass} and \textsection 14 of~\cite{LPW} for a more thorough treatment of plabic graphs.

Continuing with our running example in bijection with the same positroid, we present the reduced plabic graph in bijection with the \Le-diagram $D$ from Example~\ref{ex:le}.

\tpointn{Example}\label{ex:plabic}
\statement{
  In Figure~\ref{fig:ex-plabic}, we draw the reduced plabic graph $G$ associated to Examples~\ref{ex:perm},~\ref{ex:le} in two different ways.
  The top representation hints at the bijection with \Le-diagrams, and the bottom representation hints at the connection to on-shell diagrams.
  \begin{figure}[ht]
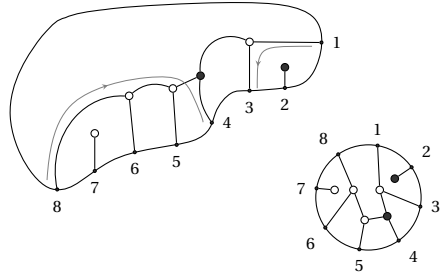

    \centering
    \[ \explabic \]
    \caption{A reduced plabic graph $G$, drawn in two ways.}
    \label{fig:ex-plabic}
  \end{figure} \\
  Very briefly to see that the decorated permutation does indeed line up with the previous examples, follow the labelled boundary vertices turning (maximally) right at every black vertex and left at every white vertex.
  Some of these paths are drawn on $G$ with the arrows indicating the direction to follow.
  Lollipops correspond to fixed points and we colour the fixed points accordingly, that is, black lollipops are loops, and white lollipops are co-loops. \\\\
  Here $G$ is of type $(3,8)$ as $G$ has $8$ boundary vertices, and $k=3$ is computed via
  \[ k = \#\text{edges} - \#\text{black vertices} - \sum_{\text{white vertices }v} (\deg(v) - 1) = 11 - 2 - 6 = 3. \]
  Thus the equivalence class of $G$ indexes the cell $S_{G}$ of $\:\Gr_{3,8}^{\geq 0}$, which corresponds to $S_{\pi}$ and $S_{D}$ from Examples~\ref{ex:perm},~\ref{ex:le}.
  The dimension of the cell $S_G$ indexed by $G$ is then $\#\text{faces} - 1 = 5$, which lines up with the dimension of $S_{D}$.
}

From our three examples, we have that the decorated permutation $\pi$ from Example~\ref{ex:perm}, the \Le-diagram $D$ from Example~\ref{ex:le}, and the reduced plabic graph $G$ from Example~\ref{ex:plabic} (and its equivalence class) are all in bijection with each other and with the same positroid.
Thus they index the same positroid cell, which we label by $S_{\pi}$, $S_{D}$, or $S_{G}$ according to which combinatorial object is of interest! \\

\section{T-duality on decorated permuations}\label{S:t-duality}

Through the recent extensive study of the combinatorial objects related to positroids, the positive Grassmannian (like its full counterpart) has been shown to have deep connections to many different areas of mathematics and physics.

To name a few, first purely combinatorially, we can study positroids themselves as a special class of matroids.
We can also study the numerous related combinatorial objects that arise directly from positroids and those arising from the many applications. Some examples include oriented matroids, polytopes and polyhedral subdivisions, non-crossing partitions, and lattice paths.
On the more algebraic and geometric side, there are connections to cluster algebras and quantum algebras, and naturally we could study the positive Grassmannian in the spirit of Schubert calculus, through varieties, flags, and symmetric functions.
Finally, on the physics side there have been applications to KP solitons, types of asymmetric exclusion processes, and scattering amplitudes (via the amplituhedron and Wilson loop diagrams).
We refer the reader to ~\cite{tpgrass, post, kps, asep, LPW, PSBW} and the references therein for further details.

We will highlight and focus on the combinatorics related to one specific application to physics, namely to scattering amplitudes in $\mathcal{N}=4$ SYM theory.
In particular, the story from ~\cite{LPW} starts with looking at the images of the positive Grassmannian under two different maps. \\

\bpoint{The hypersimplex}

The first map of interest is the classic moment map, which was used in~\cite{GGMS} to study the Grassmannian and its decomposition into matroids.
Let $A$ be a $k \times n$ matrix representing an element of $\Gr_{k,n}$.
The \textbf{moment map} $\mu$ from $\Gr_{k,n}$ to $\mathbb{R}^n$ is defined as
\[ \mu(A) \coloneqq  \frac{\sum_{I \in \binom{[n]}{k}} \abs{p_I(A)}^2 e_I}{\sum_{I \in \binom{[n]}{k}} \abs{p_I(A)}^2}. \]

Under this moment map, the object that arises as the image of both the Grassmannian and the positive Grassmannian is the well-known hypersimplex.

\tpointn{Definition}
\statement{
For $I \in \binom{[n]}{k}$, let $e_I \coloneqq \sum_{i \in I} e_i \in \mathbb{R}^n$ where the $e_i$ are the standard basis vectors of $\ \mathbb{R}^n$.
  The \textbf{(k,n)-hypersimplex} $\:\Delta_{k,n}$ is the convex hull of all such points $e_I$,
  \[ \Delta_{k,n} \coloneqq \conv\left\{e_I \ :\ I \in \binom{[n]}{k}\right\}. \]
}

Note that $\Delta_{k,n}$ lives inside $\mathbb{R}^n$ and is of dimension $n-1$.

As the hypersimplex is a convex polytope, a natural question that arises is how can one triangulate or subdivide this polytope using the moment map and in relation with the positroid cells that decompose the positive Grassmannian?

The start to an answer to this question lies in objects called \textbf{positroid polytopes}, which are matroid polytopes where the matroid is a positroid.
Analogous to how hypersimplex is defined, given a matroid we can naturally define a polytope using its bases.
That is, for a positroid $M$ on $[n]$ as defined by its bases, its positroid polytope $\Gamma_M$ is
\[ \Gamma_M \coloneqq \conv\left\{ e_I \ :\ I \in M \right\} \subset \mathbb{R}^n. \]
In particular, why we care about these polytopes is because it turns out that positroid polytopes and the closure of the image of positroid cells under the moment map coincide!

\tpointn{Theorem} (Proposition 7.10 of~\cite{bruhat}) \\
\statement{
  Let $M$ be a positroid and $S_M$ its associated positroid cell.
  Then $\Gamma_M = \overline{\mu(S_M)}$.
}

Now the question becomes, when does a collection of positroid polytopes (with positroid cells) triangulate the hypersimplex?
Here we want the positroid polytopes to be \textbf{generalized triangles} of $\Delta_{k,n}$, which restricts to positroid cells of dimension $n-1$ and where $\mu$ is injective on the cells.
This is one of the questions that~\cite{LPW} looked at. \\

\bpoint{The amplituhedron}

On the other side of the story, as part of an ongoing program on reformulating the Feynman diagram approach in quantum field theory which culminated in the book~\cite{grassampl}, a surprising connection arose between the work of physicists in the context of scattering amplitudes in $\mathcal{N}=4$ SYM theory and work on the positive Grassmannian.
What the physicists described as \textbf{on-shell diagrams}, which are analogs to Feynman diagrams, turned out to be exactly the \textbf{plabic graphs} as described by Postnikov~\cite{tpgrass} (as defined in Section~\ref{SS:plabic}) and in relation to positroid cells.

Motivated by this connection as well as observations on the BCFW recursion relations in tree-level amplitudes, Arkani-Hamed and Trnka~\cite{ampl} introduced a new geometric object called the amplituhedron in hopes of encoding these tree-level $\mathcal{N}=4$ SYM amplitudes through its "volume".

\tpointn{Definition} (\textsection 9 of~\cite{ampl})
\statement{
  For $k+m \leq n$, the (tree) \textbf{amplituhedron} $\mathcal{A}_{n,k,m}(Z)$ is defined as the image of the positive Grassmannian under the map
  \begin{align*}
    \tilde{Z} \ :\ \Gr_{k,n}^{\geq 0} \quad&\longrightarrow\quad \Gr_{k,k+m} \\
    A \quad&\longmapsto\quad AZ
  \end{align*}
  where $A$ is a $k \times n$ matrix representing an element in $\:\Gr_{k,n}^{\geq 0}$, and $Z$ is a totally positive $n \times (k + m)$ matrix.
}

Note that $\mathcal{A}_{n,k,m}$ lives inside the Grassmannian $\Gr_{k,k+m}$ and has full dimension $km$.
We will focus on the $m=2$ case which is considered a toy-model for the more physically relevant $m=4$ case.
However, we note that there are other special cases which also have interesting connections to combinatorics.
As an example, when $Z$ is a square $n\times n$ matrix, we recover the positive Grassmannian $\Gr_{k,n}^{\geq 0}$.
When $k=1$ and $m=2$, $\mathcal{A}_{n,1,2}$ is a polygon in $\mathbb{P}^2$, and more generally $\mathcal{A}_{n,1,m}$ is a cyclic polytope living inside $\mathbb{P}^{m}$.

While not a polytope, the same question from the hypersimplex can be asked about the amplituhedron: based on the positroid stratification of the positive Grassmannian, how can one "triangulate" the amplituhedron?
Relating back to scattering amplitudes, these triangulations should then be related to the different decompositions of the amplitude in the on-shell representation.
Some aspects of these questions were previously looked at in~\cite{decomps}.

For the amplituhedron, the objects of interest are $\textbf{Grasstopes}$, which are defined to be the closure of the image of positroid cells under the $\tilde{Z}$ map.
The term as used in~\cite{PSBW} is the abbreviation of "Grassmann polytopes",
Then we can analogously define \textbf{generalized triangles} of $\mathcal{A}_{n,k,m}$ as Grasstopes restricted to $km$-dimensional positroid cells and where $\tilde{Z}$ is injective on the cells.
Once again, one could ask: when does a collection of Grasstopes triangulate the amplituhedron?

\bpoint{T-duality}\label{SS:t-duality}

For us, one of the interesting results from~\cite{LPW} was based on a different question, one connecting the two stories about triangulations.
As the positive Grassmannian has this nice positroid cell decomposition, another natural question to ask is if there is any relationship between the images of the positroid cells under the previous two maps, the moment map and $\tilde{Z}$, and how they would relate back to the full objects.
Following these lines of questioning and based on numerical evidence, Lukowski, Parisi, and Williams~\cite{LPW} found a surprising connection between the hypersimplex and the $m=2$ amplituhedron, which was encapsulated by a very simple map directly on positroid cells that they called T-duality.

The name "T-duality" comes from its correspondence to the T-duality in string theory, which very simply put is a type of equivalence of theories where as one example momentum in a certain string theory is interchanged with winding numbers in its T-dual theory.
In particular, this map is named as a nod to being a manifestation of T-duality in the context of the geometry of $\mathcal{N}=4$ SYM scattering amplitudes and the $m=2$ amplituhedron.

\tpointn{Definition} (Definition 5.1 of~\cite{LPW})\label{def:t-duality}
\statement{
  The \textbf{T-duality} map from loopless decorated permutations on $[n]$ to co-loopless decorated permutations on $[n]$ is defined as
  \begin{align*}
    \pi \quad&\longmapsto\quad \hat\pi \\
    (a_1, a_2, \ldots, a_n) \quad&\longmapsto\quad (a_n, a_1, \ldots, a_{n-1})
  \end{align*}
  where the permutations are written in one-line notation, and any fixed points in $\hat\pi$ are declared to be loops.
  That is, for a given loopless $\pi$, we have $\hat\pi(i) = \pi(i-1)$ where all fixed points are loops, and we call $\hat\pi$ the \textbf{T-dual} decorated permutation.
}

\tpointn{Lemma} (Lemma 5.2, Proposition 5.17 of~\cite{LPW})
\statement{
  The T-duality map $\pi \mapsto \hat\pi$ is a bijection between loopless decorated permutations on $[n]$ with $k+1$ anti-excedances and co-loopless decorated permutations on $[n]$ with $k$ anti-excedances. \\\\
  Equivalently, the T-duality map is a bjiection between loopless positroid cells $S_{\pi}$ of $\:\Gr_{k+1,n}^{\geq 0}$ and co-loopless positroid cells $S_{\hat\pi}$ of $\:\Gr_{k,n}^{\geq 0}$.
  Furthermore, we have the following dimensional relationship
  \[ \dim(S_{\hat\pi}) - 2k = \dim(S_{\pi}) - (n-1), \]
  where in particular, if $\dim(S_{\pi}) = n-1$ then $\dim(S_{\hat\pi}) = 2k$.
}

While this dimensional relationship was proven to exist between T-dual loopless cells of $\Gr_{k+1,n}^{\geq 0}$ and co-loopless cells of $\Gr_{k,n}^{\geq 0}$, it was not clear where it was coming from.
The need to understand this relationship motivated the question of viewing T-duality as a map on \Le-diagrams, in which one very easily sees the dimension of the positroid cells they index.
This problem is the focus of Chapter~\ref{C:lebijection}.

The magic in this very simple shift map on decorated permuations is that it was exactly what was needed to go between triangulations of the ($n-1$)-dimensional hypersimplex $\Delta_{k+1,n}$ and of the $2k$-dimensional $m=2$ amplituhedron $\mathcal{A}_{n,k,2}$!
This duality then explained the mysterious appearance of $\binom{n-2}{k}$ in the number of cells needed to triangulate both these objects.
In~\cite{LPW}, the connection was first shown for a specific type of triangulation.

\tpointn{Proposition} (Theorem 6.5, Proposition 6.6 of~\cite{LPW})
\statement{
  The T-duality map provides a bijection between BCFW triangulations (resp. dissection) of the hypersimplex $\:\Delta_{k+1,n}$ and BCFW triangulations (resp. dissection) of the amplituhedron $\mathcal{A}_{n,k,2}$. \\\\
  More specifically, a collection $\{S_{\pi}\}$ of cells of $\:\Gr_{k+1,n}^{\geq 0}$, as constructed in~\cite{LPW}, gives a triangulation (resp. dissection) of $\:\Delta_{k+1,n}$ if and only if the collection $\{S_{\hat\pi}\}$ of cells of $\:\Gr_{k,n}^{\geq 0}$ gives a triangulation (resp. dissection) of $\mathcal{A}_{n,k,2}$. \\\\
  Furthermore, if $S_{\pi}$ is a generalized triangle of $\:\Delta_{k+1,n}$, then $S_{\hat{\pi}}$ is a generalized triangle of $\mathcal{A}_{n,k,2}$.
}

Very recently, Parisi, Sherman-Bennett, and Williams~\cite{PSBW} extended the work done in~\cite{LPW}, building on the theory of generalized triangles, and were able to prove the main conjecture in~\cite{LPW} which strengthens the above statement to all triangulations and dissections.
Their proof utilized looking at the T-duality map via plabic graphs and plabic tilings.

\tpointn{Theorem} (Conjecture 6.9 of~\cite{LPW}, Theorem 11.6 of~\cite{PSBW})\label{triang-conj}
\statement{
  Let $\{S_{\pi}\}$ be a collection of cells of $\:\Gr_{k+1,n}^{\geq 0}$.
  Then this collection gives a triangulation (resp. dissection) of $\:\Delta_{k+1,n}$ if and only if the collection $\{S_{\hat\pi}\}$ of cells of $\:\Gr_{k,n}^{\geq 0}$ gives a triangulation (resp. dissection) of $\mathcal{A}_{n,k,2}$.
}\\

\chapter{Reconstructing a bijection on the level of Le diagrams}\label{C:lebijection}

A natural question that arises from T-duality (Definition~\ref{def:t-duality}) is: what does this map look like diagrammatically, that is, on the level of \Le-diagrams? How are the \Le-diagrams related?

Recall from Section~\ref{SS:le-def} that given a decorated permutation $\pi$ on $[n]$ with $k$ anti-excedances, we get an associated \reflectbox{L}-diagram $D$ of type $(k,n)$ where:
\begin{itemize}
  \item $D$ has shape $\lambda$ which fits in a $k \times (n-k)$ rectangle,
  \item co-loops are rows of all $0$'s, while loops are columns of all $0$'s,
  \item $\dim(S_D)$ is the number of $\leplus$'s in $D$,  where $S_D$ is the positroid cell indexed by $D$, and
  \item anti-excedances of $\pi$ label the rows of the pipe dream associated to $D$ (including empty rows).
\end{itemize}

The direction of T-duality we will look at is $\hat\pi \mapsto \pi$, where $\hat\pi$ is a co-loopless permutation on $[n]$ with $k$ anti-excedances and $\pi$ is a loopless permutation on $[n]$ with $k+1$ anti-excedances.
On \Le-diagrams, we are looking at the map from $\hat{D} \mapsto D$, where $\hat{D}$ is the \Le-diagram associated to $\hat\pi$ and $D$ is the \Le-diagram associated to $\pi$.
As cells, recall from Section~\ref{SS:t-duality} that $S_{\hat\pi}$ is a co-loopless cell of $\Gr^{\geq 0}_{k,n}$ and $S_{\pi}$ is a loopless cell of $\Gr^{\geq 0}_{k+1,n}$ where their dimensions relate via $\dim(S_{\hat\pi}) - 2k = \dim(S_{\pi}) - (n-1)$.
One advantage of looking at this map on \Le-diagrams, as we will see, is how simple it becomes to relate the dimensions through counting the $\leplus$'s on either side.

We start with an explicit algorithm for constructing $D$ given a $\hat{D}$, with an example given at the end in Section~\ref{ex:algorithm}, and then move to a more visual perspective that showcases the structure of the $\leplus$'s of $D$.
We end with proving that the algorithm does give the correct \Le-diagram, as given by $\pi$.\\

\section{The bijection on the level of Le diagrams}\label{S:algorithm}

\bpoint{Notation}

Throughout this chapter, we will refer to the rows and columns of a \Le-diagram $D$ by the same labelling as that which gives the associated decorated permutation (i.e from its pipe dream).
Recall, we are viewing the SE border of the shape $\lambda$ as a lattice path of $n$ steps from the top-right corner to the bottom-left corner of the $k \times (n-k)$ rectangle.
The edges of this path are then labelled from $1$ to $n$.
Boxes in $D$, and the $k \times (n-k)$ rectangle, will then be referred to by their coordinates $(i,j)$ under this labelling.
A box is considered as "existing" if it is a valid box to be filled within the shape $\lambda$.
Note that all such valid boxes have $i < j$.
An example of this notation is given in Example~\ref{ex:notation}.
Finally since we are dealing with two \Le-diagrams, one on each side of the map, we refer to the "corresponding" row/column as the row/column with the same label on the opposite side of the map.

We will also order boxes in columns from top to bottom and in rows from right to left, in accordance with labels going from smallest to largest.
For rows/columns, having $i < j$ would mean $i$ is to the left and/or above $j$.
For boxes, "first" refers to right/top-most in a row/column, "last" refers to left/bottom-most in a row/column and "next" refers to the next box to the left/below in a row/column.

\tpointn{Example}\label{ex:notation}
\statement{
  Going back to Example~\ref{ex:le}, we have the following labelling of the \Le-diagram $D$:\\
  \[ \exnotation \]\\
  where there are $\leplus$'s in boxes $(1,3)$, $(1,5)$, $(4,5)$, $(4,6)$, and $(4,8)$.
  Notice box $(4,3)$ does not exist since it is outside of the shape $\lambda$ of $D$.
  By the nature of the labelling, we also never have boxes of the form $(i,i)$.
  The arrows indicate the direction of the ordering of boxes, going from first to next to last.
}\\

\bpoint{Preliminaries and shape of \texorpdfstring{$D$}{D}}\label{SS:shape}

We start with a co-loopless decorated permutation $\hat\pi$ on $[n]$ with $k$ anti-excedances, which we denote in two-line notation as
\[\hat\pi = \left(\begin{array}{cccc} 1 & 2 & \cdots & n \\ a_n & a_1 & \cdots & a_{n-1}\end{array}\right). \]\\
Let $\{b_1, \ldots, b_k\}$ be the $k$ anti-excedances of $\hat\pi$ ordered such that
$b_1 < \cdots < b_k$ (recall that $b_u$ is an anti-excedance if $\hat\pi^{-1}(b_u) > b_u$, where note we can't have $\hat\pi(b_u) = \bar{b}_u$ since $\hat\pi$ is co-loopless).
In particular, let $i_u = \hat\pi^{-1}(b_u)$ be the position of $b_u$.
Then we have $b_u = \hat\pi(i_u) = a_{i_u -1}$ for $1 \leq u \leq k$.

\begin{figure}[ht]
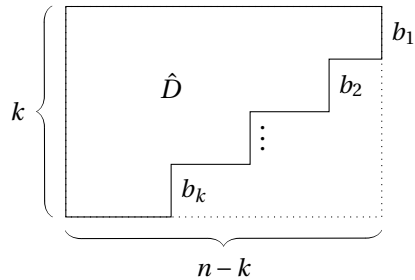

  \centering
  \[ \hatD \]
  \caption[The shape and row labels of \Le-diagram $\hat{D}$ associated to $\hat\pi$.]{The shape and row labels of \Le-diagram $\hat{D}$ associated to $\hat\pi$. Here we have illustrated when $b_1 = 1$. When $b_1 \neq 1$, there are columns of size 0 labelled $1, \ldots, b_1 - 1$ (these would be loops), indicated by a horizontal line extending the top (north) border of the diagram. Note that there are no rows of size 0 as $\hat\pi$ is co-loopless.}
  \label{fig:hat-D}
\end{figure}

We denote the associated $\Le$-diagram by $\hat{D}$ with shape $\hat\lambda$.
Since $\hat\pi$ is co-loopless, we have that every row in $\hat{D}$ has $\geq 1 \;\leplus$'s.
In particular,
\[ \hat\lambda = \left(\hat\lambda_1, \cdots, \hat\lambda_k\right), \quad\text{where }\hat\lambda_u = (n-k) - (b_{u} - u) \text{ for } 1 \leq u \leq k \]
and every $\hat\lambda_u \geq 1$.
We also have that the rows of $\hat{D}$ are thus labelled by these $b_{u}$, as in Figure~\ref{fig:hat-D}.

We want to give an algorithm that produces a $\Le$-diagram $D$ of type $(k+1, n)$ and dimension\\ ${\dim(S_{\hat{D}}) -2k + (n-1)}$ with associated decorated permutation $\pi$ of $[n]$ given by
\[\pi = \left(\begin{array}{cccc} 1 & 2 & \cdots & n \\ a_1 & a_2 & \cdots & a_{n}\end{array}\right), \]
which has $k+1$ anti-excedances and is loopless.
The loopless condition means that $D$ has to have $\geq 1 \;\leplus$'s in each of the $n-(k+1)$ columns.

\begin{figure}[th]
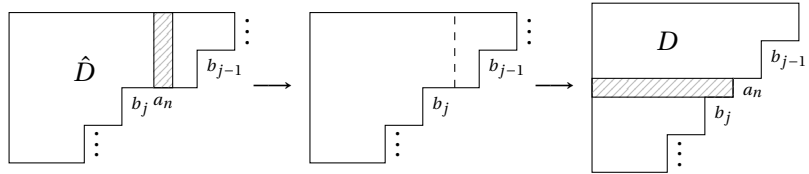

  \centering
  \[ \constructD \]
  \caption[Steps to construct the shape of $D$ from $\hat{D}$.]{Steps to construct the shape of $D$ from $\hat{D}$. Remove column $a_n$ from $\hat{D}$ and insert a row labelled $a_n$ where the dashed line is, making sure that the new boundary path is in the correct order. Here $j$ is the index such that $b_{j-1} < a_n < b_j$.}
  \label{fig:construct-D}
\end{figure}

First, by looking at T-duality as a map on permutations, we can already determine the shape $\lambda$ of $D$ which only cares about the anti-excedances of $\pi$.
Consider the following:
\begin{itemize}[itemsep=4pt]
  \item $a_n$ is not an anti-excedance of $\hat\pi$ since $\hat\pi$ is co-loopless, $\hat\pi(1) = a_n$ and $1 \not> a_n$.
    In particular, $a_n$ is the label of a column in $\hat{D}$.
  \item Under T-duality, $b_u$ for $1 \leq u \leq k$ stays an anti-excedance as $i_u > b_u = \hat\pi(i_u) = \pi(i_u-1)$.
    Since $\pi$ is loopless, we have either $\pi(i_u - 1) < i_u-1$ or $\pi(i_u-1) = i_u-1$ where $i_u-1$ must be a co-loop.
    In either case, $b_u$ is an anti-excedance of $\pi$.
  \item $a_n$ is always an anti-excedance of $\pi$ since $\pi$ is loopless.
  \item There are no other anti-excedances of $\pi$ since $\pi(i) = \hat\pi(i+1) \geq i+1 > i$ for all $1 \leq i < n$ with $i+1 \not\in \{i_1, \ldots, i_k\}$.
\end{itemize}
Based on these observations, we have that $\{b_1, \ldots, b_k\} \cup \{a_n\}$ are the $k+1$ anti-excedances of $\pi$.
In particular, the labels of the rows of $D$ (including rows of length 0) are thus exactly the same as $\hat{D}$, with the addition of $a_n$.
Then, the shape of $D$ can be constructed from $\hat{D}$ by removing the column labelled $a_n$ and inserting in a row labelled $a_n$ in the appropriate position, maintaining the order of the labels of the new boundary lattice path, see Figure~\ref{fig:construct-D}.

\begin{figure}[ht]
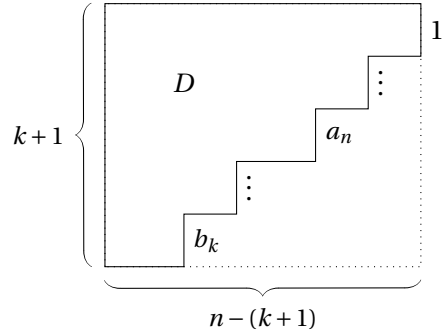

  \centering
  \[ \shapeD \]
  \caption[The shape and row labels of \Le-diagram $D$ associated to $\pi$.]{The shape and row labels of \Le-diagram $D$ associated to $\pi$. Here we have illustrated when $b_1 = 1, a_n \neq 1,n$. When $a_n = n$, we would have a row of size 0 labelled by $n$ (and thus is a co-loop) which would be indicated by a vertical line extending the west border of the diagram. Note that there are no columns of size 0 since $\pi$ is loopless.}
  \label{fig:D}
\end{figure}

The shape $\lambda$ of $D$ is then (including 0 sized parts)
\[ \lambda = \left(\lambda_1, \cdots, \lambda_{k+1}\right) \quad\text{where}\quad
  \lambda_{u} = \begin{cases} \hat\lambda_u - 1 & 1 \leq u \leq j-1, \\
  n - (k+1) - (a_n - j) & u = j,\\
  \lambda_{u} = \hat\lambda_{u-1} & j+1 \leq u \leq k+1, \end{cases} \]
and $j$ is the index such that $b_{j-1} < a_n < b_j$.
If $a_n > b_u$ for all $1 \leq u \leq k$, then let $j = k+1$.
Here we have $0 \leq \lambda_u \leq n - (k+1)$ for all $1 \leq u \leq k+1$, and at most $k+1$ non-zero parts, as needed.
Thus the order of rows (and anti-excedances) of $D$ is $\{b_1, \ldots, b_{j-1}, a_n, b_j, \ldots, b_k\}$, as in Figure~\ref{fig:D}.
Note that we always have either $b_1 = 1$ or $a_n = 1$, and thus the first row of $D$ will always be labelled with $1$ and is a full row of size $\lambda_1 = n - (k+1)$ (which also must hold as $\pi$ is loopless).

\bpoint{The algorithm}\label{SS:algorithm}

Now that we have the shape of $D$, we give an explicit algorithm for filling in $D$ such that we get the \Le-diagram with associated decorated permutation $\pi$.
We will fill $D$ row by row, from right to left, based on the corresponding row and rows below in $\hat{D}$.
For row $a_n$, since $a_n$ was a column in $\hat{D}$, we will consider its "corresponding row" in $\hat{D}$ as a row of all 0's of the same length as in $D$ and placed in the same position as $D$ (in-between rows $b_{j-1}$ and $b_j$ so that the order of labels of the SE border is maintained).

Let $L_u$ be the column containing the leftmost $\leplus$ in each row labelled $b_u$ of $\hat{D}$ (which is in box $(b_u, L_u)$) for $1 \leq u \leq k$.
Let $W_u$ be the row directly below row $b_u$, and let $W_n$ be the row directly below $a_n$, if they exist.
Otherwise, for the last row $b_k$, define $W_k=n+1$, or if $a_n$ is the last row, then define $W_n=n+1$.
Concretely we have
\[ W_u = \begin{cases} b_{u+1} & 1 \leq u < k, u \neq j-1 \\ a_n & u = j-1 \\ b_{j+1} & u = n, j \neq k+1 \\ n+1 & u=k, j \neq k+1 \text{ or } u=n, j=k+1 \end{cases} \]
In particular, $W_u - 1$ is the column right before the next row below, or for the last row, $W_u-1=n$.
Finally, we say that a $0$ is \textbf{restricted} if there is a $\leplus$ in a box to its left, in the same row, and \textbf{unrestricted} otherwise.
In Example~\ref{ex:notation}, the only restricted $0$ is in box $(1,2)$.

There are three different row types to consider for $D$:
\begin{enumerate}[label=(\Roman*)]
  \item Rows $b_u$ such that the box $(b_u, a_n)$ either does not exist, or does not have a $\leplus$ in $\hat{D}$.
  \item Rows $b_u$ such that the box $(b_u, a_n)$ has a $\leplus$ in $\hat{D}$. \\
    Note: These rows will always be above row $a_n$.
  \item Row $a_n$.
\end{enumerate}

We start by outlining the steps for row type (I); the algorithm for row types (II) and (III) will be slight modifications based on these steps.
Note that row types (I) and (III) actually follow the same steps if we consider the "leftmost" $\leplus$ in the "corresponding row" to $a_n$ in $\hat{D}$ as occurring in column $L_n > n$.

\tpointn[n]{Algorithm for filling in rows of $D$}\label{le-alg-main}
Main steps, and for row type (I)

Fill boxes in row $b_u$ from right to left, if they exist, as follows:
\begin{enumerate}[itemsep=4pt,label=Step \arabic*.]
  \item For columns $b_u+1$ to $W_u-1$, fill boxes in row $b_u$ with $\leplus$'s,\\
    i.e. boxes $(b_u, b_u+1), \ldots, (b_u, W_u-1)$ are filled with $\leplus$'s.
  \item For columns $W_u+1$ to $L_u-1$, fill boxes in row $b_u$ with $\leplus$'s under the conditions defined below.
    Otherwise, fill the boxes with 0's.\\
    i.e. fill boxes $(b_u, W_u+1), \ldots, (b_u, L_u-1)$ with $\leplus$'s if they satisfy the conditions.
  \item Starting from column $L_u$, fill the rest of the boxes in row $b_u$ with $0$'s,\\
    i.e. fill boxes $(b_u, L_u)$ and leftwards to the end of the row with $0$'s. \\
    Note: if $W_u  > L_u$, then Step 1 only fills boxes until column $L_u-1$ and Step 2 is skipped.
\end{enumerate}

Note that the algorithm can be run in parallel for each row.

\begin{figure}[h]
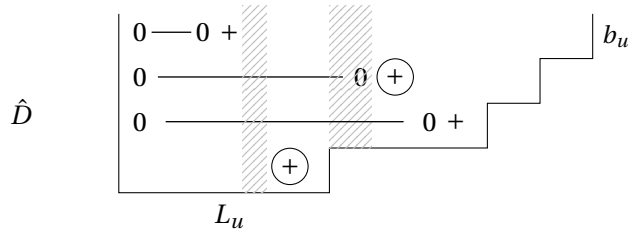

  \centering
  \[ \steptwoii \]
  \caption[Example of conditions in Step 2(ii) and (iii).]{Example of conditions in Step 2(ii) and (iii). For row $b_u$, the columns of the two circled $\leplus$'s in this $\hat{D}$ satisfy condition (ii) as all boxes above the $\leplus$ (until row $b_u$) are filled with unrestricted $0$'s. Thus in $D$, a $\leplus$ would be filled in row $b_u$ in those two columns.
  The column of the non-circled $\leplus$ not in row $b_u$ would fail this condition because of the row directly above it.
  The shaded columns indicate the columns that satisfy condition (iii) since all the boxes below row $b_u$ are filled with unrestricted $0$'s.\\}
  \label{fig:step2-conds}
\end{figure}

\tpointn[n]{Step 2 Conditions}\label{le-alg-conds}
Fill a box in row $b_u$ at $(b_u, \ell)$ with a $\leplus$ if:
\begin{enumerate}[itemsep=4pt]
  \item In $\hat{D}$, there is a $\leplus$ in box $(b_u, \ell)$.
  \item In $\hat{D}$, there is a $\leplus$ in some row below, say at $(b_m, \ell)$ where $m > u$, such that there are only unrestricted $0$'s in column $\ell$ in-between rows $b_u$ and $b_m$.
    For example, the circled $\leplus$'s in Figure~\ref{fig:step2-conds} satisfy this condition.
    Note: In the case that $m = u+1$, the condition holds trivially.\\
    Another way to phrase this condition: Look for $\leplus$'s in column $\ell$ in rows below $b_u$ where all rows in-between have no $\leplus$'s to the left of column $\ell$ (their leftmost $\leplus$ is before column $\ell$).

  \item In $\hat{D}$, column $\ell$ only has unrestricted $0$'s below row $b_u$.
    For example, the shaded columns in Figure~\ref{fig:step2-conds} satisfy this condition.
      In other words in column $\ell$, all the rows below $b_u$ has their leftmost $\leplus$ before column $\ell$, i.e. their last $\leplus$ has already passed. \\
      Note: This is a special case of (ii).\\
\end{enumerate}

\tpointn[n]{Modifications for row types (II) and (III)}\label{le-alg-mods}

For row type (II):
\begin{enumerate}[itemsep=4pt]
  \item[Steps 1 \& 2.] For columns $b_u+1$ until row $a_n$, follow Steps 1 \& 2 of row type (I), \\
    i.e. fill boxes $(b_u, b_u+1), \ldots, (b_u, a_n-1)$ as in row type (I).
  \item[Step 3.] Starting from column $a_n+1$, fill the rest of the boxes in row $b_u$ following only Step 2(i), \\
    i.e. fill boxes $(b_u, a_n+1)$ and leftwards until the end of the row using Step 2(i). \\
    Note: These boxes have the exact same filling as in the corresponding row in $\hat{D}$.
\end{enumerate}

For row type (III):
\begin{enumerate}[itemsep=4pt]
  \item[Steps 1 \& 2.] Follow Steps 1 \& 2 of row type (I) for the whole row $a_n$. \\
    Recall, the "corresponding row" to $a_n$ in $\hat{D}$ is a row of all $0$'s of the same length, in the same position as in $D$.
    Think of the leftmost $\leplus$ in this row as occurring at $(a_n, L_n)$ for some $L_n > n$.
    Under this convention, replace all instances of $b_u$ with $a_n$ and $L_u$ with $L_n$ in the main steps. \\
    Note: Step 2(i) never applies for this row.
\end{enumerate}

From these steps, we can already start to see some structure to these diagrams $D$ filled through this algorithm.
Later in Section~\ref{S:visual}, we will give a different perspective on how to build these diagrams that showcases the structure of the $\leplus$'s.

Some first observations are:
\begin{itemize}[itemsep=4pt]
  \item Step 1 tells us that there will always be a string of $\leplus$'s at the beginning of a row (until the next row below) in one of two ways:
    \begin{figure}[h]
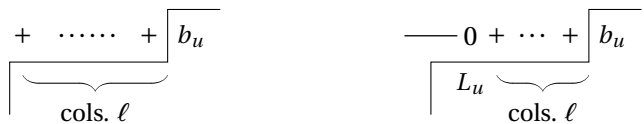

      \centering
      \[ \beginrow \]
      \caption[The beginning of a row using Step 1.]{The beginning of a row using Step 1. On the left, we have the case when $L_u > W_u$ and so $\leplus$'s are filled in columns $b_u+1 \leq \ell \leq W_u-1$. On the right, we have the case when $L_u < W_u$, in which Step 2 is skipped, and thus $b_u+1 \leq \ell \leq L_u-1$.}
    \end{figure}
  \item In particular, Step 1 also tells us that the bottom-most (non-empty) row of $D$ will be either:
    \begin{figure}[h]
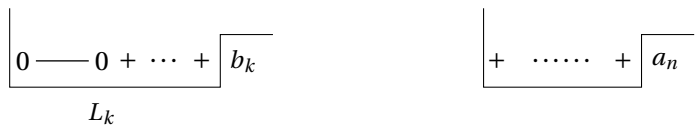

      \centering
      \[ \lastrow \]
      \caption[The last row using Step 1.]{The last row using Step 1. On the left, we have the case when $b_k > a_n$ is the last row, and so $\leplus$'s are filled until column $L_u-1$, with $0$'s filling the rest from Step 3. On the right, we have the case when $a_n > b_k$, and Step 1 fills the whole row with $\leplus$'s.}
    \end{figure}
  \item Step 2(ii) and (iii) tells us that for row type (III), there will always be a string a $\leplus$'s at the left end of this row:
    \[ \anend \]
  \item Step 2(iii) also tells us that the string of $\leplus$'s filled by satisfying this condition always comes after (is to the left of) a $\leplus$ that is filled based on Step 2(ii).
  \item Step 3 of row types (I) and (II) tells us that there will always be a string of $0$'s at the left end of row $b_u$ after column $L_u$:
      \[ \rowend \]
\end{itemize}

Finally, what's interesting about this algorithm is that for any particular row, only a part of $\hat{D}$ is looked at, namely it looks recursively at the row itself and the first rows below for which the leftmost $\leplus$ has not yet been passed. \\

\bpoint{\texorpdfstring{$D$}{D} is indeed a \texorpdfstring{\BLe-}{Le }diagram}

To even be able to consider this algorithm, we first verify that under this filling $D$ avoids the \Le-configuration and thus is a valid \Le-diagram.

\tpointn{Theorem}\label{Dlediag}
\statement{
  Under this filling, $D$ is a \Le-diagram of type $(k+1, n)$.
}

\begin{proof}
  We show that if a \Le-configuration occurs in $D$, then one must occur in $\hat{D}$ as well!

  Without loss of generality, suppose in $D$ we have the following
  \[ \leconfiglabels \]
  where the boxes indicated by the dots are filled with $0$'s.
  There are two cases for row $j$,
  \begin{enumerate}
    \item either there's a corresponding row $j$ in $\hat{D}$, or
    \item $j= a_n$.
  \end{enumerate}

  First, based on our observations from Steps 1 and 3, the three boxes in this configuration must have been filled in Step 2 for rows $i$ and $j$, regardless of row type.
  In particular, that means row $j$ is not the last row, and furthermore there must be boxes in columns $\ell, \ldots, m$ below row $j$ in both $D$ and $\hat{D}$.

  Looking at the $\leplus$ at $(j, m)$ in $D$, since it was filled based on a condition in Step 2, in each case we have:
  \begin{enumerate}[itemsep=4pt]
    \item The last $\leplus$ in the corresponding row in $\hat{D}$ must be to the left of or at column $m$, and thus in $\hat{D}$ there is a $\leplus$ at some $(j, L)$ where $L \geq m$.
    \item In this case, the box must have been filled by Step 2(ii), and thus in $\hat{D}$, there is some first row $b$ below $j$ such that its last $\leplus$ is to the left of or at column $m$, and thus there is a $\leplus$ at $(b, L)$ where $L \geq m$.
      That is, all rows in-between $j$ and $b$ must have boxes in column $m$ filled with unrestricted $0$'s i.e. their last $\leplus$ is to the right of $m$.
  \end{enumerate}

  \begin{figure}[h]
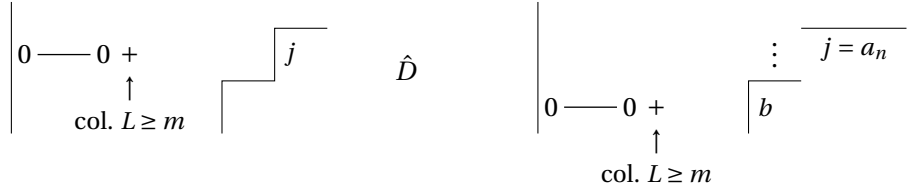

    \centering
    \[ \lediagproof \]
    \caption[The two cases for $j$ in $\hat{D}$ based on a $\leplus$ at $(j,m)$ in $D$.]{The two cases for $j$ in $\hat{D}$ based on a $\leplus$ at $(j,m)$ in $D$. On the left, we have case (i), and on the right, we have case (ii). Here $b$ is the first row below $j=a_n$ where all boxes in column $m$ in-between row $b$ and column $j$ are filled with unrestricted $0$'s. }
  \end{figure}

  Next, looking at the $0$ at $(j, \ell)$, in each case we have:
  \begin{enumerate}[itemsep=4pt]
    \item Because of Step 2(i) and the last $\leplus$ of row $j$ in $\hat{D}$ is at some column $L \geq m >  \ell$, we know that there is no $\leplus$ in $\hat{D}$ at $(j, \ell)$, and thus it is filled with $0$.
    \item Since we are assuming boxes $(j, \ell+1), \ldots, (j, m-1)$ are also filled with $0$'s in $D$, we must have that all rows in-between $j$ and $b$ in $\hat{D}$ have unrestricted $0$'s filled in those columns.
      As otherwise Step 2(ii) would have placed a $\leplus$ in $D$ somewhere in those boxes.
      But then this also means that in $\hat{D}$, $(b, \ell)$ must be filled with a $0$ since otherwise Step 2(ii) would be satisfied, placing a $\leplus$ at $(j, \ell)$ in $D$.
  \end{enumerate}

  Finally we look at the $\leplus$ at $(i, \ell)$.
  The $\leplus$ at $(j, m)$ in $D$ tells us that Step 2(iii) could not have been the reason for this $\leplus$ since in $\hat{D}$, the last $\leplus$ of row $j$ or row $b$, for each case respectively, have not passed yet.
  Thus this $\leplus$ must have come from satisfying Step 2(i) or (ii) for row $i$ (regardless of type) which means in $\hat{D}$ in column $\ell$, there must be a $\leplus$ in rows $i$ or below.
  The $0$ at $(j, \ell)$ in $D$ then tells us that in $\hat{D}$, there cannot be any $\leplus$'s in rows below $j$ in column $\ell$ or else Step 2(ii) would have been satisfied for $(j, \ell)$.
  Then, regardless of which case we are in, there must be a $\leplus$ in $\hat{D}$ in column $\ell$ above row $j$.

  Putting all these together, notice in either case we have constructed a \Le-configuration occurring in $\hat{D}$.
  But this cannot happen as $\hat{D}$ is a \Le-diagram!
  Thus the \Le-configuration cannot occur in $D$ and $D$ is indeed a valid \Le-diagram.
  Furthermore, from the shape $\lambda$ of $D$, we get that $D$ is a \Le-diagram of type $(k+1, n)$ as needed.
\end{proof}

\bpoint{An example}\label{ex:algorithm}

To illustrate this algorithm, we work through an example step by step.

Suppose we had the following \Le-diagram $\hat{D}$ of type $(4,n)$ and dimension $2k = 8$, where all unmarked boxes are filled with $0$'s.
The circled $\leplus$'s indicate the last $\leplus$'s in each row.\\
\[ \algexhatD \]

First, since in $\hat{D}$ we have that $\hat\pi(1) = a_n = m$, the shape of $D$ is\\
\[\algexDshape \]

Here we have the following variables,
\begin{align*}
  b_1 &= 1 & b_2 &= i & a_n &= m  & b_3 &= o & b_4 &= q \\
  L_1 &= u & L_2 &= \ell & L_n &> n & L_3 &= s & L_4 &= t  \\
  W_1 &= i & W_2 &= a_n & W_n &= o & W_3 &= q & W_4 &= n + 1
\end{align*}
where recall that $b_u$ are the rows of $\hat{D}$, $L_u$ is the column of the last $\leplus$ in row $b_u$ (with the last $\leplus$ for $a_n$ thought of as at some $L_n > n$) and $W_u$ is the next row below $b_u$ (or $n+1$ if its the last row).

To fill $D$, we follow the algorithm row by row, starting with row $1$.
As row $1$ has a $\leplus$ in column $a_n=m$ in $\hat{D}$, it is of row type (II) and we obtain the filling:
\[ \exrowone \]
\begin{enumerate}[itemsep=4pt,label=\arabic*.]
  \item From Step 1, we fill boxes $(1,2), \ldots, (1,i-1)$ with $\leplus$'s as the last $\leplus$ has not yet passed ($L_1 = u$) and the next row below is row $i$.
  \item From Step 2, we fill boxes $(1,i+1), \ldots, (1,m-1)$ with $\leplus$'s according to the conditions, and $0$'s otherwise.
    We stop at column $m-1$ since $a_n = m$, and we are following row type (II).\\
    \begin{enumerate}[itemsep=4pt,label=(\roman*)]
      \item Row $1$ has no $\leplus$'s in these boxes in $\hat{D}$ so Step 2(i) is not satisfied.
      \item Step 2(ii) fills a $\leplus$ in box $(1, \ell)$, as the $\leplus$ at $(i, \ell)$ in $\hat{D}$ satisfies this condition trivially since there are no rows in-between rows $1$ and $i$.
      \item Step 2(iii) then fills $\leplus$'s in boxes $(1, \ell+1), \ldots, (1,m-1)$, since the $\leplus$ at $(i, \ell)$ in $\hat{D}$ was in fact the last $\leplus$ of row $i$.
    \end{enumerate}
  \item Finally, Step 3 only fills a $\leplus$ in box $(1,u)$, with the rest filled with $0$'s, since in $\hat{D}$ there is only one $\leplus$ after column $a_n =m$ which is in column $u$.
\end{enumerate}

Next, row $i$ does not have a $\leplus$ in column $m$ in $\hat{D}$ and thus is of row type (I).
As the last (and only) $\leplus$ of row $i$ is at $(i, \ell)$ in $\hat{D}$, we have that $W_2 = m > L_2 = \ell$ and so Step 2 is skipped.
Then, Step 1 fills boxes $(i, i+1), \ldots (i, \ell - 1)$ with $\leplus$'s, and Step 3 fills the rest of the boxes with $0$'s, to obtain:
\[ \exrowi \]\\
Moving on to row $m=a_n$ which is of row type (III), recall that we are considering "row" $a_n$ in $\hat{D}$ to be a row of all $0$'s with the same length and position as in $D$, and having its last $\leplus$ at some column $L_n > n$.
\[ \exrowm \]
\begin{enumerate}[itemsep=4pt,label=\arabic*.]
  \item From Step 1, we fill boxes $(m,m+1), \ldots, (m,o-1)$ with $\leplus$'s as the last $\leplus$ has not yet passed ($L_n > n$) and the next row below is row $o$.
  \item From Step 2, we fill the rest of the boxes in the row with $\leplus$'s according to the conditions and $0$'s everywhere else.
    As we are in row type (III), we skip Step 3.\\
    \begin{enumerate}[itemsep=4pt,label=(\roman*)]
      \item Step 2(i) never applies for this row.
      \item Step 2(ii) fills a $\leplus$ in boxes $(m, p)$, $(m, s)$ and $(m,t)$.
        These columns correspond to $\leplus$'s in $\hat{D}$ in rows below $m$ such that all the boxes above are unrestricted $0$'s.
      \item Step 2(iii) then fills $\leplus$'s in boxes from $(m, t+1)$ to the end of the row, since the $\leplus$ at $(q, t)$ in $\hat{D}$ is the leftmost $\leplus$ of all rows below $m$ and thus all columns after $t$ are filled with unrestricted $0$'s below $m$.
    \end{enumerate}
\end{enumerate}

After row $a_n$, all the rows below will be of type (I).
For row $o$ we get the filling:
\[ \exrowo \]
Here, Step 1 fills the beginning of row $o$ with $\leplus$'s until the next row below, which is row $q$, as the last $\leplus$ in row $o$ in $\hat{D}$ is in column $L_3 = s$.
Note that we don't consider the $\leplus$ in box $(q,p)$ in $\hat{D}$.
Step 2 fills a $\leplus$ only at box $(o,r)$ through condition (ii), and $0$'s everywhere else until column $s$.
Step 3 then fills the rest of the row, starting from column $s$, with $0$'s.

Finally, we have row $q$.
As its the last row, we skip Step 2 (recall we are considering $W_4 = n+1$) and since the last $\leplus$ of row $q$ is in box $(q,t)$ in $\hat{D}$, we have:
\[ \exrowq \]

Putting everything together, we get the following \Le-diagram $D$:

\[ \algexD \]

To briefly check that this $D$ is indeed the right diagram for the given $\hat{D}$:
\begin{itemize}[itemsep=4pt]
  \item $D$ is a \Le-diagram of type $(5,n) = (k+1,n)$ (notice it avoids the \Le-configuration and has $5$ rows).
  \item Every column has at exactly one $\leplus$ except columns $p,r,s,u$, which have two $\leplus$'s.
    Thus $D$ is loopless and we have dimension (recall this is the number of $\leplus$'s)
    \[ \dim(S_D) = (n-5) + 4 = n - (k+1) + k = n-1, \]
    as there are $n-5$ columns and 4 columns with an extra $\leplus$.
    This gives the correct relation since we wanted $\dim{S_{\hat{D}}} - 2k = \dim(S_D) - (n-1)$, recalling that $\dim(S_{\hat{D}}) = 8 = 2k$.
  \item One can check that the associated decorated permutation $\pi$ to $D$, in two-line notation, is
    \[\indent\indent\left(\arraycolsep=1pt\begin{array}{*{23}{c}} 1 & {\cdots} & i-1 & {\cdots} & \ell-1 & {\cdots} & m-1 & {\cdots} & o-1 & {\cdots} & p-1 & {\cdots} & q-1 & {\cdots} & r-1 & {\cdots} & s-1 & {\cdots} & t-1 & {\cdots} & u-1 & {\cdots} & n  \\ 2 & & \ell & & \circled{i} & & u & & p & & s & & r & & \circled{o} & & t & & \circled{q} & & \circled{1} & & \circled{m} \end{array}\right), \]
      where the $\;\cdots\;$ denotes $\pi(a) = a+1$ for all $a$ in-between the explicitly written values.
      Notice $\pi$ is loopless and has $5$ anti-excedances (circled).
      Shifting the bottom-line one to the right (and wrapping around) exactly gives $\hat\pi$, the decorated permutation associated to $\hat{D}$.
\end{itemize}

To compare with $\hat{D}$, we redraw $D$ to illustrate the structure of the $\leplus$'s (as we will discuss in the next section).
All the unmarked boxes are filled with $0$'s.

\begin{figure}[h]
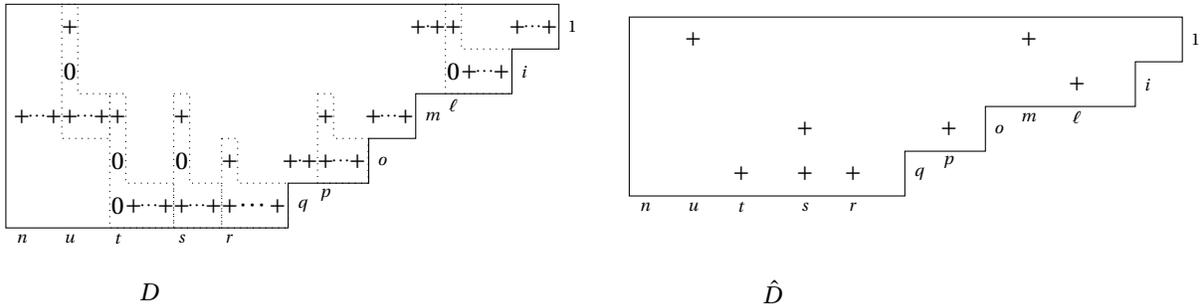

  \centering
  \[ \algexDshapes \]
  \caption[An example of ${D}$ vs. $\hat{D}$, where $D$ (left) is filled by the algorithm for the given $\hat{D}$ (right).]{An example of ${D}$ vs. $\hat{D}$, where $D$ (left) is filled by the algorithm for the given $\hat{D}$ (right). In $D$, the dotted lines outlining shapes of $\leplus$'s and $0$'s are to illustrate the structure of the $\leplus$'s.
  Comparing to the original $\hat{D}$, notice that there is one of these shapes per column of $\leplus$'s in $\hat{D}$, except column $a_n = m$.}
  \label{fig:algex}
\end{figure}

\section{A more visual perspective}\label{S:visual}

While an algorithm is nice to have, it doesn't illuminate the underlying structure of the $\leplus$'s in $D$, which we are starting to see (Figure~\ref{fig:algex}).
Perhaps a more enlightening, and surprising, viewpoint is to look at how the $\leplus$'s in $\hat{D}$ give rise to certain shapes, which when glued together fully characterize $D$ (see Theorem~\ref{glueing}).
That is, without going through the algorithm, we can reconstruct $D$ using these shapes.
Whereas the algorithm can be seen as viewing $D$ through each row, here we are viewing $D$ through each column.
Under this perspective, it becomes clearer how a reverse map on \Le-diagrams would work. \\

\bpoint{Building blocks}

The first aspect of the surprising structure that arose is that $D$ is composed of 2 distinct shapes; the main one being what we will call $\mathbb{L}$-shapes and also strings of $\leplus$'s.
In fact, the strings of $\leplus$'s can be considered as a special type of $\;\mathbb{L}$-shape where there is only the horizontal part.

We will use the same notation as in Section~\ref{S:algorithm} for \Le-diagrams, see Example~\ref{ex:notation}.
Recall that boxes in a \Le-diagram are referred to by their coordinates $(i,j)$, which denotes row $i$ and column $j$.
We also interchangeably use "first" to mean top/right-most, "last" to mean bottom/left-most, and "passed" to mean in some column to the right of the current column in consideration, with "not yet passed" meaning to the left of or at the current column.

Recall that $\hat{D}$ is a co-loopless \Le-diagram of type $(k,n)$, shape $\hat\lambda$, and with decorated permutation $\hat{\pi}$.
We are denoting its rows by $b_1, \ldots, b_k$ and distinguishing the column $a_n$, which is defined by $\hat{\pi}(1) = a_n$.
The algorithm in Section~\ref{S:algorithm} is constructing a loopless \Le-diagram $D$ of type $(k+1,n)$, shape $\lambda$, and with decorated permutation $\pi$, as given by T-duality (Definition~\ref{def:t-duality}) for the given $\hat\pi$.
The rows of $D$ are $\{b_1, \ldots, b_k, a_n\}$, with row $a_n$ inserted as in Figure~\ref{fig:construct-D}.
From the algorithm, there are three row types to consider for ${D}$: (I) rows with no $\leplus$'s in column $a_n$, (II) rows with $\leplus$'s in column $a_n$ and (III) row $a_n$.

To create these $\mathbb{L}$-shapes, we will look at columns $\ell \neq a_n$ in $\hat{D}$ with at least one $\leplus$.
As such, there are some facts we need to know about these columns:
\begin{itemize}
  \item In $\hat{D}$, since every row contains at least one $\leplus$ (as co-loopless) and we have the \Le-condition, all $0$'s below a $\leplus$ in column $\ell$ must be unrestricted, i.e. these rows have their last $\leplus$ to the right of $\ell$:
      \[ \lastplustoright \]
  \item For $\leplus$'s in column $\ell$ at $(b, \ell)$:
    \begin{itemize}
      \item If $a_n < b$ then $b$ is of row type (I) as $(b, a_n)$ is not a box within the shape $\hat\lambda$.
      \item If $b < a_n < \ell$ then $b$ is of row type (II) since there is a $\leplus$ at $(1, a_n)$ and so the \Le-condition forces there to be a $\leplus$ at $(b, a_n)$.
      \item If $\ell < a_n$, then $b$ can be of row type (I) or (II).
    \end{itemize}
\end{itemize}
We refer to $\leplus$'s in $\hat{D}$ that are not a leftmost (last) $\leplus$ in its row as \textbf{non-last}.

\tpointn{Theorem}\label{shapes}
\statement{
  For every column $\ell \neq a_n$ in $\hat{D}$ with at least one $\leplus$, in $D$ there is a corresponding shape
    \[ \lshapelabels \]
  which we call a \textbf{$\;\mathbb{L}$-shape}, and which has the following properties:
  \begin{itemize}[itemsep=4pt]
    \item There is always a $\leplus$ in box $(b_T, \ell)$, with some sequence of 0's and $\leplus$'s below.
    \item There is a string of $\leplus$'s in boxes $(b_B, m), \ldots, (b_B, \ell-1)$, possibly empty (if $m = \ell$).
    \item In $\hat{D}$, there are $s_{\ell}$ non-last $\leplus$'s in column $\ell$.
    \item In $\hat{D}$, there are $t_{\ell}$ last $\leplus$'s in rows of type (II) in column $\ell$.
    \item There is one more $\leplus$, not counted by $s_{\ell}$ or $t_{\ell}$, that is in row $a_n$ or in row $b_T \neq f, a_n$ in column $\ell$, where $f$ is the row of the topmost $\leplus$ in column $\ell$ in $\hat{D}$.
    \item In $\hat{D}$, $m - 1$ is either the label of the first column with at least one $\leplus$ to the right of $\ell$, or the first row above $\ell$, whichever comes first.
  \end{itemize}
  Moreover, we can exactly characterize the $\leplus$'s in the $\;\mathbb{L}$'s, and these are the only $\leplus$'s in these columns in $D$.
}

\begin{proof}
  Let $\ell \neq a_n$ be a column in $\hat{D}$ with at least one $\leplus$.
  Let $(f, \ell)$ be the first (topmost) $\leplus$ in this column and $(g, \ell)$ be the last (bottom-most).

  First, to determine the horizontal string of $\leplus$'s in the $\mathbb{L}$, we look at $\hat{D}$ around the $\leplus$ at $(g, \ell)$.
  Using the first fact from above, if there are boxes below $(g,\ell)$, there must be a first column call it $m-1$, to the right of $\ell$ with at least one $\leplus$.
  Furthermore, in this column there must be $\leplus$'s in rows $g$ or below.
  That is all columns in-between $m-1$ and $\ell$ and in rows $g$ and below, are filled with only 0's.
  By the \Le-condition, these columns are then actually columns of all $0$'s.
    \[ \lproofone \]

  If there are no rows below $(g, \ell)$, then we let $m-1 = g$ if $(g, \ell)$ is the first $\leplus$ in its row, or let $m-1$ be the column such that $(g,m-1)$ is the first $\leplus$ to the right of $(g, \ell)$.
    \[ \lprooftwo \]

  If row $g$ is of type (I) or of type (II) with $\ell < a_n$, then the boxes $(g, \ell -1), \ldots, (g, m)$ are filled with a string of $\leplus$'s in $D$ by Step 2(iii) if there are rows below $(g,\ell)$, or by Step 1 if there are no rows below (recall that row type (II) follows the same steps as row type (I) until the column $a_n$).
  Thus we let $b_B = g$.

  If row $g$ is of type (II) with $g < a_n < \ell$, then notice in $D$, the row $a_n$ will be inserted in-between $g$ and $\ell$.
  In this case, we have that there will always be a first column $m-1$ to the right of $\ell$ with at least one $\leplus$ in $\hat{D}$, as $m-1$ is at least $a_n$.
  Note that we can have $m-1 = a_n$.
  \[ \lproofthree \]
  Now we repeat a similar argument as before, except for boxes $(a_n, \ell -1), \ldots, (a_n, m)$.
  These are filled with a string of $\leplus$'s in $D$ by Step 1 for row type (III), if $a_n$ becomes the row right above $\ell$, or by Step 2(iii) otherwise.
  Thus we let $b_B = a_n$.

  Furthermore, we argue that these are the only $\leplus$'s in columns $m$ to $\ell - 1$ in $D$.
  Notice that from the algorithm, for columns like $m$ in $\hat{D}$ filled only with 0's, the only way that a box in this column in $D$ can be filled with a $\leplus$ is through some Step 1 or Step 2(iii) for any row type.
  So we just need to show that if say box $(b_B, m)$ is filled with a $\leplus$, then no other rows in column $m$ satisfy these two conditions.
  First note that Step 1 can only affect the very last box of column $m$ and that to have $(b_B, m)$ filled with a $\leplus$ means the last $\leplus$ in row $b_B$, say at $(b_B, L)$, must not have passed yet ($L > m$).

  Right away we see that for rows $b$ below $b_B$ in column $m$, if $(b_B, m)$ was filled through Step 1, then there's nothing to check as $b_B$ was the last row.
  If $(b_B, m)$ was filled through Step 2(iii), then all rows $b$ must have their last $\leplus$ to the right of $m$.
  But then we cannot be in Step 1 nor Step 2(iii) for $(b, m)$ since the last $\leplus$ in row $b$ has already passed.
    \[ \lprooffour \]

    For any rows $b$ above $b_B$, we can never be in Step 1 since these are not the last row in the column.
  We also can never satisfy Step 2(iii) for $(b, m)$ since this requires all rows below $b$ to have their last $\leplus$'s to the right of $m$, which does not hold for row $b_B$, and now we are done.

  Thus we get for the horizontal part of the $\mathbb{L}$ in $D$:
  \begin{itemize}[itemsep=4pt]
    \item $b_B = g$, or $b_B = a_n$ if row $g$ is of type (II) and $g < a_n < \ell$.
    \item $m$, possibly equal to $\ell$, corresponds to the column such that all columns in-between $\ell$ and $m-1$ in $\hat{D}$ are filled with only 0's.
      That is $m-1$ is either the first column with at least one $\leplus$ to the right of $\ell$, or $m-1$ is a row.
    \item For $m \neq \ell$, boxes $(b_B, \ell -1), \ldots, (b_B, m)$ are filled with $\leplus$'s, and these are the only $\leplus$'s in columns $\ell - 1$ to $m$.
      Otherwise if $m = \ell$, there are no boxes in the horizontal part.
  \end{itemize}

  Next we look at the vertical part of the $\mathbb{L}$, including box $(g, \ell)$, which corresponds to the column $\ell$.
  First, using the first fact from above once again, this time for $\leplus$'s in column $\ell$, we see that any row $b$ below $f$ with a 0 filled at box $(b, \ell)$ in $\hat{D}$ must have its last $\leplus$ in some column to the right of $\ell$.
  Note here we are not considering $b=a_n$.
  Thus for these rows, Step 3 for row types (I) and (II) will fill box $(b, \ell)$ in $D$ with a 0 as well.
  For the rows $b$ with a $\leplus$ filled at $(b, \ell)$ in $\hat{D}$, in $D$ box $(b, \ell)$ is filled with a $\leplus$ through Step 2(i) for row type (I) if it is a non-last $\leplus$, and Step 2(i) and Step 3 for row type (II), regardless if it is a last $\leplus$ or not.
  Note, if $(b, \ell)$ is a last $\leplus$ in row $b$ of type (I), then Step 3 fills this box in $D$ with a 0.
  This gives $s_{\ell} + t_{\ell}$ $\leplus$'s in column $\ell$ in $D$, determined exactly by the $\leplus$'s in column $\ell$ in $\hat{D}$, where $s_{\ell}$ is the number of non-last $\leplus$'s in column $\ell$ of $\hat{D}$ and $t_{\ell}$ is the number of last $\leplus$'s in rows of type (II) in column $\ell$ of $\hat{D}$.

  To determine row $b_T$, we look at $\hat{D}$ above the $\leplus$ at $(f, \ell)$.
  As all boxes directly above $(f, \ell)$ are filled with $0$'s in $\hat{D}$, from the algorithm, the only way to fill these boxes $(b, \ell)$ in $D$ with a $\leplus$ is through Step 2(ii).

  If row $f$ is type (I) or of type (II) with $\ell < a_n$ then we look for the first row above $f$ for which in $\hat{D}$ its last $\leplus$ has not yet passed.
  Note that such a row must always exist if $f$ is of type (II) with $\ell < a_n$, and this row will also be of type (II) with $\ell < a_n$.
  If such a row $b$ exists and it is of type (I) or type (II) with $\ell < a_n$, then Step 2(ii) for row $b$  would place a $\leplus$ at $(b, \ell)$ in $D$ and we let $b_T = b$.
  All rows in-between $b$ and $f$ would have had their last $\leplus$ to the right of $\ell$ and so these rows are filled with a 0 in column $\ell$ in $D$.

  If no such row exists, then either there are no rows above $f$ and thus $a_n = 1$ (if $f=1$ then $f$ is of row type (II) with $a_n < \ell$) or all rows above $f$ have their last $\leplus$ to the right of column $\ell$ and thus $a_n < \ell$.
  In both cases, $f$ must be of row type (I) and with $a_n < f$ since in the second case, if $f < a_n < \ell$, then by the \Le-condition there would be a $\leplus$ at $(f, a_n)$ contradicting our assumption of $f$'s row type.
  If a row $b$ exists but is of type (II) with $\ell > a_n$, in which case $a_n > b$, then once again because of the \Le-condition we actually have $a_n < f$ (as $f$ must be of type (I)).
  In either case, we get that Step 2(ii) for row $a_n$ will fill box $(a_n, \ell)$ with a $\leplus$ since all rows in-between $a_n$ and $f$ would have all their last $\leplus$'s already passed and furthermore all these rows are filled with a 0 in column $\ell$ in $D$.
  Thus we let $b_T = a_n$.

  If $f$ is of type (II) with $\ell > a_n$, in which case we also have $a_n > f$, then by the \Le-condition any rows $b$ above $f$ with $\leplus$'s to the left of column $\ell$ would have to have a $\leplus$ at $(b, a_n)$ making them of row type (II) as well.
  The rest of the rows, would have their last $\leplus$ to the right of $\ell$.
  Then Step 3 for row types (I) and (II) would fill all the rows $b$ above $f$ with a 0 at $(b, \ell)$.
  Since $f$ is of type (II), Step 2(i) would fill $(f, \ell)$ with a $\leplus$, and thus we let $b_T = f$.

  One detail left over is that for boxes in-between $g$ and $b_B$ in column $\ell$, by the \Le-condition all these rows must have their last $\leplus$ already passed (or else it contradicts $g$) in $\hat{D}$ and thus Step 3 for row types (I) and (II) would fill these boxes with a 0.

  Thus we get for the vertical part of $\mathbb{L}$ in $D$:
  \begin{itemize}[itemsep=3pt]
    \item There are three cases for $b_T$. Either $b_T$ is the first row above $f$, and below $a_n$ if $f > a_n$, with a $\leplus$ to the left of column $\ell$ in $\hat{D}$, or $b_T = f$ if row $f$ is of type (II) with $f < a_n < \ell$, or $b_T = a_n < f$.
    \item There is always a $\leplus$ at $(b_T, \ell)$.
    \item There is a $\leplus$ in box $(b, \ell)$ if there is a non-last $\leplus$ at $(b, \ell)$ in $\hat{D}$, or if there is a last $\leplus$ but $b$ is of row type (II), for rows $b$ from $f$ to $g$ in $\hat{D}$.
      Thus there is a total of $s_{\ell} + t_{\ell}$ $\leplus$'s in boxes $(f, \ell), \ldots, (g, \ell)$ in $D$, not considering row $a_n$.
    \item We always have $b_T \leq f$ and $b_B \geq g$, and thus these $\leplus$'s are all within the vertical part of the $\mathbb{L}$.
    \item Any rows in-between $b_T$ and $f$ in column $\ell$ must be filled with 0's, and similarly for $g$ and $b_B$.
  \end{itemize}

  Finally, we show that there are exactly $s_{\ell} + t_{\ell} + 1$ $\leplus$'s in column $\ell$, in the vertical part of the $\mathbb{L}$.
  To do this, we look at the different cases of $b_B$, $b_T$ and how row $a_n$ affects the $\mathbb{L}$.
  We will show that in fact $b_T < b_B$ and the extra $\leplus$ will be either in $a_n$ or we have $b_T \neq f, a_n$.
  Doing case analysis on $a_n$:
  \begin{itemize}[itemsep=4pt]
      \item If $a_n < f \leq g < \ell$, then $b_B = g$, and $b_T = a_n$ or $b_T$ is the first row above $f$, below $a_n$ with a $\leplus$ to the left of column $\ell$.\\
        In either case the extra $\leplus$ is at $(b_T, \ell)$ and we have $b_T < b_B$.
      \item If $f < a_n < g < \ell$, then $b_B = g$ and $b_T = f$.\\
        In this case, Step 2(ii) for row $a_n$ fills $(a_n, \ell)$ with a $\leplus$ since there is at least one more $\leplus$ in rows below $a_n$ in column $\ell$ in $\hat{D}$, with all rows in-between filled with unrestricted $0$'s by the \Le-condition.
      \item if $f \leq g < a_n < \ell$, then $b_B = a_n$ and $b_T = f$. \\
        Here the extra $\leplus$ is at $(b_B, \ell) = (a_n, \ell)$.
      \item Lastly if $f \leq g < \ell < a_n$, then $b_B = g$ and $b_T \neq f, a_n$.\\
        Here the extra $\leplus$ is at $(b_T, \ell)$.
  \end{itemize}
  Thus we are done since these are the only cases, and in each we have $b_T < b_B$ with exactly $s_{\ell} + t_{\ell} + 1$ $\leplus$'s in the vertical part of the $\mathbb{L}$ with a $\leplus$ at $(b_T, \ell)$.

  One final detail to check is that these are the only $\leplus$'s in column $\ell$.
  Since we always have that $b_T \leq f$ and $b_T \geq g$ with $b_T < b_B$, if there was a $\leplus$ at $(a, \ell)$ in $D$, for rows $a < b_T$ or $a > b_B$, we can show this contradicts the choice of $b_T$ or $b_B$.

  For $b < b_T \leq f$, the $\leplus$ at $(b, \ell)$ must have been filled through Step 2(ii) for row $b$ which means that row $b$'s last $\leplus$ in $\hat{D}$ has not yet passed.
  This also means that any rows $c$ below $b < c < f$ must have their last $\leplus$'s to the right of column $\ell$.
  But this contradicts when $b_T \neq f$ since in these cases, $b_T$ has their last $\leplus$ to the left of column $\ell$ (for $b_T=a_n$ remember we are thinking of row $a_n$ as having last $\leplus$ in a column $> n$).
  Thus $b_T = f$, which means that $f$ is of row type (II) with $a_n > f$.
  But because of the \Le-condition, all rows above $f$, in particular $b$, must also be of type (II) with $a_n > b$, in which case the algorithm would not have filled $(b, \ell)$ with a $\leplus$.
  Thus all rows $b$ above $b_T$ must be filled with a 0 at $(b, \ell)$.

  For $b > b_B \geq g$, the $\leplus$ at $(b, \ell)$ must have been filled through Step 2(iii) for row $b$.
  Note, it could not have been filled through Step 1 since by the \Le-condition, that would mean there was a $\leplus$ at $(b, \ell)$ in $\hat{D}$ originally, contradicting $g$ (since row $b$'s last $\leplus$ in $\hat{D}$ must have been to the left of column $\ell$ and there's a $\leplus$ at $(g, \ell)$).
  Once again, Step 2(iii) means that row $b$'s last $\leplus$ in $\hat{D}$ has not yet passed.
  By the \Le-condition for $\hat{D}$, this would force $b = a_n$ or else we would contradict $g$, and thus $b_B = g$.
  However, as there is a box at $(b, \ell)$, we actually have $\ell > b = a_n > b_B = g$.
  Looking at our case analysis above, this contradicts the choice of $b_B$, which should have been $a_n$.
  Thus all rows $b$ below $b_B$ must be filled with a 0 at $(b, \ell)$.

  This concludes the proof.
\end{proof}

Notice these $\mathbb{L}$-shapes cover all columns from $\hat{D}$, ignoring $a_n$, with at least one $\leplus$ and the columns of all 0's to the right of these.
What's left are the (consecutive) columns of all 0's to the left of a column with at least one $\leplus$ until the next row below or the end of the diagram, or if in-between two rows there are only columns of all 0's.
These columns are covered by strings of $\leplus$'s, which is actually just the horizontal part of an $\mathbb{L}$-shape.

\tpointn{Corollary}\label{strings}
\statement{
  For consecutive columns of all 0's in $\hat{D}$, in $D$ there is a corresponding string of $\leplus$'s\\
    \[ \stringplus \]
  and these are the only $\leplus$'s in these columns.
}

\begin{proof}
  From the proof of Theorem~\ref{shapes}, we saw that this statement holds for consecutive columns of all $0$'s in $\hat{D}$ in-between two columns with at least one $\leplus$, or in-between the first column to the left of a row with a least one $\leplus$ and the beginning of that row.
  What's left are consecutive columns of all $0$'s in-between two consecutive rows, and there are no columns with at least one $\leplus$ in-between these rows, and those in-between the next row below or the end of the diagram and a column with at least one $\leplus$.
  Both of these are really special cases of what we proved in Theorem~\ref{shapes}.

  Say we are considering columns with last row $a$ and the next row below is row $b$ or if we are at the end of the diagram we can think of the "next row below" as at $n+1$.
  Say $m-1 = a$ for the first case or $m-1 > a$ is the last column with at least one $\leplus$ in $\hat{D}$ for the second case, and we are considering columns $c$ of all 0's for all $m-1 < c < b$, or $m-1 < c < n+1$.
  Then we can imagine a column associated to $b$ or $n+1$, where the last row in the column is $a$ and it extends to the top of the diagram.
  We place $\leplus$'s in this "column" in the rows in $\hat{D}$ where there are $\leplus$'s to left of column $b-1$, that is their last $\leplus$'s have not yet passed and these rows are not of type (II) with $a_n < m-1$.
  Then we can think of this "column" as the column $\ell$ in Theorem~\ref{shapes}.

  If there are no $\leplus$'s in this "column", i.e. all last $\leplus$'s are before $m$, then we must have $a_n < m$ in which case Step 2(iii) for row $a_n$ would place $\leplus$'s at $(a_n, c)$ for all $c$.
  In particular, this always occurs if we are at the end of the diagram.

  If there are $\leplus$'s in this "column", then once again we look at $g$, the row of the last $\leplus$ i.e. the first row above $b$ where it's last $\leplus$ has not yet passed.
  If $g=a$, then Step 1 for row $a$ fills $(a, c)$ with $\leplus$'s for all $c$.
  If $g < a$, then Step 2(iii) for row $g$ fills $(g, c)$ with $\leplus$'s for all $c$.

  Once again, we can argue in the same way as Theorem~\ref{shapes} that these are the only $\leplus$'s in columns $c$.
\end{proof}

From Theorem~\ref{shapes} and Corollary~\ref{strings}, we immediately get the following statement.

\tpointn{Corollary}\label{number}
\statement{
  There is at least one $\leplus$ in every column of $D$.
  More specifically,
  \begin{enumerate}
    \item there is exactly one $\leplus$ in columns corresponding to those in $\hat{D}$ with no $\leplus$'s, and
    \item there are $s_{\ell} + t_{\ell} + 1$ $\leplus$'s in columns $\ell$ corresponding to those in $\hat{D}$ with at least one $\leplus$, with $s_{\ell}$ and $t_{\ell}$ as defined in Theorem~\ref{shapes}.
  \end{enumerate}
}

As an example of the different cases for the $\mathbb{L}$-shapes and strings of $\leplus$'s, we encourage the reader to look at and verify the above statements for the example given in Figure~\ref{fig:algex}.

Now that we have our building blocks, what's left is to piece them together. \\

\bpoint{Glueing these shapes together}

The other surprising aspect is the fact that there is even a structure to $D$ and how these shapes glued together.
Furthermore, in different parts of $D$ which we will call sections, these shapes consistently glued together in the same way.

\tpointn{Definition}\label{sections}
\statement{
  Given a Young diagram $D$ of shape $\lambda = (\lambda_1, \cdots, \lambda_k)$, we partition $D$ into $k$ rectangles called \textbf{sections} where each section has dimension $j \times (\lambda_j - \lambda_{j+1})$ for $1 \leq j \leq k$, and is bounded by two rows (one possibly empty). We let $\lambda_{k+1} = 0$.\\
    \[ \sectionrow \]
  We name each section by its last row.\\
  Here we include empty sections into the count of $k$ (when $\lambda_j = \lambda_{j+1}$ for some $j$).
}

\tpointn{Theorem}\label{glueing}
\statement{
  $D$ is made up of glueing together these $\;\mathbb{L}$-shapes and strings of $\leplus$'s where in each (non-empty) section $b$, we have a chain of shapes in the following form:\\
    \[ \chain \]
  More precisely, we have the following properties for each section:
  \begin{itemize}
    \item There is at least one shape, either an $\;\mathbb{L}$ or a string of $\leplus$'s, in each section.
    \item There is either none or exactly one string of $\leplus$'s.
    \item There can be any number of $\;\mathbb{L}$-shapes, including 0.
    \item $\mathbb{L}$-shapes are glued on the left to other $\;\mathbb{L}$-shapes at the last $\leplus$ in the vertical part of a previous $\;\mathbb{L}$, except the first $\;\mathbb{L}$ in the chain which is glued to row $b$.
    \item The string of $\leplus$'s is either glued on the left to the last $\;\mathbb{L}$-shape in the chain at the last $\leplus$ in the vertical part, or glued to row $b$, in which case it fills the whole width of a section as such:\\
      \[ \sectiononlystring \]
  \end{itemize}
  The rest of the boxes in $D$ (the shaded regions) are filled with $0$'s.
}

\begin{proof}
  Putting together Theorem~\ref{shapes} and Corollary~\ref{strings}, $D$ is made up of the following two types of blocks, one consisting of $\mathbb{L}$-shapes where $\ell$ corresponds to a column in with at least one $\leplus$ in $\hat{D}$, and the other consists of strings of $\leplus$'s which corresponds to consecutive columns of all 0's in $\hat{D}$.\\
  \[ \blocks \]
  In the blocks, the shaded regions represent boxes filled with all 0's, which extends to fill the rest of the rows in columns $m$ to $\ell$ in the shape $\lambda$.
  The bolded lines represent the border of $\lambda$.

  By construction we cannot have a row $b$ occurring in-between columns $m$ to $\ell$ for either shape, as that would violate Step 1 for row $b$ regardless of type since there are no $\leplus$'s in columns $m$ to $\ell-1$ in $\hat{D}$ and thus the last $\leplus$ for row $b$ must be to the left of column $\ell-1$.
  For example, we cannot have:
     \[ \rowinbetweenl \]

  In fact, we actually know more about these blocks since they cannot span across multiple rows.
  That is, each section of $D$ is made up of a chain of $\;\mathbb{L}$-shapes and a string of $\leplus$'s as such
  \[ \chain \]
  since we know the following properties:
  \begin{itemize}
    \item $\mathbb{L}$-shapes start from a column with at least one $\leplus$ in $\hat{D}$ and extend to the right to the previous column of $\leplus$'s or the start of a row.
      In particular, if there is a column with at least one $\leplus$ in section $b$ of $\hat{D}$, then there will be a first (right-most) $\mathbb{L}$ in section $b$ of $D$ which has $m-1=b_B =b$:
        \[ \firstl \]
      Here $\ell_1$ is the first column of $\leplus$'s in $\hat{D}$ in section $b$.
    \item Strings of $\leplus$'s (not in the $\mathbb{L}$-shapes) extend to the left until the start of the next row below, and to the right until a column from $\hat{D}$ with at least one $\leplus$.
      In the special case that section $b$ in $\hat{D}$ consists of columns of all $0$'s, then the string of $\leplus$'s extends to fill the whole bottom row $b$ of the section (by Step 1 of the algorithm):
        \[ \sectiononlystring \]
  \end{itemize}

  What's left is to show is where the shapes glue together.
  Consider a particular section $b$.
  First, we already saw the special case when a section contains only columns of all $0$'s, we can assume there is at least one $\mathbb{L}$ in section $b$.
  We also saw that the first $\mathbb{L}$ in the section will always glue to the beginning of the row.
  As we saw in Corollary~\ref{strings} that string's of $\leplus$'s can be thought of as special cases of $\mathbb{L}$'s, all we need to show is how $\mathbb{L}$'s glue on the left to other $\mathbb{L}$'s.

  Consider an $\mathbb{L}$-shape with vertical part in column $\ell$ (or with no vertical part) and has right-most and bottom-most box at $(b_B, m)$.
  And we want to glue to this the previous $\mathbb{L}$ with vertical part in column $m-1$, in rows $i$ to $j$.
  We want to see where $b_B$ is in relation this previous $\mathbb{L}$:

  \[ \glueproofone \]

  If $(b_B, m)$ is filled with a $\leplus$, since $D$ is a \Le-diagram and $(i, m-1)$ is always filled with a $\leplus$, we either have $b_B < i$, or $i \leq b_B \leq j$ and $(b_B, m-1)$ is also filled with a $\leplus$.
  We also know that if $g$ is the last $\leplus$ in column $\ell$ in $\hat{D}$, then $b_B \geq g$.
  Since $b_B$ is either $g$ or $a_n$ if $g < a_n < \ell$, all rows $c$ below $b_B$ (in section $b$) must have their last $\leplus$ to the right of column $m$ and in particular no such $c$ is $a_n$.
  Furthermore, no such row $c$ is of type (II) since this would mean $a_n$ is below row $c$ with $c < a_n < m-1$.
  But then this would contradict what $b_B$ is, as this inequality suggests $b_B = a_n$ since $g \leq b_B < b$ and $m-1 < \ell$.
  Thus all rows $b$ below $b_B$ must be of type (I), with last $\leplus$ in a column $< m$.
  By Step 3 of the algorithm for row $c$, regardless if the last $\leplus$ is in column $m-1$ or not, this means $(b, m-1)$ is filled with a $0$.
  Since this holds for all rows $c$ below $b_B$ in section $b$, we must have that $b_B$ is the row of the last  $\leplus$ in column $m-1$ i.e. in the previous $\mathbb{L}$:

  \[ \glueprooftwo \]

  If $(b_B, m)$ is filled with a 0, we are in the special case where $m = \ell$.
  We also have $b_B \neq a_n$, since in the proof of Theorem~\ref{shapes} we saw that if $a_n$ was in the vertical part of a $\mathbb{L}$ then $(a_n, \ell)$ would be filled with a $\leplus$.
  In particular, this means row $b_B$ is of type (I) with its last $\leplus$ in $\hat{D}$ exactly at $(b_B, m)$.
  But now we have the same analysis for the rows $c$ below $b_B$ as the $\leplus$ case, since by the \Le-condition, these rows must have their last $\leplus$'s to the right of column $m$.
  Once again, this gives that $b_B$ is the row of the last $\leplus$ in column $m-1$.

  \[ \glueproofthree \]

  Thus we get that $\mathbb{L}$-shapes, and therefore also strings of $\leplus$'s, glue to the left of $\mathbb{L}$-shapes at the last (bottom-most) $\leplus$ in the vertical part of the $\mathbb{L}$.

  Putting everything together, we get that every section $b$ of $D$ looks like one of the following, giving the properties stated in the statement, making up the whole diagram $D$!
  \begin{figure}[h]
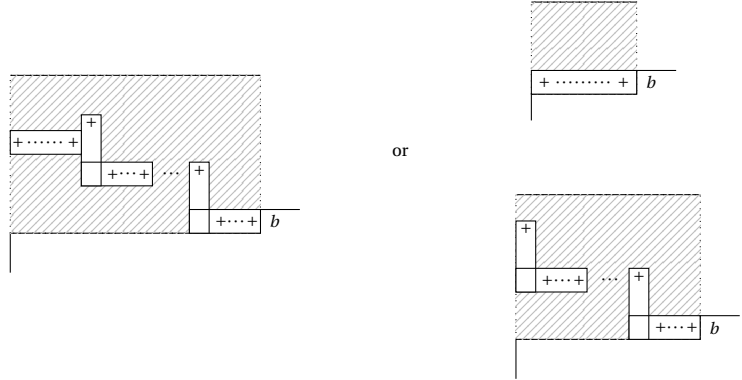

    \centering
    \[ \glueproof \]
    \caption[The possible chains of shapes in each section of $D$.]{The possible chains of shapes in each section of $D$. Generally (on the left), we have a chain of $\mathbb{L}$-shapes with a string of $\leplus$'s on the very left. In the case that section $b$ in $\hat{D}$ consists of columns of all $0$'s, we only have a string of $\leplus$'s filling the whole bottom row of the section (top-right figure). In the case that there is a column of $\leplus$'s in $\hat{D}$ right before the next row below, there is only a chain of $\mathbb{L}$-shapes (bottom-right figure). The shaded regions are filled with $0$'s.}
    \label{fig:glueproof}
  \end{figure}

\end{proof}

\bpoint{Iterating the algorithm}

If $D$ is co-loopless, that is $\pi$ has no fixed points, then the algorithm can be applied again, this time starting with $D$, to obtain a \Le-diagram $\bar{D}$.
Assuming that our algorithm is correct (which will be proven in Section~\ref{S:proof}), the associated decorated permutation to $\bar{D}$ would be:
\[\bar\pi = \left(\begin{array}{cccc} 1 & \cdots & n-1 & n \\ a_2 & \cdots & a_{n} & a_{1} \end{array}\right) \]
that is from the original $\hat\pi$, the permutation is shifted to the left twice.

First, we know that the shape of $\bar{D}$ corresponds to removing the column $a_1 = \pi(1)$ from $D$ and adding in a row labelled $a_1$ in the appropriate place.
For the filling of $\bar{D}$, from Theorem~\ref{glueing} we know that we really only need to look at what happens the $\mathbb{L}$-shapes.
Applying Theorem~\ref{shapes} to an $\mathbb{L}$-shape gives:\\
\[ \doublel \]
where now every $\leplus$ in the horizontal string of $\leplus$'s turns into its own $\mathbb{L}$, just with no horizontal part. \\

\section{Proof of the algorithm}\label{S:proof}

Finally, we show that the algorithm in Section~\ref{SS:algorithm} actually gives us the correct \Le-diagram.
We already have that $D$ is a \Le-diagram of type $(k+1, n)$ from Theorem~\ref{Dlediag}, so what's left is to show that $D$ is loopless, has the right dimension, and has $\pi$ as its associated decorated permutation. \\

\bpoint{Dimension and looplessness of \texorpdfstring{$D$}{D}}

Other than verifying that the algorithm is correct, using the visual perspective of the algorithm, we can clearly see how the dimension of $D$ arises in relation to the dimension of $\hat{D}$.
Namely, it comes from having at least one $\leplus$ in every column of $D$, with the number of additional $\leplus$'s being exactly the number of non-last $\leplus$'s in $\hat{D}$.

\tpointn{Theorem}\label{dimension}
\statement{
  Under this filling, $D$ is loopless and has dimension $\dim(S_D) = \dim(S_{\hat{D}}) - 2k + (n-1)$.
}

\begin{proof}
  For this proof, we only need Corollary~\ref{number} and refer to (i) and (ii) from that statement.

  First, it follows automatically from having $\geq 1$ $\leplus$ in every column of $D$ that $D$ is loopless.
  For the dimension $S_D$, (i) and the $+1$ in (ii) gives one $\leplus$ in every column of $D$ for a total of $n - (k+1)$ $\leplus$'s.

  Summing up what's left over in (ii) gives:
  \begin{align*}
    \sum_{\ell} s_{\ell} + t_{\ell}
      &= \sum_{\ell} s_{\ell} + \sum_{\ell} t_{\ell} \\
      &= \left( \# \text{ of non-last $\leplus$'s in $\hat{D}$ } - \# \text{ of non-last $\leplus$'s in col. $a_n$ of $\hat{D}$ }\right) \\
      &\quad\quad + \left( \# \text{ of non-last $\leplus$'s in col. $a_n$ of $\hat{D}$ }\right)  \\
      &= \# \text{ of non-last $\leplus$'s in $\hat{D}$ } \\
      &= \dim(S_{\hat{D}}) - k
  \end{align*}
  where the sums run over columns $\ell \neq a_n$ in $\hat{D}$ with at least one $\leplus$.
  The equalities come from:
  \begin{itemize}[itemsep=4pt]
    \item Summing over $s_{\ell}$ gives the number of non-last $\leplus$'s in each column in $\hat{D}$ excluding column $a_n$.
    \item Summing over $t_{\ell}$ gives the number of rows of type (II) in $\hat{D}$ with last $\leplus$'s not in column $a_n$, which in other words is the number of non-last $\leplus$'s in column $a_n$.
    \item As $\hat{D}$ is co-loopless and there are $k$ rows all with at least one $\leplus$, there are $k$ last $\leplus$'s out of a total of $\dim(S_{\hat{D}})$, which gives us the last equality.
  \end{itemize}

  Putting the two together gives:
  \begin{align*}
    \# \text{ of $\leplus's$ in $D$}
      &= \sum_{\text{cols in (i)}} 1 + \sum_{\text{cols in (ii)}} s_{\ell} + t_{\ell} + 1 \\
      &= n - (k+1) + \dim(S_{\hat{D}}) - k \\
      &= \dim(S_{\hat{D}}) - 2k + (n-1)
  \end{align*}
  which verifies that $\dim(S_D) = \dim(S_{\hat{D}}) - 2k + (n-1)$ as needed.
\end{proof}

\bpoint{\texorpdfstring{$D$}{D} corresponds to the correct decorated permutation}

Recall that $\hat{D}$ is a co-loopless \Le-diagram of type $(k,n)$ associated to the co-loopless decorated permutation $\hat\pi$ with $k$ anti-excedances, where
\[\hat\pi = \left(\begin{array}{cccc} 1 & 2 & \cdots & n \\ a_n & a_1 & \cdots & a_{n-1}\end{array}\right). \]
Applying the algorithm in Section~\ref{S:algorithm} constructs a loopless \Le-diagram $D$ of type $(k+1,n)$, with a particular dimensional relationship to $\hat{D}$.
We finally prove that $D$, as constructed, is the \Le-diagram associated to the decorated permutation $\pi$, as given by the T-duality map from $\hat\pi \mapsto \pi$ (see Section~\ref{SS:t-duality}).

\tpointn{Theorem}\label{permutation}
\statement{
  Under this filling, the associated decorated permutation to $D$ is
  \[\pi = \left(\begin{array}{cccc} 1 & 2 & \cdots & n \\ a_1 & a_2 & \cdots & a_{n}\end{array}\right), \]
  where $\pi$ is loopless and has $k+1$ anti-excedances.
}

\begin{proof}
Since $D$ is loopless and of type $(k+1,n)$, its associated decorated permutation $\pi$ is also loopless and will have $k+1$ anti-excedances.
Now, what we want to prove is that
  \[ \begin{cases} \pi(i) = \hat\pi(i+1) = \begin{cases} a_i \neq i+1 & \text{ for non-fixed points $i+1$ of $\hat\pi$ } \\ i+1 & \text{ for fixed points $i+1$ of $\hat\pi$ } \end{cases} & \text{ for } i \neq n \\ \pi(n) = \hat\pi(1) = a_n & \end{cases} \]

Recall that to get from $D$ to its decorated permutation, we go through its pipe dream (see Section~\ref{SS:le-def}).
  Starting from the label $i$ on the SE border, we follow its pipe until it reaches a label $j$ on the NW border of $D$, indicating that $\pi(i) = j$.
To read off the pipes from $D$, we start with pipes going leftwards if $i$ is a row and upwards if $i$ is a column.
  Then boxes with $0$'s correspond to the pipes continuing straight in the direction it was going, while boxes with $\leplus$'s correspond to the pipes turning.
The pipes turn leftwards if it was originally going up, or upwards if originally going to the left.
See Figure~\ref{fig:perm-proof-ex} for an example of this process.

\begin{figure}[h]
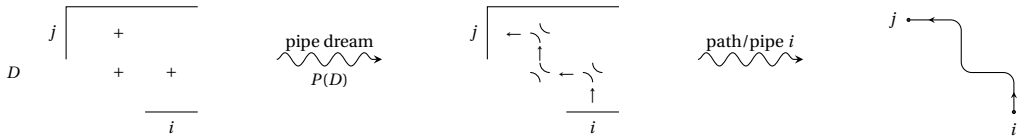

  \centering
  \[ \permproofexpath \]
  \caption[An example of reading off the decorated permutation $\pi$ from a \Le-diagram $D$.]{An example of reading off the decorated permutation $\pi$ from a \Le-diagram $D$. On the left, we have the portion of $D$ that affects label $i$, where unmarked boxes are filled with $0$'s. Translating to its pipe dream $P(D)$, we get the middle figure, where unmarked boxes are filled with crosses $\smallcross$ and the arrows indicate the direction of reading the pipe. On the right, we draw the pipe, also referred to as path, corresponding to $i$. That is, we have $\pi(i) = j$. Here we have depicted when $j < i$.}
  \label{fig:perm-proof-ex}
\end{figure}

We will also refer to the pipe starting at $i$ as the corresponding path for $\pi(i)$, or just $i$, where the turns are indicating where the $\leplus$'s are.
In the following figures, the shaded areas refer to boxes filled with 0's and are used to highlight explicit strings of 0's.

First, for the easiest case of $\pi(n)$, we look at the possibilities of what $n$ corresponds to in $D$:
  \begin{itemize}[itemsep=4pt]
    \item If $n$ is a row in $D$ (with no boxes), and in particular $\pi(n) = n$:
      \[ \permproofn \]
      we must have that $a_n = n$ since $\hat{D}$ being co-loopless (no empty rows or rows of all $0$'s) means that $n$ was a column in $\hat{D}$ with at least one box.
      Since the only change in rows and columns is through $a_n$, we get $\pi(n) = n = a_n$ as needed.
    \item If $n$ is a column in $D$, where previously in $\hat{D}$ $n$ was a column of all 0's i.e. $\hat\pi(n) = n$, then Corollary~\ref{strings} and its proof tells us there is a horizontal string of $\leplus$'s at the end of the diagram and in particular they will be in row $a_n$.
      Thus we get $n$ is a column in $D$ with exactly one $\leplus$ at $(a_n, n)$, which gives $\pi(n) = a_n$ as needed.
      \[ \permproofncolzero \]
    \item Finally, if $n$ is a column in $D$ (and thus $n \neq a_n$) and previously $n$ was a column in $\hat{D}$ with at least one $\leplus$ (i.e. $\hat\pi(n) \neq n$), then by Theorem~\ref{shapes} there is an $\mathbb{L}$-shape whose vertical part is in column $n$.
      Since we are at the last column of the diagram, we have that the $\leplus$'s in the vertical part of the $\mathbb{L}$ exactly comes from last $\leplus$'s in rows of type (II) in $\hat{D}$, with an extra $\leplus$ in row $a_n$ (since otherwise there needs to be some row whose last $\leplus$ is in some column $> n$).
      However for a row $b$ to be of type (II), i.e. with a $\leplus$ at $(b, a_n)$ in $\hat{D}$, we must have $a_n > b$.
      Thus the extra $\leplus$ at $(a_n, n)$ must be the bottom-most $\leplus$ in column $n$ in $D$.
      In particular, this means $\pi(n) = a_n$ as needed.
        \[ \permproofncol \]
  \end{itemize}

  Next, for the case with $i+1$ is a fixed point of $\hat\pi$ for $i \neq n$, that is $i+1 \neq a_n$ is a column of all 0's in $\hat{D}$, we look at what $i$ corresponds to in $D$:
  \begin{itemize}
    \item If $i$ is a row in $D$, regardless of row type for $i$, Step 1 places a $\leplus$ at $(i, i+1)$ and this is the only $\leplus$ in column $i+1$ in $D$.
      This gives $\pi(i) = i+1$ as needed.
      \[ \permprooffixedrow \]
    \item If $i$ is a column in $D$ where previously in $\hat{D}$ $i$ was a column of all 0's, then Corollary~\ref{strings} tells us there is a horizontal string of $\leplus$'s in columns $i, i+1$ in the same row, and these are the only $\leplus$'s, which gives $\pi(i) = i+1$ as needed.
        \[ \permprooffixedcolzero \]
    \item Finally, if $i$ is a column in $D$ (and thus $i \neq a_n$) where previously in $\hat{D}$ $i$ was a column with at least one $\leplus$, then by Theorem~\ref{shapes} there is an $\mathbb{L}$-shape with vertical part in column $i$.
      Since $i+1$ is a column of all 0's in $\hat{D}$, we know that the one $\leplus$ in column $i+1$ of $D$ corresponds to the rightmost $\leplus$ of either a horizontal string of $\leplus$'s or an $\mathbb{L}$-shape.
      In either case, by Theorem~\ref{glueing}, this shape must glue to the left of the $\mathbb{L}$ in column $i$ at the last $\leplus$ in column $i$, say at $(b_L, i)$.
      Thus we get $\leplus$'s at $(b_L, i)$ and $(b_L, i+1)$ with all 0's below in column $i$ and above in column $i+1$, which gives $\pi(i) = i+1$ as needed.
      \[ \permprooffixedcol \]
  \end{itemize}

  Lastly, we have the hardest case of when $i+1$ is not a fixed point of $\hat\pi$ for $i \neq n$, that is in $\hat{D}$ either $i+1$ is a column with at least one $\leplus$ or $i+1$ is a row.
  Say $\hat{\pi}(i+1) = j$ for some $j \neq i+1$ (we also have $j \neq a_n$).
  We want to show that $\pi(i) = j$.
  This is trickier as now there can be a complicated path in $\hat{D}$ to get to $j$ and it is not immediately clear how this translates in $D$ for $i$.
  In the subsequent proofs, Theorem~\ref{shapes}, Corollary~\ref{strings}, and Theorem~\ref{glueing} will be used without explicit citation.

  Figure~\ref{fig:perm-proof-ex}, when viewed as for \Le-diagram $\hat{D}$ with decorated permutation $\hat\pi$ and for $i+1$ instead of $i$, gives an example of a path when $i+1$ is not a fixed point.
  In general, every path for a non-fixed point $i+1$ must start with one $\leplus$ in row/column $i+1$ and then is built from alternating between two $\leplus$'s in the same row and two in the same column (where they share a $\leplus$), or vice versa, until row/column $j$ where the path ends with one $\leplus$ in row/column $j$.

  Consider the following general path in $\hat{D}$ which goes through $m+1$ columns of $\leplus$'s,
  \[ \permproofgeneralpath \]
  where $\ell_j$ denotes the columns of $\leplus$'s and $b_j$ denotes the rows of $\leplus$'s.
  Here, note that $i+1 > b_1 > \cdots > b_m$ and $\ell_1 < \cdots < \ell_{m+1}$.
  If instead $i+1$ is a row, there is an additional $\leplus$ in the first column of $\leplus$'s at $(i+1, \ell_1)$ where $i+1 > b_1$.
  If instead $j$ is a row, there is an additional $\leplus$ in the last column of $\leplus$'s at $(j, \ell_{m+1})$ where $j < b_m$.

  First, notice that by the \Le-condition for $\hat{D}$, there are two areas for $\leplus$'s where $\hat{D}$ must be filled with all 0's, depicted below by the shaded region around one $\leplus$.
  \[ \permproofpathzeros \]
  One area extends below to the next row of $\leplus$'s or the border of the diagram, and to the left until the end of the diagram.
  The other area extends on the right to the next column of $\leplus$'s or the border of the diagram, and above until the top of the diagram.
  In particular, the shaded areas in-between two $\leplus$'s in a path must be filled with all 0's, where for the column of $\leplus$'s these 0's extend to the left until the end of the diagram, and for the row of $\leplus$'s they extend to the top of the diagram.

  Now to tackle the problem at hand, we split the paths $\hat\pi(i+1)$ in $\hat{D}$ based on its relation to $a_n$ and consider them separately.
  As an example of the different types of paths, see Figure~\ref{fig:permpathex}.
  \begin{enumerate}[itemsep=4pt]
    \item $i+1 < j < a_n$: the path in $\hat{D}$ is completely to the right of $a_n$ and $j$ is a column.
    \item $j < i+1 \leq a_n$ or $i+1 \leq a_n < j$: the path passes through column $a_n$ and thus either the last $\leplus$ in the path is before $a_n$ (and $j$ is a row), or the path contains $\leplus$'s in column $a_n$.
    \item $j < a_n < i+1$ or $a_n < i+1 < j$: the path passes through where row $a_n$ in $D$ will be.
    \item $a_n < j < i+1$: the path is completely below where row $a_n$ in $D$ will be and $j$ is a row.
  \end{enumerate}

  \begin{figure}[h]
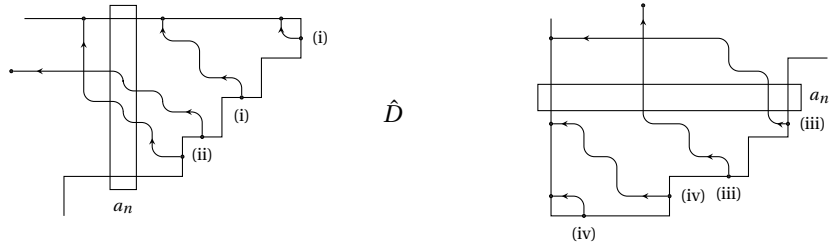

    \centering
    \[ \permproofpathex \]
    \caption[Examples of the four types of paths in $\hat{D}$ considered in the proof of $\pi$.]{Examples of the four types of paths in $\hat{D}$ considered in the proof of $\pi$. Each path is labelled by its type. On the left, the rectangle indicates where the column $a_n$ is in $\hat{D}$. On the right, the rectangle indicates where the row $a_n$ will be in $D$.}
    \label{fig:permpathex}
  \end{figure}

  Note all the paths in (ii) with $\leplus$'s to the left of $a_n$ must contain $\leplus$'s in column $a_n$ because of the \Le-condition for $\hat{D}$, since we know there is a $\leplus$ at $(1, a_n)$ if a box exists there.
  In the special case of $a_n = 1$, paths in (i) and (ii) do not exist.

  Starting with (i), we are in the case where:
  \[ \permproofihat \]
  Notice that the top $\leplus$ in each column in the path, except possibly $j$, is not the last $\leplus$ in their respective rows in $\hat{D}$.
  In the figure above, these are the circled $\leplus$'s.
  In particular this means that in the $\mathbb{L}$'s corresponding to these columns in $D$, there is also a $\leplus$ in these positions, $(b_1, \ell_1), \ldots (b_m, \ell_m)$; in the figures below these are the circled $\leplus$'s.

  We also notice that since $a_n > j$, all $\mathbb{L}$-shapes corresponding to columns $\ell_1, \ldots, \ell_m, j$ each have their $b_T \neq f, a_n$ where $f$ is the row of the first $\leplus$ in $\hat{D}$ in that column.
  Now let's look at two consecutive columns $\ell_s, \ell_{s+1}$, for $s \neq m$.
  \[ \permprooficols \]

  We know that any $\leplus$'s in column $\ell_s$ in-between rows $b_s, b_{s+1}$ must be last $\leplus$'s in their row because of the \Le-condition for $\hat{D}$.
  Thus, these rows are of type (I) and any last $\leplus$ at $(b, \ell_s)$ for $b_s > b > b_{s+1}$ corresponds to a 0 at $(b, \ell_s)$ in $D$.
  If there are $\leplus$'s in column $\ell_s$ above $b_{s+1}$ in $\hat{D}$, then by the \Le-condition there must be a $\leplus$ at $(b_{s+1},\ell_s)$ in $\hat{D}$.
  In particular, this is not the last $\leplus$ in row $b_{s+1}$ and thus there is a $\leplus$ at $(b_{s+1}, \ell_s)$ in $D$.
  If $(b_s, \ell_s)$ is the first $\leplus$ in column $\ell_s$, since there are no $\leplus$'s in rows $b$ for $b_s > b > b_{s+1}$ after column $\ell_s$ by the \Le-condition, then row $b_{s+1}$ is the first row above $b_s$ that has a $\leplus$ to the left of column $\ell_s$.
  That is, we have that $b_T = b_{s+1}$ for the $\mathbb{L}$ in column $\ell_s$ and thus there is a $\leplus$ at $(b_{s+1}, \ell_s)$.
  In either case, the $\mathbb{L}$ in column $\ell_s$ in $D$ always has a $\leplus$ at $(b_{s+1}, \ell_s)$ and at $(b_s, \ell_s)$, with 0's in-between.
  Since in $D$ the $\mathbb{L}$ in column $\ell_{s+1}$ glues to the left of the $\mathbb{L}$ in column $\ell_{s}$ at its bottom-most $\leplus$, there must be $0$'s in row $b_{s+1}$ in-between columns $\ell_s$ and $\ell_{s+1}$.
  Therefore, we get the above figure on the right.

  For column $j=\ell_{m+1}$, notice that in $\hat{D}$, there must be a row $b < b_m$ with a $\leplus$ to the left of column $j$ since we know at the least there's a $\leplus$ at $(1, a_n)$.
  Since we know that $(b_m,j)$ is the first $\leplus$ in column $j$ and that $b_T \neq j, a_n$ for the $\mathbb{L}$ in column $j$ in $D$, we must have that $b_T = b$ and thus there is a $\leplus$ at $(b, j)$ with 0's all above in $D$.
  For column $\ell_m$, using the same argument as for $\ell_s$ with row $b$ as $b_{s+1}$, in the $\mathbb{L}$ in column $\ell_m$, there must be a $\leplus$ at $(b, \ell_m)$ and at $(b_m, \ell_m)$ with 0's in-between.

  Since we know that $\mathbb{L}$'s glue to the left of previous $\mathbb{L}$'s at the first $\leplus$ from the bottom, we get the following in $D$, regardless if $i+1$ is a column or a row:
    \[ \permproofi \]

  Finally to look at what happens with $i$ in $D$.
  \begin{itemize}[itemsep=4pt]
    \item If $i+1 = \ell_1$ and $i$ are both columns, then regardless of how many $\leplus$'s there are in column $i$ in $\hat{D}$, we must have in $D$ a $\leplus$ at $(b_1, i)$ and 0's below in column $i-1$ since we know how $\mathbb{L}$'s glue.
      \[ \permproofibothcols \]
    \item If $i+1 = \ell_1$ and $i$ is a row, then we must have $b_1 = i$ (since $\hat{D}$ is co-loopless).
      \[ \permprooficolrow \]
    \item If $i+1$ is a row and $i$ is a row, then similarly we must have $b_1 = i$ (since $\hat{D}$ is co-loopless and by the \Le-condition as there are $\leplus$'s at $(i+1,\ell_1), (b_1,\ell_1)$).
      \[ \permproofibothrows \]
    \item If $i+1$ is a row and $i$ is a column, then notice $b_1 < i$ and is a row in $\hat{D}$ with its last $\leplus$ to the left of column $i$.
      If $i$ is a column of all 0's in $\hat{D}$, then Step 2(iii) for row $b_1$ would place the only $\leplus$ in column $i$ in $D$ at $(b_1, i)$.
      Otherwise, we can use the same argument for $i$ as for $\ell_s$, taking $\ell_1$ as $\ell_{s+1}$ (and only looking at rows above $i$), to show that there must be a $\leplus$ at $(b_1, i)$ and 0's all below in $D$.
      \[ \permproofirowcol \]
  \end{itemize}

  In all the cases, we can connect these path starts at the circled point to the rest of the path, and thus we get that the corresponding path to $i$ in $D$ leads to $j$, giving us $\pi(i) = j$ as needed for paths in (i).

  For paths in (iv), the proof is essentially the same as for (i) except now we need to consider when $j$ is a row.
  The other difference is we no longer have the fact that $b_T \neq f, a_n$ for columns $\ell_s$, which was only used when $(b_s, \ell_s)$ was the first $\leplus$ in its column in $\hat{D}$.
  However for paths in (iv), we have that if $(b_s, \ell_s)$ is the first $\leplus$ in the column, then $b_{s+1}$ is the first row above $b_s$ and below $a_n$ with a $\leplus$ to the left of $\ell_s$ (since $a_n < j$) and so we still have $b_T = b_{s+1}$ for the $\mathbb{L}$ in column $\ell_s$.

  When $j$ is a row, we have an additional $\leplus$ at $(j, \ell_{m+1})$ in column $\ell_{m+1}$ above $b_m$ in $\hat{D}$.
  Now we have $b = j$ for column $\ell_{m}$, and looking at column $\ell_{m+1}$, we notice that the $\leplus$ at $(j, \ell_{m+1})$ is the last $\leplus$ in row $j$.
  We also know that row $j$ is of type (I) since $j > a_n$, and thus in the $\mathbb{L}$ in column $\ell_{m+1}$ in $D$, there must be a 0 at $(j, \ell_{m+1})$.

  Thus we have, where the dotted line refers to where row $a_n$ will be in $D$:
  \[ \permproofiv \]

  The rest of the analysis is exactly the same for $i+1$ and $i$, and thus we also get that $\pi(i) = j$ as needed for paths in (iv).

  For paths in (ii), first notice we already have the case when the last $\leplus$ in the path is before $a_n$ and $j$ is a row.
  In this case, all the arguments from paths in (i) stay the same, and we can use the same argument for column $\ell_{m+1}$ as paths in (iv), since if $(j, \ell_{m+1})$ is the last $\leplus$ in row $j$ with $\ell_{m+1} < a_n$, then row $j$ is also of type (I).

  When the path contains $\leplus$'s in column $a_n$, we have (since $j \neq a_n$) either:
  \[ \permproofiihat \]
  where the dots are indicating $\leplus$'s in column $a_n$ that are not part of the path.

  Then notice for columns $\ell_1, \ldots, \ell_{s-1}$, the analysis from paths in (i) stays exactly the same.
  For columns $\ell_{s+1}, \ldots, \ell_{m+1}$, we note that all rows $b_{s}, \ldots, b_m$, and $j$ if it is a row, must be of type (II) since by \Le-condition, they must have $\leplus$'s in column $a_n$ as there is a $\leplus$ at $(1, a_n)$.
  In particular, this means that the $\mathbb{L}$'s corresponding to columns $\ell_{s+1}, \ldots, \ell_{m+1}$, which are all $> a_n$, has $\leplus$'s exactly in the same positions as in $\hat{D}$.
  However now we are done as to get $D$, column $a_n$ is removed and we know how the $\mathbb{L}$'s are glued, so we have:
  \[ \permproofii \]

  The analysis on $i+1$ and $i$ stays exactly the same as before (since it only depends on $\hat{D}$), even if $i+1 = \ell_1 = a_n$.
  In this case, the only difference is that the circled $\leplus$ in the analysis for $i$ in $D$ is now the $\leplus$ at $(b_1, \ell_2)$, and we consider $i+1$ to be a row.
  Thus we get that $\pi(i) = j$ as needed for paths in (ii).

  Lastly for paths in (iii), where the dotted line refers to where row $a_n$ will be in $D$, we have:
  \[ \permproofiiihat \]

  Here we essentially have the same proof as for (ii).
  Once again for columns $\ell_1, \ldots, \ell_{s-2}$, the analysis from paths in (iv) stays exactly the same.
  For columns $\ell_{s+1}, \ldots, \ell_{m+1}$, once again the rows $b_{s}, \ldots, b_{m}, b_{m+1},$ must be of type (II) since they are all $< a_n$ and have $\leplus$'s to the left of column $a_n$, and thus the $\mathbb{L}$'s corresponding to these columns have $\leplus$'s in the exact same positions as in $\hat{D}$.

  The only difference is we need to look at what happens to column $\ell_s$, which row $a_n$ will pass through, and column $\ell_{s-1}$.
  For column $\ell_s$, since row $a_n$ will occur in-between rows $b_s$ and $b_{s-1}$, there will be a $\leplus$ at $(a_n, \ell_s)$ in the $\mathbb{L}$ corresponding to column $\ell_s$.
  Since row $b_s$ is of type (II) with $\ell_s > a_n$, there is also a $\leplus$ at $(b_s, \ell_s)$ in the $\mathbb{L}$.

  For column $\ell_{s-1}$, since the $\leplus$ at $(b_{s-1}, \ell_{s-1})$ is not a last $\leplus$ in row $b_{s-1}$, there is a $\leplus$ in the same position in $D$.
  Now, if this $\leplus$ at $(b_{s-1}, \ell_{s-1})$ is the first one in column $\ell_{s-1}$ in $\hat{D}$, then we know that the corresponding $\mathbb{L}$ will have its $b_T = a_n$, in which case there is a $\leplus$ at $(a_n, \ell_{s-1})$.
  Otherwise if there are $\leplus$'s at $(b_{s-1}, \ell_{s-1})$ in $\hat{D}$, for $b_{s-1} > b > a_n$ these must be last $\leplus$'s in its row with $b$ of type (I) and thus corresponds to a 0 at $(b, \ell_{s-1})$ in $D$.
  For $b < a_n$, this means that $a_n$ is in-between the first and last $\leplus$'s in column $\ell_{s-1}$ and thus there must be a $\leplus$ at $(a_n, \ell_{s-1})$.
  In either case, the $\mathbb{L}$ in column $\ell_s-1$ always has $\leplus$'s at $(a_n, \ell_{s-1})$ and $(b_{s-1}, \ell_{s-1})$ with 0's in-between and so we have:\\
  \[ \permproofiii \]

  In the special case of $a_n < b_m$ where $j = \ell_{m+1} = \ell_s$, then since we know that the $\leplus$ at $(b_m, j)$ is the first in its column in $\hat{D}$, the $\mathbb{L}$ in column $j$ in $D$ has its $b_T = a_n$ with a $\leplus$ at $(a_n, j)$.
  Column $\ell_{m} = \ell_{s-1}$ stays the same as before with $\leplus$'s at $(a_n, \ell_m)$ and $(b_m, \ell_m)$, and 0's in-between, in $D$:\\
  \[ \permproofiiispecialone \]

  In the special case of $a_n > b_1$ where $i+1 = \ell_1$, now the $\mathbb{L}$ in column $\ell_1 = \ell_s$ in $D$ will have a $\leplus$ at $(a_n, \ell_1)$, since $\ell_1 > a_n$ and thus its $b_B = a_n$.
  There is also still a $\leplus$ at $(b_1, \ell_1)$, since $b_1 < a_n$ and thus is of type (II), which aligns with the $\ell_s$ as in the general case.\\
  \[ \permproofiiispecialtwo \]

  The analysis on $i+1$ and $i$ stays exactly the same as before, with the slight exception of $a_n > b_1$.
  Since we have $i \geq a_n$ with a $\leplus$ at $(a_n, \ell_1)$ in $D$, now $a_n$ plays the role of $b_1$ and we use the same arguments as before, even if $i = a_n$ in which case we consider $i$ to be a row.
  In all cases, we get that $\pi(i) = j$ as needed for paths in (iv).

  Putting (i)--(iv) together, we are now done since for all paths $\hat\pi(i+1) = j$ we get $\pi(i) = \hat\pi(i+1) = j$!
\end{proof}

\bpoint{T-duality on the level of \texorpdfstring{\BLe-}{Le }diagrams}

With the proof that $D$ gives the correct decorated permutation complete, we finally conclude by putting everything together.
Theorems~\ref{Dlediag},~\ref{dimension} and~\ref{permutation}, tell us that in fact, the algorithm given in Section~\ref{S:algorithm} is telling us what the T-duality map (see Section~\ref{SS:t-duality}) looks like on the level of \Le-diagrams.

\tpointn{Theorem}\label{letheorem}
\statement{
  Given a co-loopless \Le-diagram $\hat{D}$ of type $(k,n)$, the algorithm given in Section~\ref{S:algorithm} constructs a loopless \Le-diagram $D$ of type $(k+1,n)$ such that
  \begin{itemize}[itemsep=4pt]
    \item the dimensions of the positroid cells indexed by $D$ and $\hat{D}$ relate via
      \[ \dim(S_D) = \dim(S_{\hat{D}}) - 2k + (n-1), \]
    \item and the associated decorated permutations to $D$ and $\hat{D}$ relate via T-duality.
  \end{itemize}
  That is, the algorithm gives the T-duality map from $\hat{\pi} \mapsto \pi$ on the level of \Le-diagrams.
}

\newpage
\newrefsegment
\part{Final thoughts}\label{P:final}
\chapter{Conclusion}\label{C:conclusion}

\section{Further directions} \label{S:further}

In Part~\ref{P:c2}, we started with the tale of the $c_2$-invariant as an avatar for the Feynman period and saw how combinatorial techniques involving enumerating edge bipartitions through swapping arguments were able to resolve the $p=2$ case of the long sought after completion conjecture for the $c_2$-invariant.
However, the story is not over as of course what we would really like is a proof of the full conjecture.

Thus, the pivotal next direction to look at is: \\
Using this idea of swapping edges around a vertex, or vertices, across the edge bipartition, can we generalize the counting arguments seen in Chapter~\ref{C:completion} for higher primes $p$?

For $p > 2$, the combinatorial interpretation of the $c_2$-invariant becomes much more intricate as now we are enumerating the ways of distributing $p-1$ copies of each edge across $2p-2$ polynomials, where $p-1$ of them come from spanning trees and the other $p-1$ arises from particular spanning 2-forests.
Comparing to the $p=2$ case, recall that we were bipartitioning one copy of each edge into a spanning tree, spanning 2-forest pair.
The difficulty in higher $p$ is that now there no longer needs to be a nice pairing up of trees and forests that exactly partition full sets of edges from the graph.
Optimistically, Yeats~\cite{higherp} has been able to first answer this question in the $T$-case where we only needed to swap around 2-valent vertices.
The obstacle, naturally, in the $S$ and $R$ cases is in how to generalize this to include control vertices or perhaps instead how to circumvent the need for control vertices.

Finally, one could also ask about how much farther these enumerative techniques can take us.
Can they answer other questions about the $c_2$-invariant or compute the $c_2$'s for new families of graphs?

In Part~\ref{P:le}, we looked at the tale of the combinatorial T-duality map as a bridge between triangulations of the hypersimplex and triangulations of the amplituhedron and what this map looked like on \Le-diagrams.
One result, other than the nice structure of the \Le-diagrams arising from T-duality, was that we could directly see where the dimensional relationship between the positroid cells on either side of the map manifests from.

Continuing with this thread, a natural general question to further investigate is what other properties of positroids or triangulations can be determined through looking at T-duality via \Le-diagrams?
Two specific directions that could be of interest are:

\begin{itemize}[itemsep=4pt]
  \item The T-duality map turns out to be a special case of the $\rho_A(\pi)$ map defined by Benedetti, Chavez and Tamayo~\cite{BCT} when characterizing quotients of uniform positroids.
    Thus, can we extend our algorithm from Chapter~\ref{C:lebijection} to this more general cyclic shift map?
  \item In~\cite{PSBW}, Parisi, Sherman-Bennett, and Williams, were able to prove Theorem~\ref{triang-conj} through looking at the T-duality map via plabic graphs.
    Thus, are there any interesting connections between the two perspectives on T-duality outside of the canonical bijection between \Le-diagrams and plabic graphs from~\cite{tpgrass}? Furthermore, is there any structure to the \Le-diagrams underlying the triangulations on either side of T-duality?
\end{itemize}

\section{Final thoughts} \label{S:final}

Through the underlying theme of emerging combinatorics in the world of scattering amplitudes, one could see the parallels between the stories of the $c_2$-invariant and the T-duality map.
Both branched off of the need to better understand the scattering amplitude, one through Feynman diagrams in scalar $\phi^4$-theory and one through on-shell diagrams in $\mathcal{N}=4$ SYM.
Both then looked at related but simpler, and in some sense nicer, objects with interesting connections to many areas of mathematics, one reducing to the $c_2$-invariant and the other reducing to positroid cells of the positive Grassmannian.
In the end, both then converged to the idea that looking at these objects through a combinatorial perspective can be a powerful tool in answering questions about them, one being the completion conjecture for the $c_2$-invariant and the other being the relationship between triangulations of the amplituhedron and the hypersimplex.
At the heart of it all, they showcased a wonderful tale of beautiful mathematics at the intersection of combinatorics and quantum field theory!

\newpage
\makeatletter
\renewcommand\leftmark{\MakeUppercase{References}}
\renewcommand\rightmark{\MakeUppercase{References}}
\makeatother
\bookmarksetupnext{level=-1}
\addtocontents{toc}{\protect\setcounter{tocdepth}{0}}
\chapter*{References}
\printbibliography[segment=1,heading=subbibliography,title={Chapter 1. Introduction}]
\printbibliography[segment=2,heading=subbibliography,title={Part I. The \texorpdfstring{$c_2$}{c2}-invariant}]
\printbibliography[segment=3,heading=subbibliography,title={Part II. Le diagrams}]
\printbibliography[segment=4,heading=subbibliography,title={Part III. Final thoughts}]

\end{document}